\documentclass[12pt, english]{amsart}

\usepackage{amsmath}

\usepackage{amssymb}

\usepackage{enumerate}

\usepackage{amscd} \usepackage{tabularx}
\usepackage[all,cmtip,line]{xy}
\usepackage{graphicx}
\usepackage{tikz-cd}

\newcommand{\ra}{\rightarrow}
\newcommand{\lra}{\longrightarrow}
\newcommand{\la}{\leftarrow}
\newcommand{\lla}{\longleftarrow}
\newcommand{\into}{\hookrightarrow}

\newcommand{\ua}{\uparrow}
\newcommand{\da}{\downarrow}
\newcommand{\iso}{\stackrel{\sim}{\ra}}
\newcommand{\liso}{\stackrel{\sim}{\lra}}

\newcommand{\pfbegin}{{{\em Proof:}\;}}
\newcommand{\pfend}{$\Box$ \medskip}

\newlength{\ownl}

\newcommand{\norm}{{\mbox{\bf N}}}

\newcommand{\Art}{{\operatorname{Art}\,}}
\newcommand{\Aut}{{\operatorname{Aut}\,}}

\newcommand{\Frob}{{\operatorname{Frob}}}
\newcommand{\conju}{{\operatorname{conj}}}
\newcommand{\res}{{\operatorname{res}}}
\newcommand{\ind}{{\operatorname{ind}}}
\newcommand{\Gal}{{\operatorname{Gal}\,}}

\newcommand{\Hom}{{\operatorname{Hom}\,}}

\newcommand{\Id}{{\operatorname{Id}}}

\newcommand{\Ind}{{\operatorname{Ind}\,}}

\newcommand{\Map}{{\operatorname{Map}\,}}

\newcommand{\ad}{{\operatorname{ad}\,}}

\newcommand{\tr}{{\operatorname{tr}\,}}

\newcommand{\wt}{{\operatorname{wt}}}

\newcommand{\corr}{{\operatorname{cor}\,}}

\newcommand{\nak}{{\operatorname{nak}}}

\newcommand{\spl}{{\operatorname{sp}}}
\newcommand{\can}{{\operatorname{can}}}

\newcommand{\loc}{{\operatorname{loc}}}
\newcommand{\glob}{{\operatorname{glob}}}
\newcommand{\der}{{\operatorname{der}}}

\newcommand{\ab}{{\operatorname{ab}}}

\newcommand{\nr}{{\operatorname{nr}}}

\newcommand{\alg}{{\operatorname{alg}}}
\newcommand{\sco}{{\operatorname{SC}}}
\newcommand{\cts}{{\operatorname{cts}}}

\newcommand{\basic}{{\operatorname{basic}}}
\newcommand{\compact}{{\operatorname{compact}}}

\newcommand{\Tan}{{\operatorname{Tan}}}
\newcommand{\univ}{{\operatorname{univ}}}

\newcommand{\A}{{\mathbb{A}}}

\newcommand{\C}{{\mathbb{C}}}

\newcommand{\F}{{\mathbb{F}}}
\newcommand{\G}{{\mathbb{G}}}

\newcommand{\Q}{{\mathbb{Q}}}
\newcommand{\R}{{\mathbb{R}}}

\newcommand{\Z}{{\mathbb{Z}}}

\newcommand{\CH}{{\mathcal{H}}}

\newcommand{\ga}{{\mathfrak{a}}}

\newcommand{\gz}{{\mathfrak{z}}}

\newcommand{\barD}{\overline{{D}}}

\newcommand{\barF}{{\overline{{F}}}}

\newcommand{\barL}{\overline{{L}}}

\newcommand{\barb}{\overline{{b}}}

\newcommand{\bare}{{\overline{{e}}}}
\newcommand{\barf}{\overline{{f}}}
\newcommand{\barg}{\overline{{g}}}

\newcommand{\bark}{\overline{{k}}}

\newcommand{\bart}{\overline{{t}}}
\newcommand{\baru}{\overline{{u}}}

\newcommand{\barw}{\overline{{w}}}

\newcommand{\tN}{\widetilde{{N}}}

\newcommand{\tS}{\widetilde{{S}}}

\newcommand{\tb}{\widetilde{{b}}}

\newcommand{\tg}{{\widetilde{{g}}}}

\newcommand{\tv}{{\widetilde{{v}}}}

\newcommand{\baralpha   }{{\overline{\alpha  }}}   
   
\newcommand{\bargamma   }{{\overline{\gamma}}}

\newcommand{\barkappa     }{{\overline{\kappa}}}   
    
 \newcommand{\barmu    }{{\overline{\mu}}}   
 \newcommand{\barnu    }{{\overline{\nu}}}

\newcommand{\tgamma   }{{{\widetilde{\gamma}}}}   
   
\newcommand{\tepsilon    }{{\widetilde{\epsilon}}}

\newcommand{\teta         }{{{\widetilde{\eta}}}}

 \newcommand{\tmu    }{{\widetilde{\mu}}}

 \newcommand{\trho   }{{\widetilde{\rho}}}   
    
 \newcommand{\tsigma   }{{{\widetilde{\sigma}}}}   
    
 \newcommand{\ttau     }{{\widetilde{\tau}}}   
    
 \newcommand{\tphi    }{{\widetilde{\phi}}}

 \newcommand{\tGamma     }{{\widetilde{\Gamma}}}   
   
\newcommand{\tTheta       }{{\widetilde{\Theta}}}

 \newcommand{\hatlambda   }{{\widehat{\lambda}}}

\newcommand{\bnu}{{\boldsymbol \nu}}
\newcommand{\bzeta}{{\boldsymbol \zeta}}
\newcommand{\balpha}{{\boldsymbol \alpha}}
\newcommand{\bgamma}{{\boldsymbol \gamma}}
\newcommand{\bphi}{{\boldsymbol \phi}}
\newcommand{\bpsi}{{\boldsymbol \psi}}
\newcommand{\bhatlambda}{{{\boldsymbol {\widehat{\lambda}}}}}

\newcommand{\barFv}{{\overline{F_v}}}

\newcommand{\tcE}{{\widetilde{\cE}}}

\def\RCS$#1: #2 ${\expandafter\def\csname RCS#1\endcsname{#2}}
\RCS$Revision: 364 $
\RCS$Date: 2012-09-24 15:39:31 -0400 (Mon, 24 Sep 2012) $

\newcommand{\onto}{\twoheadrightarrow}

\newcommand{\cB}{\mathcal{B}}

\newcommand{\cE}{\mathcal{E}}

\newcommand{\cH}{\mathcal{H}}

\newcommand{\cN}{\mathcal{N}}
\newcommand{\cO}{\mathcal{O}}

\newcommand{\cR}{\mathcal{R}}

\newcommand{\cZ}{\mathcal{Z}}

\usepackage{hyperref}

\usepackage{amsthm}

\newtheorem{thm}{Theorem}[section]

\newtheorem{cor}[thm]{Corollary}
 
\newtheorem{lem}[thm]{Lemma} 
 \theoremstyle{definition}
 \theoremstyle{definition}
 \theoremstyle{remark}

\numberwithin{equation}{section}

\theoremstyle{definition}

\begin{document}

\title{Cocycles for Kottwitz cohomology}

\author{Jack Sempliner}\email{j.sempliner@imperial.ac.uk}\address{Imperial College, London}
\thanks{J.S. was supported partially by the National Science Foundation Graduate Research Fellowship (under Grant No. DGE-1656466), and partially by funding from the European Research Council (ERC) under the European Union's Horizon 2020 research and innovation program (grant agreement No. 884596).}

\author{Richard Taylor}\email{rltaylor@stanford.edu}\address{Stanford University}
\thanks{R.T. was supported in part by NSF grant DMS-1902265}
\maketitle

\section{Introduction}

Kottwitz \cite{kotbg} introduced certain extensions of local and global Galois groups and considered their `algebraic' cohomology with coefficients in an algebraic group.  However his extensions are only canonical up to conjugation. This is enough to make their cohomology canonical, but their spaces of cocyles are not canonical. The purpose of this paper is to explain how one can work with such cocyles. This will play a crucial role in our work on the formalism of Shimura varieties \cite{st}.

 If $E/F$ is a finite Galois extension of local or global fields Kottwitz considers certain extensions $\cE(E/F)$ of $\Gal(E/F)$ by certain abelian groups $\cE(E/F)^0$. The most familiar examples are the local and global Weil groups $W_{E/F}$: extensions of $\Gal(E/F)$ by $E^\times$ in the local case and by $\A_E^\times/E^\times$ in the global case. These extensions are defined in terms of a canonical class $[\alpha_{E/F}] \in H^2(\Gal(E/F),\cE(E/F)^0)$, but not by a canonical cocycle. It turns out (because $H^1(\Gal(E/F),\cE(E/F)^0)=(0)$) that $\cE(E/F)$ is unique up to an isomorphism, which is itself unique up to composition with conjugation by an element of $\cE(E/F)^0$.

Having defined these extensions $\cE(E/F)$ we will, following Kottwitz, consider what we call the {\em algebraic} (non-abelian) cohomology $H^1_\alg(\cE(E/F),G)$, where $G$ is a group (often the $E$ or $\A_E$ points of an algebraic group) with an action of $\Gal(E/F)$. To define this cohomology we consider only the set $Z^1_\alg(\cE(E/F),G) \subset Z^1(\cE(E/F),G)$ of cocycles whose restriction to $\cE(E/F)^0$ lie in some chosen class of homomorphisms, usually a class of homomorphisms coming from certain morphisms of algebraic groups. These `algebraic' cocycles will be preserved by the usual equivalence relation and hence give rise to a cohomology group $H^1_\alg(\cE(E/F),G)$. It is easy to verify that despite the ambiguity in the definition of $\cE(E/F)$, the pointed set $H^1_\alg(\cE(E/F),G)$ is well defined up to unique isomorphism. However $Z^1_\alg(\cE(E/F),G)$ is not.

If $F$ is a local field there will be only one such extension of interest to us: $\cE(E/F)=W_{E/F}$ - the Weil group defined by the usual canonical class $[\alpha_{E/F}]\in H^2(\Gal(E/F),E^\times)$. In this case, if $G/F$ is an algebraic group with centre $Z(G)$, then $Z^1_\alg(W_{E/F},G(E))$ (resp. $Z^1_\alg(W_{E/F},G(E))_\basic$) will denote those cocycles which are given on $E^\times$ by an algebraic character $\nu:\G_m\ra G$ (resp. $\nu: \G_m \ra Z(G)$). These are sometimes called `algebraic cocycles' (resp. `basic algebraic cocycles'). 

However, when $F$ is a global field, will need to consider several examples of these groups $\cE(E/F)$, which we will now describe.
\begin{enumerate}
\item We define $\cE^\loc(E/F)^{0}= \prod_{w\in V_E} E_w^\times$, where $V_E$ denotes the set of all places of $E$. There is a unique class $[\alpha^\loc_{E/F}] \in H^2(\Gal(E/F), \prod_{w \in  V_E}E_w)$ whose image in $H^2(\Gal(E_w/F_w),E_w^\times)$ equals $[\alpha_{E_w/F_w}]$ for all $w \in V_E$. We let $\cE^\loc(E/F)$ denote the corresponding extension of $\Gal(E/F)$ by $\prod_{w \in V_E} E_w^\times$. 

For an algebraic group $G/F$, basic algebraic cohomology of 
$G(\A_E)$ will be defined in terms of those cocycles whose restriction to $\cE^\loc(E/F)^0$ are of the form $\prod_w \nu_w$, where $\nu_w:\G_m \ra Z(G)_{/F_w}$ is an algebraic character, non-trivial for only finitely many $w$.

\item We define $T_{2,E}/F$ to be the protorus with character group $\Z[V_E]$, the free abelian group with basis the set $V_E$ of places of $E$ with the natural action of $\Gal(E/F)$. Then we set $\cE_2(E/F)^0=T_{2,E}(\A_E)$ and define $\cE_2(E/F)$ as the pushout of $\cE^\loc(E/F)$ along the embedding $\prod_w E_w^\times \into \prod_w \A_E^\times \cong T_{2,E}(\A_E)$, where we identify $E_w^\times$ inside inside the copy of $\A_E^\times$ indexed by $w$.

In this case, for an algebraic group $G/F$, basic algebraic cohomology of 
$G(\A_E)$ will be defined in terms of those cocycles whose restriction to $\cE_2(E/F)^0$ come from an algebraic character $\nu:T_{2,E} \ra Z(G)_{/F}$. Thus there are natural restriction maps $Z^1_\alg(\cE_2(E/F),G(\A_E))_\basic \ra Z^1_\alg(\cE^\loc(E/F),G(\A_E))_\basic$.

\item $W_{E/F}$ will denote the global Weil group, i.e. the extension of $\Gal(E/F)$ by $\A_E^\times/E^\times$ coming from the usual canonical class $[\alpha_{E/F}]\in 
H^2(\Gal(E/F),\A_E^\times/E^\times)$.

\item We will write $\cE^\glob(E/F)^0$ for the subgroup of elements of $T_{2,E}(\A_E)$ whose image in $\A_E^\times/E^\times$ under any of the characters $\pi_w$ corresponding to $w \in V_E$ is independent of $w$.  It turns out (as observed by Nakayama and Tate), that  there is a unique class $[\alpha_{E/F}^\glob] \in H^2(\Gal(E/F),\cE^\glob(E/F)^0)$ which pushes forward to $[\alpha_{E/F}] \in H^2(\Gal(E/F),\A_E^\times/E^\times)$ and to $[\alpha_{E/F}^\loc] \in H^2(\Gal(E/F),T_{2,E}(\A_E))$. We write $\cE^\glob(E/F)$ for the corresponding extension of $\Gal(E/F)$ by $\cE^\glob(E/F)^0$.

In this case, for an algebraic group $G/F$, basic algebraic cohomology of 
$G(\A_E)$ will be defined in terms of those cocycles whose restriction to $\cE_2(E/F)^0$ come from an algebraic character $\nu:T_{2,E} \ra Z(G)_{/F}$.

There are embeddings of extensions $\loc_\ga: \cE^\glob(E/F) \into \cE_2(E/F)$ giving rise to isomorphisms 
\[ \loc_\ga=(\loc_\ga^*)^{-1}: Z^1_\alg(\cE^\glob(E/F),G(\A_E))_\basic \liso Z^1_\alg(\cE_2(E/F),G(\A_E))_\basic.\]
The map of extensions is only defined up to composition with conjugation by an element of $T_{2,E}(\A_E)$; and the map of cocycles is canonically defined only up to composition with the map from $Z^1_\alg(\cE_2(E/F),G(\A_E))_\basic$ to itself given by $\phi \mapsto {}^{\phi(t)} \phi$ for some $t \in T_{2,E}(\A_E)$.

\item Finally we will write $T_{3,E}$ for the protorus over $F$ with character group $\Z[V_E]_0$, the subabelian group of $\Z[V_E]$ consisting of elements $\sum m_w w$ for which $\sum m_w=0$. We will write $\cE_3(E/F)$ for the pushout of $\cE^\glob(E/F)$ along $\cE^\glob(E/F)^0 \onto T_{3,E}(E)$. 

In this case, for an algebraic group $G/F$, basic algebraic cohomology of 
$G(E)$ will be defined in terms of those cocycles whose restriction to $T_{3,E}(E)$ come from an algebraic character $\nu:T_{2,E} \ra Z(G)_{/F}$.  Thus there is a natural morphism $Z^1_\alg(\cE_3(E/F),G(E))_\basic \ra Z^1_\alg(\cE^\glob(E/F),G(\A_E))_\basic$. 

\end{enumerate}

We have a diagram of morphisms of extensions:
\[ \begin{array}{rcccccl} && \cE_3(E/F) &\twoheadleftarrow & \cE^\glob(E/F) &\twoheadrightarrow& W_{E/F} \\ &&&& \loc_\ga \da && \\  & & \cE^\loc(E/F) & \into & \cE_2(E/F) && \\ && \bigcup &&&& \\ W_{E_w/F_v} &\twoheadleftarrow & \cE^\loc(E/F)|_{\Gal(E/F)_w} &&&& \end{array} \]
for any places $w|v$ of $E$ and $F$.
A key observation is that, although individually the extensions $\cE(E/F)$ we consider here have automorphisms, the diagram as a whole does not. Thus if we fix such a diagram it makes sense to consider algebraic cocycles and not just algebraic cohomology classes.  

However we have not specified such a diagram uniquely. There are many choices for the localization map $\loc_\ga$. To the best of our knowledge there is no preferred choice. The various choices form a set which we will denote $\cH(E/F)$, which comes with a transitive action of $T_{2,E}(\A_E)$. It seems to us that choosing one element of $\cH(E/F)$ is a bit like choosing one place of $\overline{E}$ above a given place of $E$: there are many choices, but for most purposes the choice is irrelevant. We will decorate the various extensions $\cE^?(E/F)$ with a subscript $\ga$ to indicate it is the one uniquely determined up to unique isomorphism by $\ga$. If $\ga,\ga' \in \cH(E/F)$ and $\ga'={}^t\ga$ for some $t\in T_{2,E}(\A_E)$, then we get isomorphisms
\[ z_t: Z^1_\alg(\cE^?(E/F)_\ga,G(A_E)) \liso Z^1_\alg(\cE^?(E/F)_{\ga'},G(A_E)), \]
where $A_E$ denotes $E$ or $\A_E$. There may not be a unique choice of $t$ and the isomorphism may depend on the $t$ chosen. However, after one passes to cohomology groups, it will no longer depend on the choice of $t$.

Given a choice $\ga \in \cH(E/F)$ we have maps
\[ \loc_\ga: Z^1_\alg(\cE_3(E/F)_\ga,G(E))_\basic \lra Z^1_\alg(\cE^\loc(E/F)_\ga,G(\A_E))_\basic \]
and
\[ \res: Z^1_\alg(\cE^\loc(E/F)_\ga,G(\A_E))_\basic \lra \prod_{v\in V_F} Z^1_\alg(W_{E_{w(v)}/F_w,\ga},G(E_w))_\basic, \]
where for each $v\in V_F$ we choose a place $w(v)$ of $E$ above it.
These induce maps in cohomology
\[ \loc: H^1_\alg(\cE_3(E/F),G(E))_\basic \lra H^1_\alg(\cE^\loc(E/F),G(\A_E))_\basic \]
and
\[ \res: H^1_\alg(\cE^\loc(E/F),G(\A_E))_\basic \lra \prod_{w\in V_E} H^1_\alg(W_{E_w/F_w},G(E_w))_\basic, \]
which are canonically independent of $\ga$. Moreover if $E'/F'$ is an extension of local fields isomorphic (but not canonically so) to $E_w/F_w$, then there is a well defined map
\[ \res_{E'/F'}: H^1_\alg(\cE^\loc(E/F),G(\A_E))_\basic \lra H^1_\alg(W_{E'/F'},G(E'))_\basic. \]

If $F$ is local and $\phi \in Z^1_\alg(W_{E/F,\ga },G(E))_\basic$ or if $F$ is global and $\phi \in Z^1_\alg(\cE_3(E/F)_\ga,G(E))_\basic$, then $\ad \phi \in Z^1(\Gal(E/F),(G/Z(G))(E))$ and we get an inner form, referred to as an `extended pure inner form',  ${}^\phi G$ of $G$ over $F$. If $Z(G)$ is a torus then Kottwitz showed that
\[ H^1_\alg(W_{E/F},G(E))_\basic\onto H^1(\Gal(E/F),(G/Z(G))(E)) \]
and
\[ H^1_\alg(\cE_3(E/F),G(E))_\basic\onto H^1(\Gal(E/F),(G/Z(G))(E)) \]
are surjective.

We recall that when $G$ is connected reductive, Kottwitz defined important maps
\[ \kappa: H^1_\alg(\cE^\loc(E/F),G(\A_E))_\basic \lra (\Z[V_E] \otimes \Lambda_{G})_{\Gal(E/F)} \]
and 
\[ \kappa: H^1_\alg(\cE_3(E/F),G(E))_\basic \lra (\Z[V_E]_0 \otimes \Lambda_{G})_{\Gal(E/F)}, \]
where $\Lambda_G$ denotes the algebraic fundamental group of $G$. They are compatible in that $\kappa \circ \loc$ equals $\kappa$ composed with the obvious map $
(\Z[V_E]_0 \otimes \Lambda_{G})_{\Gal(E/F)} \ra (\Z[V_E] \otimes \Lambda_{G})_{\Gal(E/F)}$. We will also write 
\[ \barkappa_G : H^1_\alg(\cE^\loc(E/F),G(\A_E)) \lra \Lambda_{G,\Gal(E/F)} \]
for the composition of $\kappa$ with the map
\[ \begin{array}{rcl} \Z[V_E] \otimes \Lambda_{G} &\lra& \Lambda_G \\ \sum_w x_w w & \longmapsto & \sum_w x_w. \end{array} \]

It turns out it is useful to fix slightly more than $\ga \in \cH(E/F)$. Although the Weil group $W_{E/F}$ is only determined up to inner automorphism by an element of $\A_E^\times/E^\times$, the absolute Weil group $W_{\barF/F}$ is much more rigid. 
It is determined up to conjugation by an element of $\ker(W_{\barF/F} \ra \Gal(\barF/F))$. The extra data we will add is roughly speaking a collection of isomorphisms between $W_{\barF/F}/\overline{[W_{\barF/E},W_{\barF/E}]}$ and $W_{E/F}$. One might wonder why one works with $W_{E/F}$ at all, and not $W_{\barF/F}/\overline{[W_{\barF/E},W_{\barF/E}]}$ directly. The answer seems to be that, when $D\supset E \supset F$, the way we compare $Z^1_\alg(\cE^?(E/F),G(A_E))_\basic$ and $Z^1_\alg(\cE^?(D/F),G(A_D))_\basic$ is {\em not} compatible with the natural map
\[ W_{\barF/F}/\overline{[W_{\barF/D},W_{\barF/D}]} \onto W_{\barF/F}/\overline{[W_{\barF/E},W_{\barF/E}]}. \]

More precisely by {\em complete rigidification data} for $\ga \in \cH(E/F)$ we will mean the choice for each place $v$ of $F$ and each $F$-linear embedding $\rho:E^\ab \into \barFv$ (giving rise to a place $w(\rho)$ of $E$) a $E_{w(\rho)}^\times$-conjugacy class $[\Gamma_{v,\rho}]$ of isomorphisms of extensions
\[ \begin{array}{ccccccccc} (0)& \lra & \Gal(E^\ab/E) & {\lra} & \Gal(E^\ab/F) & {\lra} & \Gal(E/F) & \lra &(0) \\ && \Art_E \ua \wr && \Gamma_{v,\rho} \da \wr && || && \\ 
(0)& \lra & \A_E^\times/\overline{(E^\times_\infty)^0E^\times} &{\lra} & W_{E/F,\ga}/\overline{(E_\infty^\times)^0E^\times} & \lra & \Gal(E/F) & \lra &(0) \end{array} \]
such that 
\begin{enumerate}
\item $\Gamma_{v,\rho}$ lifts to an isomorphism of extensions
\[ \tGamma_{v,\rho}: W_{\barF/F}/\overline{[W_{\barF/E},W_{\barF/E}]} \liso W_{E/F,\ga} \]
whose composition with the natural map
\[ \theta_\rho: W_{\barF_v/F_v}/\overline{[W_{\barFv/\rho(E)F_v},W_{\barFv/\rho(E)F_v}]} \lra W_{\barF/F}/\overline{[W_{\barF/E},W_{\barF/E}]} \]
is equal to the composition of a canonical map
\[ \iota_{w(\rho)}^\ga: W_{E_{w(\rho)}/F_v,\ga} \lra W_{E/F,\ga} \]
with some isomorphism of extensions
\[ \tTheta: W_{\barF_v/F_v}/\overline{[W_{\barFv/\rho(E)F_v},W_{\barFv/\rho(E)F_v}]} \liso W_{E_{w(\rho)}/F_v,\ga}; \]
\item and if $\sigma \in \Gal(E^\ab/F)$ then $[\Gamma_{v,\rho\sigma}]$ is determined in an explicit way by $[\Gamma_{v,\rho}]$ and $\sigma$.

\end{enumerate}

We will denote by $\cH(E/F)^+$ the set pairs $(\ga,\{ [\Gamma_{v,\rho}] \})$, where $\{ [\Gamma_{v,\rho}]\}$ is complete rigidification data for $\ga$. The action of $T_{2,E}(\A_E)$ on $\cH(E/F)$ lifts to a transitive action on $\cH(E/F)^+$.  One consequence of the choice of $\ga^+ \in \cH(E/F)^+$ is that if $T/F$ is a torus split by $E$, if $\mu \in X_*(T)(\barFv)$ and if $\tau \in \Aut(\barF_v/F)$, then we can associate an important element
\[ \barb_{\ga^+,v,\mu,\tau} \in T(\A_E)/\overline{T(F)T(F_\infty)^0}T(E)T(E_v). \]
In the case of the Serre torus, these elements are closely related to the structure of the Taniyama group.They are also crucial to the explicit description of the action of Galois on canonical models the Shimura varieties associated to tori.  Hence in \cite{st} they play a key role in fixing the action of $\Aut(\C)$ on Shimura varieties.

We have 
\[ \barb_{{}^t\ga^+,v,\mu,\tau}=\barb_{\ga^+,v,\mu,\tau} \prod_\rho ({}^{\rho^{-1}}\mu) (t_{w(\rho)}/t_{w(\tau\rho)}),\]
where $\rho$ runs over $F$-linear embeddings $E \into \barFv$.

Finally when $D\supset E \supset F$ and $\ga_E^+ \in \cH(E/F)^+$ and $\ga_D^+ \in \cH(D/F)^+$ we need to compare the sets $Z^1_\alg(\cE^?(E/F)_{\ga_E},G(A_E))_\basic$ with the sets $Z^1_\alg(\cE^?(D/F)_{\ga_D},G(A_D))_\basic$; and the elements $\barb_{\ga_E^+,v,\mu,\tau}$ and $\barb_{\ga_D^+,v,\mu,\tau}$. It turns out that $\ga_D^+$ and $\ga_E^+$ can be related by certain elements $t \in T_{2,E}(\A_D)$ and the choice of such an element both gives rise to maps
\[ {\inf}^?_{D/E,t}: Z^1_\alg(\cE^?(E/F)_{\ga_E},G(A_E)) \lra Z^1_\alg(\cE^?(D/F)_{\ga_D},G(A_D)) \]
and to equalities
\[ \barb_{\ga_D^+,v,\mu,\tau}=\barb_{\ga_E^+,v,\mu,\tau} \prod_\rho ({}^{\rho^{-1}}\mu) (t_{w(\rho)}/t_{w(\tau\rho)}),\]
where again $\rho$ runs over $F$-linear embeddings $E \into \barFv$.

This paper contains no major theorem, but we hope that the formalism we develop, which we found rather difficult to uncover, may prove useful in other settings besides \cite{st}.

In section 2 we will give a more complete summary of our results.
In section 3 we will recall algebraic facts, particularly concerning extensions and concerning Weil groups. In section 4 we discuss `algebraic cohomology', axiomatically in a general setting. In section 5 we discuss local Weil groups and local algebraic cohomology. In section 6 we discuss the various extensions of global Galois groups which we will need to consider, before discussing global algebraic cohomology in section 7. In section 8 we discuss how these extensions of section 6 compare to Weil groups, and how this leads to the construction of the important special elements $\barb_{\ga^+,v,\mu,\tau}$. Finally in section 9 we recall Langlands' construction of the Taniyama group and discuss its relationship to our special elements $\barb_{\ga^+,v,\mu,\tau}$. We will also prove some additional results which we will need in \cite{st}.

\newpage

\section{Summary of results}

Because our paper contains no `main theorem' but rather a framework, which we hope will be useful, we begin by summarizing the results, which we will establish.

\subsection{Some extensions of global Galois groups}

If $F$ is an algebraic extension of $\Q$, we will write $V_F$ for the set of its places. We will also write $F^\ab$ for the maximal abelian extension of $F$ inside a fixed algebraic closure $\barF$ of $F$.

If $E/F$ is a finite Galois extension of number fields we will write $T_{2,E}$ (resp. $T_{3,E}$) for the protorus over $F$ with cocharacter group $\Z[V_E]$ (resp. $\Z[V_E]_0$) with its natural action of $\Gal(E/F)$. Here $\Z[V_E]$ denotes the free abelian group with basis the set, $V_E$, of places of $E$, and $\Z[V_E]_0$ denotes the subgroup of elements $\sum_{w \in V_E} m_w w$ where $\sum_w m_w=0$. Thus there is a natural short exact sequence 
\[ (0) \lra \G_m \lra T_{2,E} \lra T_{3,E} \lra (0). \]
We will denote by $\pi_w:T_{2,E} \ra \G_m$ the character corresponding to $w \in V_E$.
Now suppose that $D\supset E \supset F$ are number fields Galois over $F$. Then we have
\[ \prod_{w \in V_E} \pi_w: T_{2,E} (\A_D) \liso \prod_{w \in V_E}\A_D^\times, \]
but with Galois action given by 
\[ \sigma ((x_w)_w)= (\sigma x_{\sigma^{-1} w})_w. \]
We have
\[ H^1(\Gal(D/F), T_{2,E}(\A_D)=(0) \]
and
\[ H^1(\Gal(D/F), T_{3,E}(D)=(0). \]

The map
\[ \begin{array}{ccc}  \Z[V_{D}] & \lra & \Z[V_{E}] \\
\sum_u m_u u & \longmapsto & \sum_u  m_u u|_E \end{array} \]
gives rise to a commutative diagram
\[ \begin{array}{ccccccccc}
(0) & \lra & \G_m & \lra & T_{2,E} & \lra & T_{3,E} & \lra & (0) \\
&& || && \iota^0_{D/E} \da &&\da \iota^0_{D/E} && \\
 (0) & \lra & \G_m & \lra & T_{2,D} & \lra & T_{3,D} & \lra & (0); \end{array} \]
 and the map
 \[ \begin{array}{ccc} \Z[V_{E}] &\lra & \Z[V_{D}]  \\
\sum_w m_w w & \longmapsto & \sum_u [D_u:E_{u|_E}] m_{u|_E} u \end{array} \]
gives rise to a commutative diagram
\[ \begin{array}{ccccccccc} (0) & \lra & \G_m & \lra & T_{2,D} & \lra & T_{3,D} & \lra & (0) \\
&& [D:E]\da && \eta^0_{D/E} \da &&\da \eta^0_{D/E} && \\
(0) & \lra & \G_m & \lra & T_{2,E} & \lra & T_{3,E} & \lra & (0), \end{array} \]
and 
\[ \eta^0_{D/E} \circ \iota^0_{D/E}=[D:E]. \]
Note that
\[ \iota^{0,*}_{D/E}: \Z[V_{D}]_{\Gal(D/E)} \liso \Z[V_{E}]. \]
In fact we also have
\[ \iota^{0,*}_{D/E}: \Z[V_{D}]_{0,\Gal(D/E)} \liso \Z[V_{E}]_0. \]

Continue to suppose that $D \supset E \supset F$ are finite Galois extensions of a number field $F$. We will set
\[ \cE^\loc(E/F)_D^0 = \prod_{w \in V_E} D_w^\times \subset T_{2,E}(\A_D) \]
(with $D_w^\times$ thought of inside the $w$-copy of $\A_D^\times$)
and
\[ \cE^\glob(E/F)_D^0 = \{ (x_w) \in T_{2,E}(\A_D):\,\, x_w \bmod D^\times\,\, {\rm is\,\, independent\,\, of}\,\, w\} \subset T_{2,E}(\A_D). \]
These are preserved by $\Gal(D/F)$. 
There are exact sequences 
 \begin{equation}\label{ses1} (0) \lra \A_D^\times \lra \cE^\glob(E/F)_D^0 \lra T_{3,E}(D) \lra (0), \end{equation}
and
 \begin{equation}\label{ses2}  (0) \lra \cE^\glob(E/F)_D^0 \lra T_{2,E}(\A_D) \times \A_D^\times/D^\times \lra T_{2,E}(\A_D)/T_{2,E}(D) \lra (0) \end{equation}
 and
  \begin{equation}\label{ses3}  (0) \lra T_{2,E}(D) \lra \cE^\glob(E/F)_D^0 \lra  \A_D^\times/D^\times  \lra (0). \end{equation}
 We have
 \[ H^1(\Gal(D/F),\cE^\loc(E/F)_D^0)=(0) \]
 and
  \[ H^1(\Gal(D/F),\cE^\glob(E/F)_D^0)=(0). \]
 We will be most interested in the case $D=E$, in which case we will drop it from the notation, writing for instance $\cE^\loc(E/F)^0=\cE^\loc(E/F)^0_E$ and $\cE^\glob(E/F)^0=\cE^\glob(E/F)^0_E$. 
 
 If $C \supset D \supset E \supset F$ are finite Galois extensions of a number field $F$, then
  \[ \eta_{D/E}^0: \cE^\loc(D/F)_C^0 \lra \cE^\loc(E/F)_C^0 \]
  and
   \[ \eta_{D/E}^0: \cE^\glob(D/F)_C^0 \lra \cE^\glob(E/F)_C^0 .\]
   The diagrams
 \[ \begin{array}{ccc} \cE^\loc(D/F)_C^0 & \stackrel{\eta_{D/E}^0}{\lra}& \cE^\loc(E/F)_C^0 \\
 \da &&\da \\ \prod_{u|w} C_u^\times & \lra & C_w^\times \\
 (x_u)_{u|w} & \longmapsto & \prod_{u|w} x_u^{[D_u:F_w]} \end{array} \]
 and
  \[ \begin{array}{ccc} \cE^\glob(D/F)_C^0 & \stackrel{\eta_{D/E}^0}{\lra}& \cE^\glob(E/F)_C^0 \\
 \da &&\da \\ \A_C^\times/C^\times & \stackrel{[D:E]}{\lra} & \A_C^\times/C^\times \end{array} \]
 commute. 
 
 In section \ref{somee} we will define abelian groups $\cZ(E/F)_D \supset \cB(E/F)_D$ with compatible actions of $T_{2,E}(\A_D)$ such that $T_{2,E}(\A_D)$ acts transitively on the quotient $\cH(E/F)_D=\cZ(E/F)_D/\cB(E/F)_D$. The stabilizer in $T_{2,E}(\A_D)$ of any element of $\cH(E/F)_D$ is $\cE^\loc(E/F)_D^0\cE^\glob(E/F)_D^0T_{2,E}(\A_F)$. To an element $\ga \in \cH(E/F)_D$ we associate (uniquely up to unique isomorphism):
 \begin{enumerate}
 \item Extensions
 \[ (0) \lra \cE^\loc(E/F)_D^0 \lra \cE^\loc(E/F)_{D,\ga} \lra \Gal(D/F) \lra (0) \]
 and
 \[ (0) \lra \cE^\glob(E/F)_D^0 \lra \cE^\glob(E/F)_{D,\ga} \lra \Gal(D/F) \lra (0). \]
 
 \item Writing $\cE_2(E/F)_{D,\ga}$ for the pushout of $\cE^\loc(E/F)_{D,\ga}$ along $\cE^\loc(E/F)_{D}^0 \into T_{2,E}(\A_D)$, a canonical map of extensions
 \[ \loc_\ga:  \cE^\glob(E/F)_{D,\ga}  \lra  \cE_2(E/F)_{D,\ga}. \] 

\item An extension
\[ (0) \lra T_{3,E}(D) \lra \cE_3(E/F)_{D,\ga} \lra \Gal(D/F) \lra (0) \]
defined as the pushout of $\cE^\glob(E/F)_{D,\ga}$ along $\cE^\glob(E/F)_D^0 \ra T_{3,E}(D)$. 

\item An extension
\[ (0) \lra \A_D^\times/D^\times \lra W_{E/F, D,\ga} \lra \Gal(D/F) \lra (0) \]
defined as the pushout of $\cE^\glob(E/F)_{D,\ga}$ along $\cE^\glob(E/F)_D^0 \ra \A_D^\times/D^\times$. The extension $W_{E/F,D,\ga}$ is isomorphic to the pushout $W_{E^\ab/F,D}$ of the pullback $W_{E^\ab/F}|_{\Gal(D/F)}$ of Weil group  $W_{E^\ab/F}$ from $\Gal(E/F)$ to $\Gal(D/F)$ along $\A_E^\times/E^\times \into \A_D^\times/D^\times$, i.e.
\[ \cE^\glob(E/F)_{D,\ga}/T_{2,E}(D) \cong W_{E^\ab/F,D} = (\A_D^\times/D^\times \rtimes (W_{E^\ab/F} \times_{\Gal(E/F)}\Gal(D/F)))/(\A_E^\times/E^\times). \]
This isomorphism is {\em not} canonical: it is only defined up to composition with conjugation by an element of $\A_D^\times/D^\times$. (The global Weil group
\[ W_{E^\ab/F}=W_{\barF/F}/\overline{[W_{\barF/E},W_{\barF/E}]},\]
is defined up to an isomorphism that is only unique up to composition with an element of $\overline{E^\times (E_\infty^\times)^0}/E^\times$.)

\item If $w|v$ are places of $E$ and $F$, an extension
\[ (0) \lra D_w^\times \lra W_{E_w/F_v, D,\ga} \lra \Gal(D/F)_w \lra (0) \]
defined as the pushout of $\cE^\loc(E/F)_{D,\ga}|_{\Gal(D/F)_w}$ along $\cE^\loc(E/F)_D^0 \ra D_w^\times$. (Here $\Gal(D/F)_w$ denotes the stabilizer of $w$ in $\Gal(D/F)$.) If $\rho:E \into \barF_v$ gives rise to $w$ then there is an isomorphism of extensions
\[ W_{E_w/F_v, D,\ga} \cong W_{(\rho(E)F_v)^\ab/F_v, \rho, D}, \]
where we define 
\[ W_{(\rho(E)F_v)^\ab/F_v,\rho,D} = (D_{w}^\times \rtimes (W_{(\rho(E)F_v)^\ab/F_v}\times_{\Gal(E/F)_{w}} \Gal(D/F)_{w}))/E_{w}^\times \]
where:
\begin{itemize}
\item $W_{(\rho(E)F_v)^\ab/F_v}$ denotes the local Weil group.
\item $W_{(\rho(E)F_v)^\ab/F_v} \onto \Gal(E/F)_{w}$ is the composition of $W_{(\rho(E)F_v)^\ab/F_v} \onto \Gal((\rho(E)F_v)/F_v)$ with the inverse of the isomorphism $\Gal(E/F)_{w} \iso \Gal((\rho(E)F_v)/F_v)$ induced by $\rho$.
\item $W_{(\rho(E)F_v)^\ab/F_v}\times_{\Gal(E/F)_{w}} \Gal(D/F)_{w}$ acts on $D_{w}^\times$ via its projection to $\Gal(D/F)_{w}$.
\item The map $E_{w}^\times \ra D_{w}^\times \rtimes (W_{(\rho(E)F_v)^\ab/F_v}\times_{\Gal(E/F)_{w}} \Gal(D/F)_{w})$ sends $a$ to $(a^{-1},(r_{\rho(E)F_v}(\rho(a)),1))$.
\end{itemize}
The isomorphism $W_{E_w/F_v, D,\ga} \cong W_{(\rho(E)F_v)^\ab/F_v, \rho, D}$ is {\em not} canonical, but only defined up to composition with conjugation by an element of $D_w^\times$. (In this case the group $W_{(\rho(E)F_v)^\ab/F_v,\rho,D}$ is defined up to unique isomorphism. If $D=E$ then $W_{(\rho(E)F_v)^\ab/F_v,\rho,D}$ is simply the local Weil group $W_{(\rho(E)F_v)^\ab/F_v}$.)

\item If $u|w|v$ are places of $D$, $E$ and $F$ an extension
\[ (0) \lra D_u^\times \lra W_{E_w/F_v,D_u,\ga} \lra \Gal(D/F)_u \lra (0) \]
defined as the pushout of $W_{E_w/F_v,D,\ga}|_{\Gal(D/F)_u}$ along $D_w^\times \onto D_u^\times$. The extension $W_{E_w/F_v,D_u,\ga}$ is isomorphic to the pushout $W_{E_u^\ab/F_v,D_w}$ of the pullback $W_{E_w^\ab/F_v}|_{\Gal(D_u/F_v)}$ of local Weil group  $W_{E_w^\ab/F_v}$ from $\Gal(E_w/F_v)$ to $\Gal(D_u/F_v)$ along $E_w^\times \into D_u^\times$, i.e.
\[ W_{E_w/F_v,D_u,\ga} \cong W_{E_w^\ab/F_v,D_u}=(D_u^\times \rtimes (W_{E_w^\ab/F_v} \times_{\Gal(E_w/F_v)} \Gal(D_u/F_v)))/E_w^\times .\]
This isomorphism is {\em not} canonical: it is only defined up to composition with conjugation by an element of $D_u^\times$.

\item If $w|v$ are places of $E$ and $F$ a map of extensions
\[ \iota^\ga_w: W_{E_w/F_v, D,\ga} \into W_{E/F,D,\ga} \]
compatible with $D_w^\times \into \A_D^\times/D^\times$. 

If  $\rho:E^\ab \into \barFv$, then there is a map of extensions (the `decomposition group') 
\[ \theta_\rho: W_{(\rho(E)F_v)^\ab/F_v} \into W_{E^\ab/F}, \]
which is well defined up to composition with conjugation by an element of $\overline{E^\times (E_\infty^\times)^0}/E^\times$. It induces a map of extensions
\[ \theta_\rho: W_{(\rho(E)F_v)^\ab/F_v,\rho,D} \into W_{E^\ab/F,D} \]
sending $[(a,(\sigma,\tau))]$ to $[(a,(\theta_\rho(\sigma),\tau))]$. Then $\theta_\rho$ and $\iota^\ga_w$ will differ (after making the above identifications) by composition with conjugation by an element of $\A_D^\times/D^\times$ (of course depending on $\rho$).

\end{enumerate}
(Thus, in the case $D=E$, the choice of $\ga \in \cH(E/F)$ inter alia gives rise to a preferred decomposition group in $\Gal(E^\ab/F)$ above each place $w$ of $E$. We think of the choice of $\ga$ as being analogous to the choice of such decomposition groups.) 

Diagramatically we have (for $w\in V_E^S$):

 \begin{tikzcd}  \cE_3(E/F)_{D,\ga} & \arrow[l, two heads]   \cE^\glob(E/F)_{D,\ga}  \arrow[d,"\loc_\ga" hook] \arrow[r, two heads] & W_{E/F,D,\ga}   \\ 
 \cE^\loc(E/F)_{D,\ga} \arrow[r, hook]& \cE_2(E/F)_{D,\ga}  & \\   
 \cE^\loc(E/F)_{D,\ga}|_{\Gal(D/F)_w}  \arrow[u, hook ]   \arrow[r, two heads]&
W_{E_w/F_v,D,\ga} \arrow[uur, bend right, "\iota^\ga_w" hook]  \\
 \cE^\loc(E/F)_{D,\ga}|_{\Gal(D/F)_u}  \arrow[u, hook ]   \arrow[r, two heads]& W_{E_w/F_v,D,\ga}|_{\Gal(D/F)_u} \arrow[u,hook] \arrow[r,two heads] & W_{E_w/F_v,D_u}  \end{tikzcd} 

Our principal interest will be in $\cE^\glob(E/F)_{\ga}$ and $W_{E_w/F_v,\ga}$, which is significantly simpler. Here the diagram becomes

\begin{tikzcd}  \cE_3(E/F)_{\ga} & \arrow[l, two heads]   \cE^\glob(E/F)_{\ga}  \arrow[d,"\loc_\ga" hook] \arrow[r, two heads] & W_{E/F,\ga}   \\ 
 \cE^\loc(E/F)_{\ga} \arrow[r, hook]& \cE_2(E/F)_{\ga}  & \\   
 \cE^\loc(E/F)_{\ga}|_{\Gal(E/F)_w}  \arrow[u, hook ]   \arrow[r, two heads]&
W_{E_w/F_v,\ga} \arrow[uur, bend right, "\iota^\ga_w" hook]  \end{tikzcd} 

\noindent with $W_{E/F,\ga}$ non-canonically isomorphic to $W_{E^\ab/F}$ and $W_{E_w/F_v,\ga}$ non-canonically isomorphic to $W_{E_w^\ab/F_v}$.
However the more general extensions with a $D$ will be required to compare the theory for different extensions $E/F$ and $D/F$. 

The choice of a `cocycle' $\balpha \in \ga$ gives rise to distinguished lifts $e^\glob_\balpha(\sigma) \in \cE^\glob(E/F)_{D,\ga}$ and $e^\loc_\balpha(\sigma) \in \cE^\loc(E/F)_{D,\ga}$ of $\sigma \in \Gal(D/F)$.

If $t \in T_{2,E}(\A_D)$ there are canonical isomorphisms 
\[ \gz_t: \cE^?(E/F)_{D,\ga} \liso \cE^?(E/F)_{D,{}^t\ga} \]
for each of the extensions considered above. They commute with all the arrows in the above diagram, except for the arrows that go between the first row and one of the other rows. For these we have
\[ \gz_t \circ \loc_\ga=\conju_t \circ \loc_{{}^t\ga} \circ \gz_t \]
and
\[ \conju_{t_w} \circ \gz_t \circ \iota_w^\ga =\iota_w^{{}^t\ga} \circ \gz_t. \]

The stabilizer of any element of $\cH(E/F)_D$ is $\cE^\loc(E/F)_D^0\cE^\glob(E/F)^0_DT_{2,E}(\A_F)$. Thus the choice of $t$ taking one element of $\cH(E/F)_D$ to another is not unique.

If $B \supset C \supset D \supset E \supset F$ are finite Galois extensions of a number field $F$, then there are maps
\[ \inf_{C/D}: \cZ(E/F)_D \lra \cZ(E/F)_C \]
and
\[ \eta_{D/E,*}: \cZ(D/F)_C \lra \cZ(E/F)_C, \]
which induce maps
\[ \inf_{C/D}: \cH(E/F)_D \lra \cH(E/F)_C \]
and
\[ \eta_{D/E,*}: \cH(D/F)_C \lra \cH(E/F)_C. \]
They satisfy:
\begin{itemize}
\item $\inf_{B/C} \circ \inf_{C/D} = \inf_{B/D}$.
\item $\eta_{D/E,*} \circ \eta_{C/D,*} = \eta_{C/E,*}$.
\item $\inf_{B/C}\circ \eta_{D/E,*}=\eta_{D/E,*} \circ \inf_{B/C}: \cH(D/F)_C \lra \cH(E/F)_B$.
\item $\inf_{C/D}{}^t\ga={}^t \inf_{C/D} \ga$.
\item $\eta_{D/E,*}{}^t \ga = {}^{\eta_{D/E}(t)} \eta_{D/E,*}(\ga)$.
\item There is a canonical identification of $\cE^?(E/F)_{C,\inf_{C/D}\ga}$ with the pushout of 
\[ \cE^?(E/F)_{D,\ga}|_{\Gal(C/D)}= \cE^?(E/F)_{D,\ga} \times_{\Gal(D/F)} \Gal(C/F) \]
 along
 \[ \cE^?(E/F)_{D}^0 \lra \cE^?(E/F)_{C}^0 \]
 for each of the extensions considered above. (In the case $ W_{E_w/F_v,D,\ga}$ we use $W_{E_w/F_v,D,\ga}|_{\Gal(C/F)_w}$.)
 These identifications commute with the maps $\loc_\ga$, $\iota_w^\ga$ and $\gz_t$. 
 
\item There is a canonical identification of 
 $\cE^?(E/F)_{C,\eta_{D/E,*}\ga}$ with the pushout of $\cE^?(D/F)_{C,\ga}$ along
 \[ \eta_{D/E}:\cE^?(D/F)_{C}^0 \lra \cE^?(E/F)_{C}^0 \]
 for each of the extensions considered above, except for the case $W_{D_u/F_v,C,\ga}$. These identifications commute with the maps $\loc_\ga$ and $\gz_t$. The case of $W_{D_u/F_v,C,\ga}$ and $\iota^\ga_u$ are more complicated and described in lemma \ref{iotamixed} in the case $C=D$.
\end{itemize}

\subsection{Algebraic cocycles}

 We define the following pointed sets of algebraic cocycles:
\begin{enumerate}
\item If $G/F$ is a linear algebraic group, we define $Z^1_\alg(\cE_3(E/F)_{D,\ga},G(D))$ to be the set of 1-cocycles $\phi: \cE_3(E/F)_{D,\ga} \ra G(D)$ such that there is a, necessarily unique, algebraic homomorphism $\bnu_\phi: T_{3,E} \ra G$ over $D$ with $\phi|_{T_{3,E}(D)}=\bnu_\phi$. We also define $Z^1_\alg(\cE_3(E/F)_{D,\ga},G(D))_\basic$ to be the subset on which $\bnu_\phi$ factors through the centre $Z(G)$ of $G$, in which case $\bnu_\phi$ is defined over $F$. Both these pointed sets of cocyles have natural actions of $G(D)$ via the usual coboundary map. 

\item If $G/\A_F$ is a linear algebraic group, we define $Z^1_\alg(\cE^\glob(E/F)_{D,\ga},G(\A_D))$ to be the set of 1-cocycles $\phi: \cE^\glob(E/F)_{D,\ga} \ra G(\A_D)$ such that there is a, necessarily unique, algebraic homomorphism $\bnu_\phi: T_{2,E} \ra G$ over $\A_D$, which is $G(\A_D)$-conjugate to one defined over $D$, with $\phi|_{\cE^\glob(E/F)_D^0}=\bnu_\phi$. We also define $Z^1_\alg(\cE^\glob(E/F)_{D,\ga},G(\A_D))_\basic$ to be the subset on which $\bnu_\phi$ factors through $Z(G)$, in which case $\bnu_\phi$ is defined over $F$. Both these pointed sets of cocyles have natural actions of $G(\A_D)$ via the usual coboundary map.

\item If $G/\A_F$ is a linear algebraic group, we define $Z^1_\alg(\cE_2(E/F)_{D,\ga},G(\A_D))$ to be the set of 1-cocycles $\phi: \cE_2(E/F)_{D,\ga} \ra G(\A_D)$ such that there is a, necessarily unique, algebraic homomorphism $\bnu_\phi: T_{2,E} \ra G$ over $\A_D$, which is $G(\A_D)$-conjugate to one defined over $D$, with $\phi|_{T_{2,E}(\A_D)}=\bnu_\phi$. We also define $Z^1_\alg(\cE_2(E/F)_{D,\ga},G(\A_D))_\basic$ to be the subset on which $\bnu_\phi$ factors through $Z(G)$, in which case $\bnu_\phi$ is defined over $F$. Both these pointed sets of cocyles have natural actions of $G(\A_D)$ via the usual coboundary map.

\item If $G/F_v$ is a linear algebraic group, we define $Z^1_\alg(W_{E_w/F_v, D_u,\ga},G(D_u))$ to be the set of 1-cocycles $\phi: W_{E_w/F_v, D_u,\ga} \ra G(D_u)$ such that there is a necessarily unique, algebraic homomorphisms $\bnu_{\phi}: \G_m \ra G$ over $D_u$ with $\phi|_{D_u^\times}=\bnu_{\phi}$. We also define $Z^1_\alg(W_{E_w/F_v, D_u,\ga},G(D_u))_\basic$ to be the subset on which each $\bnu_{\phi}$ factors through $Z(G)$, in which case $\bnu_{\phi}$ is defined over $F_v$. Both these pointed sets of cocyles have natural actions of $G(D_u)$ via the usual coboundary map.
\end{enumerate}
In each case we will write
\[ H^1_\alg(\cE^?(E/F)_{D,\ga},G(A_D))=G(A_D)\backslash Z^1_\alg(\cE^?(E/F)_{D,\ga},G(A_D)) \]
and
\[ H^1_\alg(\cE^?(E/F)_{D,\ga},G(A_D))_\basic =G(A_D)\backslash Z^1_\alg(\cE^?(E/F)_{D,\ga},G(A_D))_\basic, \]
for each of the groups $\cE^?(E/F)_{D,\ga}$ discussed above, with $A_D$ equal to $D$, $\A_D$, $\A_D$, $D_u$ respectively.

There are equivariant maps natural maps

{\footnotesize \begin{tikzcd} & Z^1_\alg(\cE_3(E/F)_{D,\ga},G(D)) \arrow[d,   dashed, "\loc_\ga"] \arrow[r] & Z^1_\alg(  \cE^\glob(E/F)_{D,\ga} ,G(\A_D))    \\ 
Z^1_\alg(\cE^\loc(E/F)_{D,\ga}|_{\Gal(D/F)_u},G(D_u)) &   Z^1_\alg(\cE^\loc(E/F)_{D,\ga},G(\A_D)) \arrow[l,"\res^S"] \arrow[d,dashed, "\res_u"]& \arrow[l] Z^1_\alg(\cE_2(E/F)_{D,\ga},G(\A_D))  \arrow[u, "\sim" ] \\   &  \arrow[lu,"\sim"] Z^1_\alg(W_{E_w/F_v,D_u,\ga},G(D_u)) &  \end{tikzcd} }
which preserve the basic subsets, are functorial in $G$,  and pass to cohomology.
This induces maps
\[ Z^1_\alg(\cE^\loc(E/F)_{D,\ga}^{V_F-\{v\}},G(D_v)) \lra Z^1_\alg(W_{E_w/F_v,D_u,\ga},G(D_u)) \]
taking basic elements to basic elements (where $u|w|v$ are any choices of places of $D \supset E \supset F$ above $v$),
and isomorphisms
\[ Z^1_\alg(\cE^\loc(E/F)_{D,\ga},G(\A_D)) \liso {\prod_{v \in V_F}}' Z^1_\alg(\cE^\loc(E/F)_{D,\ga}^{V_F-\{v\}},G(D_v)) \]
and
\[ Z^1_\alg(\cE^\loc(E/F)_{D,\ga},G(\A_D))_\basic \liso {\prod_{v \in V_F}}' Z^1_\alg(\cE^\loc(E/F)_{D,\ga}^{V_F-\{v\}},G(D_v))_\basic \]
where the products are restricted with respect to the subsets $Z^1(\Gal(D/F),G(\cO_{D,v}))$. This in turn induces isomorphisms
\[ H^1_\alg(\cE^\loc(E/F)_{D,\ga}^{V_F-\{v\}},G(D_v)) \liso H^1_\alg(W_{E_w/F_v,D_u,\ga},G(D_u)) \]
and 
\[ H^1_\alg(\cE^\loc(E/F)_{D,\ga}^{V_F-\{v\}},G(D_v))_\basic \liso H^1_\alg(W_{E_w/F_v,D_u,\ga},G(D_u))_\basic \]
and
\[ \begin{array}{rcl} H^1_\alg(\cE^\loc(E/F)_{D,\ga},G(\A_D)) &\liso& {\prod_{v \in V_F}}' H^1_\alg(\cE^\loc(E/F)_{D,\ga}^{V_F-\{v\}},G(D_v)) \\ &\liso& {\prod_{v \in V_F}}'H^1_\alg(W_{E_w/F_v,D_u,\ga},G(D_u))
\end{array}\]
and
\[ \begin{array}{rcl} H^1_\alg(\cE^\loc(E/F)_{D,\ga},G(\A_D))_\basic &\liso& {\prod_{v \in V_F}}' H^1_\alg(\cE^\loc(E/F)_{D,\ga}^{V_F-\{v\}},G(D_v))_\basic \\ &\liso& {\prod_{v \in V_F}}'H^1_\alg(W_{E_w/F_v,D_u,\ga},G(D_u))_\basic
\end{array}\]
where the products are restricted with respect to the subsets $H^1(\Gal(D/F),G(\cO_{D,v}))$ (resp. $H^1(\Gal(D_u/F_v),G(\cO_{D,u}))$). If $G$ is a connected linear algebraic group, then we even get
\[ \begin{array}{rcl} H^1_\alg(\cE^\loc(E/F)_{D,\ga},G(\A_D))_\basic &\liso& {\bigoplus_{v \in V_F}}' H^1_\alg(\cE^\loc(E/F)_{D,\ga}^{V_F-\{v\}},G(D_v))_\basic \\ &\liso& {\bigoplus_{v \in V_F}}'H^1_\alg(W_{E_w/F_v,D_u,\ga},G(D_u))_\basic
\end{array}\]
and
\[ \begin{array}{rcl} H^1_\alg(\cE^\loc(E/F)_{D,\ga},G(\A_D))_\basic &\liso& {\bigoplus_{v \in V_F}}' H^1_\alg(\cE^\loc(E/F)_{D,\ga}^{V_F-\{v\}},G(D_v))_\basic \\ &\liso& {\bigoplus_{v \in V_F}}'H^1_\alg(W_{E_w/F_v,D_u,\ga},G(D_u))_\basic
\end{array}\]

If $E_v^0/F_v$ is a finite Galois extension and $G/F_v$ is an algebraic group, we define $Z^1_\alg(W_{(E_v^0)^\ab/F_v},G(E_v^0))$ (resp. $Z^1_\alg(W_{(E_v^0)^\ab/F_v},G(E_v^0))_\basic$) to be the subset of $Z^1(W_{(E_v^0)^\ab/F_v},G(E_v^0))$ consisting of cocycles $\phi$ whose restriction to $W_{(E_v^0)^\ab/E_v^0}$ are of the form $\bnu_\phi \circ \Art_{E_v^0}^{-1}$ for some $\bnu_\phi \in X_*(G)(E_v^0)$ (resp. $X_*(Z(G))(F_v)$). The pointed sets $Z^1_\alg(W_{(E_v^0)^\ab/F_v},G(E_v^0))$ and $Z^1_\alg(W_{(E_v^0)^\ab/F_v},G(E_v^0))_\basic$ are preserved by the coboundary action of $G(E_v^0)$ and we denote the quotients $H^1_\alg(W_{(E_v^0)^\ab/F_v},G(E_v^0))$ and $H^1_\alg(W_{(E_v^0)^\ab/F_v},G(E_v^0))_\basic$ respectively. If $E/\Q$ is a finite Galois extension and $w|v$ is a place of $E$ such that $E_w \cong E_v^0$ over $F_v$, then the choice of an isomorphism of extensions $W_{(E_v^0)^\ab/F_v} \cong W_{E_w/F_v,\ga}$ gives rise to bijections 
\[ Z^1_\alg(W_{(E_v^0)^\ab/F_v},G(E_v^0)) \cong Z^1_\alg(W_{E_w/F_v,\ga},G(E_w)) \]
and
\[ H^1_\alg(W_{(E_v^0)^\ab/F_v},G(E_v^0)) \cong H^1_\alg(W_{E_w/F_v},G(E_w)) \]
identifying basic subsets. The latter is independent of the choices of isomorphisms $E_v^0\cong E_w$ and $W_{(E_v^0)^\ab/F_v} \cong W_{E_w/F_v,\ga}$. The composite of this map with $\res_w$ gives a map
\[ \res_{E_v^0/\Q_v}:  H^1_\alg(\cE^\loc(E/F),G(\A_{E}))\lra H^1_\alg(W_{(E_v^0)^\ab/F_v},G(E_v^0)) \]
which takes basic elements to basic elements and 
which is independent of all choices, including the choice of $w$.

If $F$ is a number field and $E/F$ is Galois will write $\ker^1(\Gal(E/F),G(E))$ for the kernel
\[  \ker (H^1(\Gal(E/F),G(E)) \ra \prod_{v\in V_F} H^1(\Gal(E_{w}/F_v),G(E_{w}))).  \]
We will sometimes write $\ker^1(F,G)$ for $\ker^1(\Gal(\barF/F),G(\barF))$. If $G$ is reductive then $\ker^1(\Gal(E/F),G(E))$ is finite. It vanishes if $G$ is semi-simple and either adjoint or simply connected. 
If now $D \supset E \supset F$ are finite Galois extensions of a number field $F$, then the kernel of 
\[ \loc : H^1_\alg(\cE_3(E/F)_D,G(D)) \lra H^1_\alg(\cE^\loc(E/F)_D,G(\A_D) \]
is exactly $\ker^1(\Gal(D/F),G(D))$.

If $G^\ad=G/Z(G)$, then there is a natural map
\[ Z^1_\alg(\cE^?(E/F)_{D,\ga}, G(A_D))_\basic \lra Z^1(\Gal(D/F), G^\ad(A_D)). \]
 (In the fifth case we of course have to replace $\Gal(D/F)$ with $\Gal(D_u/F_v)$.) Thus if $\phi \in Z^1_\alg(\cE^?(E/F)_{D,\ga}, G(A_D))_\basic$ there is a canonically defined inner form ${}^\phi G$ of $G$ over $A_F$, together with an isomorphism $\iota_\phi: G \times A_D \iso {}^\phi G \times A_D$ such that $\sigma \iota_\phi (g)= \iota_\phi ((\ad \phi(\sigma))(\sigma g))$ for all $\sigma \in\cE^?(E/F)_{D,\ga}$ and $g \in G(A_D)$. If $h \in G(A_D)$ then there is a unique isomorphism $\iota_h:{}^\phi G \iso {}^{{}^h\phi}G$ over $F$ such that $\iota_h \circ \iota_\phi=\iota_{{}^h\phi} \circ \conju_h$. If $\psi \in Z^1_\alg(\cE^?(E/F)_{D,\ga}, ({}^\phi G)(A_D))$, then $\psi \phi \in Z^1_\alg(\cE^?(E/F)_{D,\ga}, G(A_D))$ and this gives a bijection of sets
\[ Z^1_\alg(\cE^?(E/F)_{D,\ga}, ({}^\phi G)(A_D)) \liso Z^1_\alg(\cE^?(E/F)_{D,\ga}, G(A_D)), \]
but this map does not preserve neutral elements. This product induces an isomorphism of basic subsets, is functorial in $G$ and commutes with all the maps in the above diagram comparing the different cohomology theories we are considering. The composite 
\[ \iota_{\psi \phi} \circ \iota_{\phi}^{-1} \circ \iota_\psi^{-1}: {}^\psi({}^\phi G) \liso {}^{\psi \phi} G \]
is defined over $F$. We have $({}^g \psi)\phi = {}^g (\psi \phi)$ and so we get a bijection
\[ H^1_\alg(\cE^?(E/F)_{D,\ga}, ({}^\phi G)(A_D)) \liso H^1_\alg(\cE^?(E/F)_{D,\ga}, G(A_D)). \]
Moreover $(\iota_g \circ \psi){}^g \phi={}^g(\psi \phi)$, and if we use $\iota_g$ to identify ${}^\phi G$ and ${}^{{}^g \phi} G$ then the map induced in cohomology by $\phi$ only depends on $[\phi] \in H^1_\alg(\cE^?(E/F)_{D,\ga}, G(A_D))_\basic$.  We also have $\bnu_{\psi\phi}=\bnu_\psi\bnu_\psi$. 

If $t \in T_{2,E}(\A_D)$ then $\gz_t^{-1}$ induces maps
\[ z_t:  Z^1_\alg(\cE^?(E/F)_{D,\ga},G(A_D)) \liso Z^1_\alg(\cE^?(E/F)_{D,{}^t\ga},G(A_D)) \]
which preserve basic subsets, are functorial in $G$, and are equivariant for the $G(A_D)$-action, so that they pass to cohomology. We have $\bnu \circ z_t = \bnu$ and $z_{t_1t_2}=z_{t_1} \circ z_{t_2}$ and $z_t(\psi\phi)=z_t(\psi)z_t(\phi)$. The maps $z_t$ commute with all the arrows in the above commutative diagram, except those between the top row and the second row. In that case one has
\[( \loc_{{}^t\ga} \circ z_t) (\phi)= {}^{\bnu_\phi(t)} (z_t \circ \loc_\ga)(\phi). \]
The maps the $z_t$ induce in cohomology are independent of $t$, and so $H^1_\alg(\cE^?(E/F)_{D,\ga},G(A_D))$ and $H^1_\alg(\cE^?(E/F)_{D,\ga},G(A_D))_\basic$ are canonically independent of $\ga$. Thus we will denote them simply  $H^1_\alg(\cE^?(E/F)_{D},G(A_D))$ and $H^1_\alg(\cE^?(E/F)_{D},G(A_D))_\basic$. 

Now suppose that $C \supset D \supset E \supset F$ are finite Galois extensions of a number field $F$. We obtain maps
\[ \begin{array}{rcl} \inf_{C/D}: Z^1_\alg(\cE^?(E/F)_{D,\ga},G(A_D)) & \lra & Z^1_\alg(\cE^?(E/F)_{D,\ga}|_{\Gal(C/F)},G(A_C)) \\ &\stackrel{\sim}{\lla} & Z^1_\alg(\cE^?(E/F)_{C,\ga},G(A_C)) \end{array} \]
and
\[ \eta_{D/E}^*: Z^1_\alg(\cE^?(E/F)_{C,\eta_{D/E,*}\ga},G(A_C)) \lra Z^1_\alg(\cE^?(D/F)_{C,\ga},G(A_C)) . \]
These maps preserve basic subsets, are functorial in $G$, and are $G(A_D)$-equivariant so that they pass to cohomology. They are compatible with products and all the maps in the above commutative diagram. We have $\bnu \circ \inf_{C/D}=\bnu$ and $(\bnu \circ \eta_{D/E}^*)(\phi)=\bnu_\phi \circ \eta_{D/E}$ and $z_t \circ \inf_{C/D}=\inf_{C/D} \circ z_t$ and $z_t \circ \eta_{D/E}^* =\eta_{D/E}^* \circ z_{\eta_{D/E}(t)}$. 

If $\ga_D \in \cH(D/F)$ and $\ga_E \in \CH(E/F)$ the we can find a $t\in T_{2,E}(\A_D)$ with $\eta_{D/E,*}\ga_D={}^t \inf_{D/E} \ga_E$. Then we set
\[ \begin{array}{rcl}  \inf_{D/E,t}: Z^1_\alg(\cE^?(E/F)_{\ga_E},G(A_E)) & \stackrel{\inf_{D/E}}{\lra} & Z^1_\alg(\cE^?(E/F)_{D,\inf_{D/E}\ga_E},G(A_D)) \\ &\stackrel{z_t}{\lra}&
Z^1_\alg(\cE^?(E/F)_{D,\eta_{D/E,*} \ga_D},G(A_D)) \\ 
&\stackrel{\eta_{D/E}^*}{\lra} & Z^1_\alg(\cE^?(D/F)_{\ga_D},G(A_D)).
\end{array} \]
This map is functorial in $G$, preserves basic subsets, commutes with products, and passes to cohomology, where it is independent of $\ga_D$, $\ga_E$ and $t$ so we will denote it simply $\inf_{D/E}$. These maps are all injective even on the level of cohomology. They commute with all the maps in the above commutative diagram, except those between the top and second rows. We have
\[ \loc_{\ga_D}  (\inf_{D/E,t}(\phi))= {}^{\bnu_\phi(t)} \inf_{D/E,t}(\loc_{\ga_E} \phi). \]
Note that $\inf_{E/E,t}=z_t$
and
$\bnu_{\inf_{D/E,t} \phi}= \bnu_\phi \circ \eta_{D/E}$.
If $a \in T_{2,E}(\A_F)$ and $b \in \cE^\glob(E/F)_D^0$ and $c \in \cE^\glob(E/F)_D^0$, then
\[ \inf_{D/E, abct}(\phi)= {}^{\bnu_\phi(b)} \inf_{D/E,t}(\phi) \] 
if $\phi \in Z^1_\alg(\cE_3(E/F)_{\ga_E},G(E))$ or $Z^1_\alg(\cE^\glob(E/F)_{\ga_E},G(\A_E))$, while
\[ \inf_{D/E, abct}(\phi) = {}^{\bnu_\phi(c^{-1})} \inf_{D/E,t}(\phi) \]
if $\phi$ is in any of the other groups of cocycles $Z^1_\alg(\cE^?(E/F)_{\ga_E},G(A_E))$.
If $C\supset D$ is another finite Galois extension of $F$ and if $\ga_{C}\in \cH(C/F)$ and if $t' \in T_{2,D}(\A_{C})$ with $\eta_{C/D,*} \ga_{C}={}^{t'} \inf_{C/D} \ga_D$, then
\[ \inf_{C/D,t'} \circ \inf_{D/E,t} = \inf_{C/E, t \eta_{D/E}(t')} .\]

Following Kottwitz, we define 
\[ B(F,G) = \lim_{\ra, E} H^1_\alg(\cE_3(E/F),G(E))\]
 and 
 \[ B^\loc(F,G) = \lim_{\ra, E} H^1_\alg(\cE^\loc(E/F),G(\A_E)) \] 
 and
 \[ B(F_v,G)= \lim_{\ra, E_w} H^1_\alg(W_{E_w/F_v},G(E_w)), \]
 and similarly $B(F,G)_\basic$, $B^\loc(F,G)_\basic$ and $B(F_v,G)_\basic$.

 If $T/F$ is a torus split by $E$ and if $\balpha \in \ga \in \cH(E/F)$ we get correstriction maps
\[ \corr_\balpha^\glob: \Z[V_E]_0 \otimes X_*(T)=\Hom(T_{3,E},T) \lra Z^1_\alg(\cE_3(E/F)_\ga,T(E)) \]
and
\[ \corr_\balpha^\loc: \Z[V_E] \otimes X_*(T)=\Hom(T_{2,E},T) \lra Z^1_\alg(\cE^\loc(E/F)_\ga,T(\A_E)), \]
defined by
\[ \corr_\balpha^\glob(\chi)(\sigma) = \prod_{\eta \in \Gal(E/F)} \eta^{-1} \chi(e_\balpha^\glob(\eta)\sigma e_\balpha^\glob(\eta \sigma)^{-1}) \]
and
\[ \corr_\balpha^\loc(\chi)(\sigma) = \prod_{\eta \in \Gal(E/F)} \eta^{-1} \chi(e_\balpha^\loc(\eta)\sigma e_\balpha^\loc(\eta \sigma)^{-1}). \]
They induce maps
\[ \corr^\glob: (\Z[V_E]_0 \otimes X_*(T))_{\Gal(E/F)} \liso H^1_\alg(\cE_3(E/F),T(E)) \]
and
\[ \corr^\loc: (\Z[V_E] \otimes X_*(T))_{\Gal(E/F)} \liso H^1_\alg(\cE^\loc(E/F),T(\A_E)), \]
which are independent of $\balpha$ and $\ga$ and, which Kottwitz showed, are isomorphisms. They satisfy $\loc_\ga \circ \corr^\glob=\corr^\loc$. (However it is {\em not} true that $\loc_\ga \circ \corr^\glob_\balpha=\corr^\loc_\balpha$ on the level of cocycles.)

Suppose that $E$ is totally complex and that $G/F$ is a reductive group which contains a maximal torus (over $F$), which splits over $E$. Write $\Lambda_G$ for the arithmetic fundamental group of $G$. Then Kottwitz also showed that there are unique maps
\[ \kappa: H^1_\alg(\cE_3(E/F),G(E)) \lra (\Z[V_E]_0 \otimes_\Z \Lambda_G)_{\Gal(E/F)} \]
and
\[ \kappa: H^1_\alg(\cE^\loc(E/F),G(\A_{E})) \lra  (\Z[V_{E}] \otimes \Lambda_G)_{\Gal(E/F)}  \]
with the following properties:
\begin{enumerate}
\item They are functorial in $G$.
\item If $G=T$ is a torus then they equal $\corr^{-1}$.

\item $\kappa \circ \inf_{D/E}$ equals the composite of the natural isomorphisms
\[ (\Z[V_D]_0 \otimes_\Z \Lambda_G)_{\Gal(D/F)} \liso (\Z[V_E]_0 \otimes_\Z \Lambda_G)_{\Gal(E/F)}\]
or
\[ (\Z[V_{D}] \otimes \Lambda_G)_{\Gal(D/F)} \liso (\Z[V_{E}] \otimes \Lambda_G)_{\Gal(E/F)}\]
(induced by the maps $\Z[V_D] \ra \Z[V_E]$ sending $u \mapsto u|_E$)
with $\kappa$.  

\item $\kappa \circ \loc$ equals the composition of the natural map $(\Z[V_{E}]_0 \otimes \Lambda_G)_{\Gal(E/F)} \ra (\Z[V_{E}] \otimes \Lambda_G)_{\Gal(E/F)}$ with $\kappa$.
\end{enumerate}
We will denote by $\barkappa$ the composite
\[ H^1_\alg(\cE^\loc(E/F)_S,G(\A_{E})) \stackrel{\kappa}{\lra}  (\Z[V_{E}] \otimes \Lambda_G)_{\Gal(E/F)} \onto \Lambda_{G,\Gal(E/F)} \]
induced by $\sum_w w \otimes x_w \mapsto \sum_w x_w$.
Note that $\barkappa \circ \loc=0$. We warn the reader that our $\barkappa$ has a different meaning from Kottwitz's use of the same symbol.

Kottwitz proves that there is a cartesian square
\[ \begin{array}{ccc} B(F,G)_\basic & \stackrel{\oplus_{v|\infty} \res_v \circ \loc}{\lra} & \prod_{v|\infty} B(F_v,G)_\basic \\ \kappa \da && \prod_{v|\infty} \kappa \da \\
(\Z[V_E]_0 \otimes \Lambda_G)_{\Gal(E/F)} & \lra & \prod_{v|\infty} \Lambda_{G,\Gal(E_\tv/F_v)} \\
\sum_w w \otimes \lambda_w & \longmapsto & (\sum_{\sigma \in \Gal(E_\tv/F_v) \backslash \Gal(E/F)} \sigma \lambda_{\sigma^{-1}\tv})_{v|\infty},\end{array} \]
where $E/F$ is a sufficiently large finite Galois extension and where $\tv|v$ is a place of $E$.
In particular the fibres of $\kappa: B(F,G)_\basic \ra (\Z[V_E]_0 \otimes \Lambda_G)_{\Gal(E/F)}$ are finite.

The map 
\[ \loc: B(F,G)_\basic \lra \ker \barkappa \subset B^\loc(F,G)_\basic \]
is surjective with finite fibres.

If $S$ is a finite set of places of $F$, we will write $B(F,G)_{S,\basic}$ for the inverse image in $B(F,G)_\basic$ under $\kappa$ of the image of $\Z[V_{E,S}]_0 \otimes_\Z \Lambda_G$ in $(\Z[V_E]_0 \otimes_\Z \Lambda_G)_{\Gal(E/F)}$. (This is independent of $E/F$.) Given such a finite set $S$, there is a finite Galois extension $D/F$ such that $B(F,G)_{S,\basic}$ is contained in the image of $H^1_\alg(\cE_3(D/F),G(D))_\basic$.

\subsection{Rigidification}

There is a second set $\cH(E/F)_D^+$ with a transitive action of $T_{2,E}(\A_D)$ which comes with an equivariant surjection 
\[ \cH(E/F)_D^+ \onto \cH(E/F)_D. \]

If $\rho:E \into \barF_v$ is $F$-linear we will write $w(\rho)$ for the place of $E$ above $v$ induced by $\rho$.

An element of $\cH(E/F)_D^+$ is a pair $(\ga,\{ [\Gamma_\rho]\})$, where $\ga \in \cH(E/F)$ and where for each place $v$ of $F$ and each $F$-linear embedding $\rho:E^\ab D \into \barF_v$, the set $[\Gamma_\rho]$ is a $D_{w(\rho)}^\times$-conjugacy class of isomorphisms
\[ \Gamma_\rho: W_{E^\ab/F,D}/\overline{E^\times (E_\infty^\times)^0} \liso W_{E/F,D,\ga}/\overline{E^\times (E_\infty^\times)^0} \]
such that 
\begin{itemize}
\item There is an isomorphism of extensions
\[ \tGamma_\rho:W_{E^\ab/F,D} \liso W_{E/F,D,\ga}, \]
which lifts $\Gamma_\rho$, and an isomorphism of extensions
\[ \Theta:W_{(E F_v)^\ab/F_v,\rho,D} \liso W_{E_{w(\rho)}/F_v,D,\ga} \]
such that 
\[ \iota^\ga_{w(\rho)} \circ \Theta=\tGamma_\rho \circ
\theta_\rho . \]

\item If $\sigma \in \Gal(E^\ab D/F)$, then $[\Gamma_{\rho\sigma}]=[\Gamma_\rho^{\sigma,\balpha}]$ for some (and hence any) $\balpha \in \ga$, where
\[ \Gamma_\rho^{\sigma,\balpha} = \conju_{\pi_{w(\rho \sigma)}(e_\balpha^\loc(\sigma^{-1}) \loc_\ga(e_\balpha^\glob(\sigma^{-1})^{-1}))} \circ \conju_{e_\balpha^\glob(\sigma^{-1})} \circ \Gamma_\rho \circ \conju_\sigma \]
\end{itemize}
We call $\{ [\Gamma_\rho]\}$ {\em complete rigidification data} for the image $\ga$ of $\ga^+$ in $\cH(E/F)_D$. 

If $E$ is totally imaginary, the stabilizer in $T_{2,E}(\A_D)$ of any element of $\cH(E/F)_D^+$ is
\[ \cE^\loc(E/F)_D^0 \cE^\glob(E/F)_D^0 T_{2,F}(\A_F) \prod_{w \in V_E} \overline{E^\times (E_\infty^\times)^0}. \]

If $B \supset C \supset D \supset E \supset F$ are finite Galois extensions of a number field $F$ then there are maps
\[ \inf_{C/D}: \cH(E/F)_D^+ \lra \cH(E/F)_C^+ \]
and
\[ \eta_{D/E,*}: \cH(D/F)_C^+ \lra \cH(E/F)_C^+ \]
compatible with the corresponding maps for without complete rigidification data. Moreover these have the following properties:
\begin{itemize}
\item $\inf_{C/D} {}^t \ga^+ = {}^t \inf_{C/D} \ga^+$.
\item $\eta_{D/E,*} {}^t \ga^+ = {}^{\eta_{D/E}(t)} \eta_{D/E,*} \ga^+$.
\item $\inf_{B/C} \circ \inf_{C/D}=\inf_{B/D}$.
\item $\eta_{D/E,*} \circ \eta_{C/D,*}=\eta_{C/E,*}:\cH(C/F)_B^+ \ra \cH(E/F)_B^+$.
\item $\eta_{D/E,*} \circ \inf_{B/C}=\inf_{B/C} \circ \eta_{D/E,*}: \cH(D/F)_C^+ \ra \cH(E/F)_B^+$.
\item Suppose that $\ga_{C}^+ \in\cH(C/F)^+$ and $\ga_{D}^+ \in\cH(D/F)^+$ and $\ga_{E}^+ \in\cH(E/F)^+$ satisfy $\eta_{C/D,*}\ga_{C}^+={}^{t'} \inf_{C/D} \ga_D^+$ and $\eta_{D/E,*}\ga_{D}^+={}^{t} \inf_{D/E} \ga_E^+$ with $t' \in T_{2,D}(\A_{C})$ and $t \in T_{2,E}(\A_{D})$. Then
\[ \eta_{C/E,*}\ga_{C}^+={}^{t \eta_{D/E}(t')} \inf_{C/E} \ga_E^+. \]
\end{itemize}

Now suppose that $\ga^+ \in \cH(E/F)_D^+$, that $T/F$ is a torus split by $E$, that $v$ is a place of $F$, that $\tau \in \Aut(\barF_v/F_v)$ is a field theoretic automorphism of $\barF_v$ fixing $F$ (but not necessarily continuous), and that $\mu \in X_*(T)(\barF_v)$. Then there is a well defined element
\[ \barb_{\ga^+,v,\mu,\tau} \in T(\A_D)/T(D)T(D_v)\overline{T(F)T(F_\infty)_0} \]
with the following properties:
\begin{itemize}
\item If $D=E$ and $\balpha \in \ga$, the image of $\ga^+$ in $\cH(E/F)$, then there is a lift $b\in T(\A_D)$ of $\barb_{\ga^+,v,\mu,\tau}$ such that
\[ \loc_\ga \corr_\balpha^\glob (({}^{\rho
^{-1}}\mu) \circ (\pi_{w(\rho)}/\pi_{w(\tau \rho)}))={}^b \corr_\balpha^\loc (({}^{\rho
^{-1}}\mu) \circ (\pi_{w(\rho)}/\pi_{w(\tau \rho)})).\]

\item $\barb_{\ga^+,v,\mu,\tau_1\tau_2}=\barb_{\ga^+,v,{}^{\tau_2}\mu,\tau_1}\barb_{\ga^+,v,\mu,\tau_2}$.

\item If $\tau$ is continuous, then $\barb_{\ga^+,v,\mu,\tau}=1$.

\item If if $\rho:E^\ab \into \barF_v$ is $F$-linear, if $\tau$ fixes the image of $E$ in $\barFv$,  and if $a_\rho \in \A_E^\times$ with $\rho \circ \Art_E (a_\rho)=\tau \circ \rho$, then
\[  \barb_{\ga^+,v,\mu,\tau} =  \prod_{\eta \in \Gal(E/F)} \eta ({}^{\rho^{-1}} \mu)(a_\rho)^{-1}.\]
\end{itemize}

There elements enjoy the following functorialities.
\begin{itemize}

\item If $\ga^+\in \cH(E/F)_D^+$ and $\chi:T \ra T'$ is a morphism of algebraic groups over $F$, then $\barb_{\ga^+,v,\chi \circ \mu,\tau}=\chi(\barb_{\ga^+,v,\mu,\tau})$.
 
 \item If $\ga^+ \in \cH(E/F)_D^+$, then $\barb_{\inf_{C/D}\ga^+,v,\mu,\tau} = \barb_{\ga^+,v,\mu,\tau}$.
 
 \item If $\ga^+ \in \cH(D/F)_C^+$, then $\barb_{\eta_{D/E,*} \ga^+,v,\mu,\tau} = \barb_{\ga^+,v,\mu,\tau}$.
 
 \item If $\ga^+ \in \cH(E/F)_D^+$ and $t \in T_{2,E}(\A_D)$, then
\[ \begin{array}{rcl} \barb_{{}^t\ga^+,v,\mu,\tau} &=  &\barb_{\ga^+,v,\mu,\tau}\prod_{\rho}  ({}^{\rho^{-1}}\mu) \circ (\pi_{w(\rho)}/ \pi_{w(\tau \rho)})(t) \\ &=& \barb_{\ga^+,v,\mu,\tau} \prod_\rho ({}^{\rho^{-1}}(\mu/{}^\tau \mu) )(t_{w(\rho)})\end{array}  \]
where $\rho$ runs over $F$-linear embeddings $\rho:E \into \barFv$.

\item Suppose that $\ga_E^+ \in \cH(E/F)^+$ and $\ga_D^+ \in \cH(D/F)^+$. Then we can find a $t \in T_{2,E}(\A_D)$ such that ${}^t\inf_{D/E} \ga_E^+=\eta_{D/E,*} \ga_D^+ \in \cH(E/F)_D^+$, in which case 
\[ \begin{array}{rcl} \barb_{\ga_E^+,v,\mu,\tau} &=& \barb_{\ga_D^+,v,\mu,\tau}  \prod_\rho ({}^{\rho^{-1}}\mu)(t_{w(\tau\rho)}/t_{w( \rho)})) \\ &=& 
\barb_{\ga_D^+,v,\mu,\tau}  \prod_\rho ({}^{\rho^{-1}}({}^\tau \mu/\mu)(t_{w( \rho)})),\end{array} \]
where $\rho$ runs over $F$-linear embeddings $E \into \barFv$.

\end{itemize}

These elements are closely connected to, and generalize the cocycles that define the Taniyama group, which we will now explain.

Write $\C^\alg$ for the algebraic closure of $\Q$ in $\C$. If $E/\Q$ is a finite Galois extension we will write $E \cap \C$ for its image under any embedding into $\C$. Then there is a torus $S_{E,\C}/\Q$, called the Serre torus, with cocharacter group
\[ X^*(S_E)=\{ (\varphi,w)\in \Map(\Gal((E\cap\C)/\Q),\Z) \times \Z: \,\, \varphi(\sigma c\tau)+\varphi(\sigma\tau)=w \,\, \forall \sigma, \tau \in \Gal((E\cap\C)/\Q) \}, \]
with a left action of $\Gal(\C^\alg/\Q)$ given by 
\[ \sigma(\varphi,w)=(\tau \mapsto \varphi(\sigma^{-1} \tau), w). \]
The torus $S_{E,\C}$ has an action of $\Gal(E \cap \C/\Q)$ defined over $\Q$ and characterized by 
\[ X^*([\tau]) (\varphi,w)=(\tau' \mapsto \varphi(\tau'\tau^{-1}), w). \] 
There is also a canonical cocharacter $\mu^\can \in X_*(S_{E,\C})(\C)$ defined by
\[ X^*(\mu^\can)(\varphi,w)=\varphi(1) \in \Z \cong X^*(\G_m). \]
Importantly $S_{E,\C}(\Q)$ is a discrete subgroup of $S_{E,\C}(\A^\infty)$.

Langlands defined a pro-algebraic group $\tS_{E,\C}$, the {\em Taniyama group}, which is an extension of $\Gal(E^\ab \cap \C/\Q)$ by $S_{E,\C}$, where the induced action is via $\tau \mapsto [\tau]$:
\[ (0) \lra S_{E,\C} \lra \tS_{E,\C} \lra \Gal(E^\ab\cap \C/\Q) \lra (0). \]
It comes with a canonical section
\[ \spl:\Gal(E^\ab\cap \C/\Q) \ra \tS_{E,\C}(\A^\infty). \]
It has the properties that $\tS_{E,\C}(E)$ surjects onto $\Gal(E^\ab\cap \C/\Q)$, and that $\tS_{E,\C}|_{\Gal(E^\ab \cap\C/E\cap \C)}$ is abelian. We have the following equality:
\begin{itemize}
\item If $\alpha \in \tS_{E,\C}(E)$ and $\tau \in \Aut(\C/\Q)$ have the same image in $\Gal(E^\ab \cap \C/\Q)$ and if $\ga \in \cH(E/\Q)$, then 
\[ \alpha^{-1} \spl(\tau) \in S_{E,\C}(E) S_{E,\C}(\A^\infty) \subset S_{E,\C}(\A_E^\infty) \]
lifts
\[ \barb_{\ga^+,\infty,\mu^\can,\tau} \in S_{E,\C}(\A_E^\infty)/S_{E,\C}(E).\]
\end{itemize}

\newpage

\section{Algebraic background}\label{ag}

\subsection{Notations}

{\em For simplicity we will assume all fields we consider in this paper will be assumed to be perfect unless we specifically say otherwise.}

If $F$ is a field we will write $\barF$ for an algebraic closure of $F$ and $F^\ab$ for the maximal abelian Galois extension of $F$ in $\barF$. If $E/F$ is a Galois extension and $L/F$ any field extension, then we will write $E \cap L$ (resp. $EL$) for $\rho (E) \cap L$ (resp. $\rho(E)L$) for any $F$-linear embedding $\rho:E \into \barL$. The field $E \cap L$ (resp. $EL$) is a subfield of $L$ (resp. $\barL$) independent of the choice of $\rho$, but the identification of $E \cap L$ with a subfield of $E$ depends on $\rho$. If $L$ is any field of characteristic $0$, we will write $L^\alg$ for the subfield consisting of elements algebraic over $\Q$. If $L/K$ is any extension of fields we will write $\Aut(L/K)$ for the group of field theoretic automorphisms of $L$ which fix $K$ pointwise. If $L$ has characteristic $0$ (resp. $p>0$) will write $\Aut(L)$ for $\Aut(L/\Q)$ (resp. $\Aut(L/\F_p)$). If $L/K$ is Galois we will write $\Gal(L/K)$ for $\Aut(L/K)$. If $E$ is a subfield of $\C^\alg$ then
\[ \C^{\Aut(\C/E)}=E. \]

If $F$ is a local field of characteristic $0$ we will write $\Art_F: F^\times \ra \Gal(F^\ab/F)$ for the Artin map. (Normalized to take uniformizers to geometric Frobenius elements.)

If $F$ is an algebraic extension of $\Q$ we will write $V_F$ for the set of places of $F$ and $\A_F$ for the ring of adeles of $F$. (In the case that $F$ is an infinite extension of $\Q$ then $\A_F = \lim_{\ra E} \A_E$, where $E$ runs over subfields of $F$ finite over $\Q$.) If $v$ is a place of $F$ then by $F_v$ we will mean $\lim_{\ra E} E_v$ as $E$ runs over subextensions of $F/\Q$ which are finite over $\Q$. (So $F_v$ may not be complete, but it is algebraic over $\Q_v$.) If $F$ is a number field will write $\Art_F:\A_F^\times/\overline{F^\times (F_\infty^\times)^0} \iso \Gal(F^\ab/F)$ for the Artin map.

If $E/F$ is an algebraic extension of fields with $F$ a number field and if $S \subset V_F$ we will write $V_{E,S}$ for the set of places of $E$ above a place in $S$, and $\A_{E,S}$ for the ring of adeles of $E$ supported at the primes in $S$. (If $E$ is also a number field then $\A_{E,S}$ is the restricted product $\prod_{w:\,\, w|_F \in S}' E_w^\times$.) Moreover 
$\A_E^S=\A_{E,V_F-S}$.

We will write $\Z[V_{E,S}]$ for the free ablelian group on $V_{E,S}$ and $\Z[V_{E,S}]_0$ for the subabelian group consisting of elements $\sum_w m_w w$ with $\sum_w m_w=0$. If $E/F$ is Galois, both groups have a natural action of $\Gal(E/F)$ via $\sigma \sum_w m_w w= \sum_w m_w (\sigma w)=\sum_w m_{\sigma^{-1}w} w$.

If $F$ is an algebraic extension of $\Q$ and $K$ is a local field and $\rho:F \into \overline{K}$, then we will write $v(\rho)=v(F,\rho)$ or $w(\rho)=w(F,\rho)$ or $u(\rho)=u(F,\rho)$ for the place of $F$ induced by $\rho$. (We will tend to use $v(\rho)$ when the field is denoted $F$, $w(\rho)$ when it is denoted $E$ and $u(\rho)$ otherwise.) If moreover $F/\Q$ is Galois and $\tau \in \Aut(\overline{K})$, then we will write $\tau^\rho$ for the element of $\Gal(F/\Q)$ satisfying $\rho \circ \tau^\rho =  \tau \rho$.

If $E/F$ is a Galois extension with $F$ a number field, and if $v$ is a real place of $F$ we will write $[c_v]$ for the conjugacy class in $\Gal(E/F)$ consisting of complex conjugations at places above $v$. If $F=\Q$ and $v=\infty$ we will simply write $[c]$.

If $G$ is an (algebraic) group then $Z(G)$ will denote its centre and $G^\ad$ will denote $G/Z(G)$. Moreover $G^\der$ will denote its commutator subgroup and $C(G)=G^\ab$ will denote it co-center/abelianization $G/G^\der$. If $H\subset G$ is a subgroup we will write $N_G(H)$ for its normalizer and $Z_G(H)$ its centralizer. If $H$ has finite index, we will also write $\tr_{G/H}:G^\ab \ra H^\ab$ for the transfer map. If $H$ is normal in $G$, then the image of $\tr_{G/H}$ is contained in $(H^\ab)^{G/H}$. If $G$ acts on $X$ we will write $[x]_G$ for the $G$ orbit of $x \in X$ and $Z_G(x)$ or $G_x$ for the stabilizer of $x$ in $G$. If $G$ is an algebraic group acting on a variety $X$ over a field $F$ and $x \in X(F)$, then $[x]_G$ is a variety, and $[x]_G(F) \supset [x]_{G(F)}$, but these two sets may not be equal.

If $G$ is an affine algebraic group over $F$ then there a scheme $X_*(G)$, smooth and separated over $F$, and a homomorphism $\mu^\univ: \G_m \times_F X_*(G) \ra G \times_F X_*(G)$, such that if $S$ is any $F$-scheme and $\mu:\G_{m,S} \ra G_S$ is a homomorphism, then there is a unique morphism $S \ra X_*(G)$ under which $\mu^\univ$ pulls back to $\mu$.  
Moreover 
\[ \begin{array}{ccc}
G \times X_*(G) & \lra & X_*(G) \times X_*(G) \\ (g,\mu) & \longmapsto & (\conju_g \circ \mu , \mu) \end{array} \]
is smooth; and
\[ X_*(G)_{\barF} = \coprod_{[\mu] \in G(\barF)\backslash X_*(G)(\barF)} G / Z_G(\mu) . \]
(See sections 4 and 5 of expos\'{e} XI in \cite{sga3}.) 

{\em We will require all our reductive groups to be geometrically connected, i.e. by the term `reductive group' we will mean what is often referred to as `connected reductive group'.} 
We will write $G^\sco$ for the simply connected semi-simple cover of $G^\der$. If $T$ is a maximal torus of $G$ we will write $T^\ad$ for the image of $T$ in $G^\ad$ (a maximal torus in $G^\ad$) and $T^\der=(G^\der \cap T)$ (which is a maximal torus in $G^\der$, see remark 3.5 of \cite{cbrs}) and $T^\sco$ for the preimage of $T$ in $G^\sco$ (which is a maximal torus in $G^\sco$, see for instance proposition 4.1 of \cite{cbrs}). We have $T=Z_G(T)$. We will also write $W_T$ for the Weyl group $N_G(T)/T$, which we think of a finite algebraic group. It acts faithfully on $T$. We will also write $W_{T,F}$ for $N_G(T)(F)/T(F) \subset W_T(F)$.

We remark that if $T \subset G$ is a maximal torus and $\mu_1,\mu_2 \in X_*(T)$ are conjugate under $G(\barF)$ then they are conjugate under $W_T(\barF)$. (This is probably well known, but as we don't know a reference we will sketch the proof.
Let $H$ denote the connected component of the identity of the centralizer of $\mu_1(\G_m)$ in $G$. It is reductive. (See theorem 2.1 of \cite{cbrs}.) Suppose that $\mu_1=g\mu_2g^{-1}$. Then $\mu_1(\G_m) \subset gTg^{-1}$ so that $T$ and $gTg^{-1}$ are both maximal tori in $H$. Hence we have $gTg^{-1}=hTh^{-1}$ for some $h \in H$. Then $h^{-1}g \in N_G(T)$ and $\mu_1=h^{-1}g \mu_2 g^{-1}h$, as desired.)  

 We will let $\Lambda_G$ denote the algebraic fundamental group of $G$, i.e. $X_*(T)/X_*(T^\sco)$ for any maximal torus $T$ of $G$. Note that the Weyl group $W_T$ acts trivially on $X_*(T)/X_*(T^\sco)$. Any two maximal tori $T$ and $T'$ defined over $F$ are conjugate over the separable closure $\barF$ of $F$ by $g \in G(\barF)$ with $gN_G(T)$ uniquely defined. Then $\conju_g$ induces an isomorphism $X_*(T)/X_*(T^\sco) \iso X_*(T')/X_*(T^{\prime \sco})$. If we alter $g$ by an element $h\in N_G(T)(\barF)$ then this isomorphism changes by an element of $W_T(\barF)$, i.e. is in fact unchanged. Thus $\Lambda_G$ is canonically defined independent of the choice of $T$. In particular it has  a canonical action of $\Gal(\barF/F)$. (If $T'=\conju_g T$ and $\sigma \in \Gal(\barF/F)$, then
$\sigma \circ \conju_g = \conju_g \circ \sigma \circ \conju_{w_\sigma}$ on $X_*(T)$ for some $w_\sigma \in W_T(\barF)$, and so $\sigma \circ \conju_g = \conju_g \circ \sigma$ on $\Lambda_T$.)
If $[\mu]$ is a conjugacy class of cocharacters $\mu:\G_m \ra G$, then $[\mu]$ gives rise to well defined element $\lambda_G([\mu]) \in \Lambda_{G}$. If $\sigma \in \Gal(\barF/F)$ then $\lambda_G({}^\sigma[\mu])={}^\sigma\lambda_G([\mu])$.

If $G/F$ is an algebraic group and $E/F$ is a Galois extension then we will write $H^1(\Gal(E/F),G(E))$ for the first Galois cohomology.  More precisely, by a $1$ cocycle of $\Gal(E/F)$ we will mean a locally constant map $\phi:\Gal(E/F) \ra G(E)$ such that $\phi(\sigma_1\sigma_2)=\phi(\sigma_1){}^{\sigma_1}\phi(\sigma_2)$. We denote the set of $1$ cocycles $Z^1(\Gal(E/F),G(E))$. If $\phi \in Z^1(\Gal(E/F),G(E))$ and $g \in G(E)$ we define ${}^g\phi \in Z^1(\Gal(E/F),G(E))$ by $({}^g\phi)(\sigma)=g \phi(\sigma){}^\sigma g^{-1}$. 
We call two cocycles $\phi_1$ and $\phi_2$ equivalent if there is a $g \in G(E)$ such that $\phi_2={}^g \phi_1$. Then $H^1(\Gal(E/F),G(E))$ is the set of equivalence classes of cocycles. It is a pointed set with neutral element represented by the trivial cocycle (identically $1$). If $G$ is abelian then $H^1(\Gal(E/F),G(E))$ is an abelian group. We will sometimes write $H^1(F,G)$ for $H^1(\Gal(\barF/F),G(\barF))$.

If $\phi \in Z^1(\Gal(E/F),G^\ad(E))$ then there is a canonically defined algebraic group ${}^\phi G/F$ together with an isomorphism $\iota_\phi: G \times E \iso {}^\phi G \times E$ such that $\sigma \iota_\phi (g)= \iota_\phi ((\ad \phi(\sigma))(\sigma g))$ for all $\sigma \in \Gal(E/F)$ and $g \in G(E)$. If $h \in G^\ad(E)$ then there is a unique isomorphism $\iota_h:{}^\phi G \iso {}^{{}^h\phi}G$ over $F$ such that $\iota_h \circ \iota_\phi=\iota_{{}^h\phi} \circ \ad(h)$. If $\psi \in Z^1(\Gal(E/F), {}^\phi G^\ad (E))$, then $\psi \phi \in Z^1(\Gal(E/F),  G^\ad (E))$ and this gives a bijection of sets
\[ Z^1(\Gal(E/F), {}^\phi G^\ad (E)) \liso Z^1(\Gal(E/F), G^\ad (E)), \]
but this map does not preserve neutral elements. The composite 
\[ \iota_{\psi \phi} \circ \iota_{\phi}^{-1} \circ \iota_\psi^{-1}: {}^\psi({}^\phi G) \liso {}^{\psi \phi} G \]
is defined over $F$. 

If $F$ is a number field and $E/F$ is Galois we will write $\ker^1(\Gal(E/F),G(E))$ for the kernel
\[  \ker (H^1(\Gal(E/F),G(E)) \ra \prod_{v\in V_F} H^1(\Gal(E_{w(v)}/F_v),G(E_{w(v)}))), \]
where for $v \in V_F$ we use $w(v)$ to denote a choice of place of $E$ above $v$. If $G$ is reductive then $\ker^1(\Gal(E/F),G(E))$ is finite. It vanishes if $G$ is semi-simple and either adjoint or simply connected. 
(See for instance theorems 6.6, 6.15 and 6.22 of \cite{pr}.) 
We will sometimes write $\ker^1(F,G)$ for $\ker^1(\Gal(\barF/F),G(\barF))$.

\subsection{Group extensions}

We recall some of the general theory of group extensions. Suppose that $A$ is an abelian group and $G$ is a finite group that acts on $A$. A 2-cocyle $\alpha \in Z^2(G,A)$ is a function $G \times G \ra A$ satisfying the relation
\[ {}^{g_1} \alpha (g_2,g_3) \alpha(g_1,g_2g_3)=\alpha(g_1g_2,g_3) \alpha(g_1,g_2). \]
(We record that this implies that $\alpha(1,g)=\alpha(1,1)$ and $\alpha(g,1)={}^g\alpha(1,1)$.)
If $\alpha$ is a 2-cocycle and $\beta : G \ra A$ is any function then 
\[ {}^\beta \alpha(g_1,g_2) = \alpha(g_1,g_2) \beta(g_1g_2) \beta(g_1)^{-1} {}^{g_1}\beta(g_2)^{-1} \]
is another 2-cocycle. 
If $\alpha$ is a 2-cocycle we 
obtain an extension
\[ \begin{array}{ccccccccc} 0 &\lra& A &\lra& \cE_\alpha &\lra& G &\lra& 0 \\
&& a & \longmapsto & a\alpha(1,1)^{-1}e_\alpha(1) &&&& \\ &&&& a e_\alpha(g) & \longmapsto & g, && \end{array} \]
where $\cE_\alpha$ is the group with elements $ae_\alpha(g)$ with $a \in A$ and $g \in G$ with the multiplication rule
\[ a_1 e_\alpha(g_1) a_2 e_\alpha(g_2) = a_1 {}^{g_1}a_2 \alpha (g_1,g_2) e_\alpha(g_1g_2). \]
There is an isomorphism of extensions
\[ \begin{array}{rcl}  i_\beta: \cE_\alpha & \liso & \cE_{{}^\beta \alpha} \\ a e_\alpha(g) & \longmapsto & a \beta(g) e_{{}^\beta \alpha}(g) \end{array} \] 
for any map $\beta:G \ra A$.
Thus the isomorphism class of the extension $\cE_\alpha$ only depends on $[\alpha] \in H^2(G,A)$, but {\em not canonically}. If $a \in A$ we set $({}^a \beta)(g)=\beta(g) a/{}^ga$, and we have
\[ i_{{}^a\beta} =\conju_a \circ i_\beta .\]
Any element $\beta \in Z^1(G,A)$ gives rise to an automorphism of extensions $i_\beta: \cE_\alpha \ra \cE_\alpha$, and in fact this establishes an isomorphism between $Z^1(G,A)$ and the group of automorphisms of the extension $\cE_\alpha$. The automorphism arises as conjugation by an element of $A$ if and only if $\beta$ is a coboundary. Thus, if $H^1(G,A)=(0)$, then every automorphism of the extension $\cE_\alpha$ arises by conjugation by an element of $A$. Every extension of $G$ by $A$ arises from some $\alpha \in Z^2(G,A)$. 

If $h:G \ra G'$ and $f:A \ra A'$ are morphisms such that
\[ f({}^g a) ={}^{h(g)}f(a) \]
and if $\alpha \in Z^2(G,A)$ and $\alpha'\in Z^2(G',A')$ satisfy
\[ f(\alpha(g_1,g_2))=\alpha'(h(g_1),h(g_2)) \]
then there is a morphism
\[ \begin{array}{rcl} (f,h): \cE_\alpha &\lra& \cE_{\alpha'} \\ ae_\alpha(g) & \longmapsto & f(a)e_{\alpha'}(h(g)) \end{array} \]
such that
\[  \begin{array}{ccccccccc} 0 &\lra& A &\lra& \cE_\alpha &\lra& G &\lra& 0 \\ && f \da && (f,h) \da && h \da && \\
 0 &\lra& A' &\lra& \cE_{\alpha'} &\lra& G' &\lra& 0 \end{array} \]
 commutes.
 
 If $h: G' \ra G$ and $\cE$ is an extension of $G$ by $A$, then we can form a pull-back extension 
 \[ \cE|_{G'}=\cE|_{G',h}  = \cE \times_{G,h} G' = \{ (e,g'):\,\, e\,\, {\rm and}\,\, g'\,\, {\rm have\,\, the\,\, same\,\,image \,\, in} \,\, G\} \subset \cE \times G'  \]
 of $G'$ by $A$. If $\cE$ arises from $\alpha \in Z^2(G,A)$ then $\cE|_{G',h}$ arises from $h^*\alpha = \alpha \circ h$. Similarly if $f:A \ra A'$ is a map of $G$-modules we can form a push-out extension
 \[ f_* \cE = (A' \rtimes \cE) /A \]
 of $A'$ by $G$. Here $\cE$ acts on $A'$ via its projection to $G$ and we embed $A$ as a normal subgroup of $A' \rtimes \cE$ via $a \mapsto (f(a)^{-1},a)$. If $\cE$ arises from $\alpha \in Z^2(G,A)$ then $f_*\cE$ arises from $f_*\alpha = f \circ \alpha$. 
 
 If 
 \[ G'' \stackrel{h'}{\lra} G' \stackrel{h}{\lra} G \]
 and
 \[ A \stackrel{f}{\lra} A' \stackrel{f'}{\lra} A'' \]
 then
 \[ \begin{array}{rcl} (A'' \rtimes ((A' \rtimes \cE|_{G'}) /A)|_{G''})/A' &\liso& ((A'' \rtimes \cE|_{G''}) /A \\
 {}[(a'',([(a',(e,g'))],g''))] & \longmapsto & [(a''f'(a'),(e,g''))] . \end{array}\]
 
 If $H \lhd G$ are finite groups, if $s:G/H \ra G$ is a set theoretic section to the projection map, and if $A$ is a $G$-module, then there is a homomorphism
 \[ \nak_s: Z^2(G,A) \lra Z^2(G/H,A^H) \]
 defined by
 \[ \begin{array}{rl} & (\nak_s \alpha)(\barg_1,\barg_2) \\ =& \prod_{h \in H}\alpha(h,s(\barg_1)s(\barg_2)s(\barg_1\barg_2)^{-1}) {}^h\alpha(s(\barg_1),s(\barg_2))/{}^h\alpha(s(\barg_1)s(\barg_2)s(\barg_1\barg_2)^{-1},s(\barg_1\barg_2)) .\end{array} \]
 It induces a homomorphism
 \[ \nak^G_{G/H}:  H^2(G,A) \lra H^2(G/H,A^H) \]
 which is independent of the choice of $s$. If 
 \[ (0) \lra A \lra \cE \lra G \lra (0) \]
 is an extension with class $[\alpha] \in H^2(G,A)$, then the extension
 \[ (0) \lra A^H \lra \tr_{\cE|_H/A,*} (\cE/[\cE|_H,\cE|_H]) \lra G/H \lra (0) \]
 has class $\nak^G_{G/H} [\alpha] \in H^2(G/H,A^H)$. Here
 \[  \tr_{\cE|_H/A}: \cE|_H/[\cE|_H,\cE|_H] \lra A^H \]
 is the transfer map. We have
 \[ {\inf}^G_{G/H} \circ \nak^G_{G/H}=N_{H,*}: H^2(G,A) \lra H^2(G,A^H)  \]
 and
 \[ {\inf}^G_{G/H} \circ \nak^G_{G/H}=\# H: H^2(G,A) \lra H^2(G,A)  \]
 (the Akizuki-Witt theorem). Here
 \[ \begin{array}{rcl} N_H: A & \lra & A^H \\ a & \longmapsto & \prod_{h \in H} {}^h a. \end{array} \]
 
 \begin{lem}\label{grp} Suppose that $\pi:G \onto H$ is a surjective group homomorphism, and that $A \lhd G$ and $B \lhd H$ are normal subgroups with $\pi A \subset B$. Then there is an isomorphism of groups
\[ \begin{array}{rcl} (B \rtimes G)/A & \liso & H \times_{H/B} G/A \\ {} [(b,g)] & \longmapsto & (b \pi(g), gA). \end{array} \]
Here 
\[ \begin{array}{rcl} A & \lra & B \rtimes G \\ a & \longmapsto & (\pi(a)^{-1},a). \end{array} \]
\end{lem}

\pfbegin The map
\[ \begin{array}{rcl} \phi: B \rtimes G & \liso & H \times_{H/B} G/A \\ (b,g) & \longmapsto & (b \pi(g), gA) \end{array} \]
is easily checked to be a group homomorphism. We have $(b,g) \in \ker \phi$ if and only if $g \in A$ and $b = \pi(g)^{-1}$, i.e. $(b,g)$ is in the image of $A$. It remains to check that $\phi$ is surjective. Suppose $(h,gA) \in H \times_{H/B} G/A$. Then $h \pi(g)^{-1} \in B$ and $\phi(h \pi(g)^{-1}, g)=(h,g)$, as desired.\pfend

\subsection{Local Weil Groups}

We recall the theory of Weil groups for p-adic fields. See \cite{tatecorvallis}.

First suppose that $F$ is a p-adic field and that $\barF$ is an algebraic closure. If $k$ denotes the residue field of $F$, there is an exact sequence
\[ (0) \lra I_F \lra \Gal(\barF/F) \lra \Gal(\bark/k) \lra (0). \]
We denote by $W_{\barF/F}$ the preimage of $\Frob_k^\Z \subset \Gal(\bark/k)$ and endow it with a topology decreeing that $I_F$ should be an open subgroup with its usual topology. If $\sigma:\barF \iso \barF'$ is a continuous automorphism with $\sigma F=F'$, then there is a canonical isomorphism
\[ \begin{array}{rcl} \conju_\sigma: W_{\barF/F} & \liso & W_{\barF'/F'}\\
\tau& \longmapsto &\sigma \tau \sigma^{-1}. \end{array}\]
Note that there is a canonical map $\varphi_{\barF/F}: W_{\barF/F} \ra \Gal(\barF/F)$ with dense image.

If $E$ is an intermediate field between $F$ and $\barF$ we will write $W_{\barF/E}=\varphi_{\barF/F}^{-1}\Gal(\barF/E)$. For $E/F$ finite, we will write $E^\ab$ for the (unique) maximal Galois extension of $E$ in $\barF$ with $\Gal(E^\ab/E)$ abelian. Then
\[ W_{E^\ab/F}= W_{\barF/E}^\ab := W_{\barF/F}/\overline{[W_{\barF/E},W_{\barF/E}]}. \]
There are canonical isomorphisms
\[ r_E: E^\times \liso W_{\barF/E}^\ab \]
with the following properties:
\begin{itemize}
\item $\varphi_{\barF/F} \circ r_E=\Art_E$.
\item If $\sigma \in \Gal(\barF/F)$, then $\conju_\sigma \circ r_E = r_{{}^{\sigma}E} \circ \sigma$.
\item If $E' \subset E$ then $\tr_{W_{\barF/E'}/W_{\barF/E}} \circ r_{E'}$ equals $r_E$ composed with the inclusion $(E')^\times \into E^\times$.
\item $W_{\barF/F} \cong \lim_{\la E} W_{E^\ab/F}$ as topological groups.
\end{itemize}
These properties imply that
\begin{itemize}
\item If $E' \subset E$ then $r_{E'}\circ \norm_{E/E'}$ equals $r_E$ followed by the map $W_{\barF/E}^\ab \ra W_{\barF/E'}^\ab$ induced by the inclusion $W_{\barF/E}\subset  W_{\barF/E'}$.
\end{itemize}
There are no non-identity automorphisms of $W_{\barF/F}$ compatible with $\varphi_{\barF/F}$. 

If $F \cong \C$ as topological fields we set $W_{\barF/F}=F^\times$. If $F \cong \R$ as topological fields and $\barF$ is an algebraic closure we set 
$W_{\barF/F}=\langle \barF^\times, j:\,\, j^2=-1,\,\, jzj^{-1}={}^cz \rangle$. If $F=\R$ or $\C$ there is a natural map 
\[ \varphi_{\barF/F}: W_{\barF/F} \onto \Gal(\barF/F) \]
with kernel $\barF^\times$. If $E$ is an intermediate field between $F$ and $\barF$ we will write $W_{\barF/E}=\varphi_{\barF/F}^{-1}\Gal(\barF/E)$. For $E/F$ finite, there are canonical isomorphisms
\[ r_E: E^\times \liso W_{\barF/E}^\ab \]
which are the identity if $E \cong \C$ and, if $E \cong \R$, are induced by $-1 \mapsto j$ and $x \mapsto \sqrt{x}$ for $x>0$. These structures share the properties itemized above for p-adic fields. Again these constructions are functorial in the pair $\barF/F$. In the case $F\cong \C$ the group $W_{\barF/F}$ has no automorphisms compatible with $r_F$. On the other hand, if $F \cong \R$ then $W_{\barF/F}$ does have automorphisms compatible with $\varphi_{\barF/F}$ and $r_F$ and $r_\barF$, namely the inner automorphisms $\conju_z$ for $z \in \barF^\times$. However the only ones compatible with the functoriality $W_{\barF/F} \ra W_{\barF/F}$ induced by $c:\barF \ra \barF$ are the identity and $\conju_{\sqrt{-1}}$.

If $F$ is either a p-adic field or isomorphic to $\R$ or $\C$, and if $\barF$ is an algebraic closure of $F$ and if $E$ is an intermediate field finite and Galois over $F$, then there is a short exact sequence
\[ (0) \lra E^\times \stackrel{r_E}{\lra} W_{E^\ab/F} \lra \Gal(E/F) \lra (0) \]
which determines a class 
\[ [\alpha_{E/F}] \in H^2(\Gal(E/F),E^\times), \]
called the canonical class. It depends only on $E/F$, i.e. not on $\barF$. If $D \supset E \supset F$ are finite extensions (inside $\barF$) with $D/F$ Galois then
\begin{itemize}
\item $\res^{\Gal(D/F)}_{\Gal(D/E)} [\alpha_{D/F}]=[\alpha_{D/E}] \in H^2(\Gal(D/E),D^\times)$,
\item  and $\corr_{\Gal(D/E)}^{\Gal(D/F)} [\alpha_{D/E}]=[E:F][\alpha_{D/F}] \in H^2(\Gal(D/F),D^\times)$. 
\end{itemize}
If in addition $E/F$ is Galois then
\begin{itemize}
\item ${\inf}_{\Gal(E/F)}^{\Gal(D/F)} [\alpha_{E/F}]=[D:E][\alpha_{D/F}] \in H^2(\Gal(D/F),D^\times)$,
\item and $\nak_{\Gal(D/F)}^{\Gal(E/F)} [\alpha_{D/F}]=[\alpha_{E/F}] \in H^2(\Gal(E/F),E^\times)$.
\end{itemize}
(For the first three see for instance section XI.3 of \cite{serre}. For the final assertion use the injectivity of the map
\[ {\inf}_{\Gal(E/F)}^{\Gal(D/F)}: H^2(\Gal(E/F),E^\times) \lra H^2(\Gal(D/F),D^\times) \]
and the Akizuki-Witt theorem.) Hence
\begin{itemize}
\item $\inf_{\Gal(E/F)}^{\Gal(D/F)} [\alpha_{E/F}]= N_{D/E,*}[\alpha_{D/F}] \in H^2(\Gal(D/F),E^\times)$.
\end{itemize}
We will write
\[ [\alpha_{E/F,D}]={\inf}_{\Gal(E/F)}^{\Gal(D/F)} [\alpha_{E/F}]=[D:E][\alpha_{D/F}] \in H^2(\Gal(D/F),D^\times). \]
If $C \supset D \supset E \supset F$ are finite Galois extensions of $F$ then
\[ [\alpha_{E/F,C}]={\inf}_{\Gal(D/F)}^{\Gal(C/F)} [\alpha_{E/F,D}] \]
and
\[ [D:E] [\alpha_{D/F,C}]=[\alpha_{E/F,C}]. \]

If $D \supset E \supset F$ are finite Galois extensions of $F$ in $\barF$ then there is an obvious map
\[ W_{D^\ab/F} \onto W_{E^\ab/F} \]
which fits into a commutative diagram
\[ \begin{array}{ccccccccc} (0)& \lra & D^\times & \lra &W_{D^\ab/F} & \lra & \Gal(D/F) & \lra & (0) \\ && N_{D/E} \da && \da && \da && \\
(0)& \lra & E^\times & \lra & W_{E^\ab/F} & \lra & \Gal(E/F) & \lra & (0). \end{array} \]
We will denote this map $\sigma \mapsto \sigma|_{E^\ab}$.

If we map
\[ \begin{array}{rcl} D^\times & \lra & E^\times \rtimes W_{D^\ab/F} \\ a & \longmapsto & 
(N_{D/E}(a)^{-1}, r_D(a)), \end{array} \]
then we see that there is an isomorphism of extensions
\[ \begin{array}{rcl} (E^\times \rtimes W_{D^\ab/F})/D^\times & \liso & W_{E^\ab/F}|_{\Gal(D/F)} = W_{E^\ab/F} \times_{\Gal(E/F)} \Gal(D/F) \\ {}
 [(a,\tau)] & \longmapsto & (r_E(a)\tau|_{E^\ab},\tau|_D).\end{array} \]
(Note that $r_D(b)|_{E^\ab}=r_E(N_{D/E}(b))$ and recall lemma \ref{grp}.) 
We see that we have maps of extensions
\[ \begin{array}{ccccccccc} (0)& \lra & D^\times & \lra &W_{D^\ab/F} & \lra & \Gal(D/F) & \lra & (0) \\ && N_{D/E} \da && \da && || && \\
(0)& \lra & E^\times & \lra & W_{E^\ab/F}|_{\Gal(D/F)} & \lra & \Gal(D/F) & \lra & (0) \\ && || && \da && \da && \\
(0)& \lra & E^\times & \lra & W_{E^\ab/F} & \lra & \Gal(E/F) & \lra & (0),\end{array} \]
whose composite is the natural surjection $W_{D^\ab/F} \onto W_{E^\ab/F}$, and where the middle row can either be obtained as a pushout from the top row or a pullback from the bottom row. Moreover the maps
\[ \begin{array}{rcl} \Gal(D/E) & \lra & W_{E^\ab/F}|_{\Gal(D/F)} \\
\sigma & \longmapsto & (1,\sigma) \end{array} \]
and
\[ \begin{array}{rcl} \Gal(D/E) & \lra & (E^\times \rtimes W_{D^\ab/F})/D^\times \\
\sigma & \longmapsto & [(r_E^{-1}(\tsigma|_{E^\ab})^{-1},\tsigma)], \end{array} \]
where $\tsigma \in W_{D^\ab/E}$ is any lift of $\sigma$,
are identified. (The latter map is easily checked to be well defined, i.e. independent of the choice of lift $\tsigma$.)

Continue to suppose that $D \supset E \supset F$ are finite Galois extensions of $F$. 
We define a pushout
\[ W_{E^\ab/F,D}= (D^\times \rtimes W_{E^\ab/F}|_{\Gal(D/F)})/E^\times, \]
where
\[ \begin{array}{rcl} E^\times & \lra & D^\times \rtimes W_{E^\ab/F}|_{\Gal(D/F)} \\ a & \longmapsto & (a^{-1},(r_E(a),1)). \end{array} \]
It comes with a natural section 
\[ \begin{array}{rcl} \Gal(D/E) & \lra & W_{E^\ab/F,D} \\
\sigma & \longmapsto & [(1,(1,\sigma))]. \end{array} \]
From the above discussion we see that this has a second description as 
\[ \begin{array}{rcl} \Gal(D/E) &\into& (D^\times \rtimes W_{D^\ab/F})/D^\times \\
\sigma & \longmapsto & [(r_E^{-1}(\tsigma|_{E^\ab})^{-1},\tsigma)], \end{array}\]
 where $\tsigma \in W_{D^\ab/E}$ is any lift of $\sigma$ and
\[ \begin{array}{ccc}D^\times  & \into & D^\times \rtimes W_{D^\ab/F} \\ a & \longmapsto & ((N_{D/E}a)^{-1},r_D(a)) .\end{array}\]
The map is
\[ \begin{array}{rcl} (D^\times \rtimes W_{D^\ab/F})/D^\times & \liso & (D^\times \rtimes W_{E^\ab/F}|_{\Gal(D/F)})/E^\times \\
{}[(a,\sigma)] & \longmapsto & [(a,(\sigma|_{E^\ab},\sigma|_D))]. \end{array} \]
The extension corresponds to the class
\[ {\inf}_{\Gal(E/F)}^{\Gal(D/F)} [\alpha_{E/F}] = N_{D/E,*} [\alpha_{D/F}] \in H^2(\Gal(D/F), D^\times). \]
We have a commutative diagram of extensions
\[ \begin{array}{ccccccccc} 
(0) & \lra & E^\times & \lra & W_{E^\ab/F} & \lra & \Gal(E/F) & \lra & (0) \\
&& || && \ua && \ua && \\
(0) & \lra & E^\times & \lra & W_{E^\ab/F}|_{\Gal(D/F)} & \lra & \Gal(D/F) & \lra & (0) \\
&& \bigcap && \da && || && \\
(0) & \lra & D^\times & \lra & W_{E^\ab/F,D} & \lra & \Gal(D/F) & \lra & (0) \\
&& N_{D/E} \ua && \ua && || && \\
(0) & \lra & D^\times & \lra & W_{D^\ab/F} & \lra & \Gal(D/F) & \lra & (0) . \end{array}\]

\subsection{Global Weil groups}\label{gwg}

We now recall the theory of Weil groups for number fields fields. See \cite{tatecorvallis}.

Now suppose that $F$ is a number field and that $\barF$ is an algebraic closure of $F$. One can associate to $\barF/F$ a topological group $W_{\barF/F}$ together with:
\begin{itemize}
\item A map $\varphi_{\barF/F}: W_{\barF/F} \onto \Gal(\barF/F)$. If $E$ is an intermediate field we set $W_{\barF/E}=\varphi_{\barF/F}^{-1} \Gal(\barF/E)$.
\item For each intermediate field finite over $F$ a map
\[ r_E: \A_E^\times/E^\times \liso W_{\barF/E}^\ab \]
such that $\varphi_{\barF/F} \circ r_E=\Art_E$.
\end{itemize}
These maps also satisfy:
\begin{itemize}
\item If $\sigma \in W_{\barF/F}$, then $\conju_\sigma \circ r_E = r_{{}^{\varphi_{\barF/F}(\sigma)}E} \circ \varphi_{\barF/F}(\sigma)$.
\item If $E' \subset E$ then $\tr_{W_{\barF/E'}/W_{\barF/E}} \circ r_{E'}$ equals $r_E$ composed with the inclusion $\A_{E'}^\times/(E')^\times \into \A_E^\times/E^\times$.
\item $W_{\barF/F} \cong \lim_{\la E} W_{\barF/F}/\overline{[W_{\barF/E},W_{\barF/E}]}$ as topological groups. 
\end{itemize}
These properties imply that
\begin{itemize}
\item If $E' \subset E$ then $r_{E'}\circ \norm_{E/E'}$ equals $r_E$ followed by the map $W_{\barF/E}^\ab \ra W_{\barF/E'}^\ab$ induced by the inclusion $W_{\barF/E}\subset  W_{\barF/E'}$.
\item $W_{E^\ab/F}= W_{\barF/F}/\overline{[W_{\barF/E},W_{\barF/E}]}$. 
\end{itemize}
The only automorphisms of $W_{\barF/F}$ compatible with $\varphi_{\barF/F}$ and the $r_E$ are the inner automorphisms $\conju_\sigma$ for $\sigma \in W_{\barF/\barF}$. The structure $(W_{\barF/F}, \varphi_{\barF/F},\{ r_E\})$ is unique up to isomorphism. 
However we do not know how to make the isomorphism canonical. (If it can be made canonical.)

The image of $W_{\barF/\barF}$ in $W_{E^\ab/F}$ is $\overline{E^\times (E_\infty^\times)^0}/E^\times=\ker \Art_E$. We will denote it $\Delta_E$.

\begin{lem} \label{weillem}
\begin{enumerate}
\item For $i>0$ there is an isomorphism $H^i(\Gal(E/F), \prod_{v|\infty} E_v^\times) \liso H^i(\Gal(E/F),\Delta_E)$. Moreover 
\[  H^i(\Gal(E/F),\Delta_E) \cong H^i(\Gal(E/F), \prod_{v|\infty} E_v^\times) \cong \prod_v H^i(\Gal(E/F)_{w(v)}, E_{w(v)}^\times), \]
where $v$ runs over real places of $F$ which do not split completely in $E$, and $w(v)$ is a choice of place of $E$ above $v$. 

\item $N_{E/F} \Delta_E=\Delta_F$.

\item $H^1(\Gal(E/F),\Delta_E)=(0)$.

\item If $E$ is totally imaginary, then $\Delta_E^{\Gal(E/F)}=\overline{F^\times F_\infty^\times}/F^\times$. In particular if $F$ is totally imaginary then $\Delta_E^{\Gal(E/F)}=\Delta_F$. 

\item $\Art_E:\A_F^\times/F^\times \onto \Gal(E^\ab/E)^{\Gal(E/F)}$.

\item If $S$ is a finite set of finite places of $F$, then $\Delta_F \cap \prod_{v \in S} F_v^\times =\{1\} \subset \A_F^\times/F^\times$.

\end{enumerate} \end{lem}

\pfbegin The first four parts are proved in section III of \cite{weil}. The fifth part follows from the third. The sixth part follows from the congruence subgroup property for the group of $S$-units $\cO_F[1/S]^\times$ of $F$. (See \cite{chev}.) 

More precisely if $a \in \prod_{v \in S} F_v^\times$ maps to the identity in $\A_F^\times/F^\times \Delta_F$ then we can find $a_i \in F^\times$ with $a_i \ra a \in (\A_F^\infty)^\times$. Then for $i$ sufficiently large we have $a_i \in \cO_F[1/S]^\times$. For any $m \in \Z_{>0}$ we can find an open subgroup  $U_m \subset (\A_F^{S \cup \{v|\infty\}})^\times$ such that $U_m\cap \cO_F[1/S]^\times \subset (\cO_F[1/S]^\times)^m$. (Combine Dirichlet's and Chevalley's theorems.) Then for $i$ sufficiently large $a_i\in U_m\cap \cO_F[1/S]^\times$ so that $a_i$ is an $m^{th}$ power in $\cO_F[1/S]^\times$ and hence $a_{i,S}$ is an $m^{th}$ power in $\prod_{v \in S} F_v^\times$. Thus $a \in \bigcap_{m>0} (\prod_{v \in S} F_v^\times)^m =\{1\}$, and the sixth part follows.
\pfend

If $u$ is a place of $\barF$ then there is a continuous homomorphism
\[ \theta_u: W_{\barF_u/F_u} \lra W_{\barF/ F} \]
such that 
\begin{itemize}
\item $\varphi_{\barF/F} \circ \theta_u$ equals the composite of $\varphi_{\barF_u/F_u}$ with the inverse of the canonical map $ \Gal(\barF/F)_u\iso \Gal(\barF_u/F_u)$;
\item and, for $E$ a finite intermediate field, $\theta_u \circ r_{E_u}$ equals the composite of $r_E$ with the canonical map $E_u^\times \ra \A_E^\times/E^\times$.
\end{itemize}
The map $\theta_u$ is determined up to conjugation by an element of $W_{\barF/\barF}$. The images of $\conju_a \circ \theta_u$, for any $a \in W_{\barF/\barF}$, are referred to as {\em decomposition groups} for $u$. The closure of the image under $\varphi_{\barF/F}$ of any decomposition group for $u$ is $\Gal(\barF/F)_u$. If $\sigma \in W_{\barF/F}$, then 
\[ \theta_{\varphi_{\barF/F}(\sigma)u} \circ \conju_{\varphi_{\barF/F}(\sigma)}= \conju_\sigma \circ \theta_u,\]
 up to conjugation by an element of $W_{\barF/\barF}$. 
 
 If $v$ is a place of $F$ and if $\rho:\barF \ra \barFv$ is $F$-linear, then we define 
 \[ \theta_\rho = \theta_{u(\barF,\rho)} \circ \rho^*: W_{\barFv/F_v} \lra W_{\barF/F} \]
 (where $\rho^*: W_{\barFv/F_v} \iso W_{\barF_{u(\barF,\rho)}/F_v}$). Up to conjugation by an element of $W_{\barF/\barF}$:
 \begin{enumerate}
 \item If $\sigma \in W_{\barF/F}$, then $\theta_{\rho\sigma}=\conju_{\sigma^{-1}} \circ \theta_\rho$ up to conjugation by an element of $W_{\barF/\barF}$.
 \item If $\tau \in \Gal(\barFv/F_v)$, then $\theta_{\tau\rho}=\theta_\rho \circ \conju_{\tau^{-1}}$.
 \item If $E$ is an intermediate field finite over $F$, then $\theta_\rho \circ r_{\rho(E)F_v}$ equals the composition of $r_E$ with $(\rho(E)F_v)^\times \iso E_{u(\rho)}^\times \ra \A_E^\times/E^\times$, where the first map is the inverse of the continuous extension of $\rho$.
  \end{enumerate}
  We get an induced map 
  \[ \theta_\rho: W_{(\rho(E)F_v)^\ab/F_v} \lra W_{E^\ab/F}, \]
  defined up to conjugation by an element of $\Delta_E$. Up to this ambiguity, it only depends $\rho|_{E^\ab}$.
 
If $E$ is an intermediate field, finite and Galois over $F$, the  short exact sequence
\[ (0) \lra \A^\times_E /E^\times \stackrel{r_E}{\lra}  W_{E^\ab/F} \lra \Gal(E/F) \lra (0) \]
which determines a class 
\[ [\alpha_{E/F}^W] \in H^2(\Gal(E/F),\A_E^\times/E^\times) \]
called the {\em canonical class}. This class depends only on $E/F$. 
If $\imath_w:E_w^\times \ra \A_E^\times/E^\times$, then we have
\[ \res [\alpha_{E/F}^W] = \imath_{w,*} [\alpha_{E_w/F_v}] \in H^2(\Gal(E/F)_w,\A_E^\times/E^\times). \]
(See formula (12) in \cite{tatenakayama}.) 

If $D \supset E \supset F$ are finite extensions (inside $\barF$) with $D/F$ Galois then
\begin{itemize}
\item $\res^{\Gal(D/F)}_{\Gal(D/E)} [\alpha_{D/F}^W]=[\alpha_{D/E}^W] \in H^2(\Gal(D/E),\A_D^\times/D^\times)$,
\item  and $\corr_{\Gal(D/E)}^{\Gal(D/F)} [\alpha_{D/E}^W]=[E:F][\alpha_{D/F}^W] \in H^2(\Gal(D/F),\A_D^\times/D^\times)$. 
\end{itemize}
If in addition $E/F$ is Galois then
\begin{itemize}
\item $\inf_{\Gal(E/F)}^{\Gal(D/F)} [\alpha_{E/F}^W]=[D:E][\alpha_{D/F}^W] \in H^2(\Gal(D/F),\A_D^\times/D^\times)$,
\item and $\nak_{\Gal(D/F)}^{\Gal(E/F)} [\alpha_{D/F}^W]=[\alpha_{E/F}^W] \in H^2(\Gal(E/F),\A_E^\times/E^\times)$.
\end{itemize}
(For the first and third assertions see for instance \cite{tatecorvallis}. The second follows from the first and the fourth from the third.)
  Hence
\begin{itemize}
\item $\inf_{\Gal(E/F)}^{\Gal(D/F)} [\alpha_{E/F}^W]= N_{D/E,*}[\alpha_{D/F}^W] \in H^2(\Gal(D/F),\A_E^\times/E^\times)$.
\end{itemize}
We will write
\[ [\alpha_{E/F,D}^W]={\inf}_{\Gal(E/F)}^{\Gal(D/F)} [\alpha_{E/F}^W]=[D:E][\alpha_{D/F}^W] \in H^2(\Gal(D/F),\A_D^\times/D^\times). \]
Note that if $u|w|v$ are places of $D \supset E\supset F$, then 
\[ \res [\alpha^W_{E/F,D}] = \imath_{u,*} ([D:E]/[D_u:E_w])[\alpha_{E_w/F_v,D_u}] \in H^2(\Gal(D/F)_u, \A_D^\times/D^\times). \]
If $C \supset D \supset E \supset F$ are finite Galois extensions of $F$ then
\[ [\alpha^W_{E/F,C}]={\inf}_{\Gal(D/F)}^{\Gal(C/F)} [\alpha^W_{E/F,D}] \]
and
\[ [D:E] [\alpha^W_{D/F,C}]=[\alpha^W_{E/F,C}]. \]
 
Continue to suppose that $D\supset E\supset F$ are finite Galois extensions of $F$ in $\barF$. We have a natural map
 \[ W_{D^\ab/F} \onto W_{E^\ab/F} \]
 which we will denote $\sigma \mapsto \sigma|_{E^\ab}$. Lemma \ref{grp} tells us there is an isomorphism
\[ \begin{array}{rcl} (\A_E^\times/E^\times \rtimes W_{D^\ab/F})/(\A_D^\times/D^\times) & \liso & W_{E^\ab/F}|_{\Gal(D/F)} = W_{E^\ab/F} \times_{\Gal(E/F)} \Gal(D/F) \\ {}
 [(a,\tau)] & \longmapsto & (r_E(a)\tau|_{E^\ab},\tau|_D),\end{array} \]
 where
 \[ \begin{array}{rcl} \A_D^\times/D^\times & \lra & \A_E^\times/E^\times \rtimes W_{D^\ab/F} \\ a & \longmapsto & 
(N_{D/E}(a)^{-1}, r_D(a)). \end{array} \]
It is compatible with the maps
\[ \begin{array}{rcl} \Gal(D/E) & \lra & W_{E^\ab/F}|_{\Gal(D/F)} \\
\sigma & \longmapsto & (1,\sigma) \end{array} \]
and
\[ \begin{array}{rcl} \Gal(D/E) & \lra & (\A_E^\times/E^\times \rtimes W_{D^\ab/F})/(\A_D^\times/D^\times) \\
\sigma & \longmapsto & [(r_E^{-1}(\tsigma|_{E^\ab})^{-1},\tsigma)] \end{array} \]
for any lift $\tsigma \in W_{D^\ab/E}$ of $\sigma$. this latter map is well defined, i.e. independent of the choice of lift $\tsigma$.
We see that we have maps of extensions
\[ \begin{array}{ccccccccc} (0)& \lra & \A_D^\times/D^\times & \lra &W_{D^\ab/F} & \lra & \Gal(D/F) & \lra & (0) \\ && N_{D/E} \da && \da && || && \\
(0)& \lra & \A_E^\times/E^\times & \lra & W_{E^\ab/F}|_{\Gal(D/F)} & \lra & \Gal(D/F) & \lra & (0) \\ && || && \da && \da && \\
(0)& \lra & \A_E^\times/E^\times & \lra & W_{E^\ab/F} & \lra & \Gal(E/F) & \lra & (0),\end{array} \]
whose composite is the natural surjection $W_{D^\ab/F} \onto W_{E^\ab/F}$, and where the middle row can either be obtained as a push-out from the top row or a pullback from the bottom row.

 We define a pushout
\[ W_{E^\ab/F,D}= (\A_D^\times/D^\times \rtimes W_{E^\ab/F}|_{\Gal(D/F)})/(\A_E^\times/E^\times), \]
where
\[ \begin{array}{rcl} \A_E^\times/E^\times & \lra & \A_D^\times/D^\times \rtimes W_{E^\ab/F}|_{\Gal(D/F)} \\ a & \longmapsto & (a^{-1},(r_E(a),1)), \end{array} \]
together with a splitting
\[ \begin{array}{rcl} \Gal(D/E) & \lra & W_{E^\ab/F,D} \\
\sigma & \longmapsto & [(1,(1,\sigma))]. \end{array} \]
The extension $W_{E^\ab/F,D}$ corresponds to the class $[\alpha_{E/F,D}^W] \in H^2(\Gal(D/F), \A_D^\times/D^\times)$.
Note that if $C \supset D\supset E \supset F$ are finite Galois extensions of $F$, then 
\[ \begin{array}{rcl} (\A_C^\times/C^\times \rtimes W_{E/F,D}|_{\Gal(C/F)})/(\A_D^\times/D^\times) &\liso& W_{E/F,C} \\ {}[(b,([(a,(\sigma,\tau))],\eta)] & \longmapsto & [(ba,(\sigma,\eta))].
\end{array} \]
Also note that there is a natural embedding 
\[ \begin{array}{rcl} \Gal(E^\ab D/F) & \into & W_{E^\ab/F,D}/\Delta_E \\ \sigma & \longmapsto & [(1,(\widetilde{\sigma|_{E^\ab}},\sigma|_D))] , \end{array} \]
where $\widetilde{\sigma|_{E^\ab}}$ is any lift of $\sigma|_{E^\ab}$ to $W_{E^\ab/F}$.

From the above discussion we see that $W_{E^\ab/F,D}$ has a second description as 
\[  (\A_D^\times/D^\times \rtimes W_{D^\ab/F})/\A_D^\times \]
where
\[ \begin{array}{ccc}\A_D^\times/D^\times  & \into & \A_D^\times/D^\times \rtimes W_{D^\ab/F} \\ a & \longmapsto & ((N_{D/E}a)^{-1},r_D(a)). \end{array}\]
The map is given by
\[ \begin{array}{rcl} (\A_D^\times/D^\times \rtimes W_{D^\ab/F})/(\A_D^\times/D^\times) & \liso & (\A_E^\times/E^\times \rtimes W_{E^\ab/F}|_{\Gal(D/F)})/(\A_E^\times/E^\times)  \\ {}
 [(a,\tau)] & \longmapsto & [(a,(\tau|_{E^\ab},\tau|_D))].\end{array} \]
 In the reverse direction
\[ \begin{array}{rcl} W_{E^\ab/F}|_{\Gal(D/E)} & \lra & (\A_D^\times/D^\times \rtimes W_{D^\ab/F})/(\A_D^\times/D^\times) \\ (\sigma,\tau) & \longmapsto &  [(r_E^{-1}(\sigma\ttau|_{E^\ab}^{-1}),\ttau)] ,\end{array}\]
 where $\ttau \in W_{D^\ab/F}$ is any lift of $\tau$.

We have a commutative diagram of extensions
\[ \begin{array}{ccccccccc} 
(0) & \lra & \A_E^\times/E^\times & \lra & W_{E^\ab/F} & \lra & \Gal(E/F) & \lra & (0) \\
&& || && \ua && \ua && \\
(0) & \lra & \A_E^\times/E^\times & \lra & W_{E^\ab/F}|_{\Gal(D/F)} & \lra & \Gal(D/F) & \lra & (0) \\
&& \bigcap && \da && || && \\
(0) & \lra & \A_D^\times/D^\times & \lra & W_{E^\ab/F,D} & \lra & \Gal(D/F) & \lra & (0) \\
&& N_{D/E} \ua && \ua && || && \\
(0) & \lra & \A_D^\times/D^\times & \lra & W_{D^\ab/F} & \lra & \Gal(D/F) & \lra & (0) . \end{array}\]

 \newpage

\section{Generalities on algebraic cohomology}\label{algcoh}

We will need a modification of the Galois cohomology of  reductive group, which is described in section 12 of \cite{kotbg}. Suppose that we are given the following data:
\begin{enumerate}
\item An extension 
\[ (0) \lra \cE^0 \lra \cE \lra \Gamma \lra (0) \]
where $\cE^0$ is an abelian group and $\Gamma$ is a group. Note that $\cE^0$ has a $\Gamma$-action by conjugation.
\item A group $G$ with an action of $\Gamma$. Note that the conjugation action of $G$ and the $\Gamma$-action on $G$, piece together to give a $G \rtimes \Gamma$ action on $G$. This, together with the $\Gamma$ action on $\cE^0$, gives a $G \rtimes \Gamma$-action on $\Hom(\cE^0,G)$.
\item A set $\cN$ with an action of $G \rtimes \Gamma$.
\item  A $G \rtimes \Gamma$-equivariant map 
\[ \begin{array}{rcl} \cN &\lra &\Hom(\cE^0,G) \\ \nu & \longmapsto & \barnu, \end{array} \]
such that if $e \in \cE^0$ and $\nu \in \cN$, then
\[ {}^{\barnu(e)}\nu = \nu. \]
\end{enumerate}
We will refer to $\cN$ as an {\em pre-algebraicity condition}. In most cases we will have $\cN\subset \Hom(\cE^0,G)$. 

We define pointed sets
\[ Z^1_\cN(\cE,G) =\{(\nu,\phi) \in  \cN \times Z^1(\cE,G):\,\, \phi|_{\cE^0}=\barnu\,\, {\rm and}\,\, {}^{\phi(e)^{-1}}  \nu   = {}^e\nu\,\, \forall e \in \cE  \}. \]
We refer to elements of this pointed set as {\em algebraic cocycles}.
We write 
\[ \bnu:Z^1_\cN(\cE,G) \lra \cN \]
for the map $(\nu,\phi) \mapsto \nu$.
In cases where $\cN$ is contained $\Hom(\cE^0,G)$ we will often use $\phi$ to denote a cocycle $(\nu,\phi)$ (as $\nu=\phi|_{\cE^0}$).
The group $G$ acts on $Z^1_\cN(\cE,G)$ by
\[ {}^g(\nu,\phi) = ({}^g \nu, {}^g \phi). \]
 We define $H^1_\cN(\cE,G)$ to be the quotient of $Z^1_\cN(\cE,G)$ by $G$. We refer to this as the {\em algebraic cohomology}. 

There is a left exact sequence (of pointed sets)
\[ \begin{array}{ccccccc} (0) & \lra & H^1(\Gamma,G) & \lra & H^1_\cN(\cE,G) & \stackrel{\bnu}
{\lra} & (G \backslash \cN)^\Gamma \\ &&&& (\nu,\phi) & \longmapsto & [\nu]. \end{array} \]
If $\Gamma=\{1\}$, then 
\[ Z^1_\cN(\cE,G) \liso \cN \]
and
\[ H^1_\cN(\cE,G) \liso G\backslash \cN. \]

These sets of cocycles and cohomology sets satisfy various natural functorialities, which are a bit tedious to spell out:

\begin{enumerate}[(A)] 

\item\label{changeg} Suppose first that we have a $\Gamma$-equivariant map $h:G_1 \ra G_2$, pre-algebraicity conditions $\cN_i$ for $(\cE,G_i)$ and a $\Gamma$-equivariant map $n: \cN_1 \ra \cN_2$
such that 
\[ n({}^g\nu)={}^{h(g)}n(\nu)\]
and
\[ \overline{n(\nu)}=h \circ \barnu \]
for all $\nu \in \cN$ and $g \in G_1$. 
Then we obtain a natural map
\[ \begin{array}{rcl} h_*=(h,n)_*: Z^1_{\cN_1}(\cE,G_1) & \lra & Z^1_{\cN_2}(\cE,G_2) \\ (\nu,\phi) & \longmapsto &(n(\nu), h \circ \phi), \end{array} \]
which induces a map
\[  h_*=(h,n)_*: H^1_{\cN_1}(\cE,G_1)  \lra  H^1_{\cN_2}(\cE,G_2). \]

\item\label{changee} Second suppose that we have maps of extensions
\[ \begin{array}{ccccccccc} (0) & \lra & \cE_1^0 & \lra &\cE_1 & \lra & \Gamma_1 & \lra & (0) \\ && \da && f \da && \barf \da && \\ (0) & \lra & \cE_2^0 & \lra &\cE_2 & \lra & \Gamma_2 & \lra & (0), \end{array} \]
a group homomorphism $h:G_2 \ra G_1$, 
pre-algebraicity conditions $\cN_i$ for $(\cE_i,G_i)$ and a map $n:\cN_2 \ra \cN_1$ such that
\begin{itemize}
\item $h({}^{\barf(\sigma)}a)={}^\sigma h(a)$,
\item $n({}^g\nu)={}^{h(g)} n(\nu)$,
\item $n({}^{\barf(\sigma)}\nu)={}^\sigma n(\nu)$
\item and $\overline{n(\nu)}=h \circ \barnu \circ f$.
\end{itemize}
Then we obtain a natural map
\[ \begin{array}{rcl} f^*=(f,h,n)^*: Z^1_{\cN_2}(\cE_2,G_2) & \lra & Z^1_{\cN_1}(\cE_1,G_1) \\ (\nu,\phi) & \longmapsto &(n(\nu), h \circ  \phi \circ f), \end{array} \]
which induces a map
\[ f^* =(f,h,n)^*: H^1_{\cN_2}(\cE_2,G_2)  \lra  H^1_{\cN_1}(\cE_1,G_1) .\]
(This is a minor generalization of Kottwitz's map $\Psi(n,f)$, which he defines in the special case $G_1=G_2$ and $h=\Id_G$. In the special case $f=\Id$ we recover the map in part \ref{changeg}.)

(If we take $\cE_1=\cE_2$ and $G_1=G_2$ and $\cN_1=\cN_2$, and also take $f=\conju_e$ and $h=e^{-1}$ and $n=e^{-1}$ with $e \in \cE$; then
\[ \conju_e^* (\nu,\phi)={}^{{}^{e^{-1}}\phi(e)}(\nu,\phi). \]
In particular $\conju_e^*$ is the identity on $H^1_{\cN_1}(\cE_1,G_1)$.)

We have $(f_1,h_1,n_1)^*\circ (f_2,h_2,n_2)^*=(f_2 \circ f_1,h_1 \circ h_2,n_1 \circ n_2)^*$.

\item\label{changec} Thirdly suppose that we have maps of extensions
\[ \begin{array}{ccccccccc} (0) & \lra & \cE_1^0 & \lra &\cE_1 & \lra & \Gamma & \lra & (0) \\ && \da && f \da && || && \\ (0) & \lra & \cE_2^0 & \lra &\cE_2 & \lra & \Gamma & \lra & (0), \end{array} \]
and a $\Gamma$-equivariant group homomorphism $h:G_1 \ra G_2$ and 
pre-algebraicity conditions $\cN_i$ for $(\cE_i,G_i)$ and a  map $n:\cN_1 \ra \cN_2$ such that
\begin{itemize}
\item $n({}^g\nu)={}^{g} n(\nu)$,
\item $n({}^{\sigma}\nu)={}^\sigma n(\nu)$,
\item  and $h \circ \barnu=\overline{n(\nu)} \circ f$.
\end{itemize}
Then we obtain a natural map
\[ \begin{array}{rcl} f_*=(f,h,n)_*: Z^1_{\cN_1}(\cE_1,G_1) & \lra & Z^1_{\cN_2}(\cE_2,G_2) \\ (\nu,\phi) & \longmapsto &(n(\nu), \tphi), \end{array} \]
where $\tphi$ is defined by $\tphi(tf(e))=\overline{n (\nu)}(t) (h \circ \phi(e))$ for any $t \in \cE_2^0$ and $e \in \cE_1$. This induces
\[ f_*=(f,h,n)_* : H^1_{\cN_1}(\cE_1,G)  \lra  
H^1_{\cN_2}(\cE_2,G). \]
(Kottwitz denotes this map $\Phi(h,n,f)$. The map in part \ref{changeg} is a special case of this map in which $f$ is the identity.) 

We see that $(f_1,h_1,n_1)_* \circ (f_2,h_2,n_2)_* = (f_1 \circ f_2, h_1 \circ h_2, n_1 \circ n_2)_*$.

Suppose that we are given commutative diagrams
\[ \begin{array}{rcl} \cE_2 & \stackrel{f}{\lla} & \cE_1 \\ f_2 \da && \da f_1 \\ \cE_2' & \stackrel{f'}{\lla} & \cE_1' \end{array} \]
and
\[ \begin{array}{rcl} G_2 & \stackrel{h}{\lra} & G_1 \\ h_2 \da && \da h_1 \\ G_2' & \stackrel{h'}{\lra} & G_1' \end{array} \]
and
\[ \begin{array}{rcl} \cN_2 & \stackrel{n}{\lra} & \cN_1 \\ n_2 \da && \da n_1 \\ \cN_2' & \stackrel{n'}{\lra} & \cN_1' \end{array} \]
such that $(f,h,n)$ and $(f',h',n')$ are as in part \ref{changee}, while $(f_1,h_1,n_1)$ and $(f_2,h_2,n_2)$ are as in this part. Then
\[ \begin{array}{rcl} Z^1_{\cN_2}(\cE_2,G_2)  & \stackrel{(f,h,n)^*}{\lra} & Z^1_{\cN_1}(\cE_1,G_1) \\ (f_2,h_2,n_2)_* \da && \da (f_1,h_1,n_1)_* \\ Z^1_{\cN_2'}(\cE_2',G_2') & \stackrel{(f',h'n')^*}{\lra} & Z^1_{\cN_1'}(\cE_1',G_1') \end{array} \]
commutes.

\item\label{shapiro} Suppose that $\Delta \subset \Gamma$ is a subgroup. If $X$ is a set with an action of $\Delta$ we will write $\Ind_\Delta^\Gamma X$ for the set of functions $\varphi:\Gamma \ra X$ satisfying
\[ \varphi(\tau \sigma)= \tau \varphi(\sigma) \]
for all $\tau \in \Delta$ and $\sigma \in \Gamma$. It has an action of $\Gamma$ via
\[ (\sigma \varphi)(\sigma')=\varphi(\sigma'\sigma). \]
If $X$ is a group and $\Delta$ acts via group automorphisms, then $\Ind_\Delta^\Gamma X$ is a group via
\[ (\varphi \varphi')(\sigma)=\varphi(\sigma) \varphi'(\sigma), \]
and its $\Gamma$ action is via group automorphisms. The map
\[ \begin{array}{ccc} \epsilon:\Ind_\Delta^\Gamma X &\onto& X \\ \varphi & \longmapsto & \varphi(1) \end{array} \]
is $\Delta$-equivariant. 

Suppose we have an extension
\[ (0) \lra \cE^0 \lra \cE \lra \Delta \lra (0) \]
and a group $G$ with a $\Delta$-action and a set $\cN$ with a $G \rtimes \Delta$-action, together with a $G \rtimes \Delta$-invariant map $\overline{\,\,\mbox{}}:\cN \ra \Hom(\cE^0,G)$ such that ${}^{\barnu(e)} \nu =\nu$ for all $e \in \cE^0$ and $\nu \in \cN$. Suppose moreover that we are given a second extension
\[ (0) \lra \Ind_\Delta^\Gamma \cE^0 \lra \tcE \lra \Gamma \lra (0) \]
such that if $\tcE|_\Delta$ denotes the preimage of $\Delta$ in $\tcE$, then there is a map of extensions
\[ \begin{array}{ccccccccc} (0)& \lra &\Ind_\Delta^\Gamma \cE^0 & \lra&  \tcE|_\Delta& \lra& \Delta& \lra& (0) \\ && \epsilon \da && \tepsilon \da && || && \\ 
(0)& \lra &\cE^0 & \lra&  \cE& \lra& \Delta& \lra& (0). \end{array} \]
We will write $i$ for the natural inclusion $\tcE|_\Delta \into \tcE$.
Note that if we think of $\cN$ as a pre-algebraicity condition for  $(\tcE|_\Delta, G)$ with the new $\barnu$ equal to $\barnu \circ \epsilon$ for $\nu \in \cN$, then
\[ \tepsilon_*: Z^1_\cN(\tcE|_\Delta,G) \liso Z^1_\cN(\cE,G)  \]
and
\[ \tepsilon_*: H^1_\cN(\tcE|_\Delta,G) \liso H^1_\cN(\cE,G), \]
with inverse $\tepsilon^*$.

Note that $\Ind_\Delta^\Gamma \cN$ has an action of $(\Ind_\Delta^\Gamma G) \rtimes \Gamma$, where
\[ ({}^\varphi \nu)(\sigma)={}^{\varphi(\sigma)}(\nu(\sigma)) \]
for $\varphi \in \Ind_\Delta^\Gamma G$ and $\nu \in \Ind_\Delta^\Gamma \cN$. If $\nu \in \Ind_\Delta^\Gamma \cN$ we define $\barnu \in \Hom(\Ind_\Delta^\Gamma \cE^0, \Ind_\Delta^\Gamma G)$ by
\[ \barnu(\varphi)(\sigma)=\overline{\nu(\sigma)} (\varphi(\sigma)). \]
It is easy to check that this makes $\Ind_\Delta^\Gamma \cN$ a pre-algebraicity condition for $(\tcE,\Ind_\Delta^\Gamma G)$. Combining lemma 12.10 of \cite{kotbg} with the observation of the last paragraph we get the following result.
\begin{lem}\label{shapirolem} In the above situation the composite
\[ H^1_{\Ind_\Delta^\Gamma \cN} (\tcE,\Ind_\Delta^\Gamma G) \stackrel{(i,\epsilon,\epsilon)^*}{\lra} H^1_\cN(\tcE|_\Delta,G)  \stackrel{\tepsilon_*}{\lra} H^1_\cN(\cE,G)\]
is an isomorphism. \end{lem}

\end{enumerate}

It will be convenient for us to work with pre-algebraicity conditions $\cN$ with slightly more structure. Namely we will assume that $\cN$ is endowed with a subset $\cN_\basic \subset \cN^{G}$, an abelian group structure on $\cN_\basic$ and an action of $\cN_\basic$ on $\cN$ extending the action of $\cN_\basic$ on itself by translation; such that the following properties hold:
\begin{itemize}
\item $\Gamma$ preserves $\cN_\basic$ and acts on it via group automorphisms.
\item The action of $G \rtimes \Gamma$ on $\cN$ commutes with the action of $\cN_\basic$, i.e. ${}^g (\nu \mu)=\nu \, {}^g \mu$ and ${}^\sigma(\nu \mu)={}^\sigma \nu \, {}^\sigma \mu$ for all $\nu \in \cN_\basic$, $\mu \in \cN$, $g \in G$ and $\sigma \in \Gamma$.
\item If $\nu \in \cN_\basic$, then $\barnu$ factors through $Z(G)$.
\item $\overline{\nu \mu}=\barnu \, \barmu$.
\end{itemize}
We will refer to this additional data as an {\em algebraicity condition}. Note that in this case both $\cN$ and $\cN_\basic$ are pre-algebraicity conditions. We will often write $Z^1_\cN(\cE,G)_\basic$ and $H^1_\cN(\cE,G)_\basic$, instead of $Z^1_{\cN_\basic}(\cE,G)$ and $H^1_{\cN_\basic}(\cE,G)$. We refer to these as the set of {\em basic algebraic cocycles} and the {\em basic algebraic cohomology}.

If $G$ is abelian then $Z^1_\cN(\cE,G)_\basic$ and $H^1_\cN(\cE,G)_\basic$ are naturally abelian groups.

There is a natural map
\[ \ad: Z^1_{\cN}(\cE,G)_\basic \lra Z^1(\Gamma,G^\ad) \]
which induces a map in cohomology. If $(\nu,\zeta)\in Z^1_{\cN_\basic}(\cE, Z(G))$ and $(\mu,\phi) \in Z^1_\cN(\cE,G)$ then $(\nu\mu,\eta \phi) \in Z^1_\cN(\cE,G)$. This induces maps
\[ (\nu,\zeta): Z^1_\cN(\cE,G)_\basic \lra Z^1_\cN(\cE,G)_\basic \]
and
\[ (\nu,\zeta): H^1_\cN(\cE,G) \lra H^1_\cN(\cE,G) \]
and
\[ (\nu,\zeta): H^1_\cN(\cE,G)_\basic \lra H^1_\cN(\cE,G)_\basic . \]
This gives actions of $Z^1_{\cN_\basic}(\cE,Z(G))$ on $Z^1_\cN(\cE,G)$ and $Z^1_\cN(\cE,G)_\basic$; and of $H^1_{\cN_\basic}(\cE,Z(G))$ on $H^1_\cN(\cE,G)$ and $H^1_\cN(\cE,G)_\basic$. The map $H^1_\cN(\cE,G) \ra H^1_\cN(\cE,G^\ad)$ is constant on $H^1_{\cN_\basic}(\cE,Z(G))$-orbits.

We have the following additions to our various functorialities:
\begin{enumerate}[(A)]
\item In the situation of \ref{changeg} if $n(\cN_{1,\basic}) \subset \cN_{2,\basic}$ then $(h,n)_*$ takes basic cocycles or cohomology classes to basic ones.

\item In the situation of \ref{changee} if $n(\cN_{2,\basic}) \subset \cN_{1,\basic}$ then $(f,h,n)^*$ takes basic cocycles or cohomology classes to basic ones.

\item In the situation of \ref{changec} if $n(\cN_{1,\basic}) \subset \cN_{2,\basic}$ then $(f,h,n)_*$ takes basic cocycles or cohomology classes to basic ones.

\item In the situation of \ref{shapiro} if $(\cN,\cN_\basic)$ is an algebraicity condition for $(\cE,G)$, then $(\Ind_\Delta^\Gamma \cN, \Ind_\Delta^\Gamma \cN_\basic)$ is one for $(\tcE,\Ind_\Delta^\Gamma G)$, and lemma \ref{shapirolem} is also true for the basic algebraic cohomology.

\item\label{changecor} Suppose that $G$ is abelian and that $\Delta \subset \Gamma$ is a subgroup of finite index.
Let $\cR \subset \cE$ be a set of representatives for $\cE|_\Delta \backslash \cE$. Then we obtain a natural map
\[ \begin{array}{rcl} \corr_\cR: Z^1_{\cN}(\cE|_\Delta,G)_\basic & \lra & Z^1_{\cN}(\cE,G)_\basic \\ (\nu,\phi) & \longmapsto &(\prod_{r \in \cR} {}^{r^{-1}}\nu, \tphi), \end{array} \]
where
\[ \tphi(e) = \prod_{r \in \cR} {}^{r^{-1}}( \phi (res^{-1})) \]
with each $s \in \cR$ chosen such that $res^{-1} \in \cE|_\Delta$. It induces a map
\[ \corr: H^1_{\cN}(\cE|_\Delta,G)_\basic  \lra  H^1_{\cN}(\cE,G)_\basic \]
which is independent of the choice of $\cR$.

If $\cE=\cE_\alpha$ with $\alpha \in Z^2(\Gamma,\cE^0)$ and $\Delta=\{1\}$, then we may take $\cR=\cR_\alpha=\{ e_\alpha(\sigma):\,\, \sigma \in \Gamma\}$ and we will write
\[ \corr_\alpha=\corr_{\cR_\alpha}: \cN_\basic \lra Z^1_\cN(\cE_\alpha,G)_\basic .\]
Note that
\[ (\corr_\alpha \nu)(e_\alpha(\sigma))=\prod_{\eta \in \Gamma} \eta^{-1} \barnu(\alpha(\eta,\sigma)). \]
If $\beta: \Gamma \ra \cE^0$ and if $i_\beta: \cE_\alpha \iso \cE_{{}^\beta \alpha}$ is the canonical isomorphism sending $e_\alpha(\sigma)$ to $\beta(\sigma)e_{{}^\beta \alpha}(\sigma)$ then it is easily verified that
\[ i_\beta^*(\corr_{{}^\beta \alpha} \nu)= {}^{\prod_{\eta \in \Gamma}\eta^{-1} \barnu(\beta(\eta))^{-1}} \corr_\alpha \nu .\]

(Indeed, both sides are of the form $(\prod_{\eta \in \Gamma}{}^\eta \nu,\phi)$ for some $\phi$. Moreover
\[  \begin{array}{rl} & (i_\beta^* \corr_{{}^\beta\alpha}(\nu))(e_{\alpha}(\sigma)) \\ =&
 \corr_{{}^\beta\alpha}(\nu)(\beta(\sigma) e_{{}^\beta\alpha}(\sigma)) \\ =&
  \prod_{\eta \in \Gamma}{}^\eta \barnu(\beta(\sigma)) \prod_{\eta \in \Gamma} \eta^{-1} \barnu({}^\beta\alpha(\eta,\sigma)) \\ =&
\prod_{\eta \in \Gamma}{}^\eta \barnu(\beta(\sigma)) \prod_{\eta \in \Gamma} \eta^{-1} \barnu(\alpha(\eta,\sigma) \beta(\eta\sigma)\beta(\eta)^{-1} {}^\eta \beta(\sigma)^{-1})\\  =&
\prod_{\eta \in \Gamma}{}^\eta \barnu(\beta(\sigma)) \corr_{\alpha}(\nu)(e_{\alpha}(\sigma))
\prod_{\eta \in \Gamma} \eta^{-1} \barnu( \beta(\eta\sigma))\\ &
\prod_{\eta \in \Gamma} \eta^{-1} \barnu(\beta(\eta))^{-1}
\prod_{\eta \in \Gamma} {}^{\eta^{-1}}\barnu(  \beta(\sigma))^{-1} \\ =&
\corr_{\alpha}(\nu)(e_{\alpha}(\sigma)) \sigma \prod_{\eta \in \Gamma} \eta^{-1} \barnu(\beta(\eta))/\prod_{\eta \in \Gamma} \eta^{-1} \barnu(\beta(\eta))
.) \end{array} \]

\end{enumerate}

Now consider the case that $\Gamma=\Gal(E/F)$ and $G=H(A_E)$, where 
\begin{itemize}
\item $E/F$ is a finite Galois extension of fields;
\item $H/F$ is an algebraic group;
\item and $A_E= A \otimes_FE$ for some $F$-algebra $A$.
\end{itemize}
If $(\nu,\phi) \in Z^1_{\cN}(\cE,G(A_E))_\basic$ then we define ${}^\phi G/A$ to be the etale descent of $G \times A_E$ to $G \times A$ via the action 
\[ \sigma \longmapsto \conju_{\phi(e)} \circ \sigma ,\]
where $e\in \cE$ is any lift of $\sigma$.
Thus ${}^\phi G \times_{A_F} A_E = G \times_F A_E$. If $g \in G(A_E)$ then 
\[ \conju_g: {}^\phi G \liso {}^{{}^g\phi} G. \]
Thus ${}^\phi G$ depends only on $[(\nu,\phi)]$ up to an isomorphism that is unique {\em up to composition conjugation by an element of $({}^\phi G)(A)$}. (Note that we have ${}^\phi G(A)$ here and not ${}^\phi G^\ad (A)$. This is an important point.) When we are only concerned with properties of ${}^\phi G$ for which this ambiguity does not matter, we may write ${}^{[(\nu,\phi)]}G$. There is a bijection
\[ \begin{array}{rcl} Z^1_{\cN}(\cE,{}^\phi G(A_E)) &\liso& Z^1_\cN(\cE,G(A_E)) \\ (\mu,\psi) & \longmapsto (\mu\nu,\psi \phi) \end{array} \]
which takes basic subset to basic subset, and induces isomorphisms in cohomology. Note that 
\[ {}^{gh}((\mu,\psi)(\nu,\phi))=\conju_g ({}^h(\mu,\psi)) {}^g(\nu,\phi). \]
If $\bphi \in H^1_{\cN}(\cE,G(A_E))_\basic$ then 
\[ \{ \bzeta \in H^1_{\cN_\basic}(\cE,Z(G)(A_E)):\,\, \bzeta \bphi = \bphi \}= \ker (H^1_{\cN_\basic}(\cE,Z(G)(A_E)) \ra H^1_\cN(\cE,{}^\bphi G(A_E))) \]
is a subgroup of $H^1_{\cN_\basic}(\cE,Z(G)(A_E))$. 

\newpage

\section{Kottwitz cohomology: the local case}\label{bglocal}

\subsection{Kottwitz cohomology for local Weil groups}

Suppose that $F$ is a local field of charcateristic $0$ and $D \supset E \supset F$ are finite Galois extensions. 
If $\alpha \in [\alpha_{E/F,D}]\subset Z^2(\Gal(D/F),D^\times)$, then we get 
a well defined extension
\[ 0 \lra D^\times \lra W_{E/F,D,\alpha} \lra \Gal(D/F) \lra 0. \]
In the case $D=E$ we will write simply $W_{E/F,\alpha}$.
As $H^1(\Gal(D/F),D^\times)=(0)$ the only automorphisms of this extension are conjugation by an element of $D^\times$. If $\alpha'$ is a second element of $[\alpha_{E/F,D}]$, then there is an isomorphism of extensions
\[ \begin{array}{ccccccccc}
0 &\lra& D^\times &\lra& W_{E/F,D,\alpha} &\lra& \Gal(D/F) &\lra& 0 \\ && || && \da \wr && || && \\
0 &\lra& D^\times &\lra& W_{E/F,D,\alpha'} &\lra& \Gal(D/F) &\lra& 0. \end{array} \]
However it is only unique up to composition with conjugation by an element of $D^\times$. In particular, if $\barF$ is an algebraic closure of $F$ containing $D$, then the extension $W_{E/F,D,\alpha}$ is isomorphic to $W_{E^\ab/F,D}$, but this isomorphism is only unique up to composition with conjugation by an element of $E^\times$.

Let $G/F$ denote an algebraic group. We will consider the algebraicity conditions $\cN=X_*(G)(D)$ and $\cN_\basic=X_*(Z(G))(D)$. We will denote the corresponding algebraic cocycles, basic cocycles, cohomology and basic cohomology as $Z^1_\alg(W_{E/F,D,\alpha},G(D))$, $Z^1_\alg(W_{E/F,D,\alpha},G(D))_\basic$, $H^1_\alg(W_{E/F,D},G(D))$, and $H^1_\alg(W_{E/F,D},G(D))_\basic$ respectively. As the notation suggest the two cohomology groups are canonically independent of the choice of $\alpha \in [\alpha_{E/F,D}]$.

We will call $\bphi\in H^1_\alg(W_{E/F,D},G(D))_\basic$ {\em compact} if ${}^\bphi G^\ad(F)$ is compact, and we will write $H^1_\alg(W_{E/F,D},G(D))_\basic^\compact$ for the set of compact elements in $H^1_\alg(W_{E/F,D},G(D))_\basic$.

Choose representatives $\alpha_{E/F} \in [\alpha_{E/F}]$ and $\alpha_{D/F} \in [\alpha_{D/F}]$ and $\gamma_{D/E}: \Gal(D/F) \ra D^\times$ such that ${}^{\gamma_{D/E}}\alpha_{D/F}^{[D:E]}=\inf_{D/E} \alpha_{E/F} \in Z^2(\Gal(D/F),D^\times)$.
Then there is a commutative diagram
\[ \begin{array}{ccccccccc} (0)& \lra & D^\times & \lra & W_{D/F,\alpha_{D/F}} & \lra & \Gal(D/F) & \lra & (0) \\ 
&& [D:E] \da && \eta_{D/E,\gamma_{D/E}}\da && || && \\
(0)& \lra & D^\times & \lra & W_{E/F,D,\inf_{D/E} \alpha_{E/F}} & \lra & \Gal(D/F) & \lra & (0) \\
&& \ua && \ua  && || && \\
(0)& \lra & E^\times & \lra & W_{E/F,\alpha_{E/F}}|_{\Gal(D/F)} & \lra & \Gal(D/F) & \lra & (0) \\
&& || && \da && \da &&  \\
(0)& \lra & E^\times & \lra & W_{E/F, \alpha_{E/F}} & \lra & \Gal(E/F) & \lra & (0). \end{array} \]
Using successively functorialities (B) then (C) then (B) again from the end of section \ref{algcoh}, 
we obtain a map
\[ \begin{array}{rcl}  \inf_{D/E,\gamma_{D/E}}: Z^1_\alg(W_{E/F,\alpha_{E/F}},G(E)) &\lra& Z^1_\alg(W_{E/F,\alpha_{E/F}}|_{\Gal(D/F)},G(D)) \\ &\lra& Z^1_\alg(W_{E/F,D, \inf_{D/E} \alpha_{E/F}},G(D)) \\ &\lra& Z^1_\alg(W_{D/F,\alpha_{D/F}},G(D)). \end{array} \]
(Where we use the algebraicity conditions $\cN=X_*(G)(D)$ and $\cN_\basic = X_*(Z(G))(D)$ for the second set of cocycles.) All these maps are functorial in $G$, take basic elements to basic elements and commute with products.
The composite sends $(\nu,\phi)$ to $(\nu^{[D:E]},\tphi)$, where, if $\eta_{D/E,\gamma_{D/E}}(e)=de'$ with $d\in D^\times$ and $e' \in W_{E/F,\alpha_{E/F}}|_{\Gal(D/F)}$ then
\[ \tphi(e)=\nu(d)\phi(\bare') \]
where $\bare'$ denotes the image of $e'$ in $W_{E/F,\alpha_{E/F}}$. This composite is injective. The map $\gamma_{E/D}$ can only be replaced by 
\[ {}^d \gamma_{D/E}: \sigma \mapsto \gamma_{D/E}(\sigma) d/{}^\sigma d, \]
for some $d \in D^\times$. (As $H^1(\Gal(D/F),D^\times)=(0)$.) We 
have 
\[ \eta_{D/E,{}^d\gamma_{D/E}}= \conju_{d^{-1}} \circ \eta_{D/E,\gamma_{D/E}}, \]
and so the induced map
\[ {\inf}_{D/E,\gamma_{D/E}}: H^1_\alg(W_{E/F},G(E)) \lra  H^1_\alg(W_{D/F},G(D)) \]
is independent of the choice of $\gamma_{D/E}$ and so we will denote it simply $\inf_{D/E}$.
The maps
\[ \begin{array}{rcl}  \inf_{D/E}: H^1_\alg(W_{E/F,\alpha_{E/F}},G(E)) &\lra& H^1_\alg(W_{E/F,\alpha_{E/F}}|_{\Gal(D/F)},G(D)) \\ &\lra & H^1_\alg(W_{E/F,D,\inf_{D/E} \alpha_{E/F}},G(D)) \\ & \lra & H^1_\alg(W_{D/F,\alpha_{D/F}},G(D)) \end{array} \]
are all injective. (The first because the usual inflation map is injective on $H^1$, and the second and third immediately from the definitions.)

 Kottwitz defines 
\[ B(F,G) = \lim_{\ra E} H^1_\alg(W_{E/F},G(E)) \] 
and
\[ B(F,G)_\basic = \lim_{\ra E} H^1_\alg(W_{E/F},G(E))_\basic . \] 

If $\bphi_1, \bphi_2 \in H^1_\alg(W_{E/F},G(E))$ have the same image in $H^1(W_{E/F},G^\ad(E))$, then we can find a finite extension $D/E$ Galois over $F$ such that $\inf \bphi_i \in H^1_\alg(W_{D/F},G(D))$ can be represented by  cocycles $\phi_i$ with $\ad \phi_1 = \ad \phi_2$. (If $\phi_i'$ is a representative of $\bphi_i$ and if $g \in G^\ad(E)$ with ${}^g \ad \phi_1' =\ad \phi_2'$, then we may choose such a field $D$ and a $\tg \in G(D)$ lifting $g$. Then ${}^{\tg} \inf \phi_1'$ and $\inf \phi_2'$ will do.) Thus $\inf \bphi_1,\inf \bphi_2 \in H^1_\alg(W_{D/F},G(D))$ differ by an element of $H^1_\alg(W_{D/F},Z(G)(D))$.

\subsection{Reductive groups}

If $T/F$ is a torus split by a finite Galois extension $E/F$, then Kottwitz (see formula (9.1) of \cite{kotbg}) shows that 
\[ \corr: X_*(T)_{\Gal(E/F)} \liso H^1_\alg(W_{E/F}, T). \]
Moreover if $G$ is any reductive group which splits over $E$, then Kottwitz constructs a map
\[ \kappa: H^1_\alg ( W_{E/F},G(E)) \lra \Lambda_{G,\Gal(E/F)} \]
with the following properties
\begin{itemize}
\item it is functorial in $G$;
\item if $G=T$ is a torus, then $\kappa=\corr^{-1}$ is an isomorphism;
\item if $D/E$ is a finite extension Galois over $F$, then
\[ \begin{array}{ccc} H^1_\alg(W_{E/F},G(E)) & \stackrel{\kappa}{\lra} & \Lambda_{G,\Gal(\barF/F)} \\ \inf_{D/E} \da & \kappa \nearrow & \\ H^1_\alg(W_{D/F},G(D))  && \end{array} \]
commutes.
\end{itemize}
 (See sections 9 and 11 of \cite{kotbg}.) In the limit Kottwitz obtains a map
\[ \kappa: B(F,G) \lra \Lambda_{G,\Gal(\barF/F)}. \]
 
 \begin{lem} If $\bphi \in H^1_\alg (W_{E/F},G)_\basic$ and $\bpsi \in H^1_\alg(W_{E/F},{}^\bphi G)$ then
 \[ \kappa(\bpsi \bphi)=\kappa(\bpsi) \kappa(\bphi). \] \end{lem}
 
 \pfbegin This follows easily from the construction of $\kappa$ in section 9.3 of \cite{kotbg}: If $G$ is a torus it simply expresses the fact that $\corr$ is an abelian group  homomorphism in this case. If $G^\der$ is simply connected it follows for the corresponding fact for $C(G)$. In the general case it follows from the corresponding fact for a suitable $z$-extension of $G$. \pfend

If $T \subset G$ is a maximal torus then we define $H^1_\alg(W_{E/F},T(E))_{G-\basic}$ to be those elements $[(\nu,\phi)]$ where $\nu$ factors through $Z(G)$. If $F=\R$ suppose that $T$ is fundamental, while if $F$ is p-adic assume that $T$ is elliptic. (Recall that in the case $F=\R$ `fundamental' means that its split rank is minimal among those of all maximal tori defined over $\R$. All fundamental maximal tori are $G(\R)$-conjugate - see \cite{bw} section I.7.1.) If $E$ splits $T$ then   
\[ H^1_\alg(W_{E/F},T(E))_{G-\basic} \onto H^1_\alg(W_{E/F},G(E))_{\basic}. \]
(See proposition 13.1 and lemma 13.2 of \cite{kotbg}.) We deduce that if $E$ splits an elliptic (in the p-adic case) or fundamental (in the real case) torus, then 
\[ H^1_\alg(W_{E/F},G(E))_\basic \iso B(F,G)_\basic.\]
We further deduce that for $E$ sufficiently large the quotient of $H^1_\alg(W_{E/F},G(E))_\basic$ by its action of $H^1_\alg(W_{E/F},Z(G)(E))$ embeds into $H^1_\alg(W_{E/F},G^\ad(E))_\basic$. 
(First choose $E_0$ such that $H^1_\alg(W_{E/F},G(E))_\basic \iso B(F,G)_\basic$ for any $E \supset E_0$, and then  $E \supset E_0$ such that every element of $G^\ad(E_0)$ has a lift in $G(E)$.)

If $F$ is a p-adic field then $\kappa$ is in fact a bijection
\[ \kappa: B(F,G)_\basic \liso \Lambda_{G,\Gal(\barF/F)}. \]
and so $B(F,G)_\basic$ becomes an abelian group. (See proposition 13.1 of \cite{kotbg}.) We deduce that if $E$ splits some maximal torus of $G$ defined over $F$, then 
\[ \kappa: H^1_\alg(W_{E/F},G(E))_\basic \liso \Lambda_{G,\Gal(\barF/F)}.\]
In particular if $G$ is semi-simple, simply connected, then $H^1_\alg(W_{E/F},G(E))_\basic=\{1\}$. 
If $G$ is semisimple then 
\[ \kappa:H^1(F,G) \liso B(F,G)_\basic \liso \Lambda_{G,\Gal(\barF/F)}. \]

If $F=\C$ then $W_{\C/\C}=\C^\times$, $H^1(\C,G)=(0)$ and $B(\C,G)_\basic = X_*(Z(G))$ and $\kappa: B(\C,G)_\basic =X_*(Z(G)) \into \Lambda_G$.

\subsection{The real case} \label{real}

Now suppose that $F=\R$. Choose a representative $\alpha^0_{\C/\R}$ for $[\alpha_{\C/\R}]$ defined by
\[ \alpha^0_{\C/\R}(\sigma_1,\sigma_2)= \left\{ \begin{array}{ll} -1 & {\rm if}\,\, \sigma_1=\sigma_2=c \\ 1 & {\rm otherwise.} \end{array} \right. \]
Then 
\[ W_{\C/\R,\alpha^0_{\C/\R}} \cong \langle \C^\times, j:\,\, j^2=-1\,\, {\rm and} \,\, jzj^{-1}={}^cz\rangle, \]
with $e(1)=1$ and $e(c)=j$.
Thus an element of $Z^1_\alg(W_{\C/\R, \alpha^0_{\C/\R}},G(\C))$ is a pair $(\nu,J)$ where $\nu\in X_*(G)$ and $J \in G(\C)$ satisfy
\[ \nu = J {}^c\nu J^{-1} \]
and
\[ J{}^cJ=\nu(-1). \]
Moreover $[(\nu,J)]=[(\nu',J')] \in H^1_\alg(W_{\C/\R},G(\C))$ if there exists $g \in G(\C)$ such that 
\[ \nu'=g\nu g^{-1} \]
and
\[ J'=gJ{}^cg^{-1}. \]

If $\mu \in X_*(G)$ and the image of $\mu$ commutes with that of ${}^c\mu$, then we obtain an element $\widehat{\lambda}_{G}(\mu) \in Z^1_\alg(W_{\C/\R,\alpha^0_{\C/\R}},G(\C))$ defined by 
\[ \widehat{\lambda}_{G}(\mu) = (\mu {}^c\mu, \mu(-1)). \]
Note that $\kappa( \hatlambda_G(\mu))$ equals the image of $\mu$ in $\Lambda_{G,\Gal(\C/\R)}$. (To see this note that it suffices to treat the case that $G=T$ is a torus, in which case $\hatlambda_G(\mu)=\corr_{\alpha_{\C/\R}^0} \mu$.)
If $\mu{}^c \mu$ factors through $Z(G)$, then $\hatlambda_G(\mu)\in Z^1_\alg(W_{\C/\R,\alpha^0_{\C/\R}},G(\C))_\basic$. 
This induces maps
\[ \bhatlambda_G: G(\R) \backslash \{\mu \in X_*(G): \,\, {\rm the \,\, images\,\, of}\,\, \mu \,\, {\rm and}\,\,{}^c\mu\,\, {\rm commute}\} \lra H^1_\alg(W_{\C/\R},G(\C)) \]
and
\[ \bhatlambda_G:G(\R)\backslash \{\mu \in X_*(G): \,\, {}^c\mu=\mu^{-1} \in X_*(G^\ad) \} \onto H^1_\alg(W_{\C/\R},G(\C))_\basic . \]
(To see the latter map is surjective choose a fundamental torus $T \subset G$. Given a class in $H^1_\alg(W_{\C/\R},G(\C))_\basic$ it will be the image of a class in $H^1_\alg(W_{\C/\R},T(\C))_{G-\basic}$, which in turn is of the form $\corr_{\alpha^0_{\C/\R}} \mu$ for $\mu \in X_*(T)$. Then the class is represented by $(\mu{}^c\mu,\mu(-1))$ and $\mu{}^c\mu$ is valued in $Z(G)$.)
The image $\bhatlambda_G(\mu)$ depends only on the $G(\R)$-conjugacy class of $\mu$, so we will sometimes write $\bhatlambda_G([\mu]_{G(\R)})$. 

\newpage

\section{Some extensions of global Galois groups}\label{bgglob1}

In the case of global Kottwitz cohomology we will need to consider various different extensions of a global Galois group $\Gal(E/F)$ and the relationships between them. For now consider $D \supset E \supset F$ with $D$ and $E$ Galois over $F$. We will be primarily interested in the case $D=E$, but we will need the more general case to allow comparisons as the extension $E$ varies. {\em When we drop $D$ from the notation we imply that we are in the case $D=E$.}

\subsection{Some pro-tori and related groups}

Suppose that $D \supset E \supset F$ are number fields with $D$ and $E$ Galois over $F$. If $t\in T_{2,D}(\A_F)$ then for each place $u$ of $D$ we have $t_u\in (\A_D^\times)^{\Gal(D/F)_u}$ and $t_{\sigma u}=\sigma t_u$. Moreover if $w$ is the place of $E$ below $u$ then
\[ \begin{array}{rcl} \eta_{D/E}(t)_w&=& \prod_{u'|w} \prod_{\sigma \in \Gal(
D/E)_{u'}} {}^{\sigma} t_{u'} \\ &=& \prod_{\tau \in \Gal(D/E)/\Gal(D/E)_u} \prod_{\sigma \in \Gal(D/E)_u} {}^{\tau \sigma \tau^{-1}} t_{\tau u} \\ &=& \prod_{\tau \in \Gal(D/E)/\Gal(D/E)_u} \prod_{\sigma \in \Gal(D/E)_u} {}^{\tau \sigma } t_{u} \\ &=& N_{D/E}t_u. \end{array}\]

\begin{lem}\label{cohtori} Suppose that $D \supset E \supset F$ are number fields with $D$ and $E$ Galois over $F$. For each place $v$ of $F$ choose places $u(v)|w(v)|v$ of $D$ and $E$ respectively.
\begin{enumerate}
\item\label{ct1} $H^i(\Gal(D/F),T_{2,E}(D)) \cong \prod_{v \in V_F} H^i(\Gal(D/F)_{w(v)},D^\times)$ and \newline $H^1(\Gal(D/F),T_{2,E}(D))=(0)$.

\item\label{ct2} $H^i(\Gal(D/F),T_{2,E}(\A_D)) \cong \prod_{v \in V_F} H^i(\Gal(D/F)_{w(v)},\A_D^\times)$ and \newline $H^1(\Gal(D/F),T_{2,E}(\A_D))=(0)$.

\item\label{ct3} $H^i(\Gal(D/F),T_{2,E}(\A_D)/T_{2,E}(D)) \cong \prod_{v \in V_F} H^i(\Gal(D/F)_{w(v)},\A_D^\times/D^\times)$ and $H^1(\Gal(D/F),T_{2,E}(\A_D)/T_{2,E}(D))=(0)$.

\item\label{ct4} $H^1(\Gal(D/F),T_{3,E}(D))=(0)$.
\end{enumerate} \end{lem}

\pfbegin The first part is an application of Shapiro's lemma
\[ H^i(\Gal(D/F),T_{2,E}(D)) \cong \prod_{v \in V_F} H^i(\Gal(D/F),\prod_{w\in V_{E,v}} D^\times) \cong \prod_{v \in V_F} H^i(\Gal(D/F)_{w(v)},D^\times) \]
combined with Hilbert's theorem 90.
The second and third parts are proved similarly using the vanishing of $H^1(\Gal(D/E)_{w(v)}, \A_D^\times)$ and $H^1(\Gal(D/E)_{w(v)}, \A_D^\times/D^\times)$.

Consider the fourth part. As $D$ splits $T_{2,E}$ there is an exact sequence
\[ (0) \lra D^\times \lra T_{2,E}(D) \lra T_{3,E}(D) \lra (0), \]
and so it suffices to show that 
\[ H^2(\Gal(D/F),D^\times) \lra H^2(\Gal(D/F), T_{2,E}(D)) \cong \prod_{v \in V_F} H^2(\Gal(D/F)_{w(v)},D^\times) \]
is injective. In fact it suffices to show that the composite with the map
\[ \prod_{v \in V_F} H^2(\Gal(D/F)_{w(v)},D^\times) \lra \prod_{v \in V_F} H^2(\Gal(D/F)_{u(v)},D_{u(w)}^\times) \]
is injective. However this injectivity follows from the fact that the Brauer group of $F$ (of which $H^2(\Gal(D/F),D^\times)$ is a subgroup) injects into the product of the Brauer groups of all completions of $F$. \pfend

\begin{lem}\label{surje} If $D\supset E \supset F$ are finite Galois extensions of $F$ and if $A$ is any finitely generated ablelian group (with a trivial action of $\Gal(D/E)$), then
\[ \iota^{0,*}_{D/E}: (\Z[V_D] \otimes_\Z A)_{\Gal(D/E)} \liso \Z[V_E]\otimes_\Z A \]
and
\[ \iota^{0,*}_{D/E}: (\Z[V_D]_0 \otimes_\Z A)_{\Gal(D/E)} \liso \Z[V_E]_0 \otimes_\Z A. \]
\end{lem}

\pfbegin The first assertion is clear. For the second use the long exact sequence associated to 
\[ (0) \ra \Z[V_D]_0 \otimes_\Z A\ra \Z[V_D] \otimes_\Z A \ra A \ra (0) \]
and the surjectivity of $\bigoplus_{w \in V_E} \Gal(D_{u(w)}/E_w)^\ab \otimes_\Z A\ra \Gal(D/E)^\ab \otimes_\Z A$, where for each $w\in V_E$ we choose a place $u(w)$ of $D$ above it.
\pfend

If $D\supset E \supset F$ are finite Galos extensions of $F$ and if $S$ is a set of places of $F$, then we will write 
\[ \cE^\loc(E/F)_{D,S}^0 = \prod_{w \in V_{E,S}} D_w^\times \subset T_{2,E}(\A_D), \]
where we consider $D_w^\times \subset \A_D^\times$ inside the $w^{th}$-factor. 
 If $C \supset D$ is also finite Galois over $F$ then $\cE^\loc(E/F)_{C,S}$ is preserved by $\Gal(C/F)$ and
\[ (\cE^\loc(E/F)_{C,S}^0)^{\Gal(C/D)}=\cE^\loc(E/F)^0_{D,S}. \]
The map $\eta_{D/E}^0:T_{2,D} \ra T_{2,E}$ sends
\[ \begin{array}{rcl} \eta_{D/E}^0: \cE^\loc(D/F)_{C,S}^0 &\lra& \cE^\loc(E/F)_{C,S}^0 \\
(x_u)_{u \in V_{D,S}} & \longmapsto & (\prod_{u|w} x_u^{[D_u:E_w]})_{w\in V_{E,S}}. \end{array} \]

We also define
 \[ \begin{array}{rcl} \cE^\glob(E/F)_D^0& =& T_{2,E}(\A_D) \times_{T_{2,E}(\A_D)/T_{2,E}(D)} \A_D^\times/D^\times \\ &=& \{ a \in T_{2,E}(\A_D):\,\, \pi_w(a) \bmod E^\times \,\, {\rm is\,\,  independent\,\, of}\,\, w\} \\ & \subset & T_{2,E}(\A_D). \end{array} \]
 It is preserved by $\Gal(D/F)$ and 
 \[ (\cE^\glob(E/F)_D^0)^{\Gal(D/E)} =\cE^\glob(E/F)_E^0 \]
 which we are denoting simply $\cE^\glob(E/F)^0$).
 Moreover we have short exact sequences
 \begin{equation}\label{ses1} (0) \lra \A_D^\times \lra \cE^\glob(E/F)_D^0 \lra T_{3,E}(D) \lra (0), \end{equation}
and
 \begin{equation}\label{ses2}  (0) \lra \cE^\glob(E/F)_D^0 \lra T_{2,E}(\A_D) \times \A_D^\times/D^\times \lra T_{2,E}(\A_D)/T_{2,E}(D) \lra (0) \end{equation}
 and
  \begin{equation}\label{ses3}  (0) \lra T_{2,E}(D) \lra \cE^\glob(E/F)_D^0 \lra \A_D^\times/D^\times \lra  (0) \end{equation}
 Moreover, if $C \supset D$ in another finite Galois extension of $F$, then
 \[ \eta_{D/E}^0: \cE^\glob(D/F)_C^0 \lra \cE^\glob(E/F)_C^0 .\]
 This is compatible with the $[D:E]$-power map $\A_D^\times/D^\times$ to itself.

\begin{lem} Suppose that $D\supset E \supset F$ are finite Galos extensions of $F$ and that $S$ is a set of places of $F$. For each place $v$ of $F$ choose places $u(v)|w(v)|v$ of $D$ and $E$ respectively. 
\begin{enumerate}
\item $H^i(\Gal(D/F), \cE^\loc(E/F)_{D,S}^0) \cong \prod_{v \in S} H^i(\Gal(D/F)_{u(v)}, D_{u(v)}^\times)$. 

\item $H^1(\Gal(D/F), \cE^\loc(E/F)_{D,S}^0)=(0)$.

\item There is a unique class $[\alpha_{E/F,D,S}^\loc]\in H^2(\Gal(D/F), \cE^\loc(E/F)_{D,S}^0)$ which restricts to $[\alpha_{E_{w(v)}/F_v,D_{u(v)}}] \in H^2(\Gal
(D/F)_{u(v)},D_{u(v)}^\times)$ for all $v \in S$. This class does not depend on the choices of the $w(v)$ and $u(v)$. If $S'\subset S$ then the image of $[\alpha_{E/F,D,S}^\loc]$ in $H^2(\Gal(D/F), \cE^\loc(E/F)_{D,S}^0)$ equals $[\alpha_{E/F,D,S'}^\loc]$.

\item $H^1(\Gal(D/F), \cE^\glob(E/F)_D^0)=(0)$ and $H^1(\Gal(D/E), \cE^\glob(E/F)_D^0)=(0)$.

\item There is a left exact sequence
\[ \begin{array}{r} (0) \lra H^2(\Gal(D/F), \cE^\glob(E/F)_D^0) \lra H^2(\Gal(D/F), T_{2,E}(\A_D)) \oplus H^2(\Gal(D/F),\A_D^\times/D^\times)\\  \lra \prod_{v \in V_F} H^2(\Gal(D/F)_{w(v)},\A_D^\times/D^\times); \end{array}\]
and a unique class $[\alpha_{E/F,D}^\glob]\in H^2(\Gal(D/F),\cE^\glob(E/F)_D^0)$ which maps to 
\[ \begin{array}{rl} & ((\corr_{\Gal(D/F)_{u(v)}}^{\Gal(D/F)_{w(v)}} \imath_{u(v),*}[\alpha_{E_{w(v)}/F_v,D_{u(v)}}]),[\alpha^W_{E/F,D}]) \\ \in &\left(\prod_{v\in V_F} H^2(\Gal(D/F)_{w(v)}, \A_D^\times)\right) \oplus H^2(\Gal(D/E),\A_D^\times/D^\times). \end{array} \]
This class does not depend on the choices of the $w(v)$ and $u(v)$. 

Moreover $[\alpha_{E/F,D}^\loc]$ and $[\alpha_{E/F,D}^\glob]$ have the same image in $H^2(\Gal(D/F), T_{2,E}(\A_D))$, which we will denote $[\alpha_{2,E/F,D}]$.

\item If $C\supset D\supset E \supset F$ are finite Galois extensions of $F$, then 
\[ \eta_{D/E,*}[\alpha^\loc_{D/F,C,S}]=[\alpha^\loc_{E/F,C,S}] \]
and
\[ \eta_{D/E,*}[\alpha^\glob_{D/F,C}]=[\alpha^\glob_{E/F,C}]. \]

\item  If $C\supset D\supset E \supset F$ are finite Galois extensions of $F$, then 
\[ \inf_{C/D}[\alpha^\loc_{E/F,D,S}]=[\alpha^\loc_{E/F,C,S}] \]
and
\[ \inf_{C/D}[\alpha^\glob_{E/F,D}]=[\alpha^\glob_{E/F,C}]. \]

\end{enumerate}
\end{lem}

\pfbegin The first part is an application of Shapiro's lemma. The second and third parts follow easily. The fourth part follows from part \ref{ct4} of lemma \ref{cohtori} and the long exact sequence associated to (\ref{ses1}). 

The left exact sequence in part five comes from the long exact sequence associated to (\ref{ses2}) and part \ref{ct3} of lemma \ref{cohtori}. The existence of $[\alpha^\glob_{E/F,D}]$ follows because 
\[ \begin{array}{rl} &\corr_{\Gal(D/F)_{u(v)}}^{\Gal(D/F)_{w(v)}} \imath_{u(v),*}[\alpha_{E_{w(v)}/F_v,D_{u(v)}}] \\ = &
\corr_{\Gal(D/F)_{u(v)}}^{\Gal(D/F)_{w(v)}} [D_{u(v)}:E_{w(v)}] [\alpha^W_{D/D^{\Gal(D/F)_{u(v)}}}] \\ =  &
[D^{\Gal(D/F)_{u(v)}}:D^{\Gal(D/F)_{w(v)}}][D_{u(v)}:E_{w(v)}] [\alpha^W_{D/D^{\Gal(D/F)_{w(v)}}}] \\ =&
[D:E][\alpha^W_{D/D^{\Gal(D/F)_{w(v)}}}] \\ =&
[D:E] \res^{\Gal(D/F)}_{\Gal(D/F)_{w(v)}} [\alpha^W_{D/F}] \\ =&
\res^{\Gal(D/F)}_{\Gal(D/F)_{w(v)}} [\alpha^W_{E/F,D}].
\end{array} \]
The images of  $[\alpha_{E/F,D}^\loc]$ and $[\alpha_{E/F,D}^\glob]$ in 
\[ H^2(\Gal(D/F), T_{2,E}(\A_D)) \cong \prod_{v\in V_F} H^2(\Gal(D/F)_{w(v)}, \A_D^\times) \]
are
\[ (\imath_{w(v),*} \corr_{\Gal(D/F)_{u(v)}}^{\Gal(D/F)_{w(v)}}[\alpha_{E_{w(v)}/F_v,D_{u(v)}}])_v=( \corr_{\Gal(D/F)_{u(v)}}^{\Gal(D/F)_{w(v)}}\imath_{u(v),*}[\alpha_{E_{w(v)}/F_v,D_{u(v)}}])_v. \]

For the sixth part note that, if $t(v)$ is a place of $C$ above $u(v)$, then
\[ [\alpha_{E_{w(v)}/F_v,C_{t(v)}}]= [C_{t(v)}:E_{w(v)}][\alpha_{C_{t(v)}/F_v}] = [D_{u(v)}:E_{w(v)}] [\alpha_{D_{u(v)}/F_v,C_{t(v)}}], \]
and
\[ [\alpha_{E/F,C}^W]=[D:E][\alpha_{D/F,C}^W]. \]

For the seventh part note that
\[ \res^{\Gal(C/F)_{u(v)}}_{\Gal(C/F)_{t(v)}} \inf_{\Gal(D/F)_{u(v)}}^{\Gal(C/F)_{u(v)}} [\alpha_{E_{w(v)}/F_v,D_{u(v)}}] =\inf_{\Gal(D/F)_{u(v)}}^{\Gal(C/F)_{t(v)}} [\alpha_{E_{w(v)}/F_v,D_{u(v)}}] = [\alpha_{E_{w(v)}/F_v,C_{t(v)}}] \]
and
\[ \inf_{C/D} [\alpha^W_{E/F,D}]=[\alpha^W_{E/F,C}]. \]
\pfend

We will write $\alpha_{3,E/F,D}$ for the image of $\alpha^\glob_{E/F,D}$ in $H^2(\Gal(D/F),T_{3,E}(D))$.

 \subsection{Some extensions}\label{somee}

We will write $\cZ(E/F)_D$ for the group of 3-tuples
\[ \balpha=(\alpha^\glob,\alpha^\loc,\beta) \in [\alpha_{E/F,D}^\loc] \times [\alpha_{E/F,D}^\glob] \times C^1(\Gal(D/F),T_{2,E}(\A_D)) \times  C^1(\Gal(D/E),\cE^\glob(E/F)_D^0)\]
such that 
\[ {}^\beta \alpha^\loc = \alpha^\glob \in Z^2(\Gal(D/F),T_{2,E}(\A_D)). \]

We will also write $\cB(E/F)_D$ for the abelian group consisting of pairs $\bgamma=(\gamma^\glob,\gamma^\loc)$ where $\gamma^\glob: \Gal(D/F) \ra \cE^\glob(E/F)_D^0$ and $\gamma^\loc: \Gal(D/F) \ra \cE(E/F)_D^\loc$; with pointwise multiplication. 
Then $\cB(E/F)_D$ acts on $\cZ(E/F)_D$ via 
\[ {}^{(\gamma^\glob,\gamma^\loc)} (\alpha^\glob,\alpha^\loc, \beta)= ({}^{\gamma^\glob}\alpha^\glob, {}^{\gamma^\loc}\alpha^\loc, \gamma^\glob \beta (\gamma^\loc)^{-1}). \]
We write $\cH(E/F)_D$ for the set of orbits of $\cB(E/F)_D$ on $\cZ(E/F)_D$.  We will call two elements of $\cZ(E/F)_D$ in the same orbit {\em equivalent}.

Moreover $T_{2,E}(\A_D)$ acts on $\cZ(E/F)_D$ via
\[ {}^t (\alpha^\glob,\alpha^\loc, \beta) = (\alpha^\glob,\alpha^\loc, {}^t\beta). \]
This action commutes with the action of $\cB(E/F)_D$, and the induced action of $\cB(E/F)_D \times T_{2,E}(\A_E)_D$ is transitive (because $H^1(\Gal(D/F),T_{2,E}(\A_D))=(0)$). The stabilizer in $\cB(E/F)_D \times T_{2,E}(\A_D)$ of any $\balpha$ is the group of
\[ ( ({}^a 1, {}^b1), bc/a)  \]
with $a \in \cE^\glob(E/F)_D^0$ and $b \in \cE^\loc(E/F)_D^0$ and $c \in T_{2,E}(\A_F)$. (This is because $H^1(\Gal(D/F), \cE^\loc(E/F)_D^0)=(0)$ and $H^1(\Gal(D/F),
 \cE^\glob(E/F)_D^0)=(0)$.)

 \begin{lem} \label{eglob} \begin{enumerate}
 \item If $w_1\neq w_2$ are places of $E$, then the intersection of the images of $D_{w_1}^\times$ and $D_{w_2}^\times$ in $\A_D^\times/D^\times$ is trivial. Thus
 $\cE^\glob(E/F)_D^0 \cap \cE^\loc(E/F)_D^0 = \{ 1\}$.
 
 \item If $t \in \cE^\glob(E/F)_D^0$ and $s \in \cE^\loc(E/F)_D^0$ and $st \in T_{2,E}(\A_F)$, then $s$ and $t \in T_{2,E}(\A_F)$.
\end{enumerate} \end{lem}
 
 \pfbegin 
 For the first part, if $t_i \in D_{w_i}^\times$ have the same image in $\A_D^\times/D^\times$ for $i=1,2$, then $t_1/t_2 \in D^\times \cap D_{w_1}^\times D_{w_2}^\times = \{1\}$ and so $t_1=t_2=1$. Then $\cE^\glob(E/F)_D^0 \cap \cE^\loc(E/F)_D^0 = T_{2,E}(D) \cap \cE^\loc(E/F)_D^0 =\{ 1\}$.
 
For the second part, write $t=(t_w)$ and $s=(s_w)$ as $w$ runs over places of $E$. Also write $\bart$ for the common image of the $t_w$ in $\A_D^\times/D^\times$. We see that for all $\sigma \in \Gal(D/F)$ we have 
\[ \sigma t_{\sigma^{-1}w} \equiv t_w \bmod D_w^\times . \]
From the first part of the lemma we see that ${}^\sigma \bart/\bart=1$, and so
\[ \bart \in (\A_D^\times/D^\times)^{\Gal(D/F)}=\A_F^\times/F^\times. \]
 Thus we may write $t_w=t_0t_w'$ with $t_0 \in \A_F^\times$ independent of $w$ and $t_w' \in D^\times$. If $\sigma \in \Gal(D/F)$, then
 \[ {}^\sigma (t_w')/(t_w') = (s_w)/{}^\sigma(s_w) \in \prod_w D^\times \cap \prod_w D_w^\times \subset \prod_w \A_D^\times ,\]
 from which we can conclude that ${}^\sigma(t_w')=(t_w')$ and ${}^\sigma(s_w)=(s_w)$. The lemma follows.
 \pfend

 \begin{cor} If $\balpha$ and $\balpha_1 \in \cZ(E/F)_D$ are equivalent, then there is a {\em unique} $\bgamma=(\gamma^\glob,\gamma^\loc)\in \cB(E/F)_D$ with 
 \[ \balpha_1={}^\bgamma \balpha. \]
\end{cor}

\pfbegin Any other such triple must be of the form $({}^a \gamma^\glob,{}^b \gamma^\loc)$ with $a\in \cE^\glob(E/F)_D^0$ and $b \in \cE^\loc(E/F)_D^0$ and $a/b \in T_{2,E}(\A_F)$. Thus, by the lemma, $a,b \in T_{2,E}(\A_F)$ so that ${}^a \gamma^\glob=\gamma^\glob$ and ${}^b \gamma^\loc=\gamma^\loc$. \pfend

\begin{cor}\label{stab} The stabilizer in $T_{2,E}(\A_D)$ of a class $\ga \in \cH(E/F)_D$ is 
\[ \cE^\loc(E/F)_D^0\cE^\glob(E/F)^0_DT_{2,E}(\A_F). \]\end{cor}

\pfbegin Suppose $\balpha= (\alpha^\glob,\alpha^\loc,\beta) \in \ga$. If ${}^t\ga=\ga$, then
\[ (\alpha^\glob,\alpha^\loc,{}^t\beta)=({}^{\gamma^\glob}\alpha^\glob,{}^{\gamma^\loc}\alpha^\loc,\gamma^\glob\beta(\gamma^\loc)^{-1}) \]
for some $\bgamma=(\gamma^\glob,\gamma^\loc) \in \cB(E/F)_D$. Then $\gamma^\glob={}^a1$ for some $a \in \cE^\glob(E/F)_D^0$ (because $H^1(\Gal(D/F), \cE^\glob(E/F)_D^0)=(0)$) and $\gamma^\loc={}^b1$ for some $b\in \cE^\loc(E/F)_D^0$ (because $H^1(\Gal(D/F), \cE^\loc(E/F)_D^0)=(0)$) and ${}^t\beta={}^{a/b}\beta$, so that $tb/a \in T_{2,E}(\A_F)$ and $t \in \cE^\loc(E/F)_D^0\cE^\glob(E/F)^0_DT_{2,E}(\A_F)$. The converse is easier.
\pfend

If $C\supset D\supset E \supset F$ are finite Galois extensions of $F$, then there are maps
\[ \begin{array}{rcl} \inf_{C/D}: \cB(E/F)_D &\lra& \cB(E/F)_C \\
(\gamma^\glob,\gamma^\loc) & \longmapsto & (\inf_{C/D}\gamma^\glob,\inf_{C/D} \gamma^\loc) \end{array} \]
and
\[ \begin{array}{rcl} \inf_{C/D}: \cZ(E/F)_D &\lra& \cZ(E/F)_C \\
(\alpha^\glob,\alpha^\loc,\beta) & \longmapsto & (\inf_{C/D}\alpha^\glob,\inf_{C/D} \alpha^\loc, \inf_{C/D} \beta) ,\end{array} \]
which induce a map
\[ {\inf}_{C/D}: \cH(E/F)_D \lra \cH(E/F)_C. \]
These maps are compatible with the map $T_{2,E}(\A_D) \into T_{2,E}(\A_C)$ and the actions of these groups. There are also maps
\[ \begin{array}{rcl} \eta_{D/E,*}: \cB(D/F)_C &\lra& \cB(E/F)_C \\ (\gamma^\glob,\gamma^\loc) & \longmapsto & (\eta_{D/E} \circ \gamma^\glob, \eta_{D/E} \circ \gamma^\loc) \end{array}\]
and
\[ \begin{array}{rcl} \eta_{D/E,*}: \cZ(D/F)_C &\lra& \cZ(E/F)_C \\ (\alpha^\glob,\alpha^\loc,\beta) & \longmapsto & (\eta_{D/E} \circ \alpha^\glob, \eta_{D/E} \circ \alpha^\loc, \eta_{D/E} \circ \beta ), \end{array}\]
which induce a map
\[ \eta_{D/E,*}: \cH(D/F)_C \lra \cH(E/F)_C. \]
These maps are compatible with the map $\eta_{D/E}:T_{2,D}(\A_C) \into T_{2,E}(\A_C)$ and the actions of these groups.

If $B \supset C\supset D\supset E \supset F$ are finite Galois extensions of $F$, then
\[ {\inf}_{B/C} \circ {\inf}_{C/D} = {\inf}_{B/D} \]
and
\[ \eta_{D/E,*} \circ \eta_{C/D,*} = \eta_{C/E,*} \]
and
\[ {\inf}_{B/C}\circ \eta_{D/E,*}=\eta_{D/E,*} \circ {\inf}_{B/C}: \cH(D/F)_C \lra \cH(E/F)_B. \]
The following lemma follows immediately:
\begin{lem} Suppose that $C \supset D \supset E \supset F$ are finite Galois extensions of $F$.  Suppose also that $\balpha_{C} \in\cZ(C/F)$ and $\balpha_{D} \in\cZ(D/F)$ and $\balpha_{E} \in\cZ(E/F)$ satisfy $\eta_{C/D,*}\balpha_{C}={}^{t'} \inf_{C/D} \balpha_D$ and $\eta_{D/E,*}\balpha_{D}={}^{t} \inf_{D/E} \balpha_E$ with $t' \in T_{2,D}(\A_{D'})$ and $t \in T_{2,E}(\A_{D})$. Then
\[ \eta_{C/E,*}\balpha_{C}={}^{t \eta_{D/E}(t')} \inf_{C/E} \balpha_E. \]
\end{lem}

To an element $\balpha=(\alpha^\glob,\alpha^\loc,\beta)$ we can associate extensions, and maps between them as follows:
\begin{itemize}
\item $\cE^\glob(E/F)_{D,\balpha}=\cE^\glob(E/F)_{D,\alpha^\glob}$ the extension of $\Gal(D/F)$ by $\cE^\glob(E/F)_{D}^0$ arising from $\alpha^\glob$.

\item $W_{E/F,D,\balpha}$ an extension of $\Gal(D/F)$ by $\A_D^\times/D^\times$ obtained as the pushout of $\cE^\glob(E/F)_{D,\balpha}$ along $\cE^\glob(E/F)_{D}^0 \onto \A_D^\times/D^\times$.

\item $\cE_3(E/F)_{D,\balpha}$ an extension of $\Gal(D/F)$ by $T_{3,E}(D)$ obtained as the pushout of $\cE^\glob(E/F)_{D,\balpha}$ along $\cE^\glob(E/F)_{D}^0 \onto T_{3,E}(D)$.

\item $\cE^\loc(E/F)_{D,\balpha}=\cE^\loc(E/F)_{D,\alpha^\loc}$ the extension of $\Gal(D/F)$ by $\cE^\loc(E/F)_{D}^0$ arising from $\alpha^\loc$.

\item $\cE^\loc(E/F)_{D,S,\balpha}$ an extension of $\Gal(D/F)$ by $\cE^\loc(E/F)_{D,S}^0$ obtained as the pushout of $\cE^\loc(E/F)_{D}^0 \ra \cE^\loc(E/F)_{D,S}^0$.

\item $\cE_2(E/F)_{D,\balpha}$ an extension of $\Gal(D/F)$ by $T_{2,E}(\A_D)$ obtained as the pushout of $\cE^\loc(E/F)_{D}^0 \ra T_{2,E}(\A_D)$.

\item If $u|w|v$ are places of $D$, $E$ and $F$ respectively, then we obtain an extension $W_{E_w/F_v,D,\balpha}$ of $\Gal(D/F)_w$ by $D_w^\times$ as the pushout of $\cE^\loc(E/F)_{D,\balpha}|_{\Gal(D/F)_w}$ along $\cE^\loc(E/F)_{D}^0\onto D_w^\times$; and an extension $W_{E_w/F_v,D_u,\balpha}$ of $\Gal(D/F)_u$ by $D_u^\times$ as the pushout of $\cE^\loc(E/F)_{D,\balpha}|_{\Gal(D/F)_u}$ along $\cE^\loc(E/F)_{D}^0\onto D_u^\times$.

\item A map of extensions
\[ \begin{array}{rcl} \loc_\balpha=\loc_{\beta}=i_{\beta^{-1}}: \cE^\glob(E/F)_{D,\balpha} &\lra& \cE_2(E/F)_{D,\balpha} \\
e_{\alpha^\glob}(\sigma) & \longmapsto & \beta(\sigma)^{-1} e_{\alpha^\loc}(\sigma)\end{array} \]
extending $\cE^\glob(E/F)_D^0 \into T_{2,E}(\A_D)$.

\item A map of extensions
\[ \begin{array}{rcl} \iota^\balpha_w: W_{E_w/F_v,D,\balpha} &\lra& W_{E/F,D,\balpha}|_{\Gal(D/F)_w} \\ e_{\alpha^\loc_w}(\sigma) & \longmapsto & (\pi_w \circ \beta)(\sigma)e_{\alpha^\glob}(\sigma) \end{array} \]
extending $D_w^\times \into \A_D^\times/D^\times$.
\end{itemize}
We will write $\cE^?(E/F)_D$ for any of the above extensions and $\cE^?(E/F)_D^0$ for the kernel of $\cE^?(E/F)_D\ra \Gal(D/F)$.

If $[\balpha_1]=[\balpha] \in \cH(E/F)_D$ then we get a {\em unique} $\bgamma=(\gamma^\glob,\gamma^\loc)\in \cB(E/F)_D$ with $\balpha_1={}^\bgamma \balpha$. In this case $i_\bgamma=i_{\gamma^\glob}$ gives canonical isomorphisms of extensions from $\cE^\glob(E/F)_{D,\balpha}$ to $\cE^\glob(E/F)_{D,\balpha_1}$, and from $W_{E/F,D,\balpha}$ to $W_{E/F,D,\balpha_1}$, and from $\cE_3(E/F)_{D,\balpha}$ to $\cE_3(E/F)_{D,\balpha_1}$. Moreover $i_\bgamma=i_{\gamma^\loc}$ gives canonical isomorphisms of extensions from $\cE^\loc(E/F)_{D,S,\balpha}$ to $\cE^\loc(E/F)_{D,S,\balpha_1}$, and from $W_{E_w/F_v,D,\balpha}$ to $W_{E_w/F_v,D,\balpha_1}$, and from $W_{E_w/F_v,D_u,\balpha}$ to $W_{E_w/F_v,D_u,\balpha_1}$, and from $\cE_2(E/F)_{D,\balpha}$ to $\cE_2(E/F)_{D,\balpha_1}$. We have $i_\bgamma \circ \loc_\balpha =\loc_{\balpha_1} \circ i_\bgamma$ and $i_\bgamma \circ \iota_w^\balpha =\iota_w^{\balpha_1} \circ i_\bgamma$.
Thus to $\ga \in \cH(E/F)_D$ we can canonically associate a well defined diagram of extensions:
\[ \begin{array}{ccccccc} && \cE_3(E/F)_{D,\ga} & \twoheadleftarrow & \cE^\glob(E/F)_{D,\ga} & \onto & W_{E/F,D,\ga}   \\ &&&& \loc_\ga \da && \\
\cE^\loc(E/F)_{D,S,\ga}& \twoheadleftarrow &  \cE^\loc(E/F)_{D,\ga} & \into & \cE_2(E/F)_{D,\ga}  && \\
&& \bigcup &&&& \\  W_{E_w/F_v,D,\ga} & \twoheadleftarrow &\cE^\loc(E/F)_{D,\ga}|_{\Gal(D/F)_w} &&&& \\
&& \bigcup &&&& \\  W_{E_w/F_v,D_u,\ga} & \twoheadleftarrow &\cE^\loc(E/F)_{D,\ga}|_{\Gal(D/F)_u},&&&& \end{array} \]
together with maps
\[ \iota^\ga_w: W_{E_w/F_v,D,\ga} \lra W_{E/F,D,\ga}. \]
If $\balpha=(\alpha^\glob,\alpha^\loc,\beta) \in \ga \in \cH(E/F)_D$, we will sometimes write 
\[ e_\balpha^\glob(\sigma) = e_{\alpha^\glob}(\sigma) \in \cE_3(E/F)_{D,\ga} \]
and
\[ e_\balpha^\loc(\sigma) = e_{\alpha^\loc}(\sigma) \in \cE^\loc(E/F)_{D,\ga}. \]

Note that there is an isomorphism of extensions $W_{E/F,D,\ga} \cong W_{E^\ab/F,D}$ which is unique up to composition with conjugation by an element of $\A_D^\times/D^\times$. Moreover there is an isomorphism of extensions $W_{E_w/F_v,D_u,\ga} \cong W_{E_w^\ab/F_v,D_u}$ which is unique up to composition with conjugation by an element of $D_u^\times$.

Any map of extensions 
\[ \begin{array}{ccccccccc} 
(0) &\lra& E_w^\times &\lra& W_{E_w/F_v,\ga}& \lra &\Gal(E_w/F_v)& \lra &(0) \\ && \da && i \da && || && \\
(0) & \lra & \A_E^\times/E^\times & \lra & W_{E/F,\ga}|_{\Gal(E/F)_w} & \lra & \Gal(E/F)_w& \lra &(0) \end{array} \]
must be of the form $\conju_a \circ \iota_w^\ga$ for some $a \in \A_E^\times/E^\times$. (As $H^1(\Gal(E_w/F_v),\A_E^\times/E^\times)=(0)$.)
If $\varphi:W_{E/F,\ga} \iso W_{E^\ab/F}$ and $\varphi_w: W_{E_w/F_v,\ga}\iso W_{E_w^\ab/F_v}$ are isomorphisms of extensions, and if $u$ is a place of $E^\ab$ above $w$, then 
we conclude that $\varphi \circ \iota^\ga_w=\conju_a \circ \theta_u \circ \varphi_w$ for some $a \in \A_E^\times/E^\times$. Thus $\conju_a^{-1} \circ \varphi \circ \iota^\ga_w$ and $\theta_u$ have the same image. In particular the image of $\varphi \circ \iota^\ga_w$ is the decomposition group for some place of $E^\ab$ above $w$.
This suggests that the choice of $\ga$ is not dissimilar from the choice of a place of $E^\ab$ above each place of $E$.

If $t \in T_{2,E}(\A_D)$ then 
\[ \cE^\glob(E/F)_{D,\balpha}=\cE^\glob(E/F)_{D,\alpha^\glob}=\cE^\glob(E/F)_{D,{}^t\balpha} \]
and
\[ \cE^\loc(E/F)_{D,\balpha}=\cE^\loc(E/F)_{D,\alpha^\loc}=\cE^\loc(E/F)_{D,{}^t\balpha}. \]
Thus we get identifications
\[ \gz_t: \cE^?(E/F)_{D,\balpha} \liso \cE^?(E/F)_{D,{}^t\balpha} \]
for all the extensions considered above. They commute with all the `unnamed' maps, but satisfy
\[ \gz_t \circ \loc_\balpha=\conju_t \circ \loc_{{}^t\balpha} \circ \gz_t \]
and
\[ \conju_{t_w} \circ \gz_t \circ \iota_w^\balpha =\iota_w^{{}^t\balpha} \circ \gz_t. \]
Moreover $\gz_t \circ i_\bgamma = i_\bgamma \circ \gz_t$, and so the maps $\gz_t$ descend to well-defined maps 
\[ \gz_t: \cE^?(E/F)_{D,\ga} \liso \cE^?(E/F)_{D,{}^t\ga} \]
for all the extensions considered above, which satisfy the same compatibilities. We have $\gz_{t_1t_2}=\gz_{t_1} \circ \gz_{t_2}$.

If $a \in \cE^\loc(E/F)_D^0$ and $b\in \cE^\glob(E/F)_D^0$ and $c \in T_{2,E}(\A_F)$ then ${}^{abc}\ga=\ga$. Moreover $\gz_{abc}:\cE^?(E/F)_{D,\ga} \iso \cE^?(E/F)_{D,\ga}$ equals $\conju_b^{-1}$ for $\cE^?(E/F)_{D,\ga}=\cE^\glob(E/F)_{D,\ga}$ or $\cE_3(E/F)_{D,\ga}$ or $W_{E/F,D,\ga}$; and $\gz_{abc}=\conju_a$ for $\cE^?(E/F)_{D,\ga}=\cE^\loc(E/F)_{D,S,\ga}$ or $\cE_2(E/F)_{D,\ga}$ or $W_{E_w/F_v,D,\ga}$ or $W_{E_w/F_v,D_u,\ga}$. 

\begin{lem}\label{decog} If $\tau \in \Gal(D/F)$ and $\balpha=(\alpha^\glob,\alpha^\loc,\beta) \in \ga$ , then $\conju_{e^\loc_\balpha(\tau)}$ induces an isomorphism $W_{E_w/F_v,D,\ga} \iso W_{E_{\tau w}/F_v,D,\ga}$.
Moreover
\[ \iota^\ga_{\tau w}\circ \conju_{e_\balpha^\loc(\tau)} = \conju_{ \beta(\tau)_{\tau w}} \circ \conju_{e_\balpha^\glob(\tau)} \circ \iota^\ga_w:  W_{E_w/F_v,D,\ga} \lra W_{E/F,D,\ga}. \]
\end{lem}

\pfbegin 
Both maps send $x \in D_w^\times$ to ${}^\tau x \in D_{\tau w}^\times$. 

If $\sigma \in \Gal(D/F)_w$, then 
$\conju_{e_{\alpha^\loc}(\tau)}(e_{\alpha^\loc}(\sigma))$ is the image of
\[ \begin{array}{rl} &e_{\alpha^\loc}(\tau) e_{\alpha^\loc}(\sigma) e_{\alpha^\loc}(\tau)^{-1} \\ =&
 e_{\alpha^\loc}(\tau) e_{\alpha^\loc}(\sigma) e_{\alpha^\loc}(\tau^{-1})(e_{\alpha^\loc}(\tau)e_{\alpha^\loc}(\tau^{-1}))^{-1} \\ =&
  \alpha^\loc(\tau,\sigma) \alpha^\loc(\tau \sigma,\tau^{-1}) e_{\alpha^\loc}(\tau \sigma \tau^{-1}) (\alpha^\loc(\tau, \tau^{-1}) \alpha^\loc(1,1))^{-1} \end{array}\]
and so
\[   \conju_{e_{\alpha^\loc}(\tau)}(e_{\alpha^\loc}(\sigma)) = \pi_{\tau w}(\alpha^\loc(\tau,\sigma) \alpha^\loc(\tau \sigma,\tau^{-1}) / {}^{\tau \sigma \tau^{-1}} (\alpha^\loc(\tau, \tau^{-1}) \alpha^\loc(1,1))) e_{\alpha^\loc}(\tau \sigma \tau^{-1}). \]
Thus
\[  \begin{array}{rl} & (\iota^\ga_{\tau w}\circ \conju_{e_{\alpha^\loc}(\tau)})(e_{\alpha^\loc}(\sigma)) \\ = &
\pi_{\tau w}(\alpha^\loc(\tau,\sigma) \alpha^\loc(\tau \sigma,\tau^{-1}) \beta(\tau \sigma \tau^{-1}) / {}^{\tau \sigma \tau^{-1}} (\alpha^\loc(\tau, \tau^{-1}) \alpha^\loc(1,1))) e_{\alpha^\glob}(\tau \sigma \tau^{-1}). \end{array} \]
On the other hand 
\[ \begin{array}{rl} & (\conju_{\pi_{\tau w} \beta(\tau)} \circ \conju_{e_{\alpha^\glob(\tau)}} \circ \iota^\ga_w)(e_{\alpha^\loc}(\sigma)) \\ = &
\pi_{\tau w} (\beta(\tau)){}^\tau \pi_w(\beta(\sigma)) e_{\alpha^\glob}(\tau) e_{\alpha^\glob}(\sigma) e_{\alpha^\glob}(\tau)^{-1}
\pi_{\tau w} (\beta(\tau))^{-1} \\ =&
\pi_{\tau w} (\beta(\tau){}^\tau \beta(\sigma)/{}^{\tau \sigma \tau^{-1}} \beta(\tau)) e_{\alpha^\glob}(\tau) e_{\alpha^\glob}(\sigma) e_{\alpha^\glob}(\tau)^{-1} \\ =&

\pi_{\tau w} (\beta(\tau){}^\tau \beta(\sigma)/{}^{\tau \sigma \tau^{-1}} \beta(\tau)) e_{\alpha^\glob}(\tau) e_{\alpha^\glob}(\sigma) e_{\alpha^\glob}(\tau^{-1})(e_{\alpha^\glob}(\tau)e_{\alpha^\glob}(\tau^{-1}))^{-1} \\ =&

\pi_{\tau w} (\beta(\tau){}^\tau \beta(\sigma)\alpha^\glob(\tau,\sigma) \alpha^\glob(\tau\sigma,\tau^{-1})/{}^{\tau \sigma \tau^{-1}} \beta(\tau)) e_{\alpha^\glob}(\tau\sigma \tau^{-1})\pi_{\tau w}(\alpha^\glob(\tau,\tau^{-1})\alpha^\glob(1,1))^{-1}.
\end{array} \]
Thus to prove the lemma it suffices to check that
\[ \begin{array}{rl} & \alpha^\loc(\tau,\sigma) \alpha^\loc(\tau \sigma,\tau^{-1}) \beta(\tau \sigma \tau^{-1}) / {}^{\tau \sigma \tau^{-1}} (\alpha^\loc(\tau, \tau^{-1}) \alpha^\loc(1,1)) \\ =&
 \beta(\tau){}^\tau \beta(\sigma)\alpha^\glob(\tau,\sigma) \alpha^\glob(\tau\sigma,\tau^{-1})/{}^{\tau \sigma \tau^{-1}} (\beta(\tau) \alpha^\glob(\tau,\tau^{-1})\alpha^\glob(1,1)), \end{array} \]
 or equivalently that
 \[ \begin{array}{rl} &  \beta(\tau \sigma \tau^{-1})  \\ =&
 \beta(\tau){}^\tau \beta(\sigma) \beta(\tau\sigma) \beta(\tau \sigma \tau^{-1})  {}^{\tau \sigma \tau^{-1}}( \beta(\tau){}^\tau \beta (\tau^{-1}) \beta(1)^2)
 /\beta(\tau){}^\tau \beta (\sigma) \beta(\tau \sigma){}^{\tau \sigma}\beta(\tau^{-1}) {}^{\tau \sigma \tau^{-1}} (\beta(\tau) \beta(1)^2 ), \end{array} \]
 which is clear.
\pfend

 We have an identification of $\cE^?(E/F)_{C,\inf_{C/D}\ga}$ with the pushout of $\cE^?(E/F)_{D,\ga}|_{\Gal(C/D)}$ along
 \[ \cE^?(E/F)_{D}^0 \lra \cE^?(E/F)_{C}^0 \]
 for each of the extensions considered above. (In the case $ W_{E_w/F_v,D,\ga}$ we use $W_{E_w/F_v,D,\ga}|_{\Gal(C/F)_w}$; and in the case $ W_{E_w/F_v,D_u,\ga}$ we use $W_{E_w/F_v,D_u,\ga}|_{\Gal(C/F)_t}$.)
 These identifications commute with the maps $\loc_\ga$, $\iota_w^\ga$ and $\gz_t$. 
 
 We also have an identification of 
 $\cE^?(E/F)_{C,\eta_{D/E,*}\ga}$ with the pushout of $\cE^?(D/F)_{C,\ga}$ along
 \[ \eta_{D/E}:\cE^?(D/F)_{C}^0 \lra \cE^?(E/F)_{C}^0 \]
 for each of the extensions considered above, except for the case $W_{D_u/F_v,C,\ga}$. These identifications commute with the maps $\loc_\ga$ and $\gz_t$. 
 
 The case of $W_{D_u/F_v,C,\ga}$ is much more complicated, and will be discussed in the next section.

 \begin{lem}\label{betacocyc} If $w_1,w_2$ are places of $E$ above a place $v$ of $F$, then 
 \[ (\pi_{w_1}/\pi_{w_2})(\beta(\sigma_1\sigma_2)) \equiv (\pi_{w_1}/\pi_{w_2})(\beta(\sigma_1)){}^{\sigma_1}(\pi_{\sigma_1^{-1}w_1}/\pi_{\sigma_1^{-1}w_2})(\beta(\sigma_2)) \bmod D^\times D_v^\times. \]
 Moreover
 \[  (\pi_{w_1}/\pi_{w_2})(\beta(1)) \equiv 1 \bmod D^\times D_v^\times. \]
 \end{lem}
 
 \pfbegin For the first part we have
 \[ (\pi_{w_1}/\pi_{w_2})(\beta(\sigma_1){}^{\sigma_1}\beta(\sigma_2)/\beta(\sigma_1\sigma_2)) = (\pi_{w_1}/\pi_{w_2})(\alpha^\loc(\sigma_1,\sigma_2)/\alpha^\glob(\sigma_1,\sigma_2) \in D^\times D_v^\times .\]
 The second part follows from the first and $1^2=1$.
 \pfend

 \subsection{The extension $W_{E_w/F_v,D,\ga}$}
 
 Our first aim is to find a more canonical version of this group which does not depend on the choice of $\ga$ in the same way that $W_{E_w^\ab/F_v,D_u}$ is a more canonical group isomorphic to $W_{E_w/F_v,D_u,\ga}$. 
 
Fix an $F$-linear embedding $\rho: E \into \barFv$. We define
\[ W_{(EF_v)^\ab/F_v,\rho,D} = (D_{w(\rho)}^\times \rtimes (W_{(EF_v)^\ab/F_v}\times_{\Gal(E/F)} \Gal(D/F)))/E_{w(\rho)}^\times \]
where:
\begin{itemize}
\item $W_{(EF_v)^\ab/F_v} \onto \Gal(E/F)_{w(\rho)}$ is the composition of $W_{(EF_v)^\ab/F_v} \onto \Gal(EF_v/F_v)$ with the inverse of the isomorphism $\Gal(E/F)_{w(\rho)} \iso \Gal((EF_v)/F_v)$ induced by $\rho$.
\item $W_{(EF_v)^\ab/F_v}\times_{\Gal(E/F)_{w(\rho)}} \Gal(D/F)_{w(\rho)}$ acts on $D_{w(\rho)}^\times$ via its projection to $\Gal(D/F)_{w(\rho)} \subset \Gal(D/F)$.
\item The map $E_{w(\rho)}^\times \ra D_{w(\rho)}^\times \rtimes (W_{(EF_v)^\ab/F_v}\times_{\Gal(E/F)_{w(\rho)}} \Gal(D/F)_{w(\rho)})$ sends $a$ to $(a^{-1},(r_{EF_v}(\rho(a)),1))$.
\end{itemize}
It fits into an exact sequence
\[ (0) \lra D_{w(\rho)}^\times \lra W_{(EF_v)^\ab/F_v,\rho,D} \lra \Gal(D/F)_{w(\rho)} \lra (0), \]
which has class
\[ {\inf}_{D_u/E_{w(\rho)}} [\alpha_{E_{w(\rho)}/F_v}] \in H^1(\Gal(D/F)_u,D_u^\times) \cong H^1(\Gal(D/F)_{w(\rho)},D_{w(\rho)}^\times) \]
for any place $u$ of $D$ above $w(\rho)$.
If $\sigma \in \Gal(D/F)$ then we get an isomorphism
\[ \begin{array}{rcl} \sigma_*: W_{(EF_v)^\ab/F_v,\rho,D} &\liso& W_{(EF_v)^\ab/F_v,\rho\sigma^{-1},D} \\
{}[(a,(\tau_1,\tau_2))] & \longmapsto & [(\sigma(a), (\tau_1, \sigma \tau_2 \sigma^{-1}))].
\end{array} \]

Now suppose that $\trho:D \into \barFv$ is $F$-linear. Then we define
\[ W_{(EF_v)^\ab/F_v,D,\trho} = (D_{w(\trho)}^\times \rtimes  (W_{(DF_v)^\ab/F_v} \times_{\Gal(E/F)} \Gal(D/F)))/W_{(DF_v)^\ab/(EF_v)} \]
where:
\begin{itemize}
\item $w(\trho)$ is the place of $E$ induced by $\trho|_E$.
\item $W_{(DF_v)^\ab/F_v} \times_{\Gal(E/F)} \Gal(D/F)$ acts on $D_{w(\trho)}^\times$ via projection to $\Gal(D/F)_{w(\trho)} \subset \Gal(D/F)$.

\item The map $W_{(DF_v)^\ab/(EF_v)} \ra D_{w(\trho)}^\times \rtimes  (W_{(DF_v)^\ab/F_v} \times_{\Gal(E/F)} \Gal(D/F))$ is given by
\[ \sigma \longmapsto (\trho^{-1} r_{EF_v}^{-1}(\sigma|_{(EF_v)^\ab}^{-1}), (\sigma,1)). \]
We leave it to the reader to check that the image of this map is a normal subgroup.

\end{itemize}
Again it fits into an exact sequence
\[ (0) \lra D_{w(\trho)}^\times \lra W_{(EF_v)^\ab/F_v,D,\trho} \lra \Gal(D/F)_{w(\trho)} \lra (0), \]
which has class
\[ [D_{u(\trho)}:E_{w(\trho)}][\alpha_{D_{u(\trho)}/F_v}] \in H^1(\Gal(D/F)_{u(\trho)},D_{u(\trho)}^\times) \cong H^1(\Gal(D/F)_{w(\trho)},D_{w(\trho)}^\times). \]
If $\sigma \in \Gal(D/F)$ then we get an isomorphism
\[ \begin{array}{rcl} \sigma_*: W_{(EF_v)^\ab/F_v,D,\trho} &\liso& W_{(EF_v)^\ab/F_v,D,\trho\sigma^{-1}} \\
{}[(a,(\tau_1,\tau_2))] & \longmapsto & [(\sigma(a), (\tau_1,\sigma\tau_2\sigma^{-1}) )].
\end{array} \]
Moreover there are natural isomorphisms of extensions
\[ \begin{array}{rcl} W_{(EF_v)^\ab/F_v,D,\trho} & \liso & W_{(EF_v)^\ab/F_v,\trho|_E,D} \\
{} [(a,(\tau_1,\tau_2))] & \longmapsto & [(a,(\tau_1|_{(EF_v)^\ab},\tau_2 ))] \end{array} \]
compatible with the actions of $\Gal(D/F)$. (We leave the verification of this to the reader.)

If $\rho:E^\ab \into \barFv$ is $F$-linear, then there is a natural map
\[ \begin{array}{rcl} \theta_\rho: W_{(EF_v)^\ab/F_v,\rho|_E,D} &\lra& W_{E^\ab/F,D}=
(\A_D^\times/D^\times \rtimes W_{E^\ab/F}|_{\Gal(D/F)})/(\A_E^\times/E^\times) \\
{}[(a,(\tau_1,\tau_2))] & \longmapsto & [(a,(\theta_\rho(\tau_1),\tau_2))],
\end{array} \]
which is canonical up to composition with conjugation with an element of $\Delta_E$. 
If $\sigma \in W_{E^\ab/F}|_{\Gal(D/F)}$, then 
\[ \theta_{\rho \sigma^{-1}} \circ \sigma|_{D,*} = \conju_\sigma \circ \theta_\rho \]
up to composition with conjugation by an element of $\Delta_E$.

 If $C \supset D \supset E \supset F$ are finite Galois extensions of $F$ then
\[ W_{(EF_v)^\ab/F_v,\rho,C} \cong (C_{w(\rho)}^\times \rtimes W_{(EF_v)^\ab/F_v,\rho,D}|_{\Gal(C/F)_{w(\rho)}})/D_{w(\rho)}^\times. \]
Moreover
\[ \theta_\rho: W_{(EF_v)^\ab/F_v,\rho|_E,C} \lra W_{E^\ab/F,C} \]
is identified with
\[ \begin{array}{rcl} (C_{w(\rho)}^\times \rtimes W_{(EF_v)^\ab/F_v,\rho|_E,D}|_{\Gal(C/F)_{w(\rho)}})/D_{w(\rho)}^\times &\lra&
(\A_C^\times/C^\times \rtimes W_{E^\ab/F,D}|_{\Gal(C/F)})/(\A_D^\times/D^\times) \\ 
{} [(a,(\sigma,\tau))] & \longmapsto & [(a,(\theta_\rho(\sigma),\tau))]. \end{array} \]

If $\trho:D^\ab\into \barFv$ is $F$-linear, then there is a natural map
\[ \begin{array}{rcl} \theta_\trho: W_{(\rho(E)F_v)^\ab/F_v,D,\trho|_D} &\lra& W_{E^\ab/F,D}=(\A_D^\times/D^\times \rtimes W_{D^\ab/F})/(\A_D^\times/D^\times) \\
{}[(a,(\tau_1,\tau_2))] & \longmapsto & [(ar_E^{-1}((\theta_\trho(\tau_1)\ttau_2^{-1})|_{E^\ab}),\ttau_2)],
\end{array} \]
where $\ttau_2 \in W_{D^\ab/E}$ is any lift of $\tau_2$. This map in fact only depends on $\trho|_{E^\ab D}$, so we will sometimes write $\theta_\trho$ for $\trho:E^\ab D \into \barFv$.
This map is canonical up to composition with conjugation with an element of $\Delta_E=N_{D/E}\Delta$. If $\sigma \in W_{D^\ab/F}$, then 
\[ \theta_{\trho \sigma|_{E^\ab D}^{-1}} \circ \sigma|_{D,*} = \conju_\sigma \circ \theta_\trho \]
up to composition with conjugation by an element of $\Delta_E$. Moreover
\[ \begin{array}{rcl} W_{(EF_v)^\ab/F_v,D,\trho|_D} & \liso & W_{(EF_v)^\ab/F_v,\trho|_E,D} \\ \theta_\trho \searrow &&\swarrow \theta_{\trho|_{E^\ab}} \\ & W_{E^\ab/F,D} & \end{array} \]
commutes. We also leave these verifications to the reader.

We see that there is an isomorphism of extensions
\[ W_{E_{w(\rho)}/F_v,D,\ga} \cong W_{E_{w(\rho)}/F_v,\rho,D}, \]
which is unique up to composition with conjugation by an element of $D_{w(\rho)}^\times$. If $\varphi:W_{E/F,D,\ga} \iso W_{E^\ab/F,D}$ and $\varphi_{w(\rho)}:W_{E_{w(\rho)}/F_v,D,\ga} \iso W_{E_{w(\rho)}/F_v,\rho,D}$ are isomorphisms of extensions then $\varphi \circ \iota_{w(\rho)}^\ga$ and $\theta_\rho \circ \varphi_{w(\rho)}$ differ by composition with conjugation by an element of $\A_D^\times/D^\times$.

\begin{lem} \label{iotamixed}  Suppose that $D \supset E \supset F$ are finite Galois extensions of a number field $F$ and that $u|w|v$ are places of these fields. Suppose also that $\balpha \in \cZ(D/F)$. Then there is an isomorphism
\[ j_{D/E,\balpha}:(D_w^\times \rtimes (W_{D_u/F_v,\balpha} \times_{\Gal(E/F)} \Gal(D/F)))/W_{D_u/F_v,\balpha}|_{\Gal(D/E)_u} \liso W_{E_w/F_v,D,\eta_{D/E,*}\balpha}, \]
where
\[ \begin{array}{rcl} W_{D_u/F_v,\balpha}|_{\Gal(D/E)_u} & \lra &  D_w^\times \rtimes (W_{D_u/F_v,\balpha} \times_{\Gal(E/F)} \Gal(D/F)) \\ \sigma & \longmapsto & (\epsilon(\sigma)^{-1},(\sigma,1)), \end{array} \]
and
\[ \epsilon=\epsilon_{D/E}: W_{D_u/F_v,\balpha}|_{\Gal(D/E)_u} \lra E_w^\times,  \]
which we embed diagonally in $D_w^\times$, is the homomorphism sending
\[ \sigma \longmapsto \prod_{\eta \in \Gal(D_u/E_w)} s_\eta \sigma s_{\eta \sigma}^{-1} \]
for any section $s:\Gal(D_u/E_w) \ra W_{D_u/F_v,\balpha}|_{\Gal(D/E)_u}$. The isomorphism is given by
\[  j_{D/E,\balpha}((a,(be^\loc_\balpha(\sigma),\tau)))=a \prod_{\eta \in \Gal(D/E)} \left(\eta(b)\alpha^\loc(\eta,\sigma)/\alpha^\loc(\tau, \tau^{-1}\eta \sigma)\right)|_{\eta u}  e^\loc_{\eta_{D/E,*} \balpha}(\tau),  \]
where if $c \in \prod_x D_x^\times$ we write $c|_x \in \prod_xD_x^\times$ for the element that is $c_x$ at $x$ and $1$ elsewhere. We also have
\[ i_{\eta_{D/E,*} \bgamma} \circ j_{D/E,\balpha} = \conju_{
\prod_{\eta \in \Gal(D/E)} \gamma^\loc(\eta)|_{\eta u}} \circ j_{D/E,\eta_{D/E,*}\balpha}\circ (1 \times i_\bgamma \times 1). \]

Moreover the map
\[ \iota_w^{\eta_{D/E,*}\balpha}: (D_w^\times \rtimes ( W_{D_u/F_v,\balpha} \times_{\Gal(E/F)} \Gal(D/F))  )/W_{D_u/F_v,\balpha}|_{\Gal(D/E)_u} \lra (\A_D^\times/D^\times \rtimes W_{D/F})/(\A_D^\times/D^\times) \]
sends
\[ [(a,(1,1))]  \longmapsto [(a,1)] \]
for $a \in D_w^\times$, 
and
\[ [(1,(\sigma,\tau))] \longmapsto \conju_{\prod_{\eta \in \Gal(D/E)} \beta(\eta)_{\eta u}}  [(\prod_{\eta \in \Gal(D/E)} e^\glob_\balpha(\tau \eta \sigma^{-1}) \iota_u^\balpha(\sigma)e^\glob_\balpha(\eta)^{-1}e^\glob_\balpha(\tau)^{-1},e^\glob_\balpha(\tau))]\] 
for $(\sigma,\tau) \in W_{D_u/F_v,\balpha} \times_{\Gal(E/F)} \Gal(D/F)$.

\end{lem}

\pfbegin The map $\epsilon$ is independent of the choice of section $s$, it is a $W_{D_u/F_v,\balpha}$-equivariant homomorphism and it is valued in $E_w^\times \subset D_u^\times$. The map
\[ W_{D_u/F_v,\balpha}|_{\Gal(D/E)_u}  \lra   D_w^\times \rtimes (W_{D_u/F_v,\balpha} \times_{\Gal(E/F)} \Gal(D/F)) \]
is therefore a homomorphism with normal image. 

There is an exact sequence
\[ (0) \lra D_w^\times \lra (D_w^\times \rtimes (W_{D_u/F_v,\balpha} \times_{\Gal(E/F)} \Gal(D/F)))/W_{D_u/F_v,\balpha}|_{\Gal(D/E)_u} \lra \Gal(D/F)_w \lra (0). \]
The given map 
\[ j_{D/E,\balpha}:(D_w^\times \rtimes (W_{D_u/F_v,\balpha} \times_{\Gal(E/F)_w} \Gal(D/F)_w)) \lra W_{E_w/F_v,D,\eta_{D/E,*}\balpha} \]
is compatible with the inclusion of $D_w^\times$ and the projection to $\Gal(D/F)_w$. 
Thus to prove the first assertion of the lemma it suffices to show that $j_{D/E,\balpha}$ is a homomorphism that is trivial when restricted to $W_{D_u/F_v,\balpha}|_{\Gal(D/E)_u}$. 

Thus to prove the lemma it suffices to check that:
\begin{itemize}
\item If $(\sigma, \tau) \in \Gal(D/F)_u \times_{\Gal(E/F)} \Gal(D/F)$ and $b \in \cE^\loc(D/F)^0$, then 
\[ {}^\tau \left( \prod_{\eta \in \Gal(D/E)} \eta(b)|_{\eta u}\right) = \prod_{\eta \in \Gal(D/E)} \eta({}^\sigma b)|_{\eta u}.\]

\item If $(\sigma_1, \tau_1)$ and $(\sigma_2,\tau_2) \in Gal(D/F)_u \times_{\Gal(E/F)} \Gal(D/F)$ then
\[  \begin{array}{rl} & \prod_{\eta \in \Gal(D/E)} \left(\alpha^\loc(\eta,\sigma_1)/\alpha^\loc(\tau_1, \tau_1^{-1}\eta \sigma_1)\right)|_{\eta u} \\ & {}^{\tau_1} \left( \prod_{\eta \in \Gal(D/E)} \left(\alpha^\loc(\eta,\sigma_2)/\alpha^\loc(\tau_2, \tau_2^{-1}\eta \sigma_2)\right)|_{\eta u} \right) 
\eta_{D/E,*}(\alpha^\loc)(\tau_1,\tau_2) \\ 
=  & \prod_{\eta \in \Gal(D/E)} \left(\eta(\alpha^\loc(\sigma_1,\sigma_2))\alpha^\loc(\eta,\sigma_1\sigma_2)/\alpha^\loc(\tau_1\tau_2, \tau_2^{-1}\tau_1^{-1}\eta \sigma_1\sigma_2)\right)|_{\eta u} \end{array} \]

\item If $\sigma \in \Gal(D/E)_u$ and $b \in D_u^\times$, then
\[ \epsilon(be^\loc_\balpha(\sigma))^{-1} \prod_{\eta \in \Gal(D/E)}\left( \eta (b) \alpha^\loc(\eta,\sigma)/\alpha^\loc(1,\eta \sigma)\right)|_{\eta u}e^\loc_{\eta_{D/E,*}\alpha^\loc}(1) =1.\]

\item If $(\sigma, \tau) \in \Gal(D/F)_u \times_{\Gal(E/F)} \Gal(D/F)$ then
\[ \begin{array}{rl} & \prod_{\eta \in \Gal(D/E)} (\gamma^\loc(\eta) {}^\eta (\gamma^\loc(\sigma)|_u){}^{\gamma^\loc}\alpha^\loc(\eta,\sigma)/{}^{\gamma^\loc}\alpha^{\loc}(\tau,\tau^{-1}\eta\sigma))|_{\eta u} {}^\tau(\gamma^\loc(\eta)|_{\eta u}^{-1}) \\ =&  (\eta_{D/E}\circ \gamma^\loc)(\tau)|_w\prod_{\eta \in \Gal(D/E)} (\alpha^\loc(\eta,\sigma)/\alpha^\loc(\tau,\tau^{-1}\eta\sigma))|_{\eta u} .\end{array} \]

\item If $b\in D_u^\times$ then 
\[ [(\prod_{\eta \in \Gal(D/E)} \eta(b)|_{\eta u},1)]=\left[\left( \left(\prod_{\eta \in \Gal(D/E)} e_\balpha^\glob(\eta) b|_u e_\balpha^\glob(\eta)^{-1}\right) e_\balpha^\glob(1)^{-[D:E]},e_\balpha^\glob(1)\right)\right] .\]

\item If  $(\sigma, \tau) \in \Gal(D/F)_u \times_{\Gal(E/F)} \Gal(D/F)$ then
\[ \begin{array}{rl} &{}^{(1-\tau)}\left(\prod_{\eta \in \Gal(D/E)} \beta(\eta)_{\eta u}\right)\prod_{\eta \in \Gal(D/E)} e^\glob_\balpha(\tau \eta \sigma^{-1}) \beta(\sigma)|_u e^\glob_\balpha(\sigma)e^\glob_\balpha(\eta)^{-1}e^\glob_\balpha(\tau)^{-1} \\ = &
 \prod_{\eta \in \Gal(D/E)} \left(\alpha^\loc(\eta,\sigma)/\alpha^\loc(\tau, \tau^{-1}\eta \sigma)\right)|_{\eta u}  (\eta_{D/E} \beta)(\tau)|_w. \end{array} \]
\end{itemize}

The first assertion is equivalent to 
\[   \prod_{\eta \in \Gal(D/E)} (\tau \eta \sigma^{-1})({}^\sigma b)|_{\tau \eta \sigma^{-1} u} = \prod_{\eta \in \Gal(D/E)} \eta({}^\sigma b)|_{\eta u},\]
which is clear on changing variables in the product. 

The second assertion is equivalent to
\[  \begin{array}{rl} & \prod_{\eta \in \Gal(D/E)} \left(\alpha^\loc(\eta,\sigma_1)\alpha^\loc(\tau_1,\tau_2)/\alpha^\loc(\tau_1, \tau_1^{-1}\eta \sigma_1)\right)|_{\eta u} \\ &  \prod_{\eta \in \Gal(D/E)} \left({}^{\tau_1}\alpha^\loc(\eta,\sigma_2)/{}^{\tau_1}\alpha^\loc(\tau_2, \tau_2^{-1}\eta \sigma_2)\right)|_{\tau_1\eta \sigma_1^{-1} u}  
 \\ 
=  & \prod_{\eta \in \Gal(D/E)} \left({}^\eta\alpha^\loc(\sigma_1,\sigma_2)\alpha^\loc(\eta,\sigma_1\sigma_2)/\alpha^\loc(\tau_1\tau_2, \tau_2^{-1}\tau_1^{-1}\eta \sigma_1\sigma_2)\right)|_{\eta u} \end{array} \]
or
\[  \begin{array}{rl} & \prod_{\eta \in \Gal(D/E)} \left(\alpha^\loc(\eta,\sigma_1)\alpha^\loc(\tau_1,\tau_2){}^{\tau_1}\alpha^\loc(\tau_1^{-1}\eta\sigma_1,\sigma_2)/\alpha^\loc(\tau_1, \tau_1^{-1}\eta \sigma_1){}^{\tau_1}\alpha^\loc(\tau_2, \tau_2^{-1}\tau_1^{-1}\eta \sigma_1\sigma_2)\right)|_{\eta u}
 \\ 
=  & \prod_{\eta \in \Gal(D/E)} \left({}^\eta\alpha^\loc(\sigma_1,\sigma_2)\alpha^\loc(\eta,\sigma_1\sigma_2)/\alpha^\loc(\tau_1\tau_2, \tau_2^{-1}\tau_1^{-1}\eta \sigma_1\sigma_2)\right)|_{\eta u} \end{array} \]
or
\[  \begin{array}{rl} & \prod_{\eta \in \Gal(D/E)} \left( {}^{\tau_1}\alpha^\loc(\tau_1^{-1}\eta\sigma_1,\sigma_2)/\alpha^\loc(\tau_1, \tau_1^{-1}\eta \sigma_1)\right)|_{\eta u}
 \\ 
=  & \prod_{\eta \in \Gal(D/E)} \left({}^\eta\alpha^\loc(\sigma_1,\sigma_2)\alpha^\loc(\eta,\sigma_1\sigma_2)/\alpha^\loc(\eta,\sigma_1)\right)|_{\eta u} \\ &
 \prod_{\eta \in \Gal(D/E)} \left({}^{\tau_1}\alpha^\loc(\tau_2, \tau_2^{-1}\tau_1^{-1}\eta \sigma_1\sigma_2)/
\alpha^\loc(\tau_1\tau_2, \tau_2^{-1}\tau_1^{-1}\eta \sigma_1\sigma_2)\alpha^\loc(\tau_1,\tau_2)\right)|_{\eta u} \end{array} \]
or
\[  \begin{array}{rl}  & \prod_{\eta \in \Gal(D/E)} \left(\alpha^\loc(\eta \sigma_1,\sigma_2)/ \alpha^\loc(\tau_1, \tau_1^{-1}\eta \sigma_1\sigma_2)\right)|_{\eta u}
 \\ =  & \prod_{\eta \in \Gal(D/E)} \left(\alpha^\loc(\eta \sigma_1,\sigma_2)/ \alpha^\loc(\tau_1, \tau_1^{-1}\eta \sigma_1\sigma_2)\right)|_{\eta u},  \end{array} \]
which is clear.

For the third assertion note that for $b \in D_u^\times$ and $\sigma \in \Gal(D_u/E_w)$ we have
\[ \epsilon(b) = \prod_{\eta \in Gal(D/E)} \eta (b)|_{\eta u} \]
and
\[ \epsilon(e^\loc_\balpha(\sigma) ) = \prod_{\eta \in \Gal(D/E)_u} \alpha^\loc(\eta,\sigma)|_u \in E_w^\times \subset D_w^\times .\]
If $\cR\subset \Gal(D/E)$ is a set of representatives for $\Gal(D/E)/\Gal(D/E)_u$, the latter is equivalent to
\[ \begin{array}{rcl} \epsilon(e^\loc_\balpha(\sigma) ) &=& \prod_{\zeta \in \cR} \prod_{\eta \in \Gal(D/E)_u} ({}^\zeta \alpha^\loc(\eta,\sigma))|_{\zeta u} \\ &=&
 \prod_{\zeta \in \cR} \prod_{\eta \in \Gal(D/E)_u} ( \alpha^\loc(\zeta\eta,\sigma) \alpha^\loc(\zeta,\eta)/\alpha^\loc(\zeta, \eta\sigma))|_{\zeta u} \\ &=&
 \prod_{\zeta \in \cR} \prod_{\eta \in \Gal(D/E)_u} \alpha^\loc(\zeta \eta,\sigma)|_{\zeta \eta u} \\ &&
 \prod_{\zeta \in \cR} \prod_{\eta \in \Gal(D/E)_u} \alpha^\loc(\zeta,\eta)|_{\zeta u} / \\ &&
\prod_{\zeta \in \cR} \prod_{\eta \in \Gal(D/E)_u} \alpha^\loc(\zeta, \eta\sigma)|_{\zeta u}\\ &=&
 \prod_{\eta \in \Gal(D/E)} \alpha^\loc(\eta,\sigma)|_{\eta u}
 \\ &\in& D_w^\times .\end{array}\]
 Thus the third assertion is equivalent to
\[ \epsilon(be^\loc_\balpha(\sigma))= \left(\prod_{\eta \in \Gal(D/E)}\left( \eta (b) \alpha^\loc(\eta,\sigma)/\alpha^\loc(1,\eta \sigma)\right)|_{\eta u}\right) (\eta_{D/E,*}\alpha^\loc)(1,1)|_w ,\]
i.e. to
\[ \prod_{\eta \in \Gal(D/E)} \alpha^\loc(\eta,\sigma)|_{\eta u}= \prod_{\eta \in \Gal(D/E)}\left(  \alpha^\loc(\eta,\sigma)\alpha^\loc(1,1)/\alpha^\loc(1,\eta \sigma)\right)|_{\eta u} ,\]
which is true.

The fourth assertion is equivalent to
\[ \begin{array}{rl} & \prod_{\eta \in \Gal(D/E)} (\gamma^\loc(\eta) {}^\eta \gamma^\loc(\sigma) \gamma^\loc(\eta \sigma) \gamma^\loc(\eta)^{-1}{}^\eta \gamma^\loc(\sigma)^{-1}/\gamma^\loc(\eta \sigma)\gamma^\loc(\tau)^{-1}{}^\tau\gamma^\loc(\tau^{-1}\eta \sigma)^{-1})|_{\eta u} \\ & \prod_{\eta \in \Gal(D/E)}({}^\tau\gamma^\loc(\eta))|_{\tau\eta \sigma^{-1}u}^{-1} \\ =&  \prod_{\eta \in \Gal(D/E)} \gamma^\loc(\tau)|_{\eta u} ,\end{array} \]
i.e. to
\[ \prod_{\eta \in \Gal(D/E)}( {}^\tau\gamma^\loc(\tau^{-1}\eta \sigma)^{-1})|_{\eta u} = \prod_{\eta \in \Gal(D/E)}({}^\tau\gamma^\loc(\eta))|_{\tau\eta \sigma^{-1}u},\]
which follows immediately on changing variables.

The fifth assertion is immediate. The sixth assertion is equivalent to
\[ \begin{array}{rl} &{}^{(1-\tau)}\left(\prod_{\eta \in \Gal(D/E)} \beta(\eta)_{\eta u}\right) \\ & \prod_{\eta \in \Gal(D/E)}  ({}^{\tau\eta \sigma^{-1}}\beta(\sigma))|_{\tau \eta \sigma^{-1}u }
\alpha^\glob(\tau\eta\sigma^{-1},\sigma)e^\glob_\balpha(\tau \eta ) e^\glob_\balpha(\eta)^{-1}e^\glob_\balpha(\tau)^{-1} \\ = &
 \prod_{\eta \in \Gal(D/E)} \left(\alpha^\glob(\eta,\sigma)\beta(\eta){}^\eta \beta(\sigma)\beta(\tau)\beta(\eta \sigma)/\beta(\eta \sigma)\alpha^\glob(\tau, \tau^{-1}\eta \sigma)\beta(\tau){}^\tau \beta(\tau^{-1}\eta\sigma)\right)|_{\eta u} , \end{array} \]
i.e. to
\[ \begin{array}{rl} &\prod_{\eta \in \Gal(D/E)}  \left((\beta(\eta) {}^{\eta }\beta(\sigma))|_{ \eta u }/({}^\tau \beta(\eta))|_{\tau\eta u} \right)\prod_{\eta \in \Gal(D/E)}
(\alpha^\glob(\tau\eta\sigma^{-1},\sigma)/\alpha^\glob(\tau,\eta)) \\ = &
  \prod_{\eta \in \Gal(D/E)} (\alpha^\glob(\eta,\sigma)/\alpha^\glob(\tau, \tau^{-1}\eta \sigma))
\\ & \prod_{\eta \in \Gal(D/E)} \left(\beta(\eta){}^\eta \beta(\sigma)\right)|_{\eta u} / \prod_{\eta \in \Gal(D/E)}  ({}^\tau \beta(\eta))|_{\tau \eta u} , \end{array} \]
which is seen to be true, by change of variable.
\pfend

\subsection{Explicit cocycles}\label{explicit}

For comparison with some constructions of Langlands in \cite{langmarch}, it will be useful to us to have a more explicit form for some elements of $\cZ(E/F)$. 

Fix the following data
\begin{itemize}
\item $\alpha$ representing $[\alpha^\glob_{E/F}]$,
\item a place $w=w(v)$ above each place $v$ of $F$,
\item $\alpha_w$ representing $[\alpha_{E_w/F_v}]$ for each such $w$,
\item a section $s_w:\Gal(E/F)/\Gal(E_w/F_v) \ra \Gal(E/F)$ with image $H_w$ for each such $w$,
\item and, for each such $w$, a function $\gamma_w: \Gal(E_w/F_v) \ra \A_E^\times$ such that ${}^{\gamma_w} (\pi_w \circ \alpha)|_{\Gal(E_w/F_v)}=i_w \alpha_w \in Z^2(\Gal(E_w/F_v),\A_E^\times)$. (This is possible as $\alpha$ is equivalent to some $\alpha^\loc$ in $Z^2(\Gal(E/F),T_{2,E}(\A_E))$.)
\end{itemize}
Also set $\delta_v(\sigma)=\sigma^{-1}s_v(\sigma)$.

Then
\[ \alpha^\loc (\sigma_1,\sigma_2) = \prod_{v \in V_F} \prod_{\eta \in \Gal(E/F)/\Gal(E_w/F_v)} s_w(\eta) \alpha_w (s_w(\eta)^{-1} \sigma_1 s_w(\sigma_1^{-1} \eta), s_w(\sigma_1^{-1}\eta)^{-1} \sigma_2 s_w(\sigma_2^{-1} \sigma_1^{-1} \eta)) \]
is a representative of $[\alpha_{E/F}^\loc]$. (If one restricts the class to $\Gal(E_w/F_v)$ and projects to $E_w^\times$ one certainly recovers $\alpha_w$, so the only thing to check is that $\alpha(\sigma_1,\sigma_2)$ is a 2-cocycle. Writing out the cocycle relation and changing the variable from $\eta$ to $\sigma_1 \eta$ in one of the terms, what we need to check is that
\[ \begin{array}{rl} &(s_w(\eta)^{-1} \sigma_1 s_w(\sigma_1^{-1} \eta)) \alpha_w(s_w(\sigma_1^{-1}\eta)^{-1} \sigma_2 s_w(\sigma_2^{-1}\sigma_1^{-1}\eta), s_w(\sigma_2^{-1}\sigma_1^{-1}\eta)^{-1} \sigma_3 s_w(\sigma_3^{-1}\sigma_2^{-1}\sigma_1^{-1}\eta)) \\ &
\alpha_w(s_w(\eta)^{-1} \sigma_1 s_w(\sigma_1^{-1} \eta),s_w(\sigma_1^{-1}\eta)^{-1} \sigma_2  \sigma_3 s_w(\sigma_3^{-1}\sigma_2^{-1}\sigma_1^{-1}\eta)) \\ = &
\alpha_w(s_w(\eta)^{-1} \sigma_1   \sigma_2 s_w(\sigma_2^{-1}\sigma_1^{-1} \eta),
s_w(\sigma_2^{-1}\sigma_1^{-1} \eta)^{-1}  \sigma_3 s_w(\sigma_3^{-1}\sigma_2^{-1}\sigma_1^{-1}\eta)) \\ &
\alpha_w(s_w(\eta)^{-1} \sigma_1 s_w(\sigma_1^{-1} \eta),s_w(\sigma_1^{-1}\eta)^{-1} \sigma_2   s_w(\sigma_2^{-1}\sigma_1^{-1}\eta)), \end{array} \]
which is just the cocycle relation for $\alpha_w$.)

Define
\[ \beta: \Gal(E/F) \lra T_{2,E}(\A_E) \]
by
\[ \pi_{\eta w} \beta(\sigma) = \pi_{\eta w} (\alpha(s_w(\eta),s_w(\eta)^{-1}\sigma s_w(\sigma^{-1}\eta))/\alpha(\sigma,s_w(\sigma^{-1} \eta))) {}^{s_w(\eta)} \gamma_w(s_w(\eta)^{-1}\sigma s_w(\sigma^{-1} \eta))^{-1}. \]
We claim that $(\alpha,\alpha^\loc,\beta)\in \cZ(E/F)$, i.e. that
\[ {}^\beta \alpha^\loc = \alpha. \]
Our verification is a rather ugly cocyle computation. We need to check that
\[  \beta(\sigma_1\sigma_2) \beta(\sigma_1)^{-1} {}^{\sigma_1} \beta(\sigma_2)^{-1} \alpha^\loc(\sigma_1,\sigma_2) = \alpha(\sigma_1,\sigma_2) \]
or, after projecting under $\pi_{\eta w}$, that
\[ \begin{array}{rl} &\pi_{\eta w} \alpha(\sigma_1,\sigma_2)  \\ =& \pi_{\eta w} (\beta(\sigma_1\sigma_2) \beta(\sigma_1)^{-1}) {}^{\sigma_1} \pi_{\sigma_1^{-1}\eta w}(\beta(\sigma_2))^{-1} \\ &  {}^{s_w(\eta)} \alpha_w (s_w(\eta)^{-1} \sigma_1 s_w(\sigma_1^{-1} \eta), s_w(\sigma_1^{-1}\eta)^{-1} \sigma_2 s_w(\sigma_2^{-1} \sigma_1^{-1} \eta))
.   \end{array} \]
The right hand side equals
\[ \begin{array}{l} \pi_{\eta w} (\alpha(s_w(\eta),s_w(\eta)^{-1}\sigma_1\sigma_2 s_w(\sigma_2^{-1}\sigma_1^{-1}\eta))\alpha(\sigma_1,s_w(\sigma_1^{-1} \eta))) \\ 
\pi_{\eta w}(\alpha(s_w(\eta),s_w(\eta)^{-1}\sigma_1 s_w(\sigma_1^{-1}\eta))\alpha(\sigma_1\sigma_2,s_w(\sigma_2^{-1} \sigma_1^{-1}\eta)))^{-1} \\
{}^{\sigma_1}\pi_{\sigma_1^{-1}\eta w} (\alpha(\sigma_2,s_w(\sigma_2^{-1} \sigma_1^{-1}\eta))/\alpha(s_w(\sigma_1^{-1}\eta),s_w(\sigma_1^{-1}\eta)^{-1}\sigma_2 s_w(\sigma_2^{-1}\sigma_1^{-1}\eta)))\\
{}^{s_w(\eta)} \alpha_w (s_w(\eta)^{-1} \sigma_1 s_w(\sigma_1^{-1} \eta), s_w(\sigma_1^{-1}\eta)^{-1} \sigma_2 s_w(\sigma_2^{-1} \sigma_1^{-1} \eta)) \\
{}^{s_w(\eta)}(\gamma_w(s_w(\eta)^{-1}\sigma_1 s_w(\sigma_1^{-1} \eta)) / (\gamma_w(s_w(\eta)^{-1}\sigma_1\sigma_2 s_w(\sigma_2^{-1} \sigma_1^{-1} \eta)))   \\ {}^{\sigma_1 s_w(\sigma_1^{-1})} \gamma_w(s_w(\sigma_1^{-1}\eta)^{-1}\sigma_2 s_w(\sigma_2^{-1}\sigma_1^{-1} \eta))
\end{array}\]
or
\[ \begin{array}{l} \pi_{\eta w} (\alpha(s_w(\eta),s_w(\eta)^{-1}\sigma_1\sigma_2 s_w(\sigma_2^{-1}\sigma_1^{-1}\eta))\alpha(\sigma_1,s_w(\sigma_1^{-1} \eta)))\\
\pi_{\eta w}(\alpha(s_w(\eta),s_w(\eta)^{-1}\sigma_1 s_w(\sigma_1^{-1}\eta))\alpha(\sigma_1\sigma_2,s_w(\sigma_2^{-1} \sigma_1^{-1}\eta)))^{-1} \\
\pi_{\eta w} ({}^{\sigma_1}\alpha(\sigma_2,s_w(\sigma_2^{-1} \sigma_1^{-1}\eta))/{}^{\sigma_1}\alpha(s_w(\sigma_1^{-1}\eta),s_w(\sigma_1^{-1}\eta)^{-1}\sigma_2 s_w(\sigma_2^{-1}\sigma_1^{-1}\eta)))\\
{}^{s_w(\eta)} \pi_w \alpha (s_w(\eta)^{-1} \sigma_1 s_w(\sigma_1^{-1} \eta), s_w(\sigma_1^{-1}\eta)^{-1} \sigma_2 s_w(\sigma_2^{-1} \sigma_1^{-1} \eta)), 
\end{array}\]
which in turn equals $\pi_{\eta w}$ of
\[ \begin{array}{l} \alpha(s_w(\eta),s_w(\eta)^{-1}\sigma_1\sigma_2 s_w(\sigma_2^{-1}\sigma_1^{-1}\eta))\alpha(\sigma_1,s_w(\sigma_1^{-1} \eta))/\alpha(s_w(\eta),s_w(\eta)^{-1}\sigma_1 s_w(\sigma_1^{-1}\eta)) \\ \alpha(\sigma_1\sigma_2,s_w(\sigma_2^{-1} \sigma_1^{-1}\eta)) 
{}^{\sigma_1}\alpha(\sigma_2,s_w(\sigma_2^{-1} \sigma_1^{-1}\eta))/({}^{\sigma_1}\alpha(s_w(\sigma_1^{-1}\eta),s_w(\sigma_1^{-1}\eta)^{-1}\sigma_2 s_w(\sigma_2^{-1}\sigma_1^{-1}\eta))\\
{}^{s_w(\eta)} \alpha (s_w(\eta)^{-1} \sigma_1 s_w(\sigma_1^{-1} \eta), s_w(\sigma_1^{-1}\eta)^{-1} \sigma_2 s_w(\sigma_2^{-1} \sigma_1^{-1} \eta))). \end{array} \]
We can rewrite this
\[ \begin{array}{l} \alpha(s_w(\eta),s_w(\eta)^{-1}\sigma_1\sigma_2 s_w(\sigma_2^{-1}\sigma_1^{-1}\eta))\alpha(\sigma_1,s_w(\sigma_1^{-1} \eta))\\ 
(\alpha(s_w(\eta),s_w(\eta)^{-1}\sigma_1 s_w(\sigma_1^{-1}\eta))\alpha(\sigma_1\sigma_2,s_w(\sigma_2^{-1} \sigma_1^{-1}\eta)))^{-1} \\
(\alpha(\sigma_1,\sigma_2 s_w(\sigma_2^{-1}\sigma_1^{-1}\eta)) \alpha(\sigma_1 \sigma_2,s_w(\sigma_2^{-1} \sigma_1^{-1}\eta))\alpha(\sigma_1,\sigma_2)) \\
(\alpha(\sigma_1s_w(\sigma_1^{-1}\eta),s_w(\sigma_1^{-1}\eta)^{-1}\sigma_2 s_w(\sigma_2^{-1}\sigma_1^{-1}\eta)) \alpha(\sigma_1, s_w(\sigma_1^{-1}\eta)) \alpha(\sigma_1, \sigma_2s_w(\sigma_2^{-1} \sigma_1^{-1}\eta)))^{-1} \\
\alpha ( \sigma_1 s_w(\sigma_1^{-1} \eta), s_w(\sigma_1^{-1}\eta)^{-1} \sigma_2 s_w(\sigma_2^{-1} \sigma_1^{-1} \eta)) \alpha (s_w(\eta), s_w(\eta)^{-1} \sigma_1 s_w(\sigma_1^{-1} \eta))\\ \alpha (s_w(\eta),s_w(\eta)^{-1} \sigma_1  \sigma_2 s_w(\sigma_2^{-1} \sigma_1^{-1} \eta))^{-1}  \end{array} \]
in which almost everything cancels leaving just $\alpha(\sigma_1,\sigma_2)$. The claim follows. 

Let us be still more explicit in a special case. Assume that $F=\Q$ and that $E$ is totally imaginary. We may and will assume that
\begin{itemize}
\item $\alpha(1,1)=1$,
\item $\alpha_w(1,1)=1$ for all $w$,
\item $1 \in H_w$ for all $w$,
\item $\alpha_{w(\infty)}(\sigma_1,\sigma_2)=\left\{ \begin{array}{ll} -1 & {\rm if}\,\, \sigma_1=\sigma_2=c_{w(\infty)} \\ 1 & {\rm otherwise},
\end{array}\right.$
\item $\gamma_{w(\infty)}\equiv 1$,
\item and $\alpha(\sigma,c_{w(\infty)})=1$ if $\sigma \in H_\infty$.
\end{itemize}
(To achieve the last of these we replace $\alpha$ by ${}^\gamma \alpha$ where $\gamma(1)=\gamma(c_{w(\infty)})=1$ and $\gamma(\sigma c_{w(\infty)})=\alpha(\sigma,c_{w(\infty)})\gamma(\sigma)$ if $\sigma \in H_\infty$.) Then
\begin{itemize}
\item $\alpha(\sigma,1)=\alpha(1,\sigma)=1$;
\item $\alpha(c_{w(\infty)},c_{w(\infty)})_{w(\infty)}=-1_{w(\infty)}$;
\item $\alpha(\sigma,c_{w(\infty)})_{\sigma w(\infty)}=-1_{\sigma w(\infty)}$ if $\sigma \not\in H_\infty$;
\item $\alpha(\sigma_1,\sigma_2 c_{w(\infty)})_{\sigma_1\sigma_2 w(\infty)}=\alpha(\sigma_1,\sigma_2)_{\sigma_1\sigma_2 w(\infty)}\left\{ \begin{array}{ll} -1_{\sigma_1\sigma_2 w(\infty)} & {\rm if}\,\, \delta_\infty (\sigma_2) \neq \delta_\infty(\sigma_1\sigma_2)\\ 1 & {\rm if}\,\,  \delta_\infty (\sigma_2) = \delta_\infty(\sigma_1\sigma_2). \end{array} \right. $
\end{itemize}
(For the penultimate of these note that $\alpha(\sigma,c_{w(\infty)})={}^\sigma\alpha(c_{w(\infty)},c_{w(\infty)})\alpha(\sigma,1)/\alpha(\sigma c_{w(\infty)},c_{w(\infty)})$, and for the ultimate one use the cocycle relation.)
Thus
\[ \begin{array}{rcl} \pi_{\eta w(\infty)}\beta(\sigma) &= &\pi_{\eta w(\infty)} ( \alpha (\eta \delta_\infty(\eta), \delta_\infty(\eta)^{-1}\delta_\infty(\sigma^{-1}\eta))/ \alpha(\sigma, \sigma^{-1}\eta \delta_\infty(\sigma^{-1}\eta)))\\ &\equiv& \pi_{\eta w(\infty)} (\alpha(\sigma, \sigma^{-1}\eta))^{-1}  \left\{ \begin{array}{ll} -1_{\eta w(\infty)} & {\rm if}\,\, \delta_\infty(\sigma^{-1}\eta)\neq 1 \\
1 & {\rm if}\,\, \delta_\infty(\sigma^{-1}\eta)= 1
\end{array}\right. 
\end{array} \]
for all $\eta\in H_\infty$ and $\sigma \in \Gal(E/\Q)$. 

Although this doesn't fix $\balpha$ uniquely, we will use $\balpha_0$ to denote the element of $\cZ(E/\Q)$ arising from such a choice. Also choose $\rho_0:E_{w(\infty)} \iso \C$. Then we get a canonical identification
\[ \tTheta_0: W_{\C/\R}\liso W_{E_{w(\infty)}/\Q_\infty,\balpha_0}. \]
Moreover
\[ \iota_{w(\infty)}^{\balpha_0} (e_{\balpha}^\loc (c_{w(\infty)}))=e_{\balpha_0}^\glob(c_{w(\infty)}). \]

\newpage

\section{Kottwitz cohomology: the global case}\label{bgglob2}

\subsection{Algebraic cohomology}

We will be concerned with the following algebraicity conditions:

\begin{enumerate}
\item For the algebraic cohomology of $\cE^\loc(E/F)_{D,S,\ga}$ we will use the algebraicity conditions
\[ \cN_S=\{(\nu_w)_{w\in S_E}:\,\, \nu_w \in X_*(G)(D_w)\,\, {\rm and}\,\, \nu_w=1\,\, {\rm for\,\, all\,\, but \,\, finitely\,\, many}\,\, w\} \]
and
\[ \cN_{S,\basic}=\{ (\nu_w)\in \cN_S:\,\, \nu_w\,\, {\rm factors\,\, through}\,\, Z(G)\,\, \forall w \in S_E\}. \]

\item  For the algebraic cohomology of $\cE_2(E/F)_{D,\ga}$ we will use the algebraicity conditions
\[ \cN=\{ \nu \in \Hom(T_{2,E}, G)( \A_D):\,\, \nu\,\, {\rm is}\,\, G(\A_D)-{\rm conjugate\,\, to \,\, an\,\, element\,\, of} \,\,\Hom(T_{2,E},G)(D) \} \]
and
\[ \cN_\basic=\Hom(T_{2,E},Z(G))(D) . \]

\item For the algebraic cohomology of 
$W_{E_w/F_v,D_u,\ga}$ we will use the algebraicity conditions
\[ \cN_S=X_*(G)(D_u) \]
and
\[ \cN_{S,\basic}=X_*(Z(G))(D_u). \]

\item\label{case4} For the algebraic cohomology of $\cE^\glob(E/F)_{D,\ga}$ we will use the algebraicity conditions
\[ \cN=\{ \nu \in \Hom(T_{2,E}, G)( \A_D):\,\, \nu\,\, {\rm is}\,\, G(\A_D)-{\rm conjugate\,\, to \,\, an\,\, element\,\, of} \,\,\Hom(T_{2,E},G)(D) \} \]
and
\[ \cN_\basic=\Hom(T_{2,E},Z(G))(D) . \]

\item For the algebraic cohomology of $\cE_3(E/F)_{D,\ga}$ we will use the algebraicity conditions
\[ \cN=\Hom(T_{3,E},G)(D) \]
and 
\[ \cN_\basic=\Hom(T_{3,E},Z(G))(D). \]

\end{enumerate}
We will denote the corresponding pointed sets of cycles and cohomology classes $Z^1_\alg(\cE^?(E/F)_{D,\ga},G(A_D))$ and $Z^1_\alg(\cE^?(E/F)_{D,\ga},G(A_D))_\basic$ and $H^1_\alg(\cE^?(E/F)_D,G(A_D))$ and $H^1_\alg(\cE^?(E/F)_D,G(A_D))_\basic$, where 
\begin{enumerate}
\item For $\cE^\loc(E/F)_{D,S,\ga}$ we will use $A_D=\A_{D,S}$.

\item  For $\cE_2(E/F)_{D,\ga}$ we will use $A_D=\A_D$.

\item For 
$W_{E_w/F_v,D_u,\ga}$ we will use $A_D=D_u$.

\item For $\cE^\glob(E/F)_{D,\ga}$ we will use $A_D=\A_D$. 

\item For $\cE_3(E/F)_{D,\ga}$ we will use $A_D=D$.

\end{enumerate}
We have a canonical isomorphism
\[ H^1_\alg(W_{E_w/F_v,D_u,\ga},G(D_u)) \cong H^1_\alg(W_{E_w/F_v,D_u},G(D_u)) \]
preserving basic subsets and products.

\begin{lem} In each case $\nu$ is determined by $\barnu$. \end{lem}

\pfbegin In all cases except case \ref{case4} this follows because the the rational points in a split torus are Zariski dense. In case \ref{case4} suppose that $T$ is a split torus over a number field $D$ and that $\nu_1,\nu_2:T \ra G$ are homomorphisms, defined over $D$, to an algebraic group $G$. Suppose moreover that $g_1,g_2 \in G(\A_D)$ with $g_1\nu_1 g_1^{-1}=g_2\nu_2g_2^{-1}$ on $T(D)$. We need to prove that $g_1\nu_1 g_1^{-1}=g_2\nu_2g_2^{-1}$. We immediately reduce to the case $T=\G_m$ in which case it suffices to show that $D^\times$ is Zariski dense in $\G_{m/\A_D}$, i.e. that if $f \in \A_D[X,X^{-1}]$ and $f$ vanishes on $D^\times$, then $f=0$. As $D$ is infinite, this follows easily by using Vandermonde determinants.\pfend

The map $\loc_\ga$ induces an isomorphism
\[ \loc_\ga^*: Z^1_\alg(\cE_2(E/F)_{D,\ga},G(\A_D)) \liso Z^1_\alg(\cE^\glob(E/F)_{D,\ga},G(\A_D)) \]
which is functorial in $G$ and preserves basic subsets and passes to cohomology. We will denote its inverse simply
\[ \loc_\ga: Z^1_\alg(\cE^\glob(E/F)_{D,\ga},G(\A_D)) \liso Z^1_\alg(\cE_2(E/F)_{D,\ga},G(\A_D)). \]
We will also denote by $\loc_\ga$ the composite map
\[ \begin{array}{rcl} \loc_\ga: Z^1_\alg(\cE_3(E/F)_{D,\ga},G(\A_D)) &\lra& Z^1_\alg(\cE^\glob(E/F)_{D,\ga},G(\A_D)) \\ &\liso& Z^1_\alg(\cE_2(E/F)_{D,\ga},G(\A_D)), \end{array} \]
which is again functorial in $G$ and preserves basic subsets and passes to cohomology. We have $\loc_\ga(\phi_1\phi_2)=\loc_\ga(\phi_1)\loc_\ga(\phi_2)$. 

 If $S' \supset S$ then there are natural maps
\[ \begin{array}{rcl} Z^1_\alg(\cE^\loc(E/F)_{D, S,\ga},G(\A_{D,S})) &\lra& Z^1_\alg(\cE^\loc(E/F)_{D,S',\ga},G(\A_{D,S'})) \\ &\stackrel{\res_S}{\lra} & Z^1_\alg(\cE^\loc(E/F)_{D,S,\ga},G(\A_{D,S})) \end{array} \]
with composite the identity. 
The first arises from functoriality \ref{changee} of section \ref{algcoh} and the natural maps $\cE^\loc(E/F)_{D,S',\ga} \onto \cE^\loc(E/F)_{D,S,\ga}$ and $G(\A_{D,S}) \into G(\A_{D,S'})$ and $\cN_S \into \cN_{S'}$.
The second arises from functoriality \ref{changec} of section \ref{algcoh} and the maps $\cE^\loc(E/F)_{D,S',\ga} \onto \cE^\loc(E/F)_{D,S,\ga}$ and $G(\A_{D,S'}) \onto G(\A_{D,S})$ and $\cN_{S'} \onto \cN_S$.
These maps are functorial in $G$, preserve basic subsets, commute with products and pass to cohomology.
They give rise to isomorphisms
\[ Z^1_\alg(\cE^\loc(E/F)_{D,S,\ga},G(\A_{D,S})) \liso {\prod_{v \in S}}' Z^1_\alg(\cE^\loc(E/F)_{D,\{v\},\ga},G(D_v) ) \]
and
\[ H^1_\alg(\cE^\loc(E/F)_{D,S,\ga},G(\A_{D,S})) \liso {\prod_{v \in S}}' H^1_\alg(\cE^\loc(E/F)_{D,\{v\},\ga },G(D_v) ) \]
where the product is restricted with respect to the subsets (defined for almost all $v$) $Z^1(\Gal(D/F),G(\cO_{D,v}))$ and $H^1(\Gal(D/F),G(\cO_{D,v}))$ respectively. These isomorphisms preserve basic subsets.
If $G$ is connected, the right hand side of the second of these isomorphisms becomes simply $\oplus_{v \in S}H^1_\alg(\cE^\loc(E/F)_{D,\{v\},\ga},G(D_v))$. (See the corollary to theorem 6.8 of \cite{pr}.) 

If $u|w|v\in S$ are place of $D$, $E$ and $F$ respectively, then we get a map
\[ \begin{array}{rcl} \res_u:Z^1_\alg(\cE^\loc(E/F)_{D,S,\ga},G(\A_{D,S})) &\lra& Z^1_{X_*(G)(D_w)} (\cE^\loc(E/F)_{D,S,\ga}|_{\Gal(D/F)_u},G(\A_{D,S}))\\ & \lra& Z^1_\alg(W_{E_w/F_v, D_u,\ga},G(D_u)), \end{array}\]
the first map coming from functoriality \ref{changee} and the second functoriality \ref{changec} of section \ref{algcoh}. This map is functorial in $G$, preserves basic subsets, commutes with products and passes to cohomology. It follows from lemma \ref{shapirolem} that 
\[ \res_u: H^1_\alg(\cE^\loc(E/F)_{D,\{v\},\ga},G(D_v) \liso H^1_\alg(W_{E_w/F_v,D_u,\ga},G(D_u)), \]
and similarly for basic subsets.

 Note that $\res_{\sigma u}=\sigma_* \circ \res_u$. We deduce that if $D_v^0/E_v^0/F_v$ are finite extensions abstractly isomorphic to $D_u/E_w/F_v$ for any, and hence all, $u|w|v \in S$, then we obtain a natural map 
\[ \res_{D_v^0/F_v}: H^1_\alg(\cE^\loc(E/F)_{D,S},G(\A_{D,S})) \lra H^1_\alg(W_{E_v^0/F_v,D_u},G(D_v^0)), \]
defined as $\tau_* \circ \res_u$ for any $u|w|v$ and any $\tau:D_u \iso D_v^0$ over $F_v$. (The point being that $\res_{D_u^0/F_v}$ does not depend on the choice of $u$ or $\tau$.) This is functorial in $G$, preserves basic subsets, commutes with products and is an isomorphism if $S=\{ v\}$.
If $G$ is connected, then
\[  H^1_\alg(\cE^\loc(E/F)_{D,S},G(\A_{D,S})) \liso \bigoplus_{v \in S} H^1_\alg(W_{E_v^0/F_v,D_v^0},G(D_v^0)). \]
The composite
\[ \res_{E_v^0/F_v} \circ \loc: H^1_\alg(\cE_3(E/F),G(E)) \lra H^1_\alg(W_{E_v^0/F_v},G(E_v^0)) \]
coincides with the localization map defined by Kottwitz in \cite{kotbg}. 

If $\ga, \ga' \in \cH(E/F)_D$ and if $t \in T_{2,E}(\A_E)$ with $\ga'={}^t\ga$, then we get isomorphisms
\[ z_t=((\gz_t)^{-1})^*: Z^1_\alg(\cE^?(E/F)_{D,\ga},G(A_D)) \liso Z^1_\alg(\cE^?(E/F)_{D,\ga'},G(A_D)) \]
which preserve basic subsets. These are functorial in $G$, commute with taking products and pass to cohomology (as ${}^gz_t(\phi)=z_t({}^g\phi)$).
Moreover $z_t$ commutes with the maps $\res_S$ and $\res_w$, while
\[ \loc_{{}^t\ga} z_t(\phi) = {}^{\bnu_\phi(t)} z_t(\loc_\ga \phi).  \]
If $a \in T_{2,E}(\A_F)$ and $b \in \cE^\glob(E/F)_D^0$ and $c \in \cE^\glob(E/F)_D^0$, then
\[ z_{abct}(\phi)= {}^{\bnu_\phi(b)} z_t(\phi) \] 
if $\phi \in Z^1_\alg(\cE_3(E/F)_{D,\ga},G(D))$ or $Z^1_\alg(\cE^\glob(E/F)_{D,\ga},G(\A_D))$, while
\[ z_{abct}(\phi) = {}^{\bnu_\phi(c^{-1})} z_t(\phi) \]
if $\phi$ is in any of the other groups of cocycles $Z^1_\alg(\cE^?(E/F)_{D,\ga},G(A_D))$.
Thus on the level of cohomology $z_t$ is independent of $t$ and only depends on $\ga$ and $\ga'$. In particular
$H^1_\alg(\cE^?(E/F)_{D,\ga},G(A_D))$ and $H^1_\alg(\cE^?(E/F)_{D,\ga},G(A_D))_\basic$ are canonically independent of $\ga$, so we will drop the $\ga$ in the notation and denote them simply $H^1_\alg(\cE^?(E/F)_{D},G(A_D))$ and $H^1_\alg(\cE^?(E/F)_{D},G(A_D))_\basic$. Moreover the map 
\[ \loc_\ga:H^1_\alg(\cE_3(E/F)_{D},G(D)) \lra H^1_\alg(\cE_2(E/F)_{D },G(\A_D))\]
is independent of $\ga$ so we will denote it simply $\loc$. The same is true for $\res_S$ and $\res_u$. 

Now suppose that $C \supset D \supset E \supset F$ are finite Galois extensions of $F$. We have maps
\[ \begin{array}{rcl} \inf_{C/D}: Z^1_\alg(\cE^?(E/F)_{D,\ga},G(A_D)) &\lra&  Z^1_\alg(\cE^?(E/F)_{D,\ga}|_{\Gal(C/F)},G(A_C)) \\ &\lra&  Z^1_\alg(\cE^?(E/F)_{C,\inf_{C/D}\ga},G(A_C)) \end{array} \]
where the first map arises from functoriality \ref{changee} (and the map $\cE^?(E/F)_{D,\ga}|_{\Gal(C/F)} \onto \cE^?(E/F)_{D,\ga}$) and the second from functoriality \ref{changec} of section \ref{algcoh} (and the map $\cE^?(E/F)_{D,\ga}|_{\Gal(C/F)} \ra \cE^?(E/F)_{C,\inf_{C/D}\ga}$). These maps are functorial in $G$, preserve basic subsets, commute with products, and pass to cohomology. We also have
\[ {\inf}_{C/D} \circ \loc_\ga = \loc_{\inf_{C/D} \ga} \circ {\inf}_{C/D} \]
and $\inf_{C/D}$ also commutes with $\res_S$, $\res_u$ and $z_t$. We also have maps
\[ \eta_{D/E}^*: Z^1_\alg(\cE^?(E/F)_{C,\eta_{D/E,*}\ga},G(A_C)) \lra  Z^1_\alg(\cE^?(D/F)_{C,\ga},G(A_C)) \]
induced by functoriality \ref{changee} of section \ref{algcoh} (and the map $\cE^?(D/F)_{C,\ga} \ra \cE^?(E/F)_{C,\eta_{D/E,*}\ga}$). 
These maps are functorial in $G$, preserve basic subsets, commute with products, and pass to cohomology. We also have
\[ \eta_{D/E}^* \circ \loc_{\eta_{D/E,*}\ga} = \loc_{ \ga} \circ \eta_{D/E}^* \]
and
\[ \eta_{D/E}^* \circ z_{\eta_{D/E}(t)} = z_t\circ \eta_{D/E}^*. \]
Moreover $\eta_{D/E}^*$ commutes with $\res_S$ and $\res_u$.

If $\ga_D \in \cH(D/F)$ and $\ga_E \in \CH(E/F)$ the we can find a $t\in T_{2,E}(\A_D)$ with $\eta_{D/E,*}\ga_D={}^t \inf_{D/E} \ga_E$. Then we set
\[ \begin{array}{rcl}  \inf_{D/E,t}: Z^1_\alg(\cE^?(E/F)_{\ga_E},G(A_E)) & \stackrel{\inf_{D/E}}{\lra} & Z^1_\alg(\cE^?(E/F)_{D,\inf_{D/E}\ga_E},G(A_D)) \\ &\stackrel{z_t}{\lra}&
Z^1_\alg(\cE^?(E/F)_{D,\eta_{D/E,*} \ga_D},G(A_D)) \\ 
&\stackrel{\eta_{D/E}^*}{\lra} & Z^1_\alg(\cE^?(D/F)_{\ga_D},G(A_D)).
\end{array} \]
this map is functorial in $G$, preserves basic subsets, commutes with products, and passes to cohomology, where it is independent of $\ga_D$, $\ga_E$ and $t$ so we will denote it simply $\inf_{D/E}$. These maps are all injective even on the level of cohomology. (The first map because it is inflation, the second because it is an isomorphism, and the third from its definition.) 
If $a \in T_{2,E}(\A_F)$ and $b \in \cE^\glob(E/F)_D^0$ and $c \in \cE^\glob(E/F)_D^0$, then
\[ \inf_{D/E, abct}(\phi)= {}^{\bnu_\phi(b)} \inf_{D/E,t}(\phi) \] 
if $\phi \in Z^1_\alg(\cE_3(E/F)_{\ga_E},G(E))$ or $Z^1_\alg(\cE^\glob(E/F)_{\ga_E},G(\A_E))$, while
\[ \inf_{D/E, abct}(\phi) = {}^{\bnu_\phi(c^{-1})} \inf_{D/E,t}(\phi) \]
if $\phi$ is in any of the other groups of cocycles $Z^1_\alg(\cE^?(E/F)_{\ga_E},G(A_E))$.
Note that
\[ {\inf}_{E/E,t}=z_t \]
and
\[ \loc_{\ga_D}  ({\inf}_{D/E,t}(\phi))= {}^{\bnu_\phi(t)} {\inf}_{D/E,t}(\loc_{\ga_E} \phi) \]
and
\[ \bnu_{{\inf}_{D/E,t} \phi}= \bnu_\phi \circ \eta_{D/E} . \]
If $C\supset D$ is another finite Galois extension of $F$ and if $\ga_{C}\in \cH(C/F)$ and if $t' \in T_{2,D}(\A_{C})$ with $\eta_{C/D,*} \ga_{C}={}^{t'} \inf_{C/D} \ga_D$, then
\[ {\inf}_{C/D,t'} \circ {\inf}_{D/E,t} = {\inf}_{C/E, t \eta_{D/E}(t')} .\]

We define 
\[ B(F,G) = \lim_{\ra, E/F} H^1_\alg(\cE_3(E/F),G(E))  \]
and
\[ B^\loc(F,G)_S = \lim_{\ra,E/F} H^1_\alg(\cE^\loc(E/F)_S,G(\A_{E,S})), \]
and similarly for the basic subsets.
We have maps
\[ \loc: B(F,G) \lra B^\loc(F,G) \]
and
\[ \res_S: B^\loc(F,G) \lra B^\loc(F,G)_S \]
and, for $v \in S$,
\[ \res_v :B^\loc(F,G)_S \lra B(F_v,G) \]
for any place $v$ of $F$. These maps preserve the basic subsets. Note that 
\[ \res_v :B^\loc(F,G)_{\{v\}} \liso B(F_v,G) \]
and
\[ B^\loc(F,G)_S \liso {\prod_{v \in S}}' B(F_v,G), \]
where the product is restricted with respect to $\{ H^1(\Gal(F_v^\nr/F_v),G(\cO_{F_v^\nr}))\}$, where $F_v^\nr$ denotes the maximal unramified extension of $F_v$. If $G$ is connected, then
\[ B^\loc(F,G)_S \liso \bigoplus_{v \in S} B^\loc(F_v,G)_{v}. \]
These all preserve basic subsets.

\subsection{The algebraic cohomology of reductive groups}

{\em Now suppose that $G$ is reductive.} 

If $\balpha=(\alpha^\glob,\alpha^\loc,\beta) \in \ga \in \cH(E/F)$ and $T$ over $F$ is a torus split over $E$ then 
\[ \corr^\glob_\balpha=\corr_{\alpha^\glob}: \Z[V_E]_0 \otimes X_*(T)=\Hom(T_{3,E},T)\lra  Z^1_\alg(\cE_3(E/F)_{\ga}, T(E)) \]
and
\[ \corr_\balpha^\loc=\corr_{\alpha^\loc}: \Z[V_E] \otimes X_*(T)\lra  Z^1_\alg(\cE^\loc(E/F)_{\ga}, T(\A_E)) \]
induce bijections
\[ \corr: (\Z[V_E]_0 \otimes X_*(T))_{\Gal(E/F)} \liso  H^1_\alg(\cE_3(E/F), T(E)) \]
and
\[ \corr: (\Z[V_E] \otimes X_*(T))_{\Gal(E/F)} \liso  H^1_\alg(\cE^\loc(E/F), T(\A_E)) \]
which are independent of $\balpha$. 
(Note that $\corr_{\balpha}$ depends on $\balpha$ and not just its image in $\cH(E/F)$.) (For the global case see formula (9.1) of \cite{kotbg}. The local case follows from the corresponding result for local fields.)

We have the following special case of the general observation made in item \ref{changecor} of section \ref{algcoh}:
\begin{lem}\label{corr} Suppose that $\balpha=(\alpha^\glob,\alpha^\loc,\beta) \in \ga \in \cH(E/F)$. Suppose also that $T/F$ is a torus split by $E$ and that $\chi:T_{3,E} \ra T$ is a homomorphism (which must then be defined over $E$). Set
\[ b = \prod_{\eta \in \Gal(E/F)} \eta^{-1} \chi(\beta(\eta))^{-1}. \]
Then
\[ \loc_\ga \corr_{\alpha^\glob}(\chi) = {}^b \corr_{\alpha^\loc}(\chi). \]
\end{lem}

Suppose also that 
\begin{itemize}
\item $G$ contains a maximal torus (defined over $F$) which splits over $E$,
\item and $E$ is totally complex.
\end{itemize} 
Then there are maps
\[ \kappa: H^1_\alg(\cE^\loc(E/F)_S,G(\A_{E,S})) \lra  (\Z[V_{E,S}] \otimes \Lambda_G)_{\Gal(E/F)}  \]
and
\[ \kappa: H^1_\alg(\cE_3(E/F),G(E)) \lra (\Z[V_E]_0 \otimes_\Z \Lambda_G)_{\Gal(E/F)} \]
with the following properties:
\begin{enumerate}
\item They are functorial in $G$.
\item If $G=T$ is a torus then they equal $\corr^{-1}$.
\item $\kappa \circ \res_{S'}$ equals the composition of the natural map $(\Z[V_{E,S}] \otimes \Lambda_G)_{\Gal(E/F)} \ra (\Z[V_{E,S'}] \otimes \Lambda_G)_{\Gal(E/F)}$ with $\kappa$.
\item $\kappa \circ \res_{E_v^0/F_v}$ equals the composition of the natural map
\[ (\Z[V_{E,S}] \otimes \Lambda_G)_{\Gal(E/F)} \lra (\Z[V_{E,\{v\}}] \otimes \Lambda_G)_{\Gal(E/F)} \liso \Lambda_{G,\Gal(E_v^0/F_v)} \]
with $\kappa$.

\item $\kappa \circ \inf_{D/E}$ equals the composite of the natural isomorphisms
\[ (\Z[V_{D,S}] \otimes \Lambda_G)_{\Gal(D/F)} \liso (\Z[V_{E,S}] \otimes \Lambda_G)_{\Gal(E/F)}\]
or 
\[ (\Z[V_D]_0 \otimes_\Z \Lambda_G)_{\Gal(D/F)} \liso (\Z[V_E]_0 \otimes_\Z \Lambda_G)_{\Gal(E/F)}\]
(induced by the maps $\Z[V_D] \ra \Z[V_E]$ sending $u \mapsto u|_E$)
with $\kappa$.  

\item $\kappa \circ \loc$ equals the composition of the natural map $(\Z[V_{E}]_0 \otimes \Lambda_G)_{\Gal(E/F)} \ra (\Z[V_{E}] \otimes \Lambda_G)_{\Gal(E/F)}$ with $\kappa$.
\end{enumerate}
We will denote by $\barkappa_G$ the composite
\[ H^1_\alg(\cE^\loc(E/F)_S,G(\A_{E,S})) \stackrel{\kappa}{\lra}  (\Z[V_{E,S}] \otimes \Lambda_G)_{\Gal(E/F)} \onto \Lambda_{G,\Gal(E/F)} \]
induced by $\sum_w w \otimes x_w \mapsto \sum_w x_w$.
Note that $\barkappa_G \circ \loc=0$. (In the local case we construct $\kappa$ from the corresponding maps for local fields using the decompositions
\[ H^1_\alg(\cE^\loc(E/F)_S,G(\A_{E,S})) \cong \bigoplus_{v \in S} H^1_\alg(W_{E_v^0/F_v},G(E_v^0)) \]
and
\[ (\Z[V_{E,S}] \otimes \Lambda_G)_{\Gal(E/F)} \cong \bigoplus_{v\in S} \Lambda_{G,\Gal(E_v^0/F_v)}. \]
Then all the above properties, and the construction of $\kappa$ in the global case, can be found in \cite{kotbg}, particularly sections 9 and 11. We warn the reader that our $\barkappa$ has a different meaning from Kottwitz's use of the same symbol.)

Note that we obtain a commutative square:
\[ \begin{array}{rcl} B(F,G) &\stackrel{\kappa}{\lra} & (\Z[V_{\barF}]_0\otimes \Lambda_G)_{\Gal(\barF/F)} \\ \loc \da && \da \\ B^\loc(F,G)  &\stackrel{\kappa}{\lra} & (\Z[V_{\barF}]\otimes \Lambda_G)_{\Gal(\barF/F)}. \end{array}\]

If $E/F$ splits $G$ then Kottwitz proves that there is a cartesian square
\[ \begin{array}{ccc} B(F,G)_\basic & \stackrel{\oplus_{v|\infty} \res_v \circ \loc}{\lra} & \prod_{v|\infty} B(F_v,G)_\basic \\ \kappa \da && \prod_{v|\infty} \kappa \da \\
(\Z[V_E]_0 \otimes \Lambda_G)_{\Gal(E/F)} & \lra & \prod_{v|\infty} \Lambda_{G,\Gal(E_w/F_v)} \\
\sum_w w \otimes \lambda_w & \longmapsto & (\sum_{\sigma \in \Gal(E_\tv/F_v) \backslash \Gal(E/F)} \sigma \lambda_{\sigma^{-1}\tv})_{v|\infty},\end{array} \]
where $\tv|v$.
(See proposition 15.1 of \cite{kotbg}.) In particular the fibres of $\kappa: B(F,G)_\basic \ra (\Z[V_E]_0 \otimes \Lambda_G)_{\Gal(E/F)}$ are finite.

\begin{lem}\label{locglob}  $\loc: B(F,G)_\basic \ra \ker \barkappa \subset B^\loc(F,G)_\basic$ is surjective with finite fibres.\end{lem}

\pfbegin For the surjectivity, if $\barkappa(\bpsi)=1$, then
 we may then lift $\kappa(\bpsi) \in (\Lambda_G \otimes \Z[V_{E}])_{\Gal(E/F)}$ to an element $\lambda \in (\Lambda_G \otimes \Z[V_{E}]_0)_{\Gal(E/F)}$. The image of $\lambda$ in $\prod_{v|\infty} \Lambda_{G,\Gal(E_\tv/F_v)}$ equals the image of $\kappa(\bpsi)$ and hence the image of $\res_\infty \bpsi$. Thus we can find $\bphi \in B^(F,G)_\basic$ with $\res_\infty \loc \bphi=\res_\infty \bpsi$ and $\kappa(\bphi)=\lambda$. This implies that $\kappa \res^\infty \loc \bphi = \kappa \res^\infty\bpsi$, so that $\res^\infty \loc \bphi = \kappa \res^\infty\bpsi$. Thus $\loc \bphi = \bpsi$.
 
If $\bphi,\bphi' \in B(F,G)_\basic$ have the same image in $B^\loc(F,G)_\basic$, then $\bphi'\bphi^{-1} \in B(F,{}^\bphi G)_\basic$ has $\bnu_{\bphi'\bphi^{-1}}=1$ and so $\bphi'\bphi^{-1} \in H^1(\Gal(\barF/F),{}^\bphi G)$. In fact $\bphi'\bphi^{-1} \in \ker^1(\Gal(\barF/F),{}^\bphi G)$, which is finite.\pfend

If $S$ is a finite set of places of $F$, we will write $B(F,G)_{S,\basic}$ for the inverse image in $B(F,G)_\basic$ under $\kappa$ of the image of $\Z[V_{E,S}]_0 \otimes_\Z \Lambda_G$ in $(\Z[V_E]_0 \otimes_\Z \Lambda_G)_{\Gal(E/F)}$. (Here $E/F$ is any finite Galois extension that splits $G$. To see that the definition does not depend on the choice of $E$, use lemma \ref{surje}.) 

\begin{lem}\label{lemS} If $S$ is a finite set of places of $F$, then there is a finite Galois extension $D/F$ such that $B(F,G)_{S,\basic}$ is contained in the image of $H^1_\alg(\cE_3(D/F),G(D))_\basic$.\end{lem}

\pfbegin
Let $E/F$ be a finite Galois extension which splits $G$. Note that $\Lambda_G/X_*(Z(G)^0)$ is finite and hence $(\Z[V_{E,S}]_0 \otimes_\Z \Lambda_G)_{\Gal(E/F)}/(\Z[V_{E,S}]_0 \otimes_\Z \Lambda_{Z(G)^0})_{\Gal(E/F)}$ is finite. We conclude that $B(F,Z(G)^0)_S$ has only finitely many orbits on $B(F,G)_{S,\basic}$. Moreover 
$B(F,Z(G)^0)_S$ is finitely generated (being isomorphic to the image of $\Z[V_{E,S}] \otimes \Lambda_{Z(G)^0}$ in $(\Z[V_E]_0 \otimes_\Z \Lambda_{Z(G)^0})_{\Gal(E/F)}$). Thus there is a finite Galois extension $D/F$ containing $E$ such that $B(F,G)_{S,\basic}$ is contained in the image of $H^1_\alg(\cE_3(D/F),G(D))$.\pfend

Kottwitz also shows that there is a commutative diagram with exact rows
\[\begin{array}{ccccccccc}  (0) &\lra& \ker^1(F,G) &\lra& B(F,G)_\basic& \lra& B^\loc(F,G)& \stackrel{\barkappa_F}{\lra} & (\Lambda_{G})_{\Gal(\barF/F)} \\
&& || && \ua && \ua && || \\
 (0) &\lra &\ker^1(F,G)& \lra &H^1(F,G)& \lra & \bigoplus_{v\in V_F} H^1(F_v,G)& \stackrel{\barkappa_F}{\lra}& (\Lambda_{G})_{\Gal(\barF/F)}. \end{array} \]
(See proposition 15.6 of \cite{kotbg}.)

\newpage

\section{Rigidification data.}\label{addat}

\subsection{Technical lemmas}

Suppose that $D \supset E \supset F$ are finite Galois extensions of a number field $F$, that $u$ is a place of $D$ and that $\balpha=(\alpha^\glob,\alpha^\loc,\beta) \in \cZ(D/F)$. We define
\[ \begin{array}{rcl} \gamma_{E,u,\balpha}: W_{D/F,\balpha} &\lra& \A_D^\times/D^\times \\
ae^\glob_\balpha(\sigma) & \longmapsto & (N_{D/E}(a)/a^{[D:E]}) \prod_{\eta \in \Gal(D/E)} (\alpha^\glob(\eta,\sigma)\beta(\eta)_{\eta u})/(\alpha^\glob(\sigma,\eta){}^\sigma(\beta(\eta)_{\eta u})) \\ &=&
\prod_{\eta \in \Gal(D/E)} e^\glob_\balpha(\eta) (ae^\glob_\balpha(\sigma)) e^\glob_\balpha(\eta \sigma)^{-1})/((a e^\glob_\balpha(\sigma)) e^\glob_\balpha(\eta)e^\glob_\balpha(\sigma \eta)^{-1}) \\ && {}^{(1-\sigma)}\left( \prod_{\eta \in \Gal(D/E)} \beta(\eta)_{\eta u}\right) . \end{array} \]

\begin{lem} With the above notation and assumptions we have:
\begin{enumerate}
\item If $a \in \A_D^\times/D^\times$ then $\gamma_{E,u,\balpha}(a)=N_{D/E}(a)/a^{[D:E]}$. 

\item $\gamma_{E,u,\balpha}(\sigma_1\sigma_2)=\gamma_{E,u,\balpha}(\sigma_1){}^{\sigma_2}\gamma_{E,u,\balpha}(\sigma_2)$, i.e. $\gamma_{E,u,\balpha} \in Z^1(W_{E/F,\balpha}, \A_D^\times/D^\times)$.

\item $\gamma_{E,u,\balpha}$ descends to a map $\bargamma_{E,u,\balpha}: W_{D/F,\balpha}/\Delta_D \ra \A_D^\times/D^\times \Delta_D$.

\item $\gamma_{E,u,\balpha}(\conju_a(\sigma))= (N_{D/E} (a/{}^\sigma a)/(a/{}^\sigma a)^{[D:E]}) \gamma_{E,u,\balpha}(\sigma)$ for $a \in \A_D^\times/D^\times$.

\item ${}^{\sigma} \gamma_{E,u,\balpha}(\tau)=\gamma_{E,u,\balpha}(e^\glob_\balpha(\sigma)\tau e^\glob_\balpha(\sigma)^{-1}) {}^{(\sigma \tau \sigma^{-1}-1)}\gamma_{E,u,\balpha}(e^\glob_\balpha(\sigma))$.

\item $\gamma_{E,\sigma u, \balpha}={}^{\prod_{\eta \in \Gal(D/E)} (\beta(\eta)_{\eta \sigma u}/\beta(\eta)_{\eta u})} \gamma_{E,u,\balpha}$.

\item $\gamma_{E,u,{}^\bgamma \balpha}(i_{\gamma^\glob}(\sigma)) = \gamma_{E,u,\balpha}(\sigma) {}^{(\sigma-1)} \left(\prod_{\eta \in \Gal(D/E)} \gamma^\loc(\eta)_{\eta u} \right)$.

\item $\gamma_{E,u,{}^t\balpha}\circ \gz_t={}^{\prod_{\eta \in \Gal(D/E)} t_{\eta u}/{}^\eta (t_u)}\gamma_{E,u,\balpha}$.

\end{enumerate}\end{lem}

\pfbegin For the first part note that 
\[ \begin{array}{rcl} \gamma_{E,u,\balpha}(a) &=&  (N_{D/E}(a\alpha^\glob(1,1)^{-1})/(a \alpha^\glob(1,1)^{-1})^{[D:E]}) \prod_{\eta \in \Gal(D/E)} \alpha^\glob(\eta,1)/\alpha^\glob(1,\eta)  \\ &=& N_{D/E}(a)/a^{[D:E]}. \end{array} \]

For the second part set
\[ \gamma_{E,\balpha}(ae_\balpha^\glob(\sigma)) = (N_{D/E}(a)/a^{[D:E]}) \prod_{\eta \in \Gal(D/E)} \alpha^\glob(\eta,\sigma)/\alpha^\glob(\sigma,\eta) \]
so that
\[ \gamma_{E,u,\balpha}={}^{\prod_{\eta \in \Gal(D/E)} \beta(\eta)_{\eta u}} \gamma_{E,\balpha}. \]
It suffices to check that $\gamma_{E,\balpha}$ is a 1-cocycle. However
\[ \begin{array}{rl} & \gamma_{E,\balpha}(a_1e^\glob_\balpha(\sigma_1)a_2e^\glob_\balpha(\sigma_2)) \\ =& \gamma_{E,\balpha}(a_1{}^{\sigma_1}a_2\alpha^\glob(\sigma_1,\sigma_2)e^\glob_\balpha(\sigma_1\sigma_2)) \\ =& 
(N_{D/E}(a_1{}^{\sigma_1}a_2)/(a_1{}^{\sigma_1}a_2)^{[D:E]}) \prod_{\eta \in \Gal(D/E)} {}^\eta\alpha^\glob(\sigma_1,\sigma_2) \alpha^\glob(\eta,\sigma_1\sigma_2)/ \alpha^\glob(\sigma_1,\sigma_2)\alpha^\glob(\sigma_1\sigma_2,\eta) \\ =&
(N_{D/E}(a_1{}^{\sigma_1}a_2)/(a_1{}^{\sigma_1}a_2)^{[D:E]}) \prod_{\eta \in \Gal(D/E)} \alpha^\glob(\eta\sigma_1,\sigma_2) \alpha^\glob(\eta,\sigma_1)/ {}^{\sigma_1}\alpha^\glob(\sigma_2,\eta)\alpha^\glob(\sigma_1,\sigma_2\eta) \\ =&
(N_{D/E}(a_1{}^{\sigma_1}a_2)/(a_1{}^{\sigma_1}a_2)^{[D:E]}) \prod_{\eta \in \Gal(D/E)} \alpha^\glob(\sigma_1\eta,\sigma_2) \alpha^\glob(\eta,\sigma_1)/ {}^{\sigma_1}\alpha^\glob(\sigma_2,\eta)\alpha^\glob(\sigma_1,\eta\sigma_2) \\ =&
(N_{D/E}(a_1){}^{\sigma_1}N_{D/E}(a_2))/(a_1{}^{\sigma_1}a_2)^{[D:E]}) \prod_{\eta \in \Gal(D/E)} {}^{\sigma_1}\alpha^\glob(\eta,\sigma_2) \alpha^\glob(\eta,\sigma_1)/ {}^{\sigma_1}\alpha^\glob(\sigma_2,\eta)\alpha^\glob(\sigma_1,\eta) \\ =&
 \gamma_{E,\balpha}(a_1e^\glob_\balpha(\sigma_1)) {}^{\sigma_1}\gamma_{E,\balpha}(a_2e^\glob_\balpha(\sigma_2)).
\end{array} \]

Parts (3), (4) and (5) follow easily from the first two parts. For the fifth we have 
\[ \begin{array}{rcl} {}^{\sigma} \gamma_{E,u,\balpha}(\tau)&=&\gamma_{E,u,\balpha}(e^\glob_\balpha(\sigma))^{-1}\gamma_{E,u,\balpha}(e^\glob_\balpha(\sigma)\tau) \\
& =& \gamma_{E,u,\balpha}(e^\glob_\balpha(\sigma))^{-1} \gamma_{E,u,\balpha}(e^\glob_\balpha(\sigma)\tau e^\glob_\balpha(\sigma)^{-1}) {}^{\sigma \tau \sigma^{-1}}\gamma_{E,u,\balpha}(e^\glob_\balpha(\sigma)) 
\\ &=&  {}^{(\sigma \tau \sigma^{-1}-1)}\gamma_{E,u,\balpha}(e^\glob_\balpha(\sigma))\gamma_{E,u,\balpha}(e^\glob_\balpha(\sigma)\tau e^\glob_\balpha(\sigma)^{-1}). \end{array} \]

The final three parts are straightforward, and we leave their verification to the reader. 
\pfend

\subsection{Rigidification data}\label{sect82}

We have seen that $W_{E/F,D,\ga}$ is isomorphic as an extension to $W_{E^\ab/F,D}$, and so there exists an isomorphism of extensions
\[ \begin{array}{ccccccccc} (0)& \lra & \A_D^\times/D^\times \Delta_E &\stackrel{r_D}{\lra} & W_{E^\ab/F,D}/\Delta_E & \stackrel{\overline{{}}}{\lra} & \Gal(D/F) & \lra &(0) \\ && || && \Gamma \da \wr && || && \\ 
(0)& \lra & \A_D^\times/D^\times \Delta_E&{\lra} & W_{E/F,D,\ga} /\Delta_E & \lra & \Gal(D/F) & \lra &(0) \end{array} \]
Such an isomorphism will be called a {\em Galois rigidification} of $\ga$ if it lifts to an isomorphism of extensions
\[ \tGamma: W_{E^\ab/F,D} \liso W_{E/F,D,\ga}. \]
Because $H^1(\Gal(D/F),\A_D^\times/D^\times)=(0)$ we see that all Galois rigidifcations differ by composition with conjugation by an element of $\A_D^\times/D^\times \onto \A_D^\times /D^\times \Delta_E$. Note however that $H^1(\Gal(D/F),\A_D^\times /D^\times \Delta_E)$ may not vanish. 

Moreover if $t \in \A_D^\times/D^\times$ and $\conju_t \circ \Gamma = \Gamma$ then $t=rs$ with $r \in \A_F^\times/F^\times$ and $s \in \Delta_D$ with ${}^\sigma s/s \in \Delta_E$ for all $\sigma \in \Gal(D/F)$. (To see this note that ${}^\sigma t/t \in \Delta_E \subset \Delta_D$ for all $\sigma \in \Gal(D/F)$. As $H^1(\Gal(D/F),\Delta_D)=(0)$ we deduce that there is an $s\in \Delta_D$ with ${}^\sigma t/t={}^\sigma s/s$ for all $\sigma \in \Gal(D/F)$. Thus $t/s \in \A_F^\times/F^\times$ and the claim follows.)

We say that a Galois rigidification $\Gamma$ is {\em adapted} to an $F$-linear embedding $\rho: E^\ab D\into \barF_v$ if 
\[ \tGamma \circ \theta_{\rho|_{E^\ab}}  = \iota_{w(\rho)}^\ga \circ \Theta: W_{(\rho(E)F_v)^\ab/F_v,\rho, D} \lra W_{E/F,D,\ga} ,\]
for some isomorphism of extensions
\[ \Theta:   W_{(E F_v)^\ab/F_v,\rho, D} \liso W_{E_{w(\rho)}/F_v,D,\ga}, \]
some lifting $\tGamma: W_{E^\ab/F,D} \iso W_{E/F,D,\ga}$ of $\Gamma$ and some choice of $\theta_{\rho|_{E^\ab}}$. 

\begin{lem} Suppose that $\rho:E^\ab D\ra \barFv$ is $F$-linear.
\begin{enumerate}
\item A Galois rigidification adapted to $\rho$ exists.

\item\label{pp2} Any two Galois rigidifications adapted to $\rho$ differ by composition with conjugation by an element of $(\A_E^\times/E^\times)^{\Gal(E/F)_{w(\rho)}}D_{w(\rho)}^\times$.

\item If $\sigma \in \Gal(E^\ab D/F)$, if $\Gamma$ is a Galois rigidification of $\ga$ adapted to $\rho:E^\ab D \into \barF_v$ and if $\balpha=(\alpha^\glob, \alpha^\loc,\beta) \in \ga$ then 
\[ \Gamma^{\sigma,\balpha}= \conju_{(\beta(\sigma|_D^{-1})_{w(\rho\sigma)})} \circ \conju_{e^\glob_\balpha(\sigma|_D^{-1})} \circ \Gamma \circ \conju_\sigma\]
is a Galois rigidification of $\ga$ adapted to $\rho \circ \sigma$. (In the last term we think of 
\[ \sigma \in \Gal(E^\ab D/F) \into \Gal(E^\ab/F) \times_{\Gal(E/F)} \Gal(D/F) \into W_{E^\ab/F,D}/\Delta_E.)\]
If $\bgamma =(\gamma^\glob,\gamma^\loc) \in \cB(E/F)_D$, then we have
\[ \Gamma^{\sigma, {}^\bgamma\balpha} = \conju_{\gamma^\loc(\sigma^{-1})_{w(\rho \sigma)}}^{-1} \circ \Gamma^{\sigma,\balpha}. \]

\item $(\Gamma^{\sigma_1,\balpha})^{\sigma_2,\balpha}= \conju_{\alpha^\loc(\sigma_2^{-1},\sigma_1^{-1})_{w(\rho\sigma_1\sigma_2)} } \circ\Gamma^{\sigma_1\sigma_2,\balpha} $.

\item If $\Gamma$ is a Galois rigidification of $\ga$ adapted to $\rho$ and if $t \in T_{2,E}(\A_D)$, then
\[ {}^t \Gamma = \conju_{t_{w(\rho)}} \circ \gz_t \circ \Gamma \]
is a Galois rigidification of ${}^t\ga$ adapted to $\rho$. Moreover
\[ {}^{t_1}({}^{t_2}\Gamma)={}^{t_1t_2}\Gamma. \]

\item If $t \in T_{2,E}(\A_D)$ fixes $\ga$, then $t=abc$ with 
$a \in \cE^\loc(E/F)_D^0$, $b \in \cE^\glob(E/F)_D^0$ and $c \in T_{2,E}(\A_F)$ and
\[ {}^t\Gamma = \conju_{a_{w(\rho)}c_{w(\rho)}} \circ \Gamma. \]

Moreover if $\Gamma$ and $\Gamma'$ are two Galois rigidifications of $\ga$ adapted to $\rho$ then we can find $t \in \cE^\loc(E/F)_D^0 T_{2,E}(\A_F)$ with ${}^t\Gamma=\Gamma'$. 

\item $({}^t\Gamma)^{\sigma,{}^t\balpha}={}^t(\Gamma^{\sigma,\balpha})$.

\item Suppose that $C \supset D \supset E \supset F$ are finite Galois extensions of $F$, that $\ga \in \cH(E/F)_D$ and that $\rho:E^\ab C\into \barFv$ is $F$-linear. If $\Gamma$ is a Galois rigidification for $\ga$ adapted to $\rho|_{E^\ab D}$, then the map  $\inf_{C/D} \Gamma$ from 
\[ ((\A_C^\times/C^\times \Delta_E) \rtimes (W_{E^\ab/F,D}/\Delta_E)|_{\Gal(C/F)})/(\A_D^\times/D^\times \Delta_E) \]
to
\[ ((\A_C^\times/C^\times\Delta_E) \rtimes (W_{E/F,D,\ga}/\Delta_E)|_{\Gal(C/F)})/(\A_D^\times/D^\times\Delta_E) \]
given by
\[ {\inf}_{C/D} \Gamma:[(a,(\sigma,\tau))]  \longmapsto  [(a,(\Gamma(\sigma),\tau))]
\]
is a Galois rigidification for $\inf_{C/D} \ga$ adapted to $\rho$. 

\item $(\inf_{C/D} \Gamma)^{\sigma, \inf_{C/D}\balpha}=\inf_{C/D} (\Gamma^{\sigma|_{E^\ab D},\balpha})$ and $\inf_{C/D} {}^t \Gamma = {}^t(\inf_{C/D} \Gamma)$. If $B \supset C \supset D \supset E \supset F$ are finite Galois extensions of $F$, then $\inf_{B/C} \circ \inf_{C/D}=\inf_{B/D}$. 

\item Suppose that $\balpha \in \ga \in \cH(D/F)$ and that $\trho:D^\ab \into \barFv$ is $F$-linear. If $\Gamma$ is a Galois rigidification for $\ga$ adapted to $\trho$, then the map 
\[ \eta_{D/E,\balpha,\trho,*}\Gamma: W_{E^\ab/F,D}/\Delta_E \lra W_{E/F,D,\eta_{D/E,*}\balpha}/\Delta_E \]
given by
\[ \begin{array}{ccc} ((\A_D^\times/D^\times \Delta_E) \rtimes W_{D^\ab/F})/(\A_D^\times/D^\times) & \lra & ((\A_D^\times/D^\times\Delta_E) \rtimes W_{D/F,\ga})/(\A_D^\times/D^\times) \\ 
{} [(a,\sigma)] & \longmapsto &  [(a\gamma_{E,u(\trho),\balpha}(\Gamma(\sigma)), \Gamma(\sigma))] \end{array} \]
(where $u(\trho)$ is the place of $D$ induced by $\trho$)
is a Galois rigidification for $\eta_{D/E,*}\ga$ adapted to $\trho|_{E^\ab}$.

\item If $\bgamma=(\gamma^\glob,\gamma^\loc) \in \cB(D/F)$, then $\eta_{D/E,{}^\bgamma\balpha,\trho,*}\Gamma = \conju_{\prod_{\eta \in \Gal(D/E)} \gamma^\loc(\eta)_{\eta u(\trho)}^{-1}} \circ i_{\eta_{D/E} \circ \gamma^\glob} \circ (\eta_{D/E,\balpha,\trho,*}\Gamma)$.

\item If $t \in T_{2,D}(\A_D)$ then $\eta_{D/E,{}^t\balpha,\trho,*}({}^t\Gamma)={}^{\eta_{D/E}(t)} \eta_{D/E,\balpha,\trho,*}(\Gamma)$.

\item If $\sigma \in \Gal(D^\ab/F)$, then 
\[ (\eta_{D/E,\balpha,\trho,*}\Gamma)^{\sigma,\eta_{D/E,*}\balpha} = \conju_{\prod_{\eta \in \Gal(D/E)} \alpha^\loc(\sigma^{-1},\eta)_{\sigma^{-1}\eta u(\trho )}/ \alpha^\loc(\eta,\sigma^{-1})_{\eta u(\trho\sigma)}} \circ \eta_{D/E,\balpha,\trho\sigma,*}(\Gamma^{\sigma,\balpha}).\]

\item Suppose that $C \supset D \supset E \supset F$ are finite Galois extensions, that $\balpha \in \cZ(C/F)$ and that $\Gamma$ is a Galois rigidification of $[\balpha]$ adapted to $\trho:C^\ab \into \barFv$. Suppose also that $\balpha'\in \cZ(D/E)$ and that $\Gamma'$ is a Galois rigidification of $[\balpha']$ adapted to $\trho|_{D^\ab}$. Then we can find $t\in T_{2,D}(\A_C)$ and $\bgamma \in \cB(D/E)_C$ such that
\[ \eta_{C/D,*}\balpha={}^{t,\bgamma} \inf_{C/D} \balpha' \,\,\,\,\,\, {\rm and}\,\,\,\,\,\, \eta_{C/D,\trho,*} \Gamma = {}^t \inf_{C/D} \Gamma' .\]
Moreover 
\[ \eta_{C/E,\trho,*} \Gamma = {}^{\eta_{D/E}(t)} \inf_{C/D} \eta_{D/E,\trho,*} \Gamma' . \]

\end{enumerate}\end{lem}

\pfbegin
For the first part choose isomorphisms of extensions $\tGamma: W_{E^\ab/F,D} \iso W_{E/F,D,\ga}$ and $\Theta : W_{(\rho(E)F_v)^\ab/F_v,\rho|_{E^\ab}, D} \iso W_{E_{w(\rho)}/F_v,D,\ga}$. 
Because $H^1(\Gal(D/F)_{w(\rho)}, \A_D^\times/D^\times)=(0)$, we see that 
\[ \conju_a \circ \tGamma \circ \theta_\rho  =  \iota_{w(\rho)}^\ga \circ \Theta \]
for some $a \in \A_D^\times/D^\times$. Replacing $\Gamma$ by $\conju_{ a} \circ \Gamma$ we have the desired Galois rigidification adapted to $\rho$. 

For the second part suppose that $\Gamma_1$ and $\Gamma_2$ are two Galois rigidifications adapted to $\rho$. Suppose more explicitly that $\Gamma_i$ lifts to an isomorphism of extensions $\tGamma_i: W_{E^\ab/F,D} \iso W_{E/F,D,\ga}$ and that $\Theta_i: W_{(\rho(E)F_v)^\ab/F_v,\rho|_{E^\ab}, D} \iso W_{E_{w(\rho)}/F_v,D,\ga}$ so that
$\tGamma_i \circ \theta_\rho  =  \iota_{w(\rho)}^\ga\circ \Theta_{i}$. Replacing $\tGamma_2$ by a $\Delta_E$-conjugate, we may assume that we have the same choice of $\theta_\rho$ in both equations. Then we can find $b \in D_{w(\rho)}^\times$ such that $\Theta_2=\conju_b \circ \Theta_1$ (as $H^1(\Gal(D/F)_{w(\rho)}, D_{w(\rho)}^\times)=(0)$).  Moreover we can find $a \in \A_D^\times/D^\times$ such that $\tGamma_2 =  \conju_a \circ \tGamma_1$ (as $H^1(\Gal(D/F), \A_D^\times/D^\times)=(0)$). Then $\conju_{ab^{-1}} \circ \tGamma_1 \circ \theta_\rho = \iota_{w(\rho)}^\ga \circ \Theta_1$, so that $ab^{-1} \in (\A_D^\times/D^\times)^{\Gal(D/F)_{w(\rho)}}=(\A_E^\times/E^\times)^{\Gal(E/F)_{w(\rho)}}$, as desired.

For the third part, if $\tGamma: W_{E^\ab/F,D} \iso W_{E/F,D,\ga}$ is an isomorphism of extensions lifting $\Gamma$ and $\tsigma \in W_{E^\ab/F}|_{\Gal(D/F)}$ lifts $\sigma$, then
\[ \tGamma^{\sigma,\balpha}= \conju_{(\beta(\sigma^{-1})_{w(\rho\sigma)})} \circ \conju_{e^\glob_\balpha(\sigma^{-1})} \circ \tGamma \circ \conju_\tsigma:  W_{E^\ab/F,D} \liso W_{E/F,D,\ga}\]
is an isomorphism of extensions lifting $\Gamma^{\sigma,\balpha}$, and so $\Gamma^{\sigma, \balpha}$ is a Galois rigidification. Moreover using lemma \ref{decog} we have
\[ \begin{array}{rcl}  \tGamma^{\sigma, \balpha} \circ \theta_{\rho \sigma|_{E^\ab}} &= &
\conju_{(\beta(\sigma^{-1})_{w(\rho\sigma)})} \circ \conju_{e^\glob_\balpha(\sigma^{-1})} \circ \tGamma \circ \theta_\rho \circ (\sigma|_D)_* \\ & =& \conju_{(\beta(\sigma^{-1})_{w(\rho\sigma)})} \circ \conju_{e^\glob_\balpha(\sigma^{-1})} \circ \iota_{w(\rho)}^\ga \circ \Theta \circ (\sigma|_D)_*\\ &=&  \iota^\ga_{w(\rho\sigma)} \circ (\conju_{e_{\balpha}^\loc(\sigma^{-1})} \circ \Theta \circ (\sigma|_D)_*), \end{array}\]
and $\conju_{e_{\balpha}^\loc(\sigma^{-1})} \circ \Theta \circ (\sigma|_D)_*: W_{(\rho(E)F_v)^\ab/F_v,\rho\sigma,D} \iso W_{E_{w(\rho)}/F_v,D,\ga}$ is an isomorphism of extensions. The second assertion of this part follows immediately from the definitions.

For the fourth part we have
\[ \begin{array}{rl} &\Gamma^{\sigma_1\sigma_2,\balpha} \\ =&
\conju_{(\beta(\sigma_2^{-1}\sigma_1^{-1})_{w(\rho\sigma_1\sigma_2)})} \circ \conju_{e^\glob_\balpha(\sigma_2^{-1}\sigma_1^{-1})} \circ \Gamma \circ \conju_{\sigma_1} \conju_{\sigma_2} \\ =&
 \conju_{(\beta(\sigma_2^{-1}\sigma_1^{-1})_{w(\rho\sigma_1\sigma_2)})} \circ \conju_{\alpha^\glob(\sigma_2^{-1},\sigma_1^{-1})^{-1}}\circ \conju_{e^\glob_\balpha(\sigma_2^{-1})} \circ \conju_{e^\glob_\balpha(\sigma_1^{-1})}\circ \Gamma \circ \conju_{\sigma_1} \conju_{\sigma_2} \\ =&
 
 \conju_{(\beta(\sigma_2^{-1}\sigma_1^{-1})_{w(\rho\sigma_1\sigma_2)})} \circ \conju_{\alpha^\glob(\sigma_2^{-1},\sigma_1^{-1})^{-1}}\circ \conju_{e^\glob_\balpha(\sigma_2^{-1})} \circ \conju_{\beta(\sigma_1^{-1})_{w(\rho \sigma_1)}}^{-1} \circ \Gamma^{\sigma_1,\balpha} \circ \conju_{\sigma_2} \\ =&

\conju_{(\beta(\sigma_2^{-1}\sigma_1^{-1})_{w(\rho\sigma_1\sigma_2)})} \circ \conju_{\alpha^\glob(\sigma_2^{-1},\sigma_1^{-1})^{-1}}\circ  \conju_{({}^{\sigma_2^{-1}}\beta(\sigma_1^{-1}))_{w(\rho \sigma_1\sigma_2)}}^{-1} \circ \conju_{e^\glob_\balpha(\sigma_2^{-1})} \circ \Gamma^{\sigma_1,\balpha} \circ \conju_{\sigma_2} \\ =&

\conju_{(\beta(\sigma_2^{-1}\sigma_1^{-1})\alpha^\glob(\sigma_2^{-1},\sigma_1^{-1})^{-1}{}^{\sigma_2^{-1}}\beta(\sigma_1^{-1})^{-1} \beta(\sigma_2^{-1})^{-1})_{w(\rho \sigma_1\sigma_2)}} \circ  (\Gamma^{\sigma_1,\balpha})^{\sigma_2,\balpha}  \\ =&

\conju_{\alpha^\loc(\sigma_2^{-1},\sigma_1^{-1})^{-1}_{w(\rho \sigma_1\sigma_2)}} \circ  (\Gamma^{\sigma_1,\balpha})^{\sigma_2,\balpha} .
\end{array} \] 

For the fifth part simply recall that $\iota^{{}^t\ga}_w \circ \gz_t =\conju_{t_w} \circ \gz_t \circ \iota^\ga_w$.

The first assertion of the sixth part follows from corollary \ref{stab} and the fact that $\gz_{abc}=\conju_{b}^{-1}$ (see the paragraph before lemma \ref{decog}). The second assertion of the sixth part follows from the first assertion and part \ref{pp2}.

For the seventh part we have
\[ \begin{array}{rcl} ({}^t\Gamma)^{\sigma,{}^t\balpha} &=&
\conju_{(({}^t\beta)(\sigma^{-1})_{w(\rho\sigma)})} \circ \conju_{e^\glob_{{}^t\balpha}(\sigma^{-1})} \circ \conju_{t_{w(\rho)}} \circ \gz_t \circ \Gamma \circ \conju_\sigma \\ &=&
\gz_t \circ \conju_{t_{w(\rho\sigma)}\beta(\sigma^{-1})_{w(\rho\sigma)} /{}^{\sigma^{-1} }t_{w(\rho)} } \circ \conju_{{}^{\sigma^{-1}}t_{w(\rho)}} \circ\conju_{e^\glob_{\balpha}(\sigma^{-1})} \circ   \Gamma \circ \conju_\sigma \\ &=& {}^t(\Gamma^{\sigma,\balpha}).
\end{array} \]

For the eighth part suppose that $\tGamma: W_{E^\ab/F,D} \iso W_{E/F,D,\ga}$ lifts $\Gamma$ and that $\Theta: W_{(\rho(E)F_v)^\ab/F_v,\rho, D} \iso W_{E_{w(\rho)}/F_v,D,\ga}$ satisfies
\[ \tGamma \circ \theta_\rho =\iota^\ga_{w(\rho)}\circ \Theta. \]
Define $\inf_{C/D} \tGamma$ from 
\[ ((\A_C^\times/C^\times ) \rtimes W_{E^\ab/F,D}|_{\Gal(C/F)})/(\A_D^\times/D^\times ) \]
to
\[ ((\A_C^\times/C^\times) \rtimes W_{E/F,D,\ga}|_{\Gal(C/F)})/(\A_D^\times/D^\times) \]
by
\[ {\inf}_{C/D} \tGamma:[(a,(\sigma,\tau))]  \longmapsto  [(a,(\Gamma(\sigma),\tau))];
\]
and
\[ {\inf}_{C/D}\Theta: W_{(\rho(E)F_v)^\ab/F_v,\rho, C} \liso W_{E_{w(\rho)}/F_v,C,\ga}\]
to be
\[ \begin{array}{rcl} (C_{w(\rho)}^\times \rtimes W_{(\rho(E)F_v)^\ab/F_v,\rho, D}|_{\Gal(C/F)_{w(\rho)}})/D_{w(\rho)}^\times & \liso & (C_{w(\rho)}^\times \rtimes W_{E_{w(\rho)}/F_v,D,\ga}|_{\Gal(C/F)_{w(\rho)}})/D_{w(\rho)}^\times \\
{} [(a,(\sigma,\tau))] & \longmapsto & [(a,(\Theta(\sigma),\tau))]. \end{array} \]
Then
\[ ({\inf}_{C/D} \tGamma) \circ \theta_\rho  = \iota_{w(\rho)}^{{\inf}_{C/D} \ga} \circ (\inf_{C/D} \Theta). \]

We leave the straightforward verification of the ninth part to the reader. 

For the tenth part we recall that in $(\A_D^\times/D^\times \rtimes W_{D^\ab/F})/(\A_D^\times/D^\times)$ the implicit map is 
\[ \begin{array}{rcl} \A_D^\times/D^\times &\lra& \A_D^\times/D^\times \rtimes W_{D^\ab/F} \\ a & \lra & (N_{D/E} a^{-1},a),
\end{array} \]
while in $(\A_D^\times/D^\times \rtimes W_{D/F,\ga})/(\A_D^\times/D^\times)$ the implicit map is
\[ \begin{array}{rcl} (\A_D^\times/D^\times) &\lra &\A_D^\times/D^\times \rtimes W_{D/F,\ga} \\
a & \longmapsto & (a^{-[D:E]},a). \end{array}\]
The difficulty here is that these two maps are not directly compatible.

Let $\tGamma: W_{D^\ab/F} \iso W_{D/F,\ga}$ be an isomorphism of extensions lifting $\Gamma$, and 
\[ \Theta: W_{(\trho(D)F_v)^\ab/F_v,\rho} \liso W_{D_{u(\trho)}/F_v,\ga} \]
an isomorphism of extensions such that
\[ \tGamma \circ \theta_\trho =\iota^\ga_{u(\trho)} \circ \Theta. \]
Let
\[ \tGamma^E: W_{E^\ab/F,D} \lra W_{E/F,D,\eta_{D/E,*}\balpha} \]
be given by
\[ \begin{array}{ccc} (\A_D^\times/D^\times \rtimes W_{D^\ab/F,\ga})/(\A_D^\times/D^\times) & \lra & (\A_D^\times/D^\times \rtimes W_{D/F,\ga})/(\A_D^\times/D^\times) \\ 
{} [(a,\sigma)] & \longmapsto &  [(a \gamma_{E,u(\trho),\balpha}(\tGamma(\sigma)), \tGamma(\sigma))] .\end{array} \]
This is well defined because for $a \in \A_D^\times/D^\times$ we have $\gamma_{\balpha,u(\trho),E}(a)=N_{D/E}(a)/a^{[D:E]}$. It is a homomorphism because $\gamma_{\balpha,E}$ is a 1-cocycle. We deduce that $\tGamma^E$ is an isomorphism of extensions. Thus $\eta_{D/E,\balpha,*} \Gamma$ is a rigidification of $\eta_{D/E,*} \ga$.

We also define an isomorphism
\[ \Theta^E: W_{(E F_v)^\ab/F_v, D,\trho} \liso W_{E_{w(\trho)}/F_v,D,\eta_{D/E,*}\ga} \]
to be the map from 
\[  (D_{w(\trho)}^\times \rtimes( W_{(DF_v)^\ab/F_v} \times_{\Gal(E/F)} \Gal(D/F))) /W_{(DF_v)^\ab/(EF_v)} \]
to
\[ (D_{w(\trho)}^\times \rtimes (W_{D_{u(\trho)}/F_v,\balpha} \times_{\Gal(E/F)} \Gal(D/F)))/(W_{D_{u(\trho)}/F_v,\balpha}|_{\Gal(D_{u(\trho)}/E_{w(\trho)})})\]
given by 
\[ [(a,(\sigma,\tau))]  \mapsto  [(a,(\Theta(\sigma),\tau))]. \]
To see this is well defined we need to check that for $\sigma \in W_{(DF_v)^\ab/(EF_v)}$ we have
\[ r_{EF_v}(\epsilon(\Theta(\sigma)))=\sigma|_{(EF_v)^\ab}. \]
However if $s:\Gal(D_{u(\trho)}/E_{w(\trho)}) \ra W_{(DF_v)^\ab/EF_v}$ is a set theoretic section, we have
\[ \begin{array}{rcl} r_{EF_v}(\epsilon(\Theta(\sigma))) &=& r_{EF_v}(\prod_{\eta \in \Gal(D_{u(\trho)}/E_{w(\trho)})} \Theta(s_\eta) \Theta(\sigma) \Theta(s_{\eta \sigma}^{-1})) \\ &=& r_{EF_v}(\Theta(\tr_{W_{(DF_v)^\ab/EF_v}/W_{(DF_v)^\ab/DF_v}} \sigma|_{(EF_v)^\ab})) \\ &=& 
r_{EF_v}(\Theta\circ r_{DF_v} \circ r_{EF_v}^{-1}( \sigma|_{(EF_v)^\ab})),
\end{array} \]
and the claim follows because $\Theta\circ r_{DF_v}$ is the identity on $(DF_v)^\times$.
We also claim that
\[ \tGamma^E \circ \theta_{\trho}= \iota^{\eta_{D/E,*}\ga}_{w(\trho)} \circ \Theta^E:W_{(\trho(E)F_v)^\ab/F_v,D,\trho} \lra W_{E/F,D,\eta_{D/E,*}\ga}, \]
from which the tenth part will follow. 
We will realize these maps as maps
\[ (D_{w(\trho)}^\times \rtimes(W_{(DF_v)^\ab/F_v}\times_{\Gal(E/F)} \Gal(D/F)))/W_{(DF_v)^\ab/EF_v} \lra (\A_D^\times/D^\times \rtimes W_{D/F,\balpha})/(\A_D^\times/D^\times). \]

Both composites send $a \in D_{w(\trho)}^\times$ to $a \in \A_D^\times/D^\times$. Thus suppose $(\sigma ,\tau) \in W_{(DF_v)^\ab/F_v} \times_{\Gal(E/F)}\Gal(D/F)$. 
The first map, $\tGamma^E \circ \theta_{\trho}$, sends $(\sigma,\tau)$ to
\[ \begin{array}{rl} &[(r_E^{-1}(\theta_\trho(\sigma) \tGamma^{-1}(e^\glob_\balpha(\tau))^{-1})|_{E^\ab} \gamma_{E,u(\trho),\balpha}(e^\glob_\balpha(\tau)),e^\glob_\balpha(\tau))] \\ =&
[(r_D^{-1}(\prod_{\eta \in \Gal(D/E)} \tGamma^{-1}(e_\balpha^\glob(\eta)) \theta_\trho(\sigma) \tGamma^{-1}(e^\glob_\balpha(\tau))^{-1} \tGamma^{-1} (e_\balpha^\glob(\eta \sigma \tau^{-1}))^{-1}) \\ &\gamma_{E,u(\trho),\balpha}(e^\glob_\balpha(\tau)),e^\glob_\balpha(\tau))] \\ = &
[(\prod_{\eta \in \Gal(D/E)} e_\balpha^\glob(\eta)(\tGamma \circ \theta_\trho)(\sigma)e^\glob_\balpha(\tau)^{-1} e_\balpha^\glob(\eta \sigma \tau^{-1}))^{-1} \gamma_{E,u(\trho),\balpha}(e^\glob_\balpha(\tau)),e^\glob_\balpha(\tau))] 
.\end{array} \]
On the other hand $\iota^{\eta_{D/E,*}\ga}_{w(\trho)} \circ \Theta^E$ sends $(\sigma,\tau)$ to 
\[ \conju_{\prod_{\eta \in \Gal(D/E)} \beta(\eta)|_{\eta u(\trho)}} [( \prod_{\eta \in \Gal(D/E)} e_\balpha^\glob(\tau \eta \sigma^{-1}) \iota^\balpha_{u(\trho)}(\Theta(\sigma)) e_\balpha^\glob(\eta)^{-1} e_\balpha^\glob(\tau)^{-1}, e_\balpha^\glob(\tau))] . \]
Thus what we have to show is that for $(\sigma,\tau) \in W_{E_{w(\trho)}/F_v, D, \eta_{D/E,*} \ga} \times_{\Gal(E/F)}\Gal(D/F)$ we have
\[ \begin{array}{rl} & \prod_{\eta \in \Gal(D/E)} (e_\balpha^\glob(\eta)\iota^\balpha_{u(\trho)}(\sigma)e^\glob_\balpha(\tau)^{-1} e_\balpha^\glob(\eta \sigma \tau^{-1})^{-1} ) \gamma_{E,u(\trho),\balpha}(e^\glob_\balpha(\tau)) \\ = & {}^{(1-\tau)}\left(\prod_{\eta \in \Gal(D/E)} \beta(\eta)|_{\eta u(\trho)}\right) \prod_{\eta \in \Gal(D/E)} e_\balpha^\glob(\tau \eta \sigma^{-1}) \iota^\balpha_{u(\trho)}(\sigma) e_\balpha^\glob(\eta)^{-1} e_\balpha^\glob(\tau)^{-1},
\end{array} \]
i.e.
\[ \begin{array}{rl} & \prod_{\eta \in \Gal(D/E)} (e_\balpha^\glob(\eta)\iota^\balpha_{u(\trho)}(\sigma) e_\balpha^\glob(\eta \sigma)^{-1})
\prod_{\eta \in \Gal(D/E)} \alpha^\glob(\eta \sigma \tau^{-1}, \tau)^{-1}
\gamma_{E,u(\trho),\balpha}(e^\glob_\balpha(\tau)) \\ = & {}^{(1-\tau)}\left(\prod_{\eta \in \Gal(D/E)} \beta(\eta)|_{\eta u(\trho)}\right) \\ & \prod_{\eta \in \Gal(D/E)} (e_\balpha^\glob(\tau \eta \sigma^{-1}) \iota^\balpha_{u(\trho)}(\sigma) e_\balpha^\glob(\tau \eta)^{-1}) \prod_{\eta \in \Gal(D/E)} \alpha^\glob(\tau,\eta)^{-1}.
\end{array} \]
This in turn is equivalent to
\[  
\gamma_{E,u(\trho),\balpha}(e^\glob_\balpha(\tau))  =   \prod_{\eta \in \Gal(D/E)} (\alpha^\glob(\eta , \tau)\beta(\eta)|_{\eta u(\trho)})/(\alpha^\glob(\tau,\eta){}^\tau(\beta(\eta)|_{\eta u(\trho)})), \]
which is just the definition.

The eleventh and twelth parts are straightforward, so we leave them to the reader.

For the thirteenth part, if $\sigma \in \Gal(D^\ab/F)$ and $[(a,\tau)] \in (\A_D^\times /D^\times\Delta_E \rtimes W_{D^\ab/F})/(\A_D^\times/D^\times)$, then we have that
\[ \begin{array}{rl} &\eta_{D/E,\balpha,\trho\sigma,*}(\Gamma^{\sigma, \balpha})(a,\tau) \\ =&
 
 (a \gamma_{E,u(\trho\sigma),\balpha}(\Gamma^{\sigma,\balpha}(\tau)), \Gamma^{\sigma,\balpha}(\tau)) \\ = & 

 (a \gamma_{E,u
(\trho\sigma),\balpha}({}^{(1-\tau)}\beta(\sigma^{-1})_{u(\trho \sigma)}\conju_{e^\glob_\balpha(\sigma^{-1})} (\Gamma(\conju_\sigma(\tau)))), \\ & {}^{(1-\tau)}\beta(\sigma^{-1})_{u(\trho \sigma)}\conju_{e^\glob_\balpha(\sigma^{-1})} (\Gamma(\conju_\sigma(\tau)))) \\ = & 

 (a {}^{(1-\tau)}(N_{D/E}\beta(\sigma^{-1})_{u(\trho \sigma)}/\beta(\sigma^{-1})_{u(\trho \sigma)}^{[D:E]}) \\ & \gamma_{E,u(\trho\sigma),\balpha}(\conju_{e^\glob_\balpha(\sigma^{-1})} (\Gamma(\conju_\sigma(\tau)))) {}^{(1-\tau)}\beta(\sigma^{-1})_{u(\trho \sigma)}^{[D:E]}, \conju_{e^\glob_\balpha(\sigma^{-1})} (\Gamma(\conju_\sigma(\tau)))) \\ = & 
 
  (a  {}^{(1-\tau)}(\prod_{\eta \in \Gal(D/E)}{}^\eta(\beta(\sigma^{-1})_{u(\trho \sigma)})\beta(\eta)_{\eta u(\trho \sigma)}/\beta(\eta)_{\eta u(\trho)}) \\ & \gamma_{E,u(\trho),\balpha}(\conju_{e^\glob_\balpha(\sigma^{-1})} (\Gamma(\conju_\sigma(\tau)))) , \conju_{e^\glob_\balpha(\sigma^{-1})} (\Gamma(\conju_\sigma(\tau)))) \\ = & 
  
 \conju_{\prod_{\eta \in \Gal(D/E)}{}^\eta(\beta(\sigma^{-1})_{u(\trho \sigma)})\beta(\eta)_{\eta u(\trho \sigma)}/\beta(\eta)_{\eta u(\trho)}}   (a   \gamma_{E,u(\trho),\balpha}(\conju_{e^\glob_\balpha(\sigma^{-1})} (\Gamma(\conju_\sigma(\tau)))) , \\ & \conju_{e^\glob_\balpha(\sigma^{-1})} (\Gamma(\conju_\sigma(\tau)))).
 
\end{array} \]
On the other hand we have
\[  \begin{array}{rl} & (\eta_{E,u(\trho),\balpha,*}\Gamma)^{\sigma,\eta_{D/E,*}\balpha}(a,\tau) \\ =&

\conju_{(\eta_{D/E} \circ \beta)(\sigma^{-1})_{w(\trho \sigma)}} \circ \conju_{e^\glob_\balpha(\sigma^{-1})}( \eta_{D/E,u(\trho),\balpha,*}(\Gamma)({}^\sigma a,\conju_\sigma(\tau))) \\ =&

\conju_{(\eta_{D/E} \circ \beta)(\sigma^{-1})_{w(\trho \sigma)}} \circ \conju_{e^\glob_\balpha(\sigma^{-1})} 
({}^\sigma a \gamma_{E,u(\trho),\balpha}(\Gamma(\conju_\sigma(\tau))),\Gamma(\conju_\sigma(\tau))) \\ =&

\conju_{(\eta_{D/E} \circ \beta)(\sigma^{-1})_{w(\trho \sigma)}}  
(a {}^{\sigma^{-1}}\gamma_{E,u(\trho),\balpha}(\Gamma(\conju_\sigma(\tau))),  \conju_{e^\glob_\balpha(\sigma^{-1})}(\Gamma(\conju_\sigma(\tau)))) \\ =&

\conju_{(\eta_{D/E} \circ \beta)(\sigma^{-1})_{w(\trho \sigma)}}  
(a {}^{(\tau-1)}\gamma_{E,u(\trho),\balpha}(e^\glob_\balpha(\sigma^{-1})) \\ &\gamma_{E,\balpha}(e^\glob_\balpha(\sigma^{-1})\Gamma(\conju_\sigma(\tau)) e^\glob_\balpha(\sigma^{-1})^{-1}), \conju_{e^\glob_\balpha(\sigma^{-1})}(\Gamma(\conju_\sigma(\tau)))) \\ =&

\conju_{(\eta_{D/E} \circ \beta)(\sigma^{-1})_{w(\trho \sigma)}\gamma_{E,u(\trho),\balpha}(e^\glob_\balpha(\sigma^{-1}))^{-1}} 
(a  \gamma_{E,u(\trho),\balpha}(\conju_{e^\glob_\balpha(\sigma^{-1})}(\Gamma(\conju_\sigma(\tau))) ), \\ & \conju_{e^\glob_\balpha(\sigma^{-1})}(\Gamma(\conju_\sigma(\tau)))) \\ =&

\conju_A (\eta_{D/E,\balpha,\trho\sigma,*}(\Gamma^{\sigma, \balpha})(a,\tau)),
\end{array}\]
where
\[ \begin{array}{rl} & A \\ =& 

(\eta_{D/E} \circ \beta)(\sigma^{-1})_{w(\trho \sigma)}\gamma_{E,u(\trho),\balpha}(e^\glob_\balpha(\sigma^{-1}))^{-1} \prod_{\eta \in \Gal(D/E)}{}^\eta(\beta(\sigma^{-1})_{u(\trho \sigma)})^{-1}\beta(\eta)_{\eta u(\trho \sigma)}^{-1}\beta(\eta)_{\eta u(\trho)}\\ =&

\prod_{\eta \in \Gal(D/E)} \beta(\sigma^{-1})_{\eta u(\trho \sigma)} \alpha^\glob(\sigma^{-1},\eta){}^{\sigma^{-1}}(\beta(\eta)_{\eta u(\trho)})\alpha^\glob(\eta,\sigma^{-1})^{-1} \beta(\eta)_{\eta u(\trho)}^{-1} \\ &
 {}^\eta(\beta(\sigma^{-1})_{u(\trho \sigma)})^{-1}\beta(\eta)_{\eta u(\trho \sigma)}^{-1}\beta(\eta)_{\eta u(\trho)}\\ =&

\prod_{\eta \in \Gal(D/E)} \beta(\sigma^{-1})_{\sigma^{-1}\eta \sigma u(\trho \sigma)} \alpha^\glob(\sigma^{-1},\eta)({}^{\sigma^{-1}}\beta(\eta))_{\sigma^{-1} \eta u(\trho)})\alpha^\glob(\eta,\sigma^{-1})^{-1} \\ &
( {}^\eta\beta(\sigma^{-1}))_{\eta u(\trho \sigma)})^{-1}\beta(\eta)_{\eta u(\trho \sigma)}^{-1}\\ =&

\prod_{\eta \in \Gal(D/E)} (\beta(\sigma^{-1}) \alpha^\glob(\sigma^{-1},\eta){}^{\sigma^{-1}}\beta(\eta) \beta(\sigma^{-1}\eta)^{-1})_{\sigma^{-1}\eta  u(\trho)} \beta(\sigma^{-1}\eta)_{\sigma^{-1}\eta  u(\trho)}  \\ &
(\alpha^\glob(\eta,\sigma^{-1}) \beta(\eta) {}^\eta\beta(\sigma^{-1}) \beta(\eta \sigma^{-1})^{-1})_{\eta u(\trho \sigma)}^{-1} \beta(\eta \sigma^{-1})^{-1}_{\eta u(\trho \sigma)}\\ =&

\prod_{\eta \in \Gal(D/E)} \alpha^\loc(\sigma^{-1},\eta)_{\sigma^{-1}\eta  u(\trho)} \beta(\eta \sigma^{-1})_{\eta \sigma^{-1}  u(\trho)} 
\alpha^\loc(\eta,\sigma^{-1})_{\eta u(\trho \sigma)}^{-1} \beta(\eta \sigma^{-1})^{-1}_{\eta u(\trho \sigma)}\\ =&

\prod_{\eta \in \Gal(D/E)} (\alpha^\loc(\sigma^{-1},\eta)_{\sigma^{-1}\eta  u(\trho)}/ \alpha^\loc(\eta,\sigma^{-1})_{\eta u(\trho \sigma)}^{-1}),

 \end{array} \]
as desired.
\pfend

We call two Galois rigidifications $\Gamma_1, \Gamma_2$ of $\ga$ adapted to $\rho$ {\em equivalent}  if $\Gamma_2 = \conju_{a} \circ \Gamma_1$ for some $a \in D_{w(\rho)}^\times$. We will write $\Gamma_1 \sim \Gamma_2$ to mean $\Gamma_1$ and $\Gamma_2$ are equivalent, and $[\Gamma_1]$ for the equivalence class of $\Gamma_1$. Note that, if $\Gamma_1 \sim \Gamma_2$, then $\Gamma_1^{\sigma,\balpha} \sim \Gamma_2^{\sigma,\balpha}$ and ${}^t\Gamma_1 \sim {}^t \Gamma_2$ and $\inf_{C/D} \Gamma_1 \sim \inf_{C/D} \Gamma_2$ and $\eta_{D/E,\balpha,\trho,*} \Gamma_1 \sim \eta_{D/E,\balpha,\trho,*} \Gamma_2$. Moreover 
we have the following corollary:

\begin{cor}
\begin{enumerate}

\item If $[\Gamma]$ is an equivalence class of Galois rigidifications for $\ga$ adapted to $\rho$, then $[\Gamma^{\sigma,\balpha}]$ does not depend on the choice of $\balpha \in \ga$ and is an equivalence class of Galois rigidifications for $\ga$ adapted to $\rho\sigma$. We will denote it $[\Gamma]^\sigma$. We have 
\[ ([\Gamma]^{\sigma_1})^{\sigma_2}=[\Gamma]^{\sigma_1\sigma_2}. \]

\item If $[\Gamma]$ is an equivalence class of Galois rigidifications for $\ga$ adapted to $\rho$, then ${}^t[\Gamma]=[{}^t\Gamma]$ is an equivalence class of Galois rigidifications for ${}^t\ga$ adapted to $\rho$. We have
\[ {}^t([\Gamma]^\sigma)=({}^t[\Gamma])^\sigma. \]

If $[\Gamma]$ and $[\Gamma']$ are two equivalence classes of Galois rigidifications of $\ga$ adapted to $\rho$, then we can find $t\in T_{2,E}(\A_F)$ with $[\Gamma']={}^t[\Gamma]$. 

If $t\in T_{2,E}(\A_D)$ and ${}^t\ga=\ga$ and ${}^t [\Gamma]=[\Gamma]$, then $t=abc$ with $a \in \cE^\loc(E/F)_D^0$, $b \in \cE^\glob(E/F)_D^0$ and $c \in T_{2,E}(\A_F)$. Moreover for any such $a,b,c$ we must have $c_{w(\rho)}\in (\A_E^\times/E^\times)^{\Gal(E/F)_{w(\rho)}}$ of the form $c_{w(\rho)} =ef$ with $e \in \A_F^\times$ and $f \in \Delta_D^{\Gal(D/F)_{w(\rho)}}$ with ${}^\sigma f/f\in  \Delta_E$ for all $\sigma \in \Gal(D/F)$. If $E$ is totally imaginary, we deduce that $f\in \Delta_E^{\Gal(E/F)_{w(\rho)}}$.

\item If $C \supset D \supset E \supset F$ are finite Galois extensions, if $\ga \in \cH(E/F)_D$, and if $[\Gamma]$ is an equivalence class of Galois rigidifications for $\ga$ adapted to $\rho$, then $\inf_{C/D} [\Gamma]=[\inf_{C/D} \Gamma]$ is an equivalence class of Galois rigidifications for $\inf_{C/D} \ga$ adapted to $\rho$. We have
\[ (\inf_{C/D}[\Gamma])^\sigma = \inf_{C/D} ([\Gamma]^{\sigma|_{D^\ab}}) \;\;\;\;\;\;\; {\rm and}\;\;\;\;\;\;\; \inf_{C/D}({}^t[\Gamma]) = {}^t (\inf_{C/D} [\Gamma]), \]
for $\sigma \in \Gal(C^\ab/F)$ and $t \in T_{2,E}(\A_D)$. Moreover if $B \supset C \supset D \supset E \supset F$ are finite Galois extensions if a number field $F$, then 
\[ (\inf_{B/C} \circ \inf_{C/D}) [\Gamma]=\inf_{B/D}[\Gamma]. \]

\item If $\ga_D \in \cH(D/F)$ and if $[\Gamma]$ is an equivalence class of Galois rigidifications for $\ga$ adapted to and $F$-linear embedding $\trho:D^\ab \into \barFv$, then $[\eta_{D/E,\balpha,\trho,*}\Gamma]$ does not depend on the choice of $\balpha \in \ga_D$. We will denote it $\eta_{D/E,\trho,*}[\Gamma]$. If $t \in T_{2,D}(\A_D)$ we have 
\[ \eta_{D/E,\trho,*}({}^t[\Gamma]) = {}^{\eta_{D/E}(t)}(\eta_{D/E,\trho,*}[\Gamma]). \]
If $\sigma \in \Gal(D^\ab/F)$ we have
\[ \eta_{D/E,\trho \sigma,*}([\Gamma]^\sigma) = (\eta_{D/E,\trho,*}[\Gamma])^\sigma. \]

\end{enumerate}\end{cor}

\pfbegin
Everything but the last assertion of the second part follows immediately from the previous lemma. For the last assertion of the second part,
recall that we must have $t=abc$ with $a \in \cE^\loc(E/F)_D^0$, $b \in \cE^\glob(E/F)_D^0$ and $c \in T_{2,E}(\A_F)$. We have
${}^{abc}[\Gamma]=[\conju_{a_{w(\rho)}c_{w(\rho)}} \circ \Gamma]$
and so $\Gamma = \conju_{d^{-1}c_{w(\rho)}} \circ \Gamma$ for some $d \in D_{w(\rho)}^\times$. We deduce that $c_{w(\rho)}=def$ with $e\in \A_F^\times$ and $f\in \Delta_D$ with $f/{}^\sigma f\in \Delta_E$ for all $\sigma \in \Gal(D/F)$. If $\sigma \in \Gal(D/F)_{w(\rho)}$ then ${}^\sigma d/d = f/{}^\sigma f \in D_{w(\rho)}^\times \cap \Delta_D \subset \A_D^\times/D^\times$ and, as this intersection is trivial, we see that $d \in F_v^\times$ and $f \in \Delta_D^{\Gal(D/F)_{w(\rho)}}$. Replacing $e$ with $de$ the assertion follows.
\pfend

By {\em complete rigidification data}, $\{[\Gamma_\trho]\}$, for $\ga \in \cH(E/F)_D$ we mean the choice for each place $v$ of $F$ and for each $F$-linear embedding $\trho:E^\ab D \into \barFv$ of an equivalence class $[\Gamma_\trho]$ of rigidification data for $\ga$ adapted to $\trho|_{E^\ab}$ such that 
\begin{itemize}
\item if $\sigma \in \Gal(E^\ab D/F)$ then $[\Gamma_{\trho \sigma}]=[\Gamma_\trho]^\sigma$.
\end{itemize}

\begin{lem}\label{rcrd} \begin{enumerate}
\item If for each place $v$ of $F$ we choose an $F$-linear embedding $\trho_v:D^\ab \into \barFv$ and an equivalence class of rigidification data $[\Gamma_v]$ for $\ga$ adapted to $\trho_v$, then there is complete rigidification data $\{ [\Gamma_\trho]\}$ for $\ga$ with $[\Gamma_{\trho_v}]=[\Gamma_v]$. Moreover $\{ [\Gamma_\trho]\}$ is unique.

\item If $t \in T_{2,E}(\A_D)$ and if $\{ [\Gamma_{\trho}] \}$ is complete rigidification data for $\ga$, then $\{ {}^t[ \Gamma_{\trho}]\}$ is complete rigidification data for ${}^t\ga$.

\item If $\{ [\Gamma_{\trho}]\}$ and $\{ [\Gamma_{\trho}']\}$ are complete rigidification data for $\ga$, then there exists $t \in T_{2,E}(\A_F)$ such that $[\Gamma_{\trho}'] ={}^t[  \Gamma_{\trho}]$ for all $\trho$.

\item If $\ga_D \in \cH(D/F)$ and $\{ [\Gamma_\trho]\}$ is complete rigidification data for $\ga_D$, then $\{ \eta_{D/E,\trho,*}[\Gamma_\trho]\}$ is complete rigidification data for $\eta_{D/E,*} \ga_D \in \cH(E/F)_D$.

\item If $C \supset D \supset E \supset F$ are finite Galois extensions of $F$, if $\ga \in \cH(E/F)_D$ and if $\{ [\Gamma_\trho]\}$ is complete rigidification data for $\ga$, then $\{ \inf_{C/D} [\Gamma_{\trho|_{D^\ab}}]\}$ is complete rigidification data for $\inf_{C/D} \ga \in \cH(E/F)_C$.

\end{enumerate}\end{lem}

\pfbegin 
The first two and last two parts  are clear.
For the third part choose for each place $v$ of $F$ an $F$-linear embedding $\trho_v:D^\ab \into \barFv$.
Then we can find $t_{v} \in (\A_D^\times/D^\times)^{\Gal(D/F)_{w(\trho_v)}}$ such that
\[ \Gamma'_{\trho_v}\sim \conju_{t_{v} } \circ \Gamma_{\trho_v}. \]
Define $t \in T_{2,E}(\A_F)$ by
\[ t_{\sigma w(\trho_v)}=\sigma t_{w(\trho_v)}. \]
Then
\[[\Gamma'_{\trho_v\sigma}]= ({}^t[  \Gamma_{\trho_v}])^\sigma =  
{}^t([ \Gamma_{\trho_v)}]^\sigma)=  {}^t[  \Gamma_{\trho_v\sigma}], \]
and the third part follows.
\pfend

We will write $\cH(E/F)_D^+$ for the set of pairs $(\ga,\{[\Gamma_\trho]\})$, where $\ga \in \cH(E/F)_D$ and where $\{[\Gamma_\trho]\}$ is complete rigidification data for $\ga$. This set has a transitive action of $T_{2,E}(\A_D)$, where
\[ {}^t (\ga,\{[\Gamma_\trho]\})=({}^t\ga,\{{}^t[\Gamma_\trho]\}). \]
There is a natural $T_{2,E}(\A_D)$ equivariant map
\[ \begin{array}{rcl} \cH(E/F)_D^+ & \onto & \cH(E/F)_D \\ (\ga,\{[\Gamma_\trho]\}) & \longmapsto & \ga. \end{array} \]
{\em We will often denote an element of $\cH(E/F)_D^+$ by $\ga^+$ and then let $\ga$ denote the underlying element of $\cH(E/F)_D$.}

\begin{lem} If $\ga^+=(\ga,\{ [\Gamma_\trho]\}) \in \cH(E/F)_D^+$ and ${}^t \ga^+ = \ga^+$, then $t=abef$ with $a \in \cE^\loc(E/F)_D^0$, $b \in \cE^\glob(E/F)^0_D$, $e\in T_{2,F}(\A_F)$ and 
\[ f \in T_{2,E}(\A_F) \cap \prod_{w \in V_E} \{ d \in \Delta_D:\,\, {}^\sigma d/d \in \Delta_E \,\, \forall \sigma \in \Gal(D/F)\}.\]
 If $E$ is totally imaginary then in fact $f \in T_{2,E}(\A_F) \cap \prod_{w \in V_E} \Delta_E$. \end{lem}

\pfbegin We know that $t=abc$ with $a \in \cE^\loc(E/F)_D^0$, $b \in \cE^\glob(E/F)^0_D$ and $c \in T_{2,E}(\A_F)$. Moreover for all $w \in V_E$ we have $c_w=e_wf_w$ with $e_w \in \A_F^\times$ and $f_w \in \Delta_D^{\Gal(D/F)_w}$ with ${}^\sigma f_w/f_w \in \Delta_E$ for all $\sigma \in \Gal(D/F)$. Then for $\sigma \in \Gal(D/F)$ we have $c_{\sigma w}=e_w {}^\sigma f_w$ and so we may assume that $e_{\sigma w}=e_w$ and $f_{\sigma w} = {}^\sigma f_w$. The lemma follows. \pfend

There is a natural map
\[ \begin{array}{rcl} \eta_{D/E,*}: \cH(D/F)^+ & \onto & \cH(E/F)^+_D \\ (\ga_D,\{[\Gamma_\trho]\}) & \longmapsto & (\eta_{D/E,*} \ga_D, \{\eta_{D/E,\trho,*} [\Gamma_\trho]\}). \end{array} \]
If $t \in T_{2,D}(\A_D)$ then
\[ \eta_{D/E,*}{}^t \ga_D^+ = {}^{\eta_{D/E}(t)}(\eta_{D/E,*} \ga_D^+).\]

If $C \supset D \supset E \supset F$ are finite Galois extensions of $F$, then there is a natural map
\[ \begin{array}{rcl} \inf_{C/D}: \cH(E/F)_D^+ & \onto & \cH(E/F)^+_C \\ (\ga,\{[\Gamma_\trho]\}) & \longmapsto & (\inf_{C/D} \ga_D, \{\inf_{C/D} [\Gamma_{\trho|_{D^\ab}}]\}). \end{array} \]
If $t \in T_{2,E}(\A_D)$ then
\[ \inf_{C/D} {}^t \ga^+ = {}^t(\inf_{C/D} \ga^+). \]
If $B \supset C \supset D \supset E \supset F$ are finite Galois extensions of $F$, then
\[ \inf_{B/C} \circ \inf_{C/D}=\inf_{B/D}. \]

Now suppose that $C \supset D \supset E \supset F$ are finite Galois extensions of a number field $F$ and that $\ga^+=(\ga,\{ [\Gamma_\trho]\}) \in \cH(D/F)_C^+$. We wish to define $\eta_{D/E,*} \ga^+ \in \cH(E/F)_C^+$. To this end choose $\ga_1^+ \in \cH(D/F)^+$. Then we can find  $t_1 \in T_{2,D}(\A_C)$ with $\ga^+={}^{t_1} \inf_{C/D} \ga_1^+$. We wish to define
\[ \eta_{D/E,*} \ga^+ = {}^{\eta_{D/E}(t_1)}\inf_{C/D} \eta_{D/E} \ga_1^+. \]
To see that this is a good definition suppose also that $\ga_2^+ \in \cH(D/F)^+$ and  $t_2 \in T_{2,D}(\A_C)$ with $\ga^+={}^{t_2} \inf_{C/D} \ga_2^+$. Choose $s \in T_{2,D}(\A_D)$ with $\ga_2^+={}^s \ga_1^+$. Then ${}^{t_2st_1^{-1}} \ga^+=\ga^+$, so that 
\[ \begin{array}{l}  t_2st_1^{-1}  \in \\ \cE^\loc(D/F)_C^0 \cE^\glob(D/F)^0_C T_{2,F}(\A_F) (T_{2,D}(\A_F) \cap \prod_{v \in V_D} \{ d \in \Delta_C:\,\, {}^\sigma d/d \in \Delta_D \,\, \forall \sigma \in \Gal(C/F)\}) \end{array} \]
and
\[ \begin{array}{l}  \eta_{D/E}(t_2st_1^{-1}) \in \\ \cE^\loc(E/F)_C^0 \cE^\glob(E/F)^0_C T_{2,F}(\A_F) (T_{2,D}(\A_F) \cap \prod_{v \in V_D} \{ d \in \Delta_C:\,\, {}^\sigma d/d \in \Delta_E \,\, \forall \sigma \in \Gal(C/F)\}). \end{array} \]
(To see this suppose that 
\[ (f_u)_{u \in V_D} \in  T_{2,D}(\A_F) \cap \prod_{v \in V_D} \{ d \in \Delta_C:\,\, {}^\sigma d/d \in \Delta_D \,\, \forall \sigma \in \Gal(C/F)\}. \]
Then $f_u \in \A_D^\times$ and $\eta_{D/E}(f_u)_{u\in V_D}=(N_{D/E} f_{u(w)})_{w \in V_E}$, where for each $w \in V_E$ we choose $u(w)\in V_D$ above it. Moreover if $\tau \in \Gal(C/D)$ we have ${}^{(\tau -1)} N_{D/E} f_{u(w)} = N_{D/E} ({}^{(\tau-1)} f_{u(w)}) \in N_{D/E} \Delta_D = \Delta_E$.)
Thus
\[ \begin{array}{rcl} {}^{\eta_{D/E}(t_2)}\inf_{C/D} \eta_{D/E} \ga_2^+&=&{}^{\eta_{D/E}(t_2)}\inf_{C/D} \eta_{D/E} {}^s\ga_1^+
\\ &=& {}^{\eta_{D/E}(t_2s)}\inf_{C/D} \eta_{D/E} \ga_1^+ \\ &=& {}^{\eta_{D/E}(t_2st_1^{-1})} ({}^{\eta_{D/E}(t_1)}\inf_{C/D} \eta_{D/E} \ga_1^+) \\ &=&
{}^{\eta_{D/E}(t_1)}\inf_{C/D} \eta_{D/E} \ga_1^+, \end{array}\]
and $\eta_{D/E,*} \ga^+$ is well defined. 

\begin{lem} Suppose that $B \supset C \supset D \supset E \supset F$ are finite Galois extensions of a number field $F$.
\begin{enumerate}
\item If $\ga^+ \in \cH(D/F)_C^+$ and $t \in T_{2,D}(\A_C)$, then $\eta_{D/E,*} {}^t\ga^+ = {}^{\eta_{D/E}(t)} \eta_{D/E,*} \ga^+$.

\item If $\ga^+ \in \cH(C/F)_B^+$ then $\eta_{D/E,*}\eta_{C/D,*} \ga^+=\eta_{C/E,*} \ga^+$.

\item If $\ga^+ \in \cH(D/F)_C$ then $\eta_{D/E,*} \inf_{B/C} \ga^+ = \inf_{B/C} \eta_{D/E,*} \ga^+$.
\end{enumerate}
\end{lem}

\pfbegin These are all formal consequences of the definition, and the facts given above. \pfend

\begin{cor} Suppose that $C \supset D \supset E \supset F$ are finite Galois extensions of $F$.  Suppose also that $\ga_{C}^+ \in\cH(C/F)^+$ and $\ga_{D}^+ \in\cH(D/F)^+$ and $\ga_{E}^+ \in\cH(E/F)^+$ satisfy $\eta_{C/D,*}\ga_{C}^+={}^{t'} \inf_{C/D} \ga_D^+$ and $\eta_{D/E,*}\ga_{D}^+={}^{t} \inf_{D/E} \ga_E^+$ with $t' \in T_{2,D}(\A_{C})$ and $t \in T_{2,E}(\A_{D})$. Then
\[ \eta_{C/E,*}\ga_{C}^+={}^{t \eta_{D/E}(t')} \inf_{C/E} \ga_E^+. \]
\end{cor}

\subsection{More technical lemmas}

We have chosen to split the calculations needed to prove the results of the next section into two. In this section we begin those calculations.

Suppose that $F$ is a number field, that $v$ is a place of $F$ and that $D \supset E \supset F$ are finite Galois extensions of $F$. Also let 
\[ \ga^+=(\ga,\{ [\Gamma_\rho]\}) \in \cH(E/F)_D^+ .\]
If $\trho:E^\ab D \into \barFv$ is $F$-linear and if $\tau \in \Aut(\barFv/F)$ (a field theoretic automorphism of $\barF_v$ fixing elements of $F$, but {\em not necessarily continuous}) then we have
\[ \tau^\trho \in \Gal(E^\ab D/F) \into \Gal(E^\ab/F) \times_{\Gal(E/F)} \Gal(D/F) \subset W_{E^\ab/F,D} \]
and we can define
 \[ \barg_{\ga^+,v,\trho}(\tau) = \left( \left( \loc_\ga \widetilde{\Gamma_{\trho}(\tau^\trho)}\right)^{-1} e^\loc_\balpha(\tau^\trho)\right)_{w(\tau\trho)} \in \A_D^\times/D^\times \Delta_E D_v^\times, \]
where $\balpha=(\alpha^\glob,\alpha^\loc,\beta) \in \ga$, and $\widetilde{\Gamma_{\trho}(\tau^\trho)}$ is any lift of ${\Gamma_{\trho}(\tau^\trho)}$ (for some $\Gamma_\trho \in [\Gamma_\trho]$) to $\cE^\glob(E/F)_{D,\ga}$. This element is independent of the choices $\Gamma_\trho$ and $\widetilde{\Gamma_{\trho}(\tau^\trho)}$, and of the choice of $\balpha \in \ga$.

 \begin{lem}\label{gcocyc}\begin{enumerate}
 \item\label{gcocyc1} $\barg_{\ga^+,v,\trho\sigma}(\tau)={}^{\sigma^{-1}} \barg_{\ga^+,v,\trho}(\tau)  \beta(\sigma^{-1})_{w(\trho\sigma)} /\beta(\sigma^{-1})_{w(\tau\trho\sigma)} $. In particular $\barg_{\ga^+,v,\trho}(\tau)$ only depends on $\trho|_D$.
 
  \item\label{gcocyc2} $\barg_{{}^t\ga^+,v,\trho}(\tau)=(t_{w(\trho)}/t_{w(\tau\trho)})\barg_{\ga^+,v,\trho}(\tau)$.
  
  \item\label{gcocyc3} If $C \supset D \supset E \supset F$ are finite Galois extensions of $F$, if $\trho':E^\ab C \into \barF_v$ is $F$-linear and if $\ga^+ \in \cH(E/F)_D^+$, then $\barg_{\ind_{C/D} \ga^+,v,\trho'}(\tau)=\barg_{\ga^+,v,\trho'|_{E^\ab D}}(\tau)$.

\item $\barg_{\ga^+,v,\trho}(\tau)$ only depends on $\trho|_E$. {\em If $\rho: E \into \barF_v$ is $F$-linear, we will write $\barg_{\ga^+,v,\rho}(\tau)=\barg_{\ga^+,v,\trho}(\tau)$ for any extension $\trho: E^\ab D \into \barF_v$ of $\rho$.}

\item\label{gcocyc5} If $\sigma \in \Gal(D/E)$ then ${}^\sigma \barg_{\ga^+,v,\rho}(\tau)=\barg_{\ga^+,v,\rho} (\tau)\beta(\sigma)_{w(\tau\rho)}/\beta(\sigma)_{w(\rho)}$.

 \item If $\ga^+ \in \cH(D/F)^+$ then 
 \[ \begin{array}{rcl} \barg_{\eta_{D/E,*} \ga^+,v,\rho}(\tau)&=& (N_{D/E} \barg_{\ga^+,v,\rho}(\tau))\prod_{\eta \in \Gal(D/E)}\beta(\eta)_{\eta u(\rho)}/\beta(\eta)_{\eta u(\tau\rho)} \\ &=& \prod_{\eta \in \Gal(D/E)} \barg_{\ga^+,v,\rho \eta^{-1}}(\tau). \end{array}\]
 
 \item If $C \supset D \supset E \supset F$ are finite Galois extensions of $F$, if $\rho:E \into \barF_v$ is $F$-linear and if $\ga^+ \in \cH(D/F)_C^+$, then
 \[ \barg_{\eta_{D/E,*}\ga^+,v,\rho} = \prod_{\trho:D \into \barF_v} \barg_{\ga^+,v,\trho}(\tau), \]
 where $\trho$ runs over extensions of $\rho$.

 \end{enumerate} \end{lem}
 
 \pfbegin 
 For the first part choose $\balpha \in \ga$, then
 \[ \begin{array}{rl} &\barg_{\ga^+,v,\trho\sigma}(\tau) \\ = & 
 
 \left( \left( \loc_\balpha \widetilde{\Gamma_{\trho\sigma}(\sigma^{-1}\tau^\trho\sigma)}\right)^{-1} e^\loc_\balpha(\tau^{\trho\sigma})\right)_{w(\tau\trho\sigma)} \\ =&
 
 \left(  \loc_\balpha \widetilde{\left(\beta(\sigma^{-1})_{w(\trho\sigma)} e^\glob_\balpha(\sigma^{-1})\Gamma_{\trho}(\tau^\trho) e^\glob_\balpha(\sigma^{-1})^{-1} \beta(\sigma^{-1})_{w(\trho\sigma)}^{-1}\right)}^{-1} e^\loc_\balpha(\tau^{\trho\sigma})\right)_{w(\tau\trho\sigma)} \\ =&
 
  \left(  \beta(\sigma^{-1})_{w(\trho\sigma)} \beta(\sigma^{-1})^{-1} e^\loc_\balpha(\sigma^{-1})(\loc_\balpha \widetilde{\Gamma_{v,\trho}(\tau^\trho)})^{-1} e^\loc_\balpha (\sigma^{-1})^{-1} \beta(\sigma^{-1})\beta(\sigma^{-1})_{w(\trho\sigma)}^{-1}e^\loc_\balpha(\tau^{\trho\sigma})\right)_{w(\tau\trho\sigma)} \\ =&
  
    \left(   e^\loc_\balpha(\sigma^{-1})(\loc_\balpha \widetilde{\Gamma_{v,\trho}(\tau^\trho)})^{-1} e^\loc_\balpha(\tau^\trho) e^\loc_\balpha (\sigma^{-1})^{-1}\right)_{w(\tau\trho \sigma)} \\ &
\left(  \beta(\sigma^{-1})_{w(\trho\sigma)} \beta(\sigma^{-1})^{-1} e^\loc_\balpha(\sigma^{-1})e^\loc_\balpha(\tau^\trho)^{-1}    e^\loc_\balpha (\sigma^{-1})^{-1} \beta(\sigma^{-1})\beta(\sigma^{-1})_{w(\trho\sigma)}^{-1}e^\loc_\balpha(\sigma^{-1}\tau^\trho\sigma)\right)_{w(\tau\trho\sigma)} \\ =&

{}^{\sigma^{-1}} \barg_{\ga^+,v,\trho}(\tau)
 \beta(\sigma^{-1})_{w(\trho\sigma)} \beta(\sigma^{-1})_{w(\tau\trho\sigma)}^{-1}
 ({}^{\sigma^{-1}(\tau^\trho)^{-1}\sigma} \beta(\sigma^{-1}))_{w(\tau\trho\sigma)} {}^{\sigma^{-1}(\tau^\trho)^{-1}\sigma} (\beta(\sigma^{-1})_{w(\trho \sigma)})^{-1} \\ = &

{}^{\sigma^{-1}} \barg_{\ga^+,v,\trho}(\tau) 
 \beta(\sigma^{-1})_{w(\trho\sigma)} \beta(\sigma^{-1})_{w(\tau\trho\sigma)}^{-1}
{}^{(\tau^{\trho\sigma})^{-1}} (\beta(\sigma^{-1})_{w(\trho\sigma)}/ \beta(\sigma^{-1})_{w(\trho \sigma)}) \\ =&

{}^{\sigma^{-1}} \barg_{\ga^+,v,\trho}(\tau) \beta(\sigma^{-1})_{w(\trho\sigma)} \beta(\sigma^{-1})_{w(\tau\trho\sigma)}^{-1}.
    
 \end{array}\]
 If $\sigma \in \Gal(D^\ab/D)$, then $\barg_{\ga^+,v,\trho\sigma}(\tau)/\barg_{\ga^+,v,\trho}$ equals
 \[ \beta(1)_{w(\trho)}/\beta(1)_{w(\tau\trho)} = \alpha^\loc(1,1)_{w(\trho)}\alpha^\glob(1,1)_{w(\tau\trho)}/\alpha^\glob(1,1)_{w(\trho)}\alpha^\loc(1,1)_{w(\tau\trho)} \in D^\times D_v^\times. \]
 
 For the second part, if $\balpha \in \ga$, we have
 we have
\[ \begin{array}{rcl}  \barg_{{}^t\ga^+,v,\trho}(\tau) & =&
 \left( \left( \loc_{{}^t\balpha }\widetilde{{}^t\Gamma_{\trho}(\tau^\trho)}\right)^{-1} e^\loc_{{}^t\balpha}(\tau^\trho)\right)_{w(\tau\trho)} \\ &= &
\left( \left( t^{-1} (\loc_{\balpha }t_{w(\trho)}\widetilde{\Gamma_{\trho}(\tau^\trho)} t_{w(\trho)}^{-1}) t\right)^{-1} e^\loc_\balpha(\tau^\trho)\right)_{w(\tau\trho)} \\ &=&
 (({}^{(\tau^\trho)^{-1}}(t_{w(\trho)}^{-1}t))/(t_{w(\trho)}^{-1}t))_{w(\tau \trho)} \left(\left( \loc_{\balpha }\widetilde{\Gamma_{\trho}(\tau^\trho)}\right)^{-1} e^\loc_\balpha(\tau^\trho)\right)_{w(\tau\trho)} \\ &=&

(t_{w(\trho)}/t_{w(\tau\trho)})\barg_{\ga^+,v,\trho}(\tau).

\end{array}\]

 For the third part, under the identification
  \[ \cE_2(E/F)_{C,\inf_{C/D}\ga} \cong (T_{2,E}(\A_D) \rtimes (\cE_2(E/F)_{D,\ga} \times_{\Gal(D/F)} \Gal(C/F)))/T_{2,E}(\A_D) \]
 $e^\loc_{\inf_{C/D}\balpha}(\tau^{\trho'})$ is identified with $[(1,(e^\loc_{\balpha}(\tau^{\trho'}),\tau^{\trho'}))]$. Moreover, under the identification
  \[ W_{E/F,C,\inf_{C/D}\ga} \cong (\A_C^\times/C^\times \rtimes (W_{E/F,D,\ga} \times_{\Gal(D/F)} \Gal(C/F)))/(\A_D^\times/D^\times) \]
  $(\inf_{C/D}\Gamma_{\trho'|_{E^\ab D}})(\tau^{\trho'})$ is identified with $[(1,(\Gamma_{\trho'|_{E^\ab D}}(\tau^{\trho'}),\tau^{\trho'}))]$. Moreover
$\loc_{\inf_{C/D}\ga} \widetilde{(\inf_{C/D}\Gamma_{\trho'|_{E^\ab D}})(\tau^{\rho'})}$ can be chosen to be $[(1,(\loc_\ga \widetilde{\Gamma_{\trho'|_{E^\ab D}}(\tau^{\trho'})},\tau^{\trho'}))]$. The third part follows.

 For the fourth part choose $\ga_1^+ \in \cH(E/F)^+$ and $t \in T_{2,E}(\A_D)$ with $\ga^+={}^t \inf_{D/E} \ga_1^+$. Then
 \[ \barg_{\ga^+,v,\trho}(\tau)=(t_{w(\trho)}/t_{w(\tau \trho)})\barg_{\ga_1^+,v,\trho|_{E^\ab}}(\tau), \]
 which by part \ref{gcocyc1} only depends on $\trho|_E$. The fifth part now follows immediately from the first part.
 
 For the sixth part we will make use of the identifications 
\[ W_{E/F,D,\eta_{D/E,*}\ga} \cong (\A_D^\times/D^\times \rtimes W_{D/F,\ga})/(\A_D^\times/D^\times) \]
and
\[ \cE^\glob(E/F)_{D,\eta_{D/E,*}\ga} \cong (\cE^\glob(E/F)_D^0 \rtimes \cE^\glob(D/F)_\ga)/ \cE^\glob(D/F)^0 \]
and
\[ \cE_2(E/F)_{D,\eta_{D/E,*}\ga} \cong (T_{2,E}(\A_D) \rtimes \cE_2(D/F)_\ga)/ T_{2,D}(\A_D). \]
Under these identifications, if $\balpha \in \ga$ and $\trho:E^\ab C \into \barF_v$ extends $\rho$, then $(\eta_{D/E,\balpha,\trho,*} \Gamma_\trho)(\tau^\trho)$ corresponds to 
\[ \begin{array}{rl} &  [(\gamma_{E,u(\trho),\balpha}(\widetilde{\Gamma_{\trho}(\tau^\trho)}),\widetilde{\Gamma_{\trho}(\tau^\trho)})]\\  =&

 [(N_{D/E}(\widetilde{\Gamma_{\trho}(\tau^\trho)} e^\glob_\balpha(\tau^\trho)^{-1} )\gamma_{E,u(\trho),\balpha}(e^\glob_\balpha(\tau^\trho))/(\widetilde{\Gamma_{\trho}(\tau^\trho)} e^\glob_\balpha(\tau^\trho)^{-1} )^{[D:E]},\Gamma_{\trho}(\tau^\trho))]\\  =&

 [(N_{D/E}(\widetilde{\Gamma_{\trho}(\tau^\trho)} e^\glob_\balpha(\tau^\trho)^{-1} ) \prod_{\eta \in \Gal(D/E)} (\alpha^\glob(\eta,\tau^\trho)\beta(\eta)_{\eta u(\trho)}/\alpha^\glob(\tau^\trho,\eta){}^{\tau^\trho}(\beta(\eta)_{\eta u(\trho)}))
,\\& e^\glob_\balpha(\tau^\trho))]. \end{array}  \]
Applying $\loc_{\eta_{D/E,*}\balpha}$ we get
\[ \begin{array}{l} [(N_{D/E}(\widetilde{\Gamma_{\trho}(\tau^\trho)} e^\glob_\balpha(\tau^\trho)^{-1} ) \prod_{\eta \in \Gal(D/E)} (\alpha^\glob(\eta,\tau^\trho)\beta(\eta)_{\eta u(\trho)}/\alpha^\glob(\tau^\trho,\eta){}^{\tau^\trho}(\beta(\eta)_{\eta u(\trho)}))\\
(\eta_{D/E} \beta(\tau^\trho))^{-1}, e^\loc_\balpha(\tau^\trho))] \end{array} \]
and then taking the inverse gives
\[ \begin{array}{rl} & 
[({}^{(\tau^\trho)^{-1}}N_{D/E}( e^\glob_\balpha(\tau^\trho)\widetilde{\Gamma_{\trho}(\tau^\trho)}^{-1} ) {}^{(\tau^\trho)^{-1}}(\eta_{D/E} \beta(\tau^\trho)) \\ & \prod_{\eta \in \Gal(D/E)} \beta(\eta)_{\eta u(\trho)} {}^{(\tau^\trho)^{-1}}\alpha^\glob(\tau^\trho,\eta)/{}^{(\tau^\trho)^{-1}}(\beta(\eta)_{\eta u(\trho)}\alpha^\glob(\eta,\tau^\trho))
, e^\loc_\balpha(\tau^\trho)^{-1} )] \\ = &

[(N_{D/E}(\widetilde{\Gamma_{\trho}(\tau^\trho)}^{-1} e^\glob_\balpha(\tau^\trho) ) {}^{(\tau^\trho)^{-1}}(\eta_{D/E} \beta(\tau^\trho)) \\ & \prod_{\eta \in \Gal(D/E)} \beta(\eta)_{\eta u(\trho)} {}^{(\tau^\trho)^{-1}}(\alpha^\glob(\tau^\trho,\eta))/{}^{(\tau^\trho)^{-1}}(\beta(\eta)_{\eta u(\trho)}\alpha^\glob(\eta,\tau^\trho))
, e^\loc_\balpha(\tau^\trho)^{-1} )].

\end{array} \]
Multiplying on the right by $e^\loc_{\eta_{D/E,*}\balpha}(\tau^\trho)=[(1,e^\loc_{\balpha}(\tau^\trho))]$ and taking the $w(\tau\rho)$-component gives
\[ \begin{array}{l} N_{D/E}(\widetilde{\Gamma_{\trho}(\tau^\trho)}^{-1} e^\glob_\balpha(\tau^\trho) )  \\  \prod_{\eta \in \Gal(D/E)} \beta(\eta)_{\eta u(\trho)} {}^{(\tau^\trho)^{-1}}(\beta(\tau^\trho)_{\eta u( \trho)} \alpha^\glob(\tau^\trho,\eta))/{}^{(\tau^\trho)^{-1}}(\beta(\eta)_{\eta u(\trho)}\alpha^\glob(\eta,\tau^\trho)).\end{array} \]
However
\[ \widetilde{\Gamma_{\trho}(\tau^\trho)}^{-1} e^\glob_\balpha(\tau^\trho) = (\widetilde{\Gamma_{\trho}(\tau^\trho)}^{-1} \beta(\tau^\trho)^{-1}e^\loc_\balpha(\tau^\trho))_{u(\tau\trho)} = ({}^{(\tau^\trho)^{-1}}\beta(\tau^\trho))^{-1}_{u(\tau\trho)} \barg_{\ga^+,v,\trho}(\tau). \]
Thus $\barg_{\eta_{D/E,*}\ga^+,v,\trho}(\tau)$ equals the product of $N_{D/E} \barg_{\balpha^+,v,\trho}(\tau)$ with
\[ \begin{array}{l} (N_{D/E}{}^{(\tau^\trho)^{-1}}(\beta(\tau^\trho)^{-1}_{u(\trho)})  \\  \prod_{\eta \in \Gal(D/E)} \beta(\eta)_{\eta u(\trho)} {}^{(\tau^\trho)^{-1}}(\beta(\tau^\trho)_{\eta u( \trho)} \alpha^\glob(\tau^\trho,\eta))/{}^{(\tau^\trho)^{-1}}(\beta(\eta)_{\eta u(\trho)}\alpha^\glob(\eta,\tau^\trho)), \end{array} \]
which equals 
$(\tau^\trho)^{-1}$ applied to
\[ \begin{array}{rl} & 

 \prod_{\eta \in \Gal(D/E)}  \beta(\tau^\trho)_{\eta  u(\trho)} \alpha^\glob(\tau^\trho,\eta){}^{\tau^\trho}( \beta(\eta)_{\eta u(\trho)})/(\beta(\eta)_{\eta u(\trho)}\alpha^\glob(\eta,\tau^\trho){}^\eta(\beta(\tau^\trho)_{u(\trho)}) \\ =&

\prod_{\eta \in \Gal(D/E)}  \beta(\tau^\trho)_{\tau^\trho \eta u(\tau \trho)} \alpha^\glob(\tau^\trho,\eta) ({}^{\tau^\trho}\beta(\eta))_{\tau^\trho \eta u(\tau \trho)}
{}^{\tau^\trho}( \beta(\eta)_{\eta u(\trho)})/ \\ & (\beta(\eta)\alpha^\glob(\eta,\tau^\trho)({}^\eta\beta(\tau^\trho)))_{\eta u(\trho)}{}^{\tau^\trho}(\beta(\eta)_{ \eta u(\tau \trho)}) \\ =&

\prod_{\eta \in \Gal(D/E)}  (\beta(\tau^\trho\eta) \alpha^\loc(\tau^\trho,\eta))_{\tau^\trho \eta u(\tau \trho)}
{}^{\tau^\trho}( \beta(\eta)_{\eta u(\trho)})/ (\beta(\eta\tau^\trho)\alpha^\loc(\eta,\tau^\trho))_{\eta u(\trho)}{}^{\tau^\trho}(\beta(\eta)_{ \eta u(\tau \trho)}) \\ =&

\prod_{\eta \in \Gal(D/E)}  (\alpha^\loc(\tau^\trho,\eta)_{\tau^\trho \eta u(\tau \trho)} /\alpha^\loc(\eta,\tau^\trho)_{\eta u(\trho)})(\beta(\tau^\trho\eta)_{\tau^\trho \eta u(\tau \trho) }/\beta(\eta\tau^\trho)_{\eta u(\trho)}) \\ &
{}^{\tau^\trho}( \beta(\eta)_{\eta u(\trho)})/ \beta(\eta)_{ \eta u(\tau \trho)}) \\ =&

\prod_{\eta \in \Gal(D/E)}  (\alpha^\loc(\tau^\trho,\eta)_{\tau^\trho \eta u(\tau \trho)} /\alpha^\loc(\eta,\tau^\trho)_{\eta u(\trho)})(\beta(\eta \tau^\trho)_{\eta \tau^\trho  u(\tau \trho) }/\beta(\eta\tau^\trho)_{\eta u(\trho)}) \\ &
{}^{\tau^\trho}( \beta(\eta)_{\eta u(\trho)})/ \beta(\eta)_{ \eta u(\tau \trho)}) \\ =&

\prod_{\eta \in \Gal(D/E)}  (\alpha^\loc(\tau^\trho,\eta)_{\tau^\trho \eta u(\tau \trho)} /\alpha^\loc(\eta,\tau^\trho)_{\eta u(\trho)})
{}^{\tau^\trho}\prod_{\eta \in \Gal(D/E)} ( \beta(\eta)_{\eta u(\trho)})/ \beta(\eta)_{ \eta u(\tau \trho)}) \\ \equiv&

{}^{\tau^\trho}\prod_{\eta \in \Gal(D/E)} ( \beta(\eta)_{\eta u(\trho)})/ \beta(\eta)_{ \eta u(\tau \trho)}),

    \end{array} \]
    and the first assertion of the sixth part of the lemma follows. The second assertion of this part then follows from the first part of the lemma.
    
    For the seventh part choose $\ga_1^+ \in \cH(D/F)^+$ and $t \in T_{2,D}(\A_C)$ such that $\ga^+={}^t\inf_{C/D} \ga_1^+$. Then
    \[ \begin{array}{rl} &\barg_{\eta_{D/E,*}\ga^+,v,\rho}(\tau) \\ =& 
    \barg_{{}^{\eta_{D/E}(t)} \inf_{C/D} \eta_{D/E,*}\ga_1^+,v,\rho}(\tau) \\ =& 
(\eta_{D/E}(t)_{w(\rho)}/\eta_{D/E}(t)_{w(\tau\rho)}) \prod_{\trho} \barg_{\inf_{C/D} \ga_1^+,v,\trho}(\tau) \\ =& \prod_{\trho} (t_{w(\trho)}/t_{w(\tau \trho)}) \prod_\trho (t^{-1}_{w(\trho)}/t^{-1}_{w(\tau \trho)}) \barg_{\ga,v,\trho}(\tau), \end{array}\]
and the result follows.
\pfend

\begin{lem}\label{bargc} \begin{enumerate}
 \item If $\tau$ fixes the image $E$ in $\barF_v$ and $\trho:E^\ab D \into \barF_v$ extends $\rho$, then $\barg_{\ga^+,v,\rho}(\tau) =\Art_E^{-1}(\tau^\trho|_{E^\ab})^{-1}$.
 
\item If $\tau$ is continuous then $\barg_{\ga^+,v,\rho}(\tau)=1$.

 \item $\barg_{\ga^+,v,\rho}( \tau_1\tau_2)= \barg_{\ga^+,v,\tau_2\rho}(\tau_1) \barg_{\ga^+,v,\rho}(\tau_2)$. 
\end{enumerate} \end{lem}

\pfbegin
By parts \ref{gcocyc2} and \ref{gcocyc3} of the previous lemma we reduce the first part to the case that $D=E$. In this case, if $\trho:E^\ab \ra \barF_v$ and $\balpha \in \ga$, then we have
\[ \begin{array}{rcl} \barg_{\ga^+,v,\rho}(\tau)&=&(\widetilde{\Art_E^{-1}(\tau^\trho)}^{-1} e^\loc_\balpha(\tau^\trho))_{w(\tau \rho)}   \\ &=&
\Art_E^{-1}(\tau^\trho)^{-1}\alpha^\loc(1,1)_{w(\tau \rho)} \\ &=&\Art_E^{-1}(\tau^\trho)^{-1} \end{array}\]
as desired.

For the second part, by continuity, we may suppose that $\tau$ is the image of some $\ttau \in W_{\barF_V/F_v}$. If $\tGamma_\trho:W_{E^\ab/F,D} \iso W_{E/F,D,\ga}$ lifts $\Gamma_\trho$, then we may take $\widetilde{\Gamma_\trho(\tau^\trho)}$ to be a lift of $\tGamma_\trho(\theta_\trho(\ttau))=\iota^\ga_{w(\rho)}(\Theta(\ttau))$ for some $\Theta: W_{E F_v)^\ab/F_v,\rho,D} \iso W_{E_{w(\rho)}/F_v,D,\ga}$. Thus, if $\balpha \in \ga$, we need to prove that if $\sigma \in W_{E_{w(\rho)}/F_v,D,\ga}$, then 
\[ \left( \left( \loc_\ga \widetilde{\iota^\ga_{w(\rho)}(\sigma)} \right)^{-1} e_\balpha^\loc(\sigma) \right)_{w(\rho)} =1\in \A_D^\times/D^\times \Delta_E D_v^\times. \]
If $\sigma \in D_{w(\rho)}^\times$ then this equals $(\sigma^{-1} \alpha_\balpha^\loc(1,1))_{w(\rho)} \in D_{w(\rho)}^\times$. On the other hand if $\sigma=e^\loc_\balpha(\eta)$ with $\eta \in \Gal(D/F)_{w(\rho)}$, then this equals
\[ \left( \left( \beta(\eta)^{-1} \widetilde{\beta(\eta)_{w(\rho)}} e^\loc_\balpha(\eta)\right)^{-1} e^\loc_\balpha(\eta)\right)_{w(\rho)} = {}^{\eta^{-1}}(\beta(\eta)^{-1} \widetilde{\beta(\eta)_{w(\rho)}})_{w(\rho)}=1. \]
The second part follows.
 
 For the third part we have we will expand $\barg_{\ga^+,v,\tau_2\rho}(\tau_1) \barg_{\ga^+,v,\rho}(\tau_2)$. First note that, if $\balpha \in \ga$ we can rewrite
 \[ \begin{array}{rl} &\barg_{\ga^+,v,\tau_2\rho}(\tau_1)  \\ =&

 \left( \left( \loc_\balpha \widetilde{\Gamma_{\tau_2\rho}(\tau_1^{\tau_2\rho})}\right)^{-1} e^\loc_\balpha(\tau_1^{\tau_2\rho})\right)_{w(\tau_1\tau_2\rho)}   \\ =&
 
\left( \loc_\balpha\left( \conju_{\beta((\tau_2^\rho)^{-1})_{w(\tau_2\rho)}e^\glob_\balpha((\tau_2^\rho)^{-1})} (\widetilde{\Gamma_{\rho}(\tau_1^\rho)})\right)^{-1}e^\loc_\balpha(\tau_1^{\tau_2\rho})\right)_{w(\tau_1\tau_2\rho)}  \\ =&
 
 \left( \left( \conju_{\beta((\tau_2^\rho)^{-1})_{w(\tau_2\rho)}\beta((\tau_2^\rho)^{-1})^{-1}e^\loc_\balpha((\tau_2^\rho)^{-1})} \loc_\balpha (\widetilde{\Gamma_{\rho}(\tau_1^\rho)})^{-1}\right) e^\loc_\balpha(\tau_1^{\tau_2\rho})\right)_{w(\tau_1\tau_2\rho)} \end{array} \]
 and
 \[ \begin{array}{rl} & \barg_{\ga^+,v,\rho}(\tau_2) \\ =&

 \left( \left( \loc_\balpha \widetilde{\Gamma_{\rho}(\tau_2^\rho)}\right)^{-1} e^\loc_\balpha(\tau_2^\rho)\right)_{w(\tau_2\rho)} \\ =&
 
  \left(  \left(\widetilde{\Gamma_{\rho}(\tau_2^\rho)}\right)^{-1} e^\glob_\balpha(\tau_2^\rho)\right) ({}^{(\tau_2^\rho)^{-1}}\beta(\tau_2^\rho))_{w(\tau_2\rho)} \\ =&
 
 \left( \left( \loc_\balpha \widetilde{\Gamma_{\rho}(\tau_2^\rho)}\right)^{-1} e^\loc_\balpha(\tau_2^\rho)\right)_{w(\tau_1\tau_2\rho)} ({}^{(\tau_2^\rho)^{-1}} \beta(\tau_2^\rho))_{w(\tau_2 \rho)}/({}^{(\tau_2^\rho)^{-1}} \beta(\tau_2^\rho))_{w(\tau_1\tau_2 \rho)}. \end{array} \] 
 Thus
  \[ \begin{array}{rl} &\barg_{\ga^+,v,\tau_2\rho}(\tau_1) \barg_{\ga^+,v,\rho}(\tau_2) \\ & \\ =&

  \left(  \beta((\tau_2^\rho)^{-1})_{w(\tau_2\rho)}\beta((\tau_2^\rho)^{-1})^{-1}\left( \loc_\balpha \widetilde{\Gamma_{\rho}(\tau_2^\rho)}\right)^{-1} e^\loc_\balpha(\tau_2^\rho)e^\loc_\balpha((\tau_2^\rho)^{-1})\loc_\balpha (\widetilde{\Gamma_{\rho}(\tau_1^\rho)})^{-1} e^\loc_\balpha(\tau_1^\rho \tau_2^\rho) \right. \\ & \left. e^\loc_\balpha(\tau_1^\rho \tau_2^\rho)^{-1} 
 e^\loc_\balpha((\tau_2^\rho)^{-1})^{-1}  \beta((\tau_2^\rho)^{-1}) \beta((\tau_2^\rho)^{-1})_{w(\tau_2\rho)}^{-1} 
  e^\loc_\balpha(\tau_1^{\tau_2\rho})\right)_{w(\tau_1\tau_2\rho)} \\ &
 ({}^{(\tau_2^\rho)^{-1}} \beta(\tau_2^\rho))_{w(\tau_2 \rho)}/({}^{(\tau_2^\rho)^{-1}} \beta(\tau_2^\rho))_{w(\tau_1\tau_2 \rho)} \\ & \\ =&

 \beta((\tau_2^\rho)^{-1})_{w(\tau_2\rho)}\beta((\tau_2^\rho)^{-1})_{w(\tau_1\tau_2\rho)}^{-1} ({}^{(\tau_2^\rho)^{-1}}(\alpha^\loc(\tau_2^\rho, (\tau_2^\rho)^{-1})\alpha^\loc(1,1)))_{w(\tau_1\tau_2\rho)} \barg_{\ga^+,v,\rho}(\tau_1\tau_2)  \\ & \left( e^\loc_\balpha(\tau_1^\rho \tau_2^\rho)^{-1} 
 e^\loc_\balpha((\tau_2^\rho)^{-1} )^{-1}  \beta((\tau_2^\rho)^{-1}) \beta((\tau_2^\rho)^{-1})_{w(\tau_2\rho)}^{-1} 
  e^\loc_\balpha((\tau_2^\rho)^{-1}\tau_1^{\rho} \tau_2^\rho)\right)_{w(\tau_1\tau_2\rho)} \\ &
 ({}^{(\tau_2^\rho)^{-1}} \beta(\tau_2^\rho))_{w(\tau_2 \rho)}/({}^{(\tau_2^\rho)^{-1}} \beta(\tau_2^\rho))_{w(\tau_1\tau_2 \rho)} \\ & \\  =&

\barg_{\ga^+,v,\rho}(\tau_1\tau_2) (\beta((\tau_2^\rho)^{-1}) {}^{(\tau_2^\rho)^{-1}} \beta(\tau_2^\rho))_{w(\tau_2\rho)}(\beta((\tau_2^\rho)^{-1}){}^{(\tau_2^\rho)^{-1}} \beta(\tau_2^\rho))_{w(\tau_1\tau_2\rho)}^{-1}  \\ &
  ({}^{(\tau_2^\rho)^{-1}}(\alpha^\loc(\tau_2^\rho, (\tau_2^\rho)^{-1})\alpha^\loc(1,1)))_{w(\tau_1\tau_2\rho)} \\ &  {}^{(\tau_2^\rho)^{-1}(\tau_1^\rho)^{-1} \tau_2^\rho}(\beta((\tau_2^\rho)^{-1})_{w(\tau_1\tau_2\rho(\tau_2^\rho)^{-1}(\tau_1^\rho)^{-1} \tau_2^\rho)} /
\beta((\tau_2^\rho)^{-1})_{w(\tau_2\rho)}) \\ & ({}^{(\tau_2^\rho)^{-1}(\tau_1^\rho)^{-1} \tau_2^\rho}\alpha^\loc((\tau_2^\rho)^{-1}, \tau_1^\rho\tau_2^\rho))_{w(\tau_1\tau_2\rho)}^{-1}
 \\ & \\ =&

\barg_{\ga^+,v,\rho}(\tau_1\tau_2) (\pi_{w(\tau_2\rho)}/\pi_{w(\tau_1\tau_2\rho)} )(\alpha^\loc((\tau_2^\rho)^{-1},\tau_2^\rho) \beta(1)/\alpha^\glob((\tau_2^\rho)^{-1},\tau_2^\rho)) \\ & ({}^{(\tau_2^\rho)^{-1}}(\alpha^\loc(\tau_2^\rho, (\tau_2^\rho)^{-1})\alpha^\loc(1,1)))_{w(\tau_1\tau_2\rho)} ({}^{(\tau_2^\rho)^{-1}(\tau_1^\rho)^{-1} \tau_2^\rho}\alpha^\loc((\tau_2^\rho)^{-1}, \tau_1^\rho\tau_2^\rho))_{w(\tau_1\tau_2\rho)}^{-1}
 \\ & \\ =&
 
 \barg_{\ga^+,v,\rho}(\tau_1\tau_2) (\pi_{w(\tau_2\rho)}/\pi_{w(\tau_1\tau_2\rho)} )(\alpha^\loc((\tau_2^\rho)^{-1},\tau_2^\rho) \alpha^\loc(1,1)/\alpha^\glob((\tau_2^\rho)^{-1},\tau_2^\rho) \alpha^\glob(1,1)) \\ & ({}^{(\tau_2^\rho)^{-1}}(\alpha^\loc(\tau_2^\rho, (\tau_2^\rho)^{-1})\alpha^\loc(1,1)))_{w(\tau_1\tau_2\rho)} ({}^{(\tau_2^\rho)^{-1}(\tau_1^\rho)^{-1} \tau_2^\rho}\alpha^\loc((\tau_2^\rho)^{-1}, \tau_1^\rho\tau_2^\rho))_{w(\tau_1\tau_2\rho)}^{-1}
 \\ & \\ =&

 \barg_{\ga^+,v,\rho}(\tau_1\tau_2) (\pi_{w(\tau_1\tau_2\rho)}/\pi_{w(\tau_2\rho)} )(\alpha^\glob((\tau_2^\rho)^{-1},\tau_2^\rho) \alpha^\glob(1,1)) \\ &
\alpha^\loc((\tau_2^\rho)^{-1},\tau_2^\rho)_{w(\tau_2\rho)} \alpha^\loc(1,1)_{w(\tau_2\rho)}\alpha^\loc((\tau_2^\rho)^{-1},\tau_2^\rho)_{w(\tau_1\tau_2\rho)}^{-1} \alpha^\loc(1,1)_{w(\tau_1\tau_2\rho)}^{-1}\\ & ({}^{(\tau_2^\rho)^{-1}}(\alpha^\loc(\tau_2^\rho, (\tau_2^\rho)^{-1})\alpha^\loc(1,1)))_{w(\tau_1\tau_2\rho)} ({}^{(\tau_2^\rho)^{-1}(\tau_1^\rho)^{-1} \tau_2^\rho}\alpha^\loc((\tau_2^\rho)^{-1}, \tau_1^\rho\tau_2^\rho))_{w(\tau_1\tau_2\rho)}^{-1}
 \\ & \\ \equiv&

\barg_{\ga^+,v,\rho}(\tau_1\tau_2) .

\end{array}\] 

 \pfend

 \subsection{Comparing correstrictions}\label{corres}
 
 Suppose that $F$ is a number field, that $v$ is a place of $F$ and that $E \supset F$ is a finite Galois extensions of $F$. Also let 
\[ \ga^+=(\ga,\{ [\Gamma_\rho]\}) \in \cH(E/F)^+ ;\]
and suppose $\tau \in \Aut(\barFv/F)$. 
Now suppose further that $T/F$ is a torus which splits over $E$  and that $\mu$ a cocharacter of $T$ defined over $\barFv$. 
 
 If $\rho:E \into \barFv$ is $F$-linear and $\balpha=(\alpha^\glob,\alpha^\loc,\beta)\in \ga$, then we saw in lemma \ref{corr} how to find an element $b_\rho \in T(\A_E)$ with
 \[ \loc_\ga \corr_{\alpha^\glob} ({}^{\rho^{-1}} \mu) \circ (\pi_{w(\rho)}/\pi_{w(\tau \rho)})= {}^{b_\rho}  \corr_{\alpha^\loc} ({}^{\rho^{-1}} \mu)\circ (\pi_{w(\rho)}/\pi_{w(\tau \rho)}). \]
 However it turns out that the choice of $\ga^+$ above $\ga$ allows us to choose $b_\rho$ more canonically, so that it is largely independent of the choice of $\rho$. That is what we will explain next.
 
If $\tau \in \Aut(\barF_v/F)$ define
\[ \begin{array}{rcl} \barb_{\ga^+,v,\mu,\tau} &=& \prod_{\eta \in \Gal(E/F)} \teta^{-1} ({}^{\rho^{-1}}\mu)((\beta(\teta)_{w(\tau \rho)}/\beta(\teta)_{w(\rho)}) \barg_{\ga^+,v,\rho}(\tau)) \\ &\in &T(\A_D)/\overline{T(F)T(F_\infty)^0}T(D)T(D_v), \end{array} \]
where $\rho:E \into \barFv$ is $F$-linear, $\balpha=(\alpha^\glob,\alpha^\loc,\beta)\in \ga$ and for each $\eta \in \Gal(E/F)$ we choose a lift $\teta \in \Gal(D/F)$.
We see immediately that this does not depend on the lift of $\barg_{\ga^+,v,\rho}(\tau)$ we choose (because $N_{E/F} (\overline{T(E)T(E_\infty)^0}) \subset \overline{T(F)T(F_\infty)^0}$); nor on the choice of $\balpha \in \ga$ (because we are modding out by $T(D)T(D_v)$). It does not depend on the choices $\teta$ of lifts of elements $\eta \in \Gal(E/F)$, for if $\sigma_\eta \in \Gal(D/E)$ for each $\eta \in \Gal(E/F)$, then
\[ \begin{array}{rl} & \prod_{\eta \in \Gal(E/F)} \teta^{-1} \sigma_\eta^{-1} ({}^{\rho^{-1}}\mu)((\beta(\sigma_\eta \teta)_{w(\tau \rho)}/\beta(\sigma_\eta \teta)_{w(\rho)}) \barg_{\ga^+,v,\rho}(\tau)) \\ = &
 \prod_{\eta \in \Gal(E/F)} \teta^{-1}  ({}^{\rho^{-1}}\mu)({}^{\sigma_\eta^{-1}} (\pi_{w(\tau \rho)}/\pi_{w(\rho)})(\beta(\sigma_\eta) {}^{\sigma_\eta}\beta( \teta)) {}^{\sigma_\eta^{-1}}\barg_{\ga^+,v,\rho}(\tau)) \\ =&
  \prod_{\eta \in \Gal(E/F)} \teta^{-1}  ({}^{\rho^{-1}}\mu)( (\pi_{w(\tau \rho)}/\pi_{w(\rho)})({}^{\sigma_\eta^{-1}}\beta(\sigma_\eta) \beta( \teta)) (\pi_{w(\tau \rho)}/\pi_{w(\rho)})(\beta(\sigma_\eta^{-1}) )\barg_{\ga^+,v,\rho}(\tau)) \\ =&
  \prod_{\eta \in \Gal(E/F)} \teta^{-1}  ({}^{\rho^{-1}}\mu)( (\pi_{w(\tau \rho)}/\pi_{w(\rho)})( \beta( \teta)) \barg_{\ga^+,v,\rho}(\tau)) \\ &
  \prod_{\eta \in \Gal(E/F)} \teta^{-1}  ({}^{\rho^{-1}}\mu)( (\pi_{w(\tau \rho)}/\pi_{w(\rho)})(\beta(\sigma_\eta^{-1}){}^{\sigma_\eta^{-1}}\beta(\sigma_\eta)) ) \\ =&
  \prod_{\eta \in \Gal(E/F)} \teta^{-1}  ({}^{\rho^{-1}}\mu)( (\pi_{w(\tau \rho)}/\pi_{w(\rho)})( \beta( \teta)) \barg_{\ga^+,v,\rho}(\tau)) \\ &
  \prod_{\eta \in \Gal(E/F)} \teta^{-1}  ({}^{\rho^{-1}}\mu)( (\pi_{w(\tau \rho)}/\pi_{w(\rho)})(\beta(1)) )\\ =&
  \prod_{\eta \in \Gal(E/F)} \teta^{-1}  ({}^{\rho^{-1}}\mu)( (\pi_{w(\tau \rho)}/\pi_{w(\rho)})( \beta( \teta)) \barg_{\ga^+,v,\rho}(\tau)),
\end{array} \] 
where we have applied lemma \ref{betacocyc} and part \ref{gcocyc5} of lemma \ref{gcocyc}. Finally $\barb_{\ga^+,v,\mu,\tau}$ does not depend on the choice of $F$-linear $\rho:E \into \barFv$. Indeed if $\sigma \in \Gal(E/F)$, then
\[ \begin{array}{rl} & \prod_{\eta \in \Gal(E/F)} \teta^{-1} ({}^{(\rho\sigma)^{-1}}\mu)((\beta(\teta)_{w(\tau \rho\sigma)}/\beta(\teta)_{w(\rho\sigma)}) \barg_{\ga^+,v,\rho\sigma}(\tau)) \\ 

=& \prod_{\eta \in \Gal(E/F)} \teta^{-1} \sigma ({}^{(\rho\sigma)^{-1}}\mu)((\pi_{w(\tau\rho\sigma)}/\pi_{w(\rho \sigma)} )(\beta(\sigma^{-1}\teta))  \barg_{\ga^+,v,\rho\sigma}(\tau)) \\ 

= &\prod_{\eta \in \Gal(E/F)} \teta^{-1}  ({}^{\rho^{-1}}\mu)\left({\sigma}(\pi_{w(\tau\rho\sigma)}/\pi_{w(\rho \sigma)} )(\beta(\teta))  {}^\sigma \barg_{\ga^+,v,\rho \sigma}(\tau)\right)  \\

=& \prod_{\eta \in \Gal(E/F)} \teta^{-1}  ({}^{\rho^{-1}}\mu)\left( (\pi_{w(\tau \rho)}/\pi_{w(\rho)})({}^\sigma \beta(\sigma^{-1}\teta)) \barg_{\ga^+,v,\rho }(\tau){}^\sigma ((\pi_{w(\rho\sigma)}/\pi_{w(\tau \rho\sigma)}) ( \beta(\sigma^{-1}))) \right) \\

=& \prod_{\eta \in \Gal(E/F)} \teta^{-1}  ({}^{\rho^{-1}}\mu)\left( (\pi_{w(\tau \rho)}/\pi_{w(\rho)})({}^\sigma \beta(\sigma^{-1}) \beta(\teta)) \barg_{\ga^+,v,\rho }(\tau)
(\pi_{w(\rho)}/\pi_{w(\tau \rho)}) ({}^\sigma \beta(\sigma^{-1})) \right),

\end{array} \]
again using lemma \ref{betacocyc} and part \ref{gcocyc1} of lemma \ref{gcocyc}.

We claim that we have
\[ \prod_{\rho:E \into \barFv} ({}^{\rho^{-1}} \mu)(\barg_{\ga^+,v,\rho}(\tau)) = \barb_{\ga^+, v, \mu,\tau} \in 
T(\A_D)/\overline{T(E)T(E_\infty)^0}T(D)T(D_v),\]
where $\rho$ runs over $F$-linear embeddings $E \into \barF_v$. Although this new expression for $\barb_{\ga^+, v, \mu,\tau}$ is perhaps simpler than our original definition, it is less precise - it gives the value only in $T(\A_D)/\overline{T(E)T(E_\infty)^0}T(D)T(D_v)$ and not in $T(\A_D)/\overline{T(F)T(F_\infty)^0}T(D)T(D_v)$. 
To prove the claim note that
\[ \begin{array}{rl} & \prod_{\rho:E \into \barFv} ({}^{\rho^{-1}} \mu)(\barg_{\ga^+,v,\rho}(\tau)) \\

=& \prod_{\eta\in \Gal(E/F)} \teta^{-1} ({}^{\rho_0^{-1}} \mu)({}^{\teta}\barg_{\ga^+,v,\rho_0 \eta }(\tau)) \\

=& \prod_{\eta\in \Gal(E/F)} \teta^{-1} ({}^{\rho_0^{-1}} \mu)(\barg_{\ga^+,v,\rho_0 }(\tau) {}^{\teta}((\pi_{w(\rho_0 \eta)}/\pi_{w(\tau\rho_0\eta)})(\beta(\teta^{-1})))) \\

=& \prod_{\eta\in \Gal(E/F)} \teta^{-1} ({}^{\rho_0^{-1}} \mu)(\barg_{\ga^+,v,\rho_0 }(\tau) (\pi_{w(\tau \rho_0)}/\pi_{w(\rho_0)})(\beta(\teta))) \\ & 
\prod_{\eta\in \Gal(E/F)} \teta^{-1} ({}^{\rho_0^{-1}} \mu)((\pi_{w(\rho_0 )}/\pi_{w(\tau\rho_0)})(\beta(\teta){}^\teta\beta(\teta^{-1})))) \\

=& \barb_{\ga^+,v,\mu,\tau} \prod_{\eta\in \Gal(E/F)} \teta^{-1} ({}^{\rho_0^{-1}} \mu)((\pi_{w(\rho_0 )}/\pi_{w(\tau\rho_0)})(\beta(1))) \\

=& \barb_{\ga^+,v,\mu,\tau},

\end{array} \]
where $\rho_0 : E \into \barF_v$ is $F$-linear and $\teta \in \Gal(D/F)$ is a fixed lift of $\eta \in \Gal(E/F)$, and where we again use lemma \ref{betacocyc} and part \ref{gcocyc1} of lemma \ref{gcocyc}. 

\begin{lem}\label{bbarfun} Suppose that $C \supset D \supset E \supset F$ are finite Galois extensions of a number field $F$, that $T$ and $T'/F$ are tori split by $E$, that $v$ is a place of $F$, that $\tau \in \Aut(\barF_v/F)$ and that $\mu \in X_*(T)(\barF_v)$.

\begin{enumerate}
 \item\label{bbar5} If $\ga^+\in \cH(E/F)_D^+$ and $\chi:T \ra T'$ is a morphism of algebraic groups over $F$, then $\barb_{\ga^+,v,\chi \circ \mu,\tau}=\chi(\barb_{\ga^+,v,\mu,\tau})$.
 
 \item If $\ga^+ \in \cH(E/F)_D^+$, then $\barb_{\inf_{C/D}\ga^+,v,\mu,\tau} = \barb_{\ga^+,v,\mu,\tau}$.
 
 \item If $\ga^+ \in \cH(D/F)_C^+$, then $\barb_{\eta_{D/E,*} \ga^+,v,\mu,\tau} = \barb_{\ga^+,v,\mu,\tau}$.
 
 \item\label{bbar6} If $\ga^+ \in \cH(E/F)_D^+$ and $t \in T_{2,E}(\A_D)$, then
\[ \begin{array}{rcl} \barb_{{}^t\ga^+,v,\mu,\tau} &=  &\barb_{\ga^+,v,\mu,\tau}\prod_{\rho}  ({}^{\rho^{-1}}\mu) \circ (\pi_{w(\rho)}/ \pi_{w(\tau \rho)})(t) \\ &=& \barb_{\ga^+,v,\mu,\tau} \prod_\rho ({}^{\rho^{-1}}(\mu/{}^\tau \mu) )(t_{w(\rho)})\end{array}  \]
where $\rho$ runs over $F$-linear embeddings $\rho:E \into \barFv$.

 \end{enumerate} \end{lem}
 
 \pfbegin The first part is clear from the definition, and the second follows immediately from part \ref{gcocyc3} of lemma \ref{gcocyc}. For the third we have
 \[ \begin{array}{rl} & \barb_{\eta_{D/E,*} \ga^+,v,\mu,\tau} \\
 
 =& \prod_{\eta \in \Gal(E/F)} \teta^{-1} ({}^{\rho^{-1}}\mu)\left( (\prod_{\zeta \in \Gal(D/E)} \beta(\teta)_{\zeta u(\tau \rho)}/\prod_{\zeta \in \Gal(D/E)}\beta(\teta)_{\zeta u(\rho)})\right. \\ &
 \left. \prod_{\zeta \in \Gal(D/E)} \barg_{\ga^+,v,\rho \zeta^{-1}}(\tau)  \right) \\
 
 =& \prod_{\eta \in \Gal(E/F)} \teta^{-1} ({}^{\rho^{-1}}\mu)\left( \prod_{\zeta \in \Gal(D/E)} \left( (\pi_{\zeta u(\tau\rho)}/\pi_{\zeta u(\rho)})( \beta(\teta))  {}^\zeta \barg_{\ga^+,v,\rho }(\tau) \beta(\zeta)_{\zeta u(\rho)}/\beta(\zeta)_{\zeta u(\tau \rho)} \right) \right) \\
 
  =& \prod_{\eta \in \Gal(E/F)} \teta^{-1} ({}^{\rho^{-1}}\mu)\left( \prod_{\zeta \in \Gal(D/E)} \left( (\pi_{\zeta u(\tau\rho)}/\pi_{\zeta u(\rho)})( \beta(\teta)/\beta(\zeta))  {}^\zeta \barg_{\ga^+,v,\rho }(\tau)  \right) \right) \\
  
    =& \prod_{\eta \in \Gal(E/F)} \teta^{-1} ({}^{\rho^{-1}}\mu)\left( \prod_{\zeta \in \Gal(D/E)} \left( (\pi_{\zeta u(\tau\rho)}/\pi_{\zeta u(\rho)})( {}^\zeta\beta(\zeta^{-1}\teta))  {}^\zeta \barg_{\ga^+,v,\rho }(\tau)  \right) \right) \\
    
        =& \prod_{\eta \in \Gal(E/F)}\prod_{\zeta \in \Gal(D/E)} \teta^{-1} ({}^{\rho^{-1}}\mu)\left(  {}^\zeta \left( (\pi_{u(\tau\rho)}/\pi_{ u(\rho)})( \beta(\zeta^{-1}\teta))   \barg_{\ga^+,v,\rho }(\tau)  \right) \right) \\
        
=& \prod_{\eta \in \Gal(E/F)}\prod_{\zeta \in \Gal(D/E)} (\zeta^{-1}\teta)^{-1} ({}^{\rho^{-1}}\mu) \left( (\pi_{u(\tau\rho)}/\pi_{ u(\rho)})( \beta(\zeta^{-1}\teta))   \barg_{\ga^+,v,\rho }(\tau)  \right)  \\

=& \barb_{ \ga^+,v,\mu,\tau}.
 \end{array} \]
 
  For the fourth part we have
\[ \begin{array}{rl} & \barb_{{}^t\ga^+,v,\mu,\tau} /\barb_{\ga^+,v,\mu,\tau} \\

=& \prod_{\eta \in \Gal(E/F)} \teta^{-1} ({}^{\rho^{-1}}\mu) ( (t/{}^\teta t)_{w(\tau \rho)} ({}^\teta t/t)_{w(\rho)} (t_{w(\rho)}/t_{w(\tau\rho)})) \\

=& \prod_{\eta \in \Gal(E/F)} \teta^{-1} ({}^{\rho^{-1}}\mu) 
((\pi_{w(\rho)}/\pi_{w(\tau \rho)})({}^\teta t)) \\

=& \prod_{\eta \in \Gal(E/F)}  ({}^{(\rho\eta)^{-1}}\mu) 
((\pi_{w(\rho\eta)}/\pi_{w(\tau \rho \eta)})( t)) \\

=&  \prod_{\rho}  ({}^{\rho^{-1}}\mu) (t_{w(\rho)}/ t_{w(\tau \rho)}) \\

=& \prod_{\rho}  ({}^{\rho^{-1}}\mu) (t_{w(\rho)})/ \prod_{\rho}  ({}^{\rho^{-1}}\mu)(t_{w(\tau \rho)})  \\

=& \prod_{\rho}  ({}^{\rho^{-1}}\mu) (t_{w(\rho)})/ \prod_{\rho}  ({}^{\rho^{-1}\tau}\mu)(t_{w( \rho)})  \\

=& \prod_{\rho}  ({}^{\rho^{-1}}(\mu/{}^\tau \mu)) (t_{w(\rho)}).
\end{array} \]

\pfend

\begin{cor}\label{bbarcor} Suppose that $D \supset E \supset F$ are finite Galois extensions of a number field $F$, that $T/F$ is a torus split by $E$, that $v$ is a place of $F$, that $\tau \in \Aut(\barF_v/F)$ and that $\mu \in X_*(T)(\barF_v)$. Suppose also that $\ga_E^+ \in \cH(E/F)^+$ and $\ga_D^+ \in \cH(D/F)^+$. Then we can find a $t \in T_{2,E}(\A_D)$ such that ${}^t\inf_{D/E} \ga_E^+=\eta_{D/E,*} \ga_D^+ \in \cH(E/F)_D^+$, in which case 
\[ \begin{array}{rcl} \barb_{\ga_E^+,v,\mu,\tau} &=& \barb_{\ga_D^+,v,\mu,\tau}  \prod_\rho ({}^{\rho^{-1}}\mu)(t_{w(\tau\rho)}/t_{w( \rho)})) \\ &=& 
\barb_{\ga_D^+,v,\mu,\tau}  \prod_\rho ({}^{\rho^{-1}}({}^\tau \mu/\mu)(t_{w( \rho)})),\end{array} \]
where $\rho$ runs over $F$-linear embeddings $E \into \barFv$.
\end{cor}

\begin{lem} \label{bbar0} Suppose that $D \supset E \supset F$ are finite Galois extensions of a number field $F$, that $\ga^+\in \cH(E/F)_D^+$, that $T/F$ is a torus split by $E$, that $v$ is a place of $F$, that $\tau \in \Aut(\barF_v/F)$ and that $\mu \in X_*(T)(\barF_v)$.

\begin{enumerate}
\item\label{bbar1}  If $D=E$ and $\balpha=(\alpha^\glob,\alpha^\loc,\beta)\in \ga$, then there is a lift $b\in T(\A_E)$ of $\barb_{\ga^+,v,\mu,\tau}$ such that 
\[ \loc_\ga \corr_{\alpha^\glob}({}^{\rho^{-1}} \mu) \circ (\pi_{w(\rho)}/\pi_{w(\tau \rho)})={}^{b} \corr_{\alpha^\loc} ({}^{\rho^{-1}} \mu)\circ (\pi_{w(\rho)}/\pi_{w(\tau \rho)}).\]

\item\label{bbar2} If $\tau$ fixes the image of $E$ in $\barFv$, then 
\[  \barb_{\ga^+,v,\mu,\tau} =  \prod_{\eta \in \Gal(E/\Q)} \eta ({}^{\trho^{-1}} \mu)(\widetilde{\Art_E^{-1} \tau^{\trho}|_{E^\ab}})^{-1},\]
where $\trho: \barD \into \barFv$ is any $F$-linear embedding, and where $\widetilde{\Art_E^{-1} \tau^{\trho}|_{E^\ab}}$ is any lift of $\Art_E^{-1} \tau^{\trho}|_{E^\ab}$ to $\A_E^\times$.

\item\label{bbar2.5} If $\tau$ is continuous, then $\barb_{\ga^+,v,\mu,\tau}=1$.

\item\label{bbar3} $\barb_{\ga^+,v,\mu,\tau_1\tau_2}=\barb_{\ga^+,v,{}^{\tau_2}\mu,\tau_1}\barb_{\ga^+,v,\mu,\tau_2}$.

\item\label{bbar4} If $\sigma \in \Gal(D/F)$, then 
\[ \sigma \barb_{\ga^+,v,\mu,\tau} = \barb_{\ga^+,v,\mu,\tau} \prod_{\rho} ({}^{\rho^{-1}}({}^\tau\mu/ \mu) )(\beta(\sigma)_{w(\rho)}) \]
where $\rho$ runs over $F$-linear embeddings $\rho:E \into \barFv$.

\end{enumerate} \end{lem}

\pfbegin
The first part follows immediately from lemma \ref{corr} because $\barg_{\ga^+,v,\rho}(\tau) \in \A_E^\times/E^\times\Delta_EE_v^\times$ and so
\[ \prod_{\eta \in \Gal(E/F)} \eta^{-1}({}^{\rho^{-1}}\mu) (\barg_{\ga^+,v,\rho}(\tau)) \in T(\A_F)/\overline{T(F)T(F_\infty)^0}. \]

The second and third parts are immediate from the definition and lemma \ref{bargc}. 

For the fourth part note that
\[ \begin{array}{rl} & \barb_{\ga^+,v,\mu,\tau_1\tau_2} \\ 

=& \prod_{\eta \in \Gal(E/F)} \teta^{-1} ({}^{\rho^{-1}}\mu)((\beta(\teta)_{w(\tau_1\tau_2 \rho)}/\beta(\teta)_{w(\rho)}) \barg_{\ga^+,v,\rho}(\tau_1\tau_2))\\

=& \prod_{\eta \in \Gal(E/F)} \teta^{-1} ({}^{\rho^{-1}}\mu)((\beta(\teta)_{w(\tau_1\tau_2 \rho)}/\beta(\teta)_{w(\tau_2\rho)}) \barg_{\ga^+,v,\tau_2\rho}(\tau_1)
(\beta(\teta)_{w(\tau_2\rho)}/\beta(\teta)_{w(\rho)})\barg_{\ga^+,v,\rho}(\tau_2))\\

=& \barb_{\ga^+,v,\mu,\tau_2} \prod_{\eta \in \Gal(E/F)} \teta^{-1} ({}^{(\tau_2\rho)^{-1}\tau_2}\mu)((\beta(\teta)_{w(\tau_1\tau_2 \rho)}/\beta(\teta)_{w(\tau_2\rho)}) \barg_{\ga^+,v,\tau_2\rho}(\tau_1))\\

=& \barb_{\ga^+,v,\mu,\tau_2}\barb_{\ga^+,v,{}^{\tau_2}\mu,\tau_1}. \end{array} \]

For the fifth part we have
\[ \begin{array}{rl} & {}^\sigma \barb_{\ga^+,v,\mu,\tau} \\

=& \prod_{\eta \in \Gal(E/F)} \teta^{-1} ({}^{\rho_0^{-1}}\mu)((\beta(\teta\sigma)_{w(\tau \rho_0)}/\beta(\teta\sigma)_{w(\rho_0)}) \barg_{\ga^+,v,\rho_0}(\tau)) \\

=& \prod_{\eta \in \Gal(E/F)} \teta^{-1} ({}^{\rho_0^{-1}}\mu)(((\beta(\teta){}^\teta\beta(\sigma))_{w(\tau \rho_0)}/(\beta(\teta){}^\teta\beta(\sigma))_{w(\rho_0)}) \barg_{\ga^+,v,\rho_0}(\tau)) \\

=& \barb_{\ga^+,v,\mu,\tau} \prod_{\eta \in \Gal(E/F)} \teta^{-1} ({}^{\rho_0^{-1}}\mu)({}^\teta (\beta(\sigma)_{w(\tau \rho_0\eta)}/\beta(\sigma)_{w(\rho_0\eta)}))  \\

=& \barb_{\ga^+,v,\mu,\tau}\prod_{\rho} ({}^{\rho^{-1}} \mu) ( \beta(\sigma)_{w(\tau \rho)}/ \beta(\sigma)_{w(\rho)}) \\

=&\barb_{\ga^+,v,\mu,\tau} \prod_{\rho} ({}^{\rho^{-1}\tau} \mu) ( \beta(\sigma)_{w( \rho)})/ \prod_{\rho} ({}^{\rho^{-1}} \mu)(\beta(\sigma)_{w(\rho)}) \\

=&\barb_{\ga^+,v,\mu,\tau} \prod_{\rho} ({}^{\rho^{-1}}({}^{\tau} \mu/\mu))( \beta(\sigma)_{w( \rho)}). 
\end{array} \]
\pfend

\subsection{Explicit formulae}\label{trivdata}

Return to the situation discussed at the end of section \ref{explicit} with $F=\Q$ and $E$ totally imaginary. We chose an infinite place $w(\infty)$ of $E$ and an isomorphism$\rho_0:E_{w(\infty)} \iso \C$ and an element $\balpha_0 \in \cZ(E/\Q)$. Choose Galois rigidification data $\Gamma_{0,\rho_0}: \Gal(E^\ab/\Q) \liso W_{E/\Q,\balpha_0}/\Delta_E$ for $[\balpha_0]$ such that
\[ \Gamma_{0,\rho_0}\circ \theta_{\rho_0}=\iota_{w(\infty)}^{\balpha_0} \circ \tTheta_0. \]
Thus in particular $\Gamma_{0,\rho_0}$ is adapted to $\rho_0$.
We may extend $[\Gamma_{0,\rho_0}]$ to complete rigidification data $\{ [\Gamma_{0,\rho}]\}$ for $[\balpha_0]$. Then we will set
\[ \ga_0^+=([\balpha_0], \{ [\Gamma_{0,\rho}]\}) \in \cH(E/\Q)^+. \]
We will also choose a lift $\tGamma_{0,\rho_0}$ of $\Gamma_{0,\rho_0}$ to an isomorphism $W_{E^\ab/\Q} \iso W_{E/\Q,\ga_0}$.

If $T/\Q$ is a torus split by $E$, if $\mu \in X_*(T)(\C)$ and if $\tau \in \Aut(\C)$, then, keeping the notation of section \ref{explicit},
we have:
\[ \begin{array}{rl} & \barb_{\ga_0^+,\infty,\mu,\tau} \\ = &

 \prod_{\eta \in \Gal(E/\Q)} \eta^{-1} ({}^{\rho_0^{-1}}\mu)(
 (\beta(\eta)_{w(\tau \rho_0)}/\beta(\eta)_{w(\rho_0)})((\loc_{\balpha_0} \Gamma_{0,\rho_0}(\tau^{\rho_0}))^{-1} e^\loc_{\balpha_0}(\tau^{\rho_0}))_{w(\tau \rho_0)})
 \\ =&
 
  \prod_{\eta \in \Gal(E/\Q)} \eta^{-1} ({}^{\rho_0^{-1}}\mu)(
 (\beta(\eta)_{w(\tau \rho_0)}/\beta(\eta)_{w(\rho_0)})( \tGamma_{0,\rho_0}(\tau^{\rho_0})^{-1} \beta(\tau^{\rho_0})e^\glob_{\balpha_0}(\tau^{\rho_0}))_{w(\tau \rho_0)})
 \\ =&
 
   \prod_{\eta \in \Gal(E/\Q)} \eta^{-1} ({}^{\rho_0^{-1}}\mu)(
 ( \tGamma_{0,\rho_0}(\tau^{\rho_0})^{-1} e^\glob_{\balpha_0}(\tau^{\rho_0})) {}^{\tau^{\rho_0,-1}} (\beta(\tau^{\rho_0})_{w(\rho_0)})\beta(\eta)_{w(\tau \rho_0)}/\beta(\eta)_{w(\rho_0)})
 \\ =&
 
    \prod_{\eta \in \Gal(E/\Q)} \eta^{-1} ({}^{\rho_0^{-1}}\mu)(
 ( \tGamma_{0,\rho_0}(\tau^{\rho_0})^{-1} e^\glob_{\balpha_0}(\tau^{\rho_0})) {}^{\tau^{\rho_0,-1}} \alpha(\tau^{\rho_0},\tau^{\rho_0,-1})^{-1}\alpha(\eta,\eta^{-1} \tau^{\rho_0,-1})^{-1} \alpha(\eta,\eta^{-1}))
 \\ =&
 
     \prod_{\eta \in \Gal(E/\Q)} \eta^{-1} ({}^{\rho_0^{-1}}\mu)(
  \tGamma_{0,\rho_0}(\tau^{\rho_0})^{-1} e^\glob_{\balpha_0}(\tau^{\rho_0}) e^\glob_{\balpha_0}(\tau^{\rho_0})^{-1}  (e^\glob_{\balpha_0}(\tau^{\rho_0})e^\glob_{\balpha_0}(\tau^{\rho_0,-1}))^{-1} e^\glob_{\balpha_0}(\tau^{\rho_0}) \\ &
  (e^\glob_{\balpha_0}(\eta) e^\glob_{\balpha_0}(\eta^{-1}\tau^{\rho_0,-1}) e^\glob_{\balpha_0}(\tau^{\rho_0,-1})^{-1}   )^{-1} (e^\glob_{\balpha_0}(\eta) e^\glob_{\balpha_0}(\eta^{-1})))
 \\ =&

 \prod_{\eta \in \Gal(E/\Q)} \eta ({}^{\rho^{-1}}\mu)
(\tGamma_{0,\rho_0}(\tau^\rho)^{-1} e^\glob_{\balpha_0}(\eta (\tau^\rho)^{-1})^{-1}  e^\glob_{\balpha_0}(\eta)) .

\end{array} \]

\newpage

\section{Taniyama groups}

In a special case the elements $\barb_{\ga^+,v,\mu,\tau}$ are related to Langlands' Taniyama group. We will explain this in this section. Note however that the results of this section are not required for the statements of the main theorems in \cite{st}. They are only required to compare the results of that paper with the work of Langlands and Milne.

\subsection{The Serre torus}

Suppose that $E$ is a finite Galois extension of $\Q$. We will write $R_{E,\C}$ for the restriction of scalars of $\G_m$ from $E \cap \C$ (which we recall means $\rho(E)$ for any $\rho:E \into \C$) to $\Q$. Thus $X^*(R_{E,\C})=\Map(\Gal((E\cap \C)/\Q),\Z)$ with $\Gal(\C^\alg/\Q)$ action given by $({}^\sigma \varphi)(\tau)=\varphi(\sigma^{-1}\tau)$. 

If $\tau \in \Gal(E \cap \C/\Q)$ we define an automorphism $[\tau]$ of $R_{E,\C}/\Q$ by $X^*([\tau])(\varphi)(\tau')=\varphi(\tau'\tau^{-1})$. We have $[\tau_1\tau_2]=[\tau_1][\tau_2]$, i.e. this gives a left action of $\Gal(\C^\alg/\Q)$ on $R_{E,\C}/\Q$. 

There is also a cocharacter $\mu^\can=\mu_E^\can :\G_m \ra R_{E,\C}$ over $\C$ characterized by $X^*(\mu^\can)(\varphi)=\varphi(1)$. We have ${}^\sigma \mu^\can = [\sigma^{-1}] \circ \mu^\can$.

If $T/\Q$ is a torus split by $E$ and if $\mu\in X_*(T)(\C)$ then there is a unique map of tori $\tmu: R_{E,\C} \ra T$ over $\Q$ such that $\mu = \tmu \circ \mu^\can$. ($X^*(\tmu)(\chi)(\tau)=\chi \circ {}^\tau \mu$.) If $\sigma \in \Aut(\C)$ then $\widetilde{{}^\sigma \mu}=\tmu \circ [\sigma^{-1}]$.

If $D$ is a finite Galois extension of $\Q$ containing $E$ then $N_{D/E}=\widetilde{\mu_E^\can}: R_{D,\C} \ra R_{E,\C}$ is a homomorphism also characterized by $X^*(N_{D/E})(\varphi)(\tau)=\varphi(\tau|_E)$. If $T/\Q$ is a torus split by $E$ and if $\mu \in X^*(T)(\C)$, then we get $\tmu_E:R_{E,\C} \ra T$ and $\tmu_D:R_{D,\C} \ra T$. They satisfy:
\[ \tmu_E\circ N_{D/E}  = \tmu_D. \]

We will also write $S_{E,\C}$ for the torus over $\Q$ characterized by
\[ X^*(S_E)=\{ (\varphi,w)\in \Map(\Gal((E\cap\C)/\Q),\Z) \times \Z: \,\, \varphi(\sigma c\tau)+\varphi(\sigma\tau)=w \,\, \forall \sigma, \tau \in \Gal((E\cap\C)/\Q) \}, \]
with a left action of $\Gal(\C^\alg/\Q)$ given by 
\[ \sigma(\varphi,w)=(\tau \mapsto \varphi(\sigma^{-1} \tau), w). \]
It is called the {\em Serre torus}. 
There is an obvious injection $X^*(S_{E,\C}) \into X^*(R_{E,\C})$ (sending $(\varphi, w)$ to $\varphi$) with torsion free cokernel; and hence an epimorphism $R_{E,\C} \onto S_{E,\C}$ with connected kernel $R^1_{E,\C}$. The exact sequence
\[ (0) \lra R_{E,\C}^1 \lra R_{E,\C} \lra S_{E,\C} \lra (0) \]
splits over $E$.
The action of $\Gal(E \cap \C/\Q)$ on $R_{E,\C}$ (via $\tau \mapsto [\tau]$) over $\Q$ descends to an action on $S_{E,\C}$, which we will denote in the same way. It is also characterized by
\[ X^*([\tau]) (\varphi,w)=(\tau' \mapsto \varphi(\tau'\tau^{-1}), w). \] 
We will again denote the composite of $\mu^\can$ with the map $R_{E,\C} \onto S_{E,\C}$ by $\mu^\can \in X_*(S_{E,\C})(\C)$. It is also characterized by
\[ X^*(\mu^\can)(\varphi,w)=\varphi(1) \in \Z \cong X^*(\G_m). \]
Again we have ${}^\sigma \mu^\can=[\sigma^{-1}] \circ \mu^\can$. Also note that the $\{ [\tau] \circ \mu^\can\}_{\tau \in \Gal(\C^\alg/\Q)}$ spans $X_*(S_{E,\C})$. 
There is a second cocharacter $\wt:\G_m \ra S_{E,\C}$ over $\Q$ characterized by
\[ X^*(\wt)(\varphi,w)=w \in \Z \cong X^*(\G_m). \]
Note that $\wt={}^{(c+1)\tau}\mu^\can=([c\tau]\mu^\can)([\tau]\mu^\can)$ for any $\tau \in \Gal(\C^\alg/\Q)$. 
If $D \supset E$ is another finite Galois extension of $\Q$ then $N_{D/E}:R_{D,\C} \ra R_{E,\C}$ descends to a map
 $N_{D/E}:S_{D,\C} \ra S_{E,\C}$, also characterized by 
\[ X^*(N_{D/E})(\varphi,w)=(\tau \mapsto \varphi(\tau|_E),w). \]
If $E_0$ denotes the maximal CM subfield of $E$ then
\[ N_{E/E_0}:S_{E,\C} \liso S_{E_0,\C}. \]

We recall (for instance from \cite{msIII}) that $S_{E,\C}(\Q)$ is a discrete subgroup of $S_{E,\C}(\A^\infty)$ and that $\ker^1(\Q,S_{E,\C})=(0)$.

If $T/\Q$ is any torus split by $E$ and if $\mu \in X_*(T)(\C)$ satisfies 
\begin{itemize}
\item ${}^{\sigma c \sigma^{-1}} \mu$ is independent of $\sigma \in \Gal(\C^\alg/\Q)$,
\item and ${}^{(c+1)}\mu \in X_*(T)(\Q)$;
\end{itemize}
then there is a unique morphism $\tmu:S_{E,\C} \ra T$ over $\Q$ such that $\mu = \tmu \circ \mu^\can$. We have
\[ X^*(\tmu)(\chi) = ( \tau \mapsto \chi \circ {}^\tau \mu, \chi \circ ({}^c\mu \mu)). \]
Note that
\[ \widetilde{{}^\sigma \mu}= \tmu \circ [\sigma^{-1}]. \]

\begin{lem}\label{trivial}\label{vanish} \begin{enumerate} 
\item Suppose that $\chi \in \Z[V_{E,\infty}]_0 \otimes X_*(S_{E,\C})(E) \subset X^*(T_{3,E}) \otimes X_*(S_{E,\C})(E)$. Then
\[ \prod_{\eta \in \Gal(E/\Q)} {}^\eta \chi =1. \]

\item $\prod_{\rho:E \into \C} {}^{\rho^{-1}}({}^\tau \mu^\can/\mu^\can) \circ \pi_{w(\rho)}: T_{2,E} \lra S_{E,\C}$
is trivial.
\end{enumerate}\end{lem}

\pfbegin To prove the first part of the lemma we may replace $E$ with $E \cap \C$.
By linearity we only need to consider the case $\chi = [\tau] \mu^\can (\pi_{v_1}/\pi_{v_2})$. As $[\tau]$ is rationally defined we are further reduced to the case $\chi= \mu^\can (\pi_{v_1}/\pi_{v_2})$. Thought of as an element of $\Hom(T_{2,{E\cap\C}},S_{E,\C})(E\cap\C)$ we have
\[ \begin{array}{rcl} \prod_{\eta \in \Gal(E\cap \C/\Q)} {}^\eta (\mu^\can \circ \pi_v) &=& \prod_{\eta \in \Gal(E\cap \C/\Q)}({}^\eta \mu^\can) \circ \pi_{\eta v} \\ &
=& \prod_{\eta \in \Gal(E\cap\C/\Q)/\Gal((E\cap \C)_v/\R)}({}^\eta \mu^\can {}^{\eta c_v}\mu^\can) \circ \pi_{\eta v}  \\ &=& \prod_{\eta \in \Gal(E\cap\C/\Q)/\Gal((E\cap \C)_v/\R)}({}^\eta \wt) \circ \pi_{\eta v} \\ &=& \wt \circ \prod_{w|\infty} \pi_w. \end{array} \]
As this does not depend on $v$, the first part of the lemma follows.

For the second part note that this character can be rewritten
\[ \begin{array}{rl} &\prod_{\rho:E \into \C} {}^{\rho^{-1}}({}^\tau \mu^\can/\mu^\can) \circ \pi_{w(\rho)} \\ =& \prod_{\eta:E \into \C} {}^{\eta \rho_1^{-1}}({}^\tau \mu^\can/\mu^\can) \circ \pi_{w(\rho_1\eta^{-1})} \\ =&
\prod_{\eta:E \into \C} {}^{\eta \rho_1^{-1}}({}^\tau \mu^\can/\mu^\can) \circ \pi_{\eta w(\rho_1)} \\ =&
\prod_{\eta:E \into \C} {}^{\eta}({}^{ \rho_1^{-1}}({}^\tau \mu^\can/\mu^\can) \circ \pi_{ w(\rho_1)}) \\ =&
\prod_{\eta:E \into \C} {}^{\eta}({}^{ \rho_1^{-1}\tau} \mu^\can \circ \pi_{ w(\rho_1)}) / \prod_{\eta:E \into \C} {}^{\eta}({}^{ \rho_1^{-1}} \mu^\can \circ \pi_{ w(\rho_1)}) \\ =&
\prod_{\eta:E \into \C} {}^{\eta}({}^{ \tau^{\rho_1}\rho_1^{-1}} \mu^\can \circ \pi_{\tau^{\rho_1} w(\tau \rho_1)}) / \prod_{\eta:E \into \C} {}^{\eta}({}^{ \rho_1^{-1}} \mu^\can \circ \pi_{ w(\rho_1)}) \\ =&
\prod_{\eta:E \into \C} {}^{\eta \tau^{\rho_1}}({}^{\rho_1^{-1}} \mu^\can \circ \pi_{ w(\tau \rho_1)}) / \prod_{\eta:E \into \C} {}^{\eta}({}^{ \rho_1^{-1}} \mu^\can \circ \pi_{ w(\rho_1)}) \\ =&
\prod_{\eta:E \into \C} {}^{\eta }({}^{\rho_1^{-1}} \mu^\can \circ (\pi_{ w(\tau \rho_1)}/\pi_{w(\rho_1)})) \\=&1,
\end{array}\]
by the first part.
\pfend

\subsection{The Taniyama group}\label{taniy}

Langlands considers extensions
\[ 1 \lra S_{E,\C} \lra \tS \lra \Gal(E^\ab \cap \C/\Q) \lra 1 \]
(as a pro-algebraic group over $\Q$) such that the induced action of $\Gal(E^\ab \cap \C/\Q)$ on $S_{E,\C}$ is given by $[\,\,\,]$, together with a continuous group theoretic section 
\[ \spl:\Gal(E^\ab\cap \C/\Q) \lra \tS(\A^\infty).\]
 (We will follow \cite{msIII}, which in turn followed \cite{langmarch}. However the two articles use different conventions so it is hard to directly compare the details in the two sources.)
Note that such a pair $(\tS,\spl)$ has no automorphisms (where we consider $\tS$ with its structure of an extension of $\Gal(E^\ab\cap \C/\Q)$ by $S_{E,\C}$). Also note that $\tS(E) \onto \Gal(E^\ab\cap \C/\Q)$. (See \cite{msIII} p235.) Finally note that $\tS|_{\Gal(E^\ab\cap \C/E\cap \C)}$ is abelian. 

Langlands showed that giving such a pair is equivalent to giving an element
\[ \barb \in Z^1(\Gal(E^\ab\cap \C/\Q), (S_{E,\C}(\A_{E\cap\C}^\infty)/S_{E,\C}(E\cap \C))^{\Gal(E\cap\C/\Q)}), \]
where $\Gal(E^\ab\cap\C/\Q)$ acts via $[\,\,\,]$ and $\Gal(E\cap\C/\Q)$ acts by its Galois action on $\A_{E\cap\C}^\infty$, such that $\barb$ lifts to a continuous map
\[ b: \Gal(E^\ab\cap\C/\Q) \lra S_{E,\C}(\A^\infty_{E\cap\C}) \]
such that 
\[ \begin{array}{rcl} \Gal(E^\ab\cap\C/\Q)^2 &\lra& S_{E,\C}(E\cap\C) \\
(\tau_1,\tau_2) & \longmapsto & b(\tau_1) [\tau_1](b(\tau_2)) b(\tau_1\tau_2)^{-1} \end{array} \]
is locally constant. (See proposition 2.7 of \cite{msIII}.)
We will will write 
\[ Z^1_\cts(\Gal(E^\ab\cap\C/\Q), (S_E(\A_{E\cap\C}^\infty)/S_{E,\C}(E\cap\C))^{\Gal(E\cap\C/\Q)})\]
 for the set of elements $\barb \in Z^1(\Gal(E^\ab\cap\C/\Q), (S_{E,\C}(\A_{E\cap\C}^\infty)/S_{E,\C}(E\cap\C))^{\Gal(E\cap\C/\Q)})$ with such a lift. If $\alpha \in \tS(E \cap \C)$ has image $\baralpha \in \Gal(E^\ab \cap \C/\Q)$, then 
 \[ \barb(\baralpha) = \alpha\, \spl (\baralpha^{-1}) . \]

Langlands defines a particular element 
\[ \barb_{E}^\Tan \in Z^1_\cts(\Gal(E^\ab\cap\C/\Q), (S_{E,\C}(\A_{E\cap\C}^\infty)/S_{E,\C}(E\cap\C))^{\Gal(E\cap\C/\Q)})\]
 as follows. There is an exact sequence
\[1 \lra \A_{E\cap\C}^\times/(E\cap\C)^\times \lra W_{E^\ab\cap\C/\Q} \lra \Gal(E\cap \C/\Q) \lra 1. \]
Choose preimages $w_\eta \in W_{E^\ab\cap \C/\Q}$ of each $\eta \in\Gal(E\cap \C/\Q)$ such that the following conditions hold:
\begin{itemize}
\item $w_1=1$
\item and there is a set of representatives $1 \in H \subset \Gal(E\cap\C/\Q)$ for $\Gal(E\cap\C/\Q)/\Gal(\C/\R)$ such that $w_{\eta c}=w_\eta \theta_{\rho_0}(j)$ for all $\eta \in H$. (Here $\rho_0$ denotes the identity embedding $E^\ab\cap \C \into \C$.)
\end{itemize}
If $w \in W_{E^\ab\cap \C/\Q}$ then
\[ w_\eta w = \baru_{v,\eta, w} w_{\eta \barw} \]
where $\baru_{v,\eta, w} \in \A_{E\cap\C}^\times/(E\cap\C)^\times$ and where $\barw$ denotes the image of $w$ in $\Gal(E\cap \C/\Q)$. Then Langlands takes
\[ \barb_{E}^\Tan(w)= \prod_{\eta \in \Gal(E\cap\C/\Q)} ({}^\eta \mu^\can)(\baru_{v,\eta,w}) \in S_{E,\C}(\A_{E\cap \C}^\infty)/S_{E,\C}(E\cap \C). \]
He verifies that it lies in $(S_{E,\C}(\A_{E\cap\C}^\infty)/S_{E,\C}(E\cap \C))^{\Gal(E\cap \C/\Q)}$; that it doesn't depend on the choices of preimages $w_{\eta}$ (as long as they satisfy the above conditions); that it only depends on the image $\tau$ of $w$ in $\Gal(E^\ab\cap \C/\Q)$ (so we will often write $\barb_{E}^\Tan(\tau)$); and that $\barb_{E}^\Tan \in Z^1_\cts(\Gal(E^\ab\cap\C/\Q), (S_{E,\C}(\A_{E\cap \C}^\infty)/S_{E,\C}(E\cap \C))^{\Gal(E\cap \C/\Q)})$. It also doesn't depend on the choice of $\theta_{\rho_0}$ associated to $\rho_0$. (If $\theta_{\rho_0}$ is replaced by $a \theta_{\rho_0} a^{-1}$ with $a \in \overline{(E\cap\C)^\times ((E\cap\C)_\infty^\times)^0}$, then $w_\eta$ gets replaced by $w_\eta {}^{c}a/a$ if $\eta \not\in H$ and is unchanged if $\eta \in H$. Thus $\baru_{v,\eta,w}$ gets multiplied by
\[ \left\{ \begin{array}{ll} {}^{\eta c} a/{}^\eta a & {\rm if}\,\, \eta \not\in H \\ 1 & {\rm if}\,\, \eta \in H \end{array}\right\} \left\{ \begin{array}{ll} {}^{\eta \barw} a/{}^{\eta \barw c}a & {\rm if}\,\, \eta\barw \not\in H \\ 1 & {\rm if}\,\, \eta\barw \in H \end{array}\right\} . \]
Thus $\barb^\Tan_{E}(w)$ changes by
\[ \begin{array}{rl} & \prod_{\eta \not\in H} \eta \left(\mu^\can ({}^{c}a/a) ({}^{\barw^{-1}} \mu^\can)(a/{}^{c}a)\right) \\ = &
\left(\prod_{\eta \in H} \eta \wt(a) /\prod_{\eta \in \Gal(E\cap\C/\Q)} \eta \mu (a)\right) \left(\prod_{\eta \in \Gal(E\cap\C/\Q)} \eta {}^{\barw^{-1}}\mu (a)/ \prod_{\eta \in H} \eta \wt(a) \right) \\ = & \prod_{\eta \in \Gal(E\cap\C/\Q)} \eta ({}^{\barw^{-1}}\mu /\mu)(a)\\  \in &\overline{S_{E,\C}(\Q)S_{E,\C}(\R)}=S_{E,\C}(\Q)S_{E,\C}(\R)\subset S_{E,\C}(E\cap \C)S_{E,\C}((E\cap \C)_\infty), \end{array}\]
i.e. it is unchanged.)

We will write $(\tS_{E,\C},\spl_{E,\C})$ for the corresponding extension with a finite adelic section. It is usually referred to as the `Taniyama group'.
We also write $\tS_{E,\C,\tau}$ for the pre-image in $\tS_{E,\C}$ of $\tau \in \Gal(E^\ab\cap\C/\Q)$, a right $S_{E,\C}$-torsor.

If $D \supset E$ is another finite Galois extension there is a unique commutative diagram:
\[ \begin{array}{ccccccccc} 
(0) & \lra & S_{D,\C} & \lra & \tS_{D,\C} & \lra & \Gal(D^\ab \cap\C/\Q) & \lra & (0) \\
&& N_{D/E} \da && \tN_{D/E}\da  && \da && \\
(0) & \lra & S_{E,\C} & \lra & \tS_{E,\C} & \lra & \Gal(E^\ab \cap\C/\Q) & \lra & (0) \end{array} \]
with
\[ \tN_{D/E} \circ \spl_{D,\C} = \spl_{E/C}. \]

We remark that $\spl_{E,\C}(c)\in \tS_{E,\C}(\Q)$. Indeed, 
we have 
\[ \baru_{v,\eta,w_{c}}=\left\{ \begin{array}{ll} 1 & {\rm if}\,\, \eta \in H \\ -1_{\eta v} & {\rm if}\,\, \eta \not\in H. \end{array} \right. \]
Thus $\barb^\Tan_{E}(c)=1$ and $\spl_{E,\C}(c)\in \tS_{E,\C}(E)\cap \tS_{E,\C}(\A^\infty)=\tS_{E,\C}(\Q)$.

If $\alpha \in \tS_{E,\C}(E)$ has image $\baralpha \in \Gal(E^\ab\cap\C/\Q)$ we will write 
\[ b_{E,\Tan,\alpha} = \alpha^{-1} \spl(\baralpha) \in S_{E,C}(E)S_{E,\C}(\A^\infty) \subset S_{E,\C}(\A_E^\infty). \]
If $\rho:E \into \C$ then
\[ b_{E,\Tan,\alpha} \equiv \rho^{-1} \barb^\Tan_E(\baralpha^{-1}) \bmod S_{E,\C}(E). \]
We have the following properties:
\begin{enumerate}
\item If $\gamma \in S_{E,\C}(E)$, then $b_{E,\Tan,  \alpha\gamma}=\gamma^{-1} b_{E,\Tan,\alpha}$.
\item\label{form2} If $\rho:E^\ab \into \C$ then 
\[ b_{E,\Tan,\alpha}\equiv\prod_{\eta \in \Gal(E/\Q)} \eta ({}^{\rho^{-1}}\mu^\can)( \ttau^{-1}w_{({}^\rho\eta)\baralpha^{-1}}^{\rho,-1}w_{{}^\rho\eta}^\rho) \bmod S_{E,\C}(E), \]
where $\ttau\in W_{E^\ab/\Q}$ lifts $\baralpha^\rho\in \Gal(E^\ab/\Q)$ and where $w^\rho\in W_{E^\ab/\Q}$ denotes the pull back of $w\in W_{E^\ab\cap\C/\Q}$ along $\rho$. 

\item \[b_{E,\Tan,\alpha_1\alpha_2}= [\baralpha_2^{-1}]( b_{E,\Tan,\alpha_1}) b_{E,\Tan,\alpha_2}. \]
(Indeed $(\alpha_1\alpha_2)^{-1}\spl(\baralpha_1\baralpha_2)=
\alpha_2^{-1} (\alpha_1^{-1}\spl(\baralpha_1)) \alpha_2 (\alpha_2^{-1}\spl(\baralpha_2))$.)

\item 
If $\alpha \in \tS_{E,\C}(E)$ and $\sigma \in \Gal(E/\Q)$, then 
\[ b_{E,\Tan,\alpha} {}^\sigma (b_{E,\Tan,\alpha})^{-1} \in S_{E,\C}(E). \]
(Indeed ${}^\sigma (\alpha^{-1}\spl(\baralpha)) = (({}^\sigma \alpha)^{-1}\alpha) (\alpha^{-1}\spl(\baralpha))$.)

\item\label{bctan} If $D \supset E$ is another finite Galois extension of $\Q$, then $N_{D/E} b_{D,\Tan,\alpha} = b_{E,\Tan,\tN_{D/E}\alpha}$.

\end{enumerate}

\begin{lem}\label{ety2} Any element of $\Gal(E^\ab\cap \C/E\cap\C)$ has a lift to $\tS_{E,\C}(\Q)$. If $a \in \A_{E\cap \C}^\times$, then there is a unique element $\alpha(a) \in \tS_{E,\C}(\Q)$ such that $\alpha(a)$ maps to $\Art_{E \cap \C}(a)$ and 
\[ \alpha(a)^{-1} \spl(\Art_{E \cap \C}(a)) = \prod_{\eta \in \Gal(E/\Q)} \eta ( \mu^\can)(a)^{-1}. \]

Thus if $a \in \A_E^\times$ and $\rho:E \into \C$ then $\alpha(\rho(a))$ is the unique element of $\tS_{E,\C}(\Q)$ which lies above $\Art_{E \cap \C} \rho(a)$ and satisfies
\[ b_{E,\Tan,\alpha(a)}=\prod_{\eta \in \Gal(E\cap\C/\Q)} \eta \mu^\can(\rho(a))^{-1} = \prod_{\eta \in \Gal(E/\Q)} \eta ( {}^{\rho^{-1}} \mu^\can)(a)^{-1}. \]
\end{lem}

\pfbegin
Suppose that $a \in \A_{E\cap \C}^\times$ and that $\alpha \in \tS_{E,\C}(E\cap\C)$ lifts $\Art_{E \cap \C}(a)$. Then 
\[ \alpha \spl(\Art_{E \cap \C}(a))^{-1}=\barb^\Tan_E(\Art_{E \cap \C}(a))\in S_{E,\C}(\A_{E\cap \C}^\infty)/S_{E,\C}(E\cap \C)\]
equals the image of
\[  \prod_{\eta \in \Gal(E\cap \C/\Q)} \eta \mu^\can(a)  \in S_{E,\C}(\A^\infty). \]
Thus
\[ \prod_{\eta \in \Gal(E\cap\C/\Q)} \eta  
\mu^\can (a)=\gamma \alpha \spl(\Art_{E \cap \C}(a))^{-1} \in S_{E,\C}(\A_{E\cap \C}^\infty) \]
for some $\gamma \in S_{E,\C}(E\cap\C)$. If $\eta \in \Gal(E\cap \C/\Q)$ we see that 
\[ {}^\eta(\gamma\alpha)=\gamma \alpha, \]
so that $\gamma \alpha \in \tS_{E,\C}(\Q)$ lifts $\Art_{E \cap \C}(a)$. We set $\alpha(a)=\gamma \alpha$. Then
\[  \begin{array}{rcl} \alpha(a)^{-1} \spl(\Art_{E \cap \C}(a)) &=& \spl(\Art_{E \cap \C}(a)) \alpha(a)^{-1} \\ &=& (\gamma \alpha \spl(\Art_{E \cap \C}(a))^{-1})^{-1} \\ &=& \prod_{\eta \in \Gal(E\cap\C/\Q)} \eta  
\mu^\can (a)^{-1} . \end{array}\]
\pfend

We will call $\alpha(a)$ {\em well placed} with respect to $a$.

\subsection{Relationship between Taniyama groups and the elements $\barb_{\ga^+,\infty,\mu^\can,\tau}$}

\begin{lem}
\label{independ}
If $\alpha \in \tS_{E,\C}(E)$ and $\tau \in \Aut(\C)$ have the same image in $\Gal(E^\ab\cap\C/\Q)$, then 
 \[ \barb_{\ga^+,\infty,\mu^\can,\tau} \in S_{E,\C}(\A_E^\infty)/S_{E,\C}(E) =   R_{E,\C}(\A_E^\infty)/  R_{E,\C}(E) R^1_{E,\C}(\A_E^\infty) \]
 equals the image of $b_{E,\Tan,\alpha}$.
\end{lem}

\pfbegin 
Replacing $\ga^+$ by ${}^t \ga^+$ leaves $\barb_{\ga^+,\infty,\mu^\can,\tau}$ unchanged. (Use corollary \ref{vanish}, lemma \ref{bbar0} and note that $\overline{S_{E,\C}(\Q)S_{E,\C}(\R)}=S_{E,\C}(\Q)S_{E,\C}(\R)$.) Thus it suffices to prove the assertion with $\ga^+=\ga_0^+$, the class defined in section \ref{trivdata}. In this case we may take $\Gamma_{0,\rho_0}(w_\eta^{\rho_0})=e^\glob_{\balpha_0}(\eta^{\rho_0})$, and the result then follows from comparing 
the formula at the end of section \ref{trivdata} with the formula in item (\ref{form2}) of section \ref{taniy}.
\pfend

\begin{lem} Suppose that $D \supset E$ are finite Galois extensions of $\Q$.
\begin{enumerate}
\item $\overline{R_{E,\C}(\Q)} R_{E,\C}(D) = \overline{R_{E,\C}^1(\Q)} R_{E,\C}(D)$. 
\item $( \overline{R_{E,\C}(\Q)} R_{E,\C}(D)) \cap  R^1_{E,\C}(\A_D^\infty) =  \overline{R^1_{E,\C}(\Q)} R^1_{E,\C}(D)$. \end{enumerate}
(Here the closures are taken in $R_{E,\C}(\A_D^\infty)$. \end{lem}

\pfbegin For the first part suppose that $\gamma_i \in R_{E,\C}(\Q)$ tend to $h \in \overline{R_{E,\C}(\Q)}$. As $S_{E,\C}(\Q) \subset S_{E,\C}(\A^\infty)$ is discrete, we may suppose that all the $\gamma_i$ have the same image $\delta \in S_{E,\C}(\Q)$. Then $h= (\lim_{\ra i} \gamma_1^{-1}\gamma_i) \gamma_1\in \overline{R_{E,\C}^1(\Q)} R_{E,\C}(\Q)$.

For the second part 
\[ ( \overline{R_{E,\C}(\Q)} R_{E,\C}(D)) \cap  R^1_{E,\C}(\A_D^\infty)=( \overline{R^1_{E,\C}(\Q)} R_{E,\C}(D)) \cap  R^1_{E,\C}(\A_D^\infty)= \overline{R^1_{E,\C}(\Q)} R_{E,\C}^1(D) .\]
\pfend

\begin{cor}. If $\alpha \in \tS_{E,\C}(E)$ and $\tau \in \Aut(\C)$ have the same image $\Gal(E^\ab\cap\C/\Q)$, then
\[ \barb_{\ga^+,\infty,\mu^\can,\tau} \in    R_{E,\C}(\A_E^\infty)/ \overline{R^1_{E,\C}(\Q)} R_{E,\C}(E)  \]
and
\[ b_{E,\Tan,\alpha} \in R_{E,\C}(\A_E^\infty)/ R^1_{E,\C}(\A_E^\infty) \]
have a unique common lift
\[ b_{\ga^+,\infty,\mu^\can,\alpha} \in    R_{E,\C}(\A_E^\infty)/ \overline{R^1_{E,\C}(\Q)} R^1_{E,\C}(E).  \]
(It is independent of the lift $\tau$ of the image of $\alpha$ in $\Gal(E^\ab\cap \C/\Q)$.)
\end{cor}

\begin{cor} \label{propb} \begin{enumerate}
\item If $\gamma \in R_{E,\C}(E)$ then $b_{\ga^+,\infty,\mu^\can,\alpha\gamma}=\gamma^{-1}  b_{\ga^+,\infty,\mu^\can,\alpha}$.

\item $b_{\ga^+,\infty,\mu^\can,\alpha_1\alpha_2}= [\baralpha_2^{-1}]( b_{\ga^+,\infty,\mu^\can,\alpha_1}) b_{\ga^+,\infty,\mu^\can,\alpha_2}$.

\item $b_{\ga^+,\infty,\mu^\can,\alpha} {}^\sigma (b_{\ga^+,\infty,\mu^\can,\alpha})^{-1} \in R^1_{E,\C}(\A_E^\infty)R_{E,\C}(E)$.

\item Suppose that $D \supset E$ are finite Galois extensions of $\Q$, that $\ga^+_E \in \cH(E/\Q)^+$, that $\ga_D^+ \in \cH(D/\Q)^+$ and that $t \in T_{2,E}(\A_D)$ with $\eta_{D/E,*}\ga_D^+ = {}^t \inf_{D/E}\ga_E^+$. Suppose also that $\alpha_D \in \tS_{D,\C}(D)$ and $\alpha_E \in \tS_{E,\C}(E)$ have the same image in $\Gal(E^\ab \cap \C/\Q)$, so that $\alpha_E^{-1} \tN_{D/E} (\alpha_D) \in S_{E,\C}(D)$. 
Then
\[  \begin{array}{rcl} b_{\ga_E^+,\infty,\mu_E^\can, \alpha_E} &=
& (\alpha_E^{-1}\tN_{D/E} (\alpha_D) ) N_{D/E}(b_{\ga_D^+,\infty,\mu_D^\can, \alpha_D}) \prod_{\rho: E \into \C} {}^{\rho^{-1}} ({}^{\baralpha_E} \mu_E^\can /\mu_E^\can)(t_{w(\rho)}) \\ &\in& R_{E,\C}(\A_D^\infty)/\overline{R^1_{E,\C}(\Q)} R_{E,\C}(D), \end{array} \]
where $\baralpha_E$ denotes the image of $\alpha_E$ in $\Gal(E^\ab \cap \C/\Q)$.

\item If $\ga^+_0$ is as in section \ref{trivdata} and if $\tau \in \Aut(\C)$, then there is an element $\alpha_0(\tau) \in \tS_{E,\C}(E)$ such that  $\tau$ and $\alpha_0(\tau)$ have the same image in $\Gal(E^\ab \cap \C/\Q)$ and
\[ b_{\ga_0^+,\infty,\mu^\can,\alpha_0(\tau)}=\prod_{\eta \in \Gal(E/\Q)} \eta ({}^{\rho^{-1}}\mu^\can)
(\widetilde{\tGamma_{0,\rho_0}(\tau^\rho)}^{-1} e^\glob_{\balpha_0}(\eta (\tau^\rho)^{-1})^{-1}  e^\glob_{\balpha_0}(\eta)) ,
\]
where $\widetilde{\tGamma_{0,\rho_0}(\tau^\rho)}$ is any lift of ${\tGamma_{0,\rho_0}(\tau^\rho)}$ to $W_{E/\Q,\ga}$.

\item If $a \in \A_E^\times$ and $\rho:E \into \C$, then $\alpha(\rho(a)) \in \tS_{E,\C}(\Q)$  lies above $\Art_{E \cap \C}(\rho(a))$ and we have 
\[ b_{\ga^+,\infty,\mu^\can,\alpha(\rho(a))}= \prod_{\eta \in \Gal(E/\Q)} \eta ( {}^{ \rho^{-1}}\mu^\can )(a)^{-1} . \]
\end{enumerate}\end{cor}

\pfbegin
The first part follows from the corresponding property listed at the end of section \ref{taniy}. The second part follows from the corresponding property listed at the end of section \ref{taniy} and properties (\ref{bbar3}) and (\ref{bbar5}) of lemma \ref{bbar0} and the equality ${}^{\baralpha_2}\mu^\can=[\baralpha_2^{-1}] \circ \mu^\can$.
The third part follows from part (\ref{bbar4}) of lemma \ref{bbar0} and lemma \ref{vanish}. The fourth part follows from comment \ref{bctan} before lemma \ref{ety2}, lemma \ref{trivial}, corollary \ref{bbarcor} and part \ref{bbar5} of lemma \ref{bbarfun}. For the fifth part note that if $\alpha \in \tS_{E,\C}(E)$ has the same image as $\tau$ in $\Gal(E^\ab \cap \C/\Q)$, then
\[ \prod_{\eta \in \Gal(E/\Q)} \eta ({}^{\rho^{-1}}\mu^\can)
(\widetilde{\tGamma_{0,\rho_0}(\tau^\rho)}^{-1} e^\glob_{\balpha_0}(\eta (\tau^\rho)^{-1})^{-1}  e^\glob_{\balpha_0}(\eta)) =  \gamma b_{\ga_0^+,\infty,\mu^\can,\alpha}\]
with $\gamma \in  R_{E,\C}(E)$. (See section \ref{trivdata}.)
Thus $\alpha_0(\tau)=\alpha\gamma^{-1}$ will do. The sixth part follows from part (\ref{bbar2}) of lemma \ref{bbar0} and lemma \ref{ety2}.
\pfend

 Note the following simple remark:
\begin{lem}\label{wd} Suppose $\phi_1,\phi_2 \in Z^1_\alg(\cE_3(E/\Q)_\ga,T(E))$ for some torus $T/\Q$, and that $[\phi_1]=[\phi_2] \in H^1_\alg(\cE_3(E/\Q)_\ga,T(E))$. Suppose also that $\loc^\infty(\phi_i)={}^{b_i} 1$ and that $b_1^{-1}b_2 \in T(\A^\infty)$. Then $\phi_1=\phi_2$. \end{lem}
(The point being that $\phi_2={}^\gamma \phi_1$ with $\gamma \in T(E)$, so that $b_2^{-1} \gamma b_1  \in T(\A^\infty)$. Thus $\gamma \in T(\A^\infty) \cap T(E)=T(\Q)$.)

\begin{lem} If $\tb_{\ga^+,\infty,\mu^\can,\alpha}$ is a lift of $ b_{\ga^+,\infty,\mu^\can,\alpha}$ to $R_{E,\C}(\A_E^\infty)$, then there is a unique element $\tphi_{\ga^+,\infty,\mu^\can,\alpha} \in Z^1_\alg(\cE_3(E/\Q)_\ga,R_{E,\C}(E))$ such that:
\begin{enumerate}
\item $\res^\infty \loc_\ga \tphi_{\ga^+,\infty,\mu^\can,\alpha} = {}^{\tb_{\ga^+,\infty,\mu^\can,\alpha}}1$.

\item If $\rho: E \into \C$ and $w(\rho)$ is the corresponding place of $E$, then
\[ \begin{array}{rcl} \kappa(\tphi_{\ga^+,\infty,\mu^\can,\alpha})&=&[{}^{\rho^{-1}}\mu^\can \otimes (w(\rho)-w(\tau\rho))] \\ &=&[{}^{\rho^{-1}}(\mu^\can/{}^{\tau}\mu^\can) \otimes w(\rho)]. \end{array}\]
This is independent of the choice of $\rho$.
\end{enumerate}

Moreover 
\begin{enumerate}
\item[(3)] $\res_{\C/\R} \loc [\tphi_{\ga^+,\infty,\mu^\can,\alpha}]= [(\mu^\can ([c] \circ \mu^\can)/([\tau^{-1}] \circ \mu^\can)([\tau^{-1} c] \circ \mu^\can), (\mu^\can/[\tau]^{-1} \circ \mu^\can)(-1))]$.

\item[(4)] If we replace $\tb_{\ga^+,\infty,\mu^\can,\alpha}$ by $h\gamma \tb_{\ga^+,\infty,\mu^\can,\alpha}$ with $h \in \overline{R_{E,\C}^1(\Q)}$ and $\gamma \in R_{E,\C}^1(E)$, then $\tphi_{\ga^+,\infty,\mu^\can,\alpha}$ changes to ${}^\gamma \tphi_{\ga^+,\infty,\mu^\can,\alpha}$. 
\end{enumerate}
\end{lem}

\pfbegin
Choose $\balpha \in \ga$ and $\rho:E \into \C$.
By part (\ref{bbar1}) of lemma \ref{bbar0}, there is a lift $b \in R_{E,\C}(\A_E^\infty)$ of $\barb_{\ga^+,\infty,\mu^\can,\tau}$ such that
$\res^\infty \loc_\ga \corr_{\alpha^\glob}({}^{\rho^{-1}} \mu) \circ (\pi_{w(\rho)}/\pi_{w(\tau \rho)})={}^{b}1$. Then $\tb_{\ga^+,\infty,\mu^\can,\alpha}=h\gamma b$ for some $h \in \overline{R_{E,\C}(\Q)}$ and $\gamma \in R_{E,\C}(E)$. Set 
\[ \tphi_{\ga^+,\infty,\mu^\can,\alpha} = {}^\gamma  \corr_{\alpha^\glob}({}^{\rho^{-1}} \mu) \circ (\pi_{w(\rho)}/\pi_{w(\tau \rho)}). \]
Then $\tphi_{\ga^+,\infty,\mu^\can,\alpha}$ has the first property by construction and the second by definition of $\kappa$. Moreover
\[ \begin{array}{rcl} {}^{\rho^{-1}}\mu^\can \otimes (w(\rho)-w(\tau\rho)) &=& {}^{\rho^{-1}}\mu^\can \otimes w(\rho)-{}^{\rho^{-1}}\mu^\can \otimes \tau^{\rho,-1}w(\rho)) \\ &=&
({}^{\rho^{-1}}\mu^\can/{}^{\tau^\rho\rho^{-1}}\mu^\can) \otimes w(\rho) \\ &=& 
{}^{\rho^{-1}}(\mu^\can/{}^{\tau}\mu^\can) \otimes w(\rho) \\ &\in& (\Z[V_E]_0 \otimes X_*(R_{E,\C}))_{\Gal(E/\Q)}. \end{array} \]
Because $R_{E,\C}$ is a torus, the second property determines $[\tphi_{\ga^+,\infty,\mu^\can,\alpha}]$ uniquely. Thus the first property implies the uniqueness of $\tphi_{\ga^+,\infty,\mu^\can,\alpha}$ by lemma \ref{wd}.

To prove the third property note that $\kappa(\res_{\C/\R} \loc [\tphi_{\ga^+,\infty,\mu^\can,\alpha}])=[\mu^\can/{}^{\tau}\mu^\can]$, and so
\[ \res_{\C/\R} \loc [\tphi_{\ga^+,\infty,\mu^\can,\alpha}]= [(\mu^\can ([c] \circ \mu^\can)/([\tau^{-1}] \circ \mu^\can)([\tau^{-1} c] \circ \mu^\can), (\mu^\can/[\tau]^{-1} \circ \mu^\can)(-1))]. \]
The final assertion is clear.
\pfend

\begin{lem} \begin{enumerate}

\item If $\alpha \in \tS_{E,\C}(E)$ and $\gamma \in S_{E,\C}(E)$ then we can find $\tgamma \in R_{E,\C}(E)$ lifting $\gamma$ and $h \in \overline{R^1_{E,\C}(\Q)}$ such that $\tphi_{\ga^+,\infty,\mu^\can,\alpha\gamma}={}^{\tgamma^{-1}} \tphi_{\ga^+,\infty,\mu^\can,\alpha}$ and $\tb_{\ga^+,\infty,\mu^\can,\gamma\alpha}=h\tgamma^{-1} \tb_{\ga^+,\infty,\mu^\can,\alpha}$.

\item $\bnu_{\tphi_{\ga^+,\infty,\mu^\can,\alpha}}=\prod_{\rho:E \into \C} {}^{\rho^{-1}}\mu^\can \circ (\pi_{w(\rho)}/\pi_{w(\tau \rho)})= \prod_{\rho:E \into \C} {}^{\rho^{-1}} (\mu^\can/{}^\tau \mu^\can) \circ \pi_{w(\rho)}$, where $\tau\in \Aut(\C)$ has the same image in $\Gal(\C^\alg/\Q)$ as $\alpha$.

 \item Given $\alpha_i \in \tS_{E,\C}(E)$  for $i=1,2$, there exists $\beta \in R^1_{E,\C}(E)$ such that 
\[ h \beta \tb_{\ga^+,\infty,\mu^\can,\alpha_1\alpha_2}\equiv [\baralpha_2^{-1}]_\C( \tb_{\ga^+,\infty,\mu^\can,\alpha_1}) \tb_{\ga^+,\infty,\mu^\can,\alpha} \bmod  \overline{R_{E,\C}^1(\Q)} \]
and
\[  {}^\beta \tphi_{\ga^+,\infty,\mu^\can,\alpha_1\alpha_2}= [\baralpha_2^{-1}]_\C( \tphi_{\ga^+,\infty,\mu^\can,\alpha_1}) \tphi_{\ga^+,\infty,\mu^\can,\alpha_2}, \]
where $\baralpha_2 \in \Gal(E^\ab \cap \C/\Q)$ denotes the image of $\alpha_2$.

\item If $\alpha \in \tS_{E,\C}(E)$ and $\sigma \in \Gal(E/\Q)$, then 
\[ \tb_{\ga^+,\infty,\mu^\can,\alpha} {}^\sigma (\tb_{\ga^+,\infty,\mu^\can,\alpha})^{-1} \in R^1_{E,\C}(\A_E^\infty) R_{E,\C}(E) \subset R_{E,\C}(\A_E^\infty). \]

\item Suppose that $D \supset E$ are finite Galois extensions of $\Q$, that $\ga^+_E \in \cH(E/\Q)^+$, that $\ga_D^+ \in \cH(D/\Q)^+$ and that $t \in T_{2,E}(\A_D)$ with $\eta_{D/E,*}\ga_D^+ = {}^t \inf_{D/E}\ga_E^+$. Suppose also that $\alpha_D \in \tS_{D,\C}(D)$ and $\alpha_E \in \tS_{E,\C}(E)$ have the same image in $\Gal(E^\ab \cap \C/\Q)$, so that $\alpha_E^{-1} \tN_{D/E} (\alpha_D) \in S_{E,\C}(D)$. Choose $\tb_{\ga_E^+,\infty,\mu_E^\can, \alpha_E}$ lifting $b_{\ga_E^+,\infty,\mu_E^\can, \alpha_E}$ and $\tb_{\ga_D^+,\infty,\mu_D^\can, \alpha_D}$ lifting $b_{\ga_D^+,\infty,\mu_D^\can, \alpha_D}$.
Then there exists $\beta \in R^1_{E,\C}(D)$ with 
\[  \tb_{\ga_E^+,\infty,\mu_E^\can, \alpha_E} \equiv \beta (\alpha_E^{-1}\tN_{D/E} (\alpha_D) ) N_{D/E}(\tb_{\ga_D^+,\infty,\mu_D^\can, \alpha_D}) \prod_{\rho: E \into \C} {}^{\rho^{-1}} ({}^{\baralpha_E} \mu_E^\can /\mu_E^\can)(t_{w(\rho)}) \bmod \overline{R^1_{E,\C}(\Q)} \]
and
\[  \inf_{D/E,t} \tphi_{\ga_E^+,\infty,\mu_E^\can,\alpha_E} = {}^{\beta (\alpha_E^{-1}\tN_{D/E} (\alpha_D) )} N_{D/E} \circ \tphi_{\ga_D^+,\infty,\mu_D^\can,\alpha_D} \in Z^1_\alg(\cE_3(D/\Q)_{\ga_D^+},R_{E,\C}(D)).\]
(Again $\baralpha_E$ denotes the image of $\alpha_E$ in $\Gal(E^\ab \cap \C/\Q)$.)

\item If $a \in \A^\times_E$ and $\rho:E \into \C$, then we may take
\[ \tb_{\ga^+,\infty,\mu^\can,\alpha(\rho(a))}= \prod_{\eta \in \Gal(E/\Q)} \eta ( {}^{ \rho^{-1}}\mu^\can )(a)^{-1} \in R_{E,\C}(\A^\infty)\]
and $\tphi_{\ga^+,\infty,\mu^\can,\alpha(\rho(a))}=1$.

\item If $\ga^+_0$ is as in section \ref{trivdata} and $\tau \in \Aut(\C)$, then there is an element $\alpha_0(\tau) \in \tS_{E,\C}(E)$ above $\tau|_{E^\ab \cap\C}$ such that we may take
\[ \tb_{\ga_0^+,\infty,\mu^\can,\alpha_0(\tau)}= \prod_{\eta \in \Gal(E/\Q)} \eta ({}^{\rho^{-1}}\mu^\can)
(\widetilde{\tGamma_{0,\rho_0}(\tau^\rho)}^{-1} e^\glob_{\balpha_0}(\eta (\tau^\rho)^{-1})^{-1}  e^\glob_{\balpha_0}(\eta))_{w(\rho_0)}, \]
where $\widetilde{\tGamma_{0,\rho_0}(\tau^\rho)}$ is a lift of $\tGamma_{0,\rho_0}(\tau^\rho)$ to $\cE^\glob(E/F)^0$.

 \end{enumerate}

\end{lem}

\pfbegin
The first part follows from the first part of corollary \ref{propb}. The second part follows from the construction of $[\tphi_{\ga^+,\infty,\mu^\can,\alpha}]$ as a correstriction. The third part follows from the second part of corollary \ref{propb} and lemma \ref{wd}. The fourth part follows from the third part of corollary \ref{propb}. 

The first half of the fifth part follows from the fourth part of corollary \ref{propb}. For the second half of the fifth part note that  $[\inf_{D/E,t} \tphi_{\ga_E^+,\infty,\mu_E^\can,\alpha_E}]$ and $[ N_{D/E} \circ \tphi_{\ga_D^+,\infty,\mu_D^\can,\alpha_D}]$ in $H^1_\alg(\cE_3(D/\Q),R_{E,\C}(D))$ have the same image under $\kappa$ in
\[ (\Z[V_D]_0 \otimes X_*(R_{E,\C}))_{\Gal(D/\Q)} \liso (\Z[V_E]_0 \otimes X_*(R_{E,\C}))_{\Gal(E/\Q)}, \]
and so are equal. Thus it suffices to check that
\[ \res^\infty \loc_{\ga_D} \inf_{D/E,t} \tphi_{\ga_E^+,\infty,\mu_E^\can,\alpha_E} = \res^\infty \loc_{\ga_D} {}^{\beta (\alpha_E^{-1}\tN_{D/E} (\alpha_D) )}N_{D/E} \circ \tphi_{\ga_D^+,\infty,\mu_D^\can,\alpha_D},\]
or equivalently that
\[ {}^{\bnu_{\tphi_{\ga_E^+,\infty,\mu_E^\can,\alpha_E}}(t)} \inf_{D/E,t} \res^\infty \loc_{\ga_E} \tphi_{\ga_E^+,\infty,\mu_E^\can,\alpha_E} =  {}^{\beta (\alpha_E^{-1}\tN_{D/E} (\alpha_D) )} N_{D/E} \circ (\res^\infty \loc_{\ga_D} \tphi_{\ga_D^+,\infty,\mu_D^\can,\alpha_D}),\]
or equivalently again that
\[ {}^{\bnu_{\tphi_{\ga_E^+,\infty,\mu_E^\can,\alpha_E}}(t) \tb_{\ga_E^+,\infty,\mu_E^\can,\alpha_E}} 1 = {}^{ \beta (\alpha_E^{-1}\tN_{D/E} (\alpha_D) )N_{D/E}(\tb_{\ga_D^+,\infty,\mu_D^\can,\alpha_D}) }1.\]
This follows from the second part and the first half of this part.

The first half of the sixth part follows from the sixth part of corollary \ref{propb}. The second half then follows because $\tb_{\ga^+,\infty,\mu^\can,\alpha}\in R_{E,\C}(\A^\infty)$ and 
$[\tphi_{\ga^+,\infty,\mu^\can,\alpha}]$ is trivial (because in turn its image under $\kappa$ is). 

The seventh part follows from the fifth part of corollary \ref{propb}.
\pfend

\end{document}